\documentclass[12pt, oneside]{amsart}

\renewcommand{\thepart}{\Roman{part}}

\makeatletter
\renewcommand\part{%
  \clearpage
  \thispagestyle{plain}%
  \null\vskip 4em
  \secdef\@part\@spart}

\def\@part[#1]#2{%
  \ifnum \c@secnumdepth >\m@ne
    \refstepcounter{part}%
    \addcontentsline{toc}{part}{\protect\tocpart{}{\thepart}{#1}}%
  \else
    \addcontentsline{toc}{part}{\protect\tocpart{}{}{#1}}%
  \fi
  \markboth{}{}%
  {\centering
   \normalfont
   \large\bfseries \partname\ \thepart: #2\par}%
  \vskip 5em
  \@afterheading}

\def\@spart#1{%
  {\centering
   \normalfont
   \large\bfseries #1\par}%
  \vskip 5em
  \@afterheading}

\renewcommand\l@part[2]{%
  \par\addvspace{16pt plus 2pt}%
  \begingroup
    \parindent\z@ \leftskip\z@ \rightskip\z@
    \centering\normalfont\bfseries
    #1\par
  \endgroup
  \nobreak\addvspace{6pt plus 1pt}}
\makeatother

\usepackage[T1]{fontenc}
\usepackage{palatino}
\usepackage{microtype}

\usepackage{hyperref}
\newcommand\arXiv[1]{\href{http://arxiv.org/abs/#1}{arXiv:#1}}
\usepackage{xcolor}

\hypersetup{
    colorlinks,
    linkcolor={red!50!black},
    citecolor={blue!50!black},
    urlcolor={blue!80!black}
}

\usepackage{graphicx}
\usepackage{booktabs}
\usepackage{tikz-cd}
\usetikzlibrary{decorations.pathmorphing}
\usepackage{amsthm,amssymb}
\usepackage{amsaddr}
\usepackage{mathtools}
\makeatletter
\newcommand*{\relbarfill@}{\arrowfill@\relbar\relbar\relbar}
\newcommand*{\xsimeq}[1]{\ext@arrow 0055{\relbarfill@}{\smash{\raisebox{1pt}{$\scriptstyle#1$}}\vphantom{#1}}{\substack{\vphantom{#1}\\ \raisebox{0.7ex}{$\scriptstyle\sim$}}}}
\makeatother
\usepackage{enumitem}
\setlist[itemize]{leftmargin=1.5em}
\setlist[enumerate,1]{label=\textup{(\roman*)},ref=\roman*,leftmargin=*,align=left,widest=vii}
\newcommand{\itemref}[1]{(\ref{#1})}
\numberwithin{equation}{section}
\newtheorem{theorem}[subsubsection]{Theorem}
\newtheorem{proposition}[subsubsection]{Proposition}
\newtheorem{corollary}[subsubsection]{Corollary}
\newtheorem{definition}[subsubsection]{Definition}
\newtheorem{lemma}[subsubsection]{Lemma}
\newtheorem{conjecture}[subsubsection]{Conjecture}
\newtheorem{question}[subsubsection]{Question}
\newtheorem{construction}[subsubsection]{Construction}

\theoremstyle{remark}
\newtheorem{remark}[subsubsection]{Remark}
\newtheorem{example}[subsubsection]{Example}

\theoremstyle{plain}

\usepackage{tcolorbox}
\renewcommand{\baselinestretch}{1.22}
\renewcommand{\ge}{\geqslant}
\renewcommand{\geq}{\geqslant}
\renewcommand{\le}{\leqslant}
\renewcommand{\leq}{\leqslant}

\usepackage{xpatch}
\makeatletter
\xpatchcmd{\@tocline}
  {\hfil\hbox to\@pnumwidth{\@tocpagenum{#7}}\par}
  {\ifnum#1<0\hfill\else\dotfill\fi\hbox to\@pnumwidth{\@tocpagenum{#7}}\par}
  {}{}
\def\l@subsection{\@tocline{2}{0pt}{4pc}{6pc}{}}
\def\l@subsubsection{\@tocline{3}{0pt}{8pc}{8pc}{}}
\makeatother

\title{Equivariant elliptic cohomology via loop spaces}
\author{Emile Bouaziz \quad Adeel A. Khan}
\address{SIMIS, Shanghai \quad Academia Sinica, Taipei}

\newcommand{\Cat}{\operatorname{\mathsf{Cat}}}
\newcommand{\TenCat}{\operatorname{\mathsf{TenCat}}}

\newcommand{\ShvCat}{\operatorname{\mathsf{ShvCat}}}
\newcommand{\ShvTenCat}{\operatorname{\mathsf{ShvTenCat}}}

\newcommand{\CCC}{\TenCat}
\newcommand{\stCat}{\ShvCat}
\newcommand{\stCCC}{\ShvTenCat}

\newcommand{\CC}{\textsc{C}}
\newcommand{\DD}{\textsc{D}}
\newcommand{\stCC}{\textbf{\textsc{C}}}
\newcommand{\stDD}{\textbf{\textsc{D}}}

\newcommand{\Mod}{\textsc{Mod}}
\newcommand{\MOD}{\mathcal{M}\textit{od}}
\newcommand{\COMOD}{\mathcal{C}\textit{omod}}
\newcommand{\grMod}{\textsc{Mod}}

\DeclareMathOperator{\QCoh}{\textsc{QCoh}}
\newcommand{\stQCoh}{\textbf{\textsc{QCoh}}}
\DeclareMathOperator{\QCOH}{\mathcal{QC}\textit{oh}}
\DeclareMathOperator{\Perf}{\textsc{Perf}}
\DeclareMathOperator{\PERF}{\mathcal{P}\textit{erf}}

\DeclareMathOperator{\Coh}{\textsc{Coh}}
\DeclareMathOperator{\COH}{\mathcal{C}\textit{oh}}

\DeclareMathOperator{\IndPerf}{\textsc{IndPerf}}
\NewDocumentCommand{\INDPERF}{e{^}}{\operatorname{Ind\mathcal{P}\textit{erf}}\IfValueT{#1}{^{\,#1}}}
\newcommand{\stIndPerf}{\textbf{\textsc{IndPerf}}}

\newcommand{\Rep}{\textsc{Rep}}
\newcommand{\stRep}{\textbf{\textsc{Rep}}}

\newcommand{\TensFun}{\operatorname{Fun}^{\otimes}}
\newcommand{\TENSFUN}{\operatorname{\mathbf{Fun}}^{\otimes}}

\DeclareMathOperator{\Map}{Maps}
\DeclareMathOperator{\MAP}{\mathbf{Maps}}

\newcommand{\SpecTens}{\operatorname{Spec}^{\otimes}}
\newcommand{\Hcat}{\textup{\textsc{H}}}
\newcommand{\HE}{\Hcat_E}
\newcommand{\HA}{\Hcat_A}
\newcommand{\HAS}{\Hcat_{A/S}}

\newcommand{\stHcat}{\textup{\textbf{\textsc{H}}}}
\newcommand{\stHA}{\stHcat_A}
\newcommand{\stHAS}{\stHcat_{A/S}}

\DeclareMathOperator{\Hom}{Hom}
\DeclareMathOperator{\HOM}{\mathbf{Hom}}
\newcommand{\hhoms}{\underline{\operatorname{Hom}}}
\newcommand{\hends}{\underline{\operatorname{End}}}

\DeclareMathOperator{\CharOne}{\chi^{(1)}}

\newcommand{\R}{\mathcal{R}}
\newcommand{\RA}{\R_A}
\newcommand{\fra}{\mathfrak{a}}
\newcommand{\frav}{\mathfrak{a}^\vee}
\newcommand{\frM}{\mathfrak{M}}
\newcommand{\frL}{\mathfrak{L}}
\newcommand{\frX}{\mathfrak{X}}
\newcommand{\lam}{\Lambda}
\newcommand{\lamco}{\Lambda^\vee}

\newcommand{\Loop}{\operatorname{L}}

\newcommand{\OO}{\mathcal{O}}
\newcommand{\YY}{\mathcal{Y}}
\newcommand{\ZZ}{\mathcal{Z}}
\newcommand{\XX}{\mathcal{X}}
\DeclareMathOperator{\spc}{Spec}
\newcommand{\Specnc}{\operatorname{Spec}^{\mathrm{nc}}}
\newcommand{\modmod}{/\!\!/}

\newcommand{\GG}{\mathbf{G}}
\newcommand{\GGhat}{\widehat{\mathbf{G}}}
\newcommand{\GGa}{\mathbf{G}_{a}}
\newcommand{\GGm}{\mathbf{G}_{m}}
\newcommand{\GGaS}{\mathbf{G}_{a,S}}
\newcommand{\GGmS}{\mathbf{G}_{m,S}}
\newcommand{\GM}{\mathbf{G}_{m}}
\newcommand{\GL}{\mathbf{GL}}
\newcommand{\bZ}{\mathbf{Z}}
\newcommand{\bQ}{\mathbf{Q}}

\newcommand{\dAlgk}{\operatorname{CAlg}_{k}}
\newcommand{\CAlg}{\operatorname{CAlg}}
\newcommand{\CGrp}{\operatorname{CGrp}}
\newcommand{\GrpEoo}{\operatorname{Grp}^{E_\infty}}

\newcommand{\Ind}{\operatorname{Ind}}
\newcommand{\FM}{\Phi}
\newcommand{\Ren}{\Upsilon}

\newcommand{\St}{\mathbf{St}}
\newcommand{\Stk}{\mathbf{St}_{k}}
\newcommand{\Sch}{\mathbf{Sch}_{k}}
\newcommand{\Var}{\mathbf{Var}_{k}}

\newcommand{\eval}{\operatorname{ev}}
\newcommand{\TT}{\mathbf{T}}
\newcommand{\UU}{\mathcal{U}}
\newcommand{\cl}{{\operatorname{cl}}}
\newcommand{\id}{\mathrm{id}}
\newcommand{\aff}{\mathbf{A}^1}
\newcommand{\HH}{\operatorname{HH}}
\newcommand{\Hcoh}{\operatorname{H}}
\newcommand{\K}{\operatorname{K}}
\newcommand{\Kcoh}{\operatorname{\mathbf{K}}}
\newcommand{\LL}{\mathbf{L}}
\newcommand{\Nm}{\operatorname{Nm}}
\renewcommand{\ss}{\mathrm{ss}}
\newcommand{\pr}{\mathrm{pr}}
\newcommand{\ev}{\mathrm{ev}}
\newcommand{\rightrightrightarrows}
  {\mathrel{\substack{\textstyle\rightarrow\\[-0.6ex]
   \textstyle\rightarrow \\[-0.6ex]
   \textstyle\rightarrow}}}
\newcommand{\rightrightrightrightarrows}
  {\mathrel{\substack{\textstyle\rightarrow\\[-0.6ex]
   \textstyle\rightarrow \\[-0.6ex]
   \textstyle\rightarrow \\[-0.6ex]
   \textstyle\rightarrow}}}
\DeclareMathOperator{\Sym}{Sym}
\newcommand{\Thom}[1]{\langle #1 \rangle}

\newcommand{\eps}{\varepsilon}
\newcommand{\epsb}{\bar{\varepsilon}}
\newcommand{\varpib}{\bar{\varpi}}

\newcommand{\sfp}{\mathsf{p}}
\newcommand{\sfe}{\mathsf{e}}

\newcommand{\Mell}{\mathcal{M}_{\mathrm{ell}}}

\newcommand{\dashgr}{\text{-}\mathrm{gr}}
\newcommand{\Pic}{\operatorname{Pic}}
\newcommand{\PIC}{\operatorname{\mathbf{Pic}}}

\newcommand{\Bun}{\operatorname{Bun}}
\newcommand{\frBun}{\operatorname{\mathfrak{B}un}}
\newcommand{\uAut}{\underline{\operatorname{Aut}}}

\newcommand{\Shv}{\operatorname{Shv}}
\newcommand{\Shvper}{\operatorname{Shv}^{\mathrm{per}}}
\newcommand{\Rot}{\operatorname{Rot}}

\newcommand{\un}{\mathbf{1}}
\newcommand{\cF}{\mathcal{F}}
\newcommand{\cG}{\mathcal{G}}

\begin{document}

\markboth{\MakeUppercase{\shortauthors}}{\MakeUppercase{\shorttitle}}
\begingroup
\thispagestyle{plain}
\raggedbottom
\vspace*{16pt}
{\centering
{\fontsize{14}{20}\selectfont\bfseries EQUIVARIANT ELLIPTIC COHOMOLOGY\\VIA LOOP SPACES\par}
\vspace{25pt}
\begin{tabular}{@{}c@{\hspace{3.5em}}c@{}}
{\normalsize\scshape Emile Bouaziz} & {\normalsize\scshape Adeel A. Khan}\\[3pt]
{\small\itshape SIMIS, Shanghai} & {\small\itshape Academia Sinica, Taipei}
\end{tabular}\par
}
\vspace{34pt}
\noindent\hfill\begin{minipage}{0.90\textwidth}
\small\renewcommand{\baselinestretch}{1.16}\selectfont
\setlength{\parindent}{1em}
\noindent\textsc{Abstract.} We construct a cohomology theory $\Hcoh_A(-)$ canonically attached to any $1$-dimensional commutative group scheme $A$ over a connective $E_\infty$-ring.
The additive group in characteristic zero gives $2$-periodic de Rham cohomology, the multiplicative group gives topological periodic cyclic homology, and an elliptic curve gives elliptic cohomology.
The construction extends to the $G$-equivariant setting for an algebraic group $G$, yielding in the elliptic case a theory of equivariant elliptic cohomology on algebraic varieties with $G$-action that takes values in sheaves on the stack of $G$-bundles.
Over the complex numbers, and for $G$ connected reductive, this lifts Grojnowski's theory from the moduli \emph{space} of $G$-bundles to the stack.

The key player is the category of perfect complexes on $A$ with $0$-dimensional support, equipped with the convolution product induced by the group law of $A$.
It plays the role of a categorified circle and is used to define a loop space $\Loop_A^G(X)$, whose structure sheaf gives rise to $\Hcoh_A^G(X)$ after periodization.
Sheaves on $\Loop^G_A(X)$ equivariant under loop rotation provide a category of coefficients for $\Hcoh_A^G(-)$ in the same way that D-modules provide coefficients for de Rham cohomology.
At the point $X=*$, the theory also gives rise to an $E_\infty$-ring spectrum (or sheaf thereof) $\K_A$, which recovers complex K-theory for $A=\GGm$, and the sheaf $\OO^{\mathrm{top}}$ of topological modular forms for the universal elliptic curve.
\end{minipage}\hfill\mbox{}\par

\vfill
\clearpage
\endgroup

\tableofcontents

\section{Introduction}

We develop foundations for equivariant elliptic cohomology in algebraic geometry.
It is realized as the elliptic case of a general construction \emph{canonically} assigning to every $1$-dimensional commutative group scheme $A$ a cohomology theory $\HH_A^G(-)$ on algebraic varieties with an action of the algebraic group $G$.
We call it \emph{$G$-equivariant $A$-loop cohomology}.

\subsection{Overview}\label{ssec:intro/coh}

Fix a base commutative ring $k$, or more generally a connective $E_\infty$-ring spectrum, a $1$-dimensional commutative group scheme $A$ over $k$, and $G$ an affine fppf group scheme over $k$.\footnote{%
  We work in spectral algebraic geometry in the sense of Lurie \cite{SAG}.
  See Subsect.~\ref{ssec:intro/conv} for our conventions; for example, by \emph{scheme} we implicitly mean ``spectral scheme'', and commutative group schemes are taken in the strict sense.
}
Our main construction is a cohomology theory
\[
  \HH_A^G(-) : (\Sch^G)^{\mathrm{op}} \longrightarrow \Rot^G(A),
\]
on the category of schemes with $G$-action, taking values in the category of \emph{$G$-equivariant $A$-rotation modules} that we will introduce below (Subsect.~\ref{ssec:intro/shv}).
For the moment, let us point out that its structure allows us to extract two invariants of particular interest: its periodization $\Hcoh_A^G(-)$, and a cohomology theory $\K_A(-)$ on topological spaces.
For orientation, the table below gives the flavour of each theory in the additive, multiplicative and elliptic cases.

\medskip
\begin{center}
\small
\renewcommand{\arraystretch}{1.15}
\setlength{\tabcolsep}{6pt}
\begin{tabular}{cccc}
  \toprule
  $A$ & $\HH_A^G(X)$ & $\Hcoh_A^G(X)$ & $\K_A$ \\ \midrule
  $\GGa$ & Hodge-filtered de Rham & $2$-periodic de Rham & $k[\beta^{\pm 1}]$ \\
  $\GGm$ & Hochschild & periodic cyclic & $\mathrm{KU}$ \\
  $E$\,\footnotemark & elliptic Hodge cohomology & elliptic cohomology & $\mathcal{O}^{\mathrm{top}}$\,\footnotemark \\
  \bottomrule
\end{tabular}
\end{center}
\addtocounter{footnote}{-1}\footnotetext{for $E$ an elliptic curve, or a family of such}
\stepcounter{footnote}\footnotetext{for $E$ the universal elliptic curve over the moduli stack of elliptic curves}
\medskip

\subsubsection{Additive case, nonequivariantly.}\label{sssec:intro/coh/add}
Let $A = \GGa$ be the additive group over $k=\bQ$.
The category $\Rot(\GGa)$ of \emph{additive rotation modules} is compactly generated, with compact objects perfect complexes over $k$ with a $B\GGhat_a$-equivariant structure.
These are equivalently mixed complexes \cite{Kassel} with perfect underlying complex, since $\OO(B\GGhat_a)^\vee \simeq k[\epsilon]/(\epsilon^2)$, $|\epsilon| = 1$, and one can show $\Rot(\GGa)$ contains the category of mixed complexes as a localization.
We have
\[
  \HH_{\GGa}(X) \simeq \bigl(L\Omega^*_{X},\, \epsilon = d_{\mathrm{dR}}\bigr),
\]
i.e. the complex $L\Omega^*_X := \bigoplus_{n \ge 0} L\Omega^n_X[n]$ (where $L\Omega^n_X := \wedge^n_{\OO_X}\LL_X$), with $\epsilon$ acting by the de Rham differential.

Taking homotopy invariants for the $B\GGhat_a$-action gives a canonical functor $M \mapsto M^{hB\GGhat_a}$ from $\Rot(\GGa)$ to $k[u]$-modules, $|u| = -2$, which on mixed complexes sends $(K, d, \epsilon)$ to the negative cyclic complex $(K[u], d + u\epsilon)$.\footnote{%
  cf. \cite[Lem.~2.2]{HoyoisCirc}
}
Thus $\HH_{\GGa}(X)^{hB\GGhat_a}$ is the $u$-twisted derived de Rham complex of $X$,
\[
  \HH_{\GGa}(X)^{hB\GGhat_a} \simeq \bigl(L\Omega^*_{X}[u],\, u\, d_{\mathrm{dR}}\bigr),
\]
and $\Hcoh_{\GGa}(X)^{hB\GGhat_a}$ its $2$-periodization $(L\Omega^*_{X}[u, u^{-1}],\, u\, d_{\mathrm{dR}})$.
The scaling action of $\GGm$ on $\GGa$ puts $L\Omega^n_{X}$ in weight $n$ and $u$ in weight $-1$; after shearing, $\HH_{\GGa}(X)^{hB\GGhat_a}$ becomes the Rees construction of the Hodge filtration on the derived de Rham complex.
Thus $\HH_{\GGa}(X)$ and $\Hcoh_{\GGa}(X)$ refine Hodge-filtered and $2$-periodic de Rham cohomology, respectively.

\subsubsection{Multiplicative case, nonequivariantly.}\label{sssec:intro/coh/mult}
Take $A=\GGm$ over a connective $E_\infty$-ring spectrum $k$.
The category $\Rot(\GGm)$ of \emph{multiplicative rotation modules} is compactly generated, with compact objects perfect complexes over $k$ with an equivariant structure for $S^1 \simeq B\bZ$.
This contains $\operatorname{Fun}(BS^1, \MOD_k)$ as a localization, and for a qcqs scheme $X$, $\HH_{\GGm}(X)$ identifies with Hochschild homology and in particular recovers negative cyclic homology:\footnote{%
  These invariants are formed relative to the base ring $k$, which we leave implicit in the notation; e.g. $\operatorname{HH}(X)$ is often denoted in the literature by $\operatorname{HH}(X/k)$ or $\operatorname{HH}^k(X)$.
}
\[
  \HH_{\GGm}(X)
  \simeq \operatorname{HH}(X),
  \qquad
  \HH_{\GGm}(X)^{hS^1}
  \simeq \operatorname{HC}^-(X).
\]
For the \emph{periodization} $\Hcoh_{\GGm}(X)$, homotopy invariants recover the Tate construction, periodic cyclic homology:
\[
  \Hcoh_{\GGm}(X)^{hS^1}
  \simeq \operatorname{HH}(X)^{tS^1}
  \simeq \operatorname{HP}(X).
\]
Over $k=\mathbf{S}$, $\HH_{\GGm}(X)$ is thus topological Hochschild homology $\operatorname{THH}(X)$, and $\Hcoh_{\GGm}(X)$ is a refinement of topological periodic cyclic homology $\operatorname{TP}(X)$.
Over $k=\bQ$, we have $\bQ[S^1] \simeq \bQ[\epsilon]/(\epsilon^2)$, whence an equivalence $\Rot(\GGa) \simeq \Rot(\GGm)$, and the Hochschild--Kostant--Rosenberg isomorphism (see \cite[Cor.~3.2]{TVHKR}) identifies
\begin{equation}\label{eq:add=mult nonequiv}
  \HH_{\GGm}(X)
  \simeq \operatorname{HH}(X)
  \simeq \bigl(L\Omega^*_{X},\, d_{\mathrm{dR}}\bigr)
  \simeq \HH_{\GGa}(X)
\end{equation}
and similarly $\Hcoh_{\GGm}(X) \simeq \Hcoh_{\GGa}(X)$.

\subsubsection{Additive and multiplicative cases, equivariantly.}\label{sssec:intro/coh/equiv}
In $\Rot^G(\GGa)$, with $G$ an affine algebraic group over $k=\bQ$, the compact objects are perfect complexes on the adjoint quotient stack $\mathfrak{g}/G$, equivariant for the $B\GGhat_a$-action where $t \in \GGhat_a$ acts on $x \in \mathfrak{g}$ by the automorphism $\exp(tx) \in G$.
There is an induction functor $\operatorname{ind}_G^1 : \Rot^G(\GGa) \to \Rot(\GGa)$, given on compact objects by $B\GGhat_a$-equivariant pushforward along $\mathfrak{g}/G \to *$.

For a $G$-quasi-projective $X$, $\operatorname{ind}_G^1\HH_{\GGa}^G(X)$ is the derived de Rham complex of the quotient stack $X/G$ as a mixed complex:
\[
  \operatorname{ind}_G^1\HH^G_{\GGa}(X) \simeq \bigl(L\Omega^*_{X/G},\, \epsilon = d_{\mathrm{dR}}\bigr).
\]
In particular, $\HH^G_{\GGa}(X)$ and $\Hcoh_{\GGa}^G(X)$ are refinements of equivariant Hodge-filtered and $2$-periodic de Rham cohomology, respectively.

In the multiplicative case, with $G$ any linearly reductive group scheme over a connective $E_\infty$-ring $k$, compact objects of $\Rot^G(\GGm)$ are similarly perfect complexes on the adjoint quotient stack $G/G$, equivariant for the $S^1$-action where $n \in \bZ$ acts on $g \in G$ by the automorphism $g^n \in G$, and $\operatorname{ind}_G^1 : \Rot^G(\GGm) \to \Rot(\GGm)$ is again given on compact objects by $S^1$-equivariant pushforward to the point.
Then $\HH^G_{\GGm}(X)$ recovers $\operatorname{HH}(X/G)$ with its $S^1$-action:
\[
  \operatorname{ind}_G^1\HH^G_{\GGm}(X) \simeq \operatorname{HH}(X/G).
\]
In particular, $\HH^G_{\GGm}(X)$ and $\Hcoh^G_{\GGm}(X)$ recover $\operatorname{HC}^-(X/G)$ and $\operatorname{HP}(X/G)$, respectively, after applying $\operatorname{ind}_G^1$ and taking $S^1$-invariants.

Equivariantly, the HKR identification over $k=\bQ$ \eqref{eq:add=mult nonequiv} breaks down: the categories $\Rot^G(\GGa)$ and $\Rot^G(\GGm)$ are rarely equivalent, and $\HH_{\GGa}^G(X)$ and $\HH_{\GGm}^G(X)$ need not agree even after induction to $\Rot(\GGa) \simeq \Rot(\GGm)$.
This reflects the difference between the \emph{Borel}-type equivariance of de Rham cohomology (which satisfies descent along $X \twoheadrightarrow X/G$) versus the \emph{genuine} equivariance of Hochschild homology (which does not).

\subsubsection{Elliptic case.}\label{sssec:intro/coh/ell}
For $A = E$ an elliptic curve over $k$, the compactly generated category $\Rot(E)$ of \emph{elliptic rotation modules} has as compact objects the perfect complexes over $k$ with an equivariant structure for the dual elliptic curve $\check{E}$.
Over $k=\bQ$, there exists a noncanonical equivalence $\Rot(E) \simeq \Rot(\GGa)$, under which we have an HKR-type isomorphism
\[
  \HH_E(X) \simeq \HH_{\GGa}(X) \simeq \HH_{\GGm}(X)
\]
for every qcqs scheme $X$, and similarly for $\Hcoh_E(X)$.
Over a general $k$, $\HH_E(-)$ and $\Hcoh_E(-)$ are our definitions of \emph{elliptic Hodge cohomology} and \emph{elliptic cohomology}; these are thus cohomology theories on $k$-schemes, valued in elliptic rotation modules.

For $G$ a linearly reductive group scheme over $k$ (any connective $E_\infty$-ring), compact objects of $\Rot^G(E)$ are perfect complexes on $\Bun_G^\diamond(\check{E})$, a substack\footnote{%
  namely, the substack of $G$-bundles $F$ with $F \times^G V$ semistable of degree zero for every finite-dimensional representation $V$; see Subsect.~\ref{ssec:classifying stacks}
} of $\Bun_G(\check{E}) = \MAP(\check{E}, BG)$, equivariant for the action of $\check{E}$ by translation.
Our cohomology theories $\HH_E^G(-)$ and $\Hcoh_E^G(-)$, called $G$-equivariant \emph{elliptic Hodge cohomology} and \emph{elliptic cohomology}, respectively, take values in this category of $G$-equivariant elliptic rotation modules.

Over $k=\mathbf{C}$, let us explain how these constructions recover and generalize the approaches of Grojnowski \cite{Gr} and Ginzburg--Kapranov--Vasserot \cite{GKV}.
Let $\mathsf{G}$ be the maximal compact subgroup of $G^{\mathrm{an}}$.
The axiomatics of \cite{GKV} predict $\mathsf{G}$-equivariant elliptic cohomology sheaves $\operatorname{Ell}_{\mathsf{G}}^i(X^{\mathrm{an}})$ on the moduli space $\frM_E(G)^{\mathrm{GKV}}$ of semistable $G$-bundles on $\check{E}$ which are fibrewise of topologically trivial type.
The stack $\Bun_G^\diamond(\check{E})$ is classical and admits a good moduli space $\frM_E(G)$ (Theorem~\ref{thm:gms elliptic}), which contains $\frM_E(G)^{\mathrm{GKV}}$ as an open and closed subspace.\footnote{%
  If $G$ is connected with $\pi_1(G)$ torsion-free (e.g. $G=\GL_n$, $G=\mathrm{SL}_n$, or $G=T$ a torus), then $\frM_E(G)^{\mathrm{GKV}} = \frM_E(G)$.
  The same holds for $G$ finite.
}
Write $\mathsf{H}_E^G(X)$ for the $\check{E}$-invariants of $\Hcoh_E^G(X)$, pushed forward to $\frM_E(G)$, and consider the cohomology sheaves
\[ \mathcal{H}^i(\mathsf{H}_E^G(X)|_{\frM_E(G)^{\mathrm{GKV}}}). \]
The properties and structure we will establish for $\Hcoh_E^G(-)$ in \S\S \ref{sec:coh formal}--\ref{sec:pushforwards} will in particular imply that these sheaves satisfy most of the axioms stated in \cite{GKV}, such as periodicity, exactness, homotopy, K\"unneth, change of groups, induction, and existence of Euler classes and pushforwards (Gysin maps).

For $G = T$ a \emph{torus}, $\frM_E(T) \simeq \frM_E(T)^{\mathrm{GKV}} \simeq E \otimes_{\bZ} \mathbf{X}_*(T)$, and for $X$ proper, the $\operatorname{Ell}_{\mathsf{T}}^i(X^\mathrm{an})$ were constructed on this space by Grojnowski \cite{Gr}.
An algebraic construction was given in the inspiring work of Sibilla--Tomasini \cite{ST}.\footnote{%
  Working in characteristic zero, they construct quasi-coherent sheaves $HH_E([X/T])$ and $HP_E([X/T])$ on $\frM_E(T)$ which are pushed forward from $\Bun_T^0(\check{E}) \simeq \Bun_T^\diamond(\check{E})$, the stack of degree zero bundles.
  For example, their $HP_E([X/T])$ can be recovered from our $\Hcoh_E^T(X)$ by taking $\check{E}$-invariants and pushing forward to $\frM_E(T)$.
}

For $G$ \emph{connected}, with maximal torus $T$ and Weyl group $W$, Grojnowski \cite[\S 2.6]{Gr} and Ginzburg--Kapranov--Vasserot \cite[(1.11)]{GKV} then take $W$-invariants of the $T$-equivariant theory to produce a theory on the space $\frM_E(T) \modmod W \simeq (E \otimes_{\bZ} \mathbf{X}_*(T)) \modmod W$, which is the component of the trivial bundle in $\frM_E(G)$ (Theorem~\ref{thm:gms elliptic}(ii)).

For $G$ \emph{finite}, our $\mathsf{H}_E^G(X)$ is defined on $\frM_E(G) \simeq \frM_E(G)^{\mathrm{GKV}}$, the set of conjugacy classes of commuting pairs in $G$ (whereas $\frM_E(T) \modmod W$ is in this case a single point).

\subsubsection{The cohomology theories $\K_A$.}

The objects $\K_A$ are $E_\infty$-ring spectra, which we think of as representing cohomology theories on topological spaces.
In a word, $\K_A$ is entirely determined by the value $\HH_A(*)$ of the nonequivariant theory at the point $*=\spc(k)$; like $\Hcoh_A(-)$ it is defined by periodization, but via \emph{co}invariants.
We will see that it is essentially equivalent (before incorporating equivariance) to a construction developed by Lurie in \cite[\S 4.3]{LurieEll2}.

For $A=\GGa$ over $k=\bQ$, $\HH_{\GGa}(*) \simeq k$ is the trivial mixed complex, and its $B\GGhat_a$-coinvariants form the polynomial algebra $k[\beta]$, $|\beta|=2$, linearly dual to $\HH_{\GGa}(*)^{hB\GGhat_a} \simeq k[u]$.
Inverting $\beta$ gives $\K_{\GGa} \simeq k[\beta^{\pm 1}]$, which represents $2$-periodic rational cohomology.

For $A=\GGm$ over $k=\mathbf{S}$, $\HH_{\GGm}(*) \simeq \operatorname{THH}(\mathbf{S}) \simeq \mathbf{S}$ with trivial $S^1$-action, and its $S^1$-coinvariants are $\mathbf{S}_{hS^1} \simeq \Sigma^\infty_+\mathbf{C}\mathbf{P}^\infty$.
We have $\K_{\GGm} \simeq \Sigma^\infty_+\mathbf{C}\mathbf{P}^\infty[\beta^{-1}] \simeq \mathrm{KU}$ (Snaith's theorem \cite{Snaith}), the spectrum representing complex K-theory, where $\beta \in \pi_2$ is the class of the bottom cell $\mathbf{C}\mathbf{P}^1 \subset \mathbf{C}\mathbf{P}^\infty$.

For $A=E$ an elliptic curve, $\HH_E(*) \simeq k$ with trivial $\check{E}$-action, and its $\check{E}$-coinvariants are the functions $\OO(\Omega^2 E)$ on the double loop space of $E$.
Periodizing with respect to an analogue of the Bott class gives rise to $\K_E$, which we will see is the algebra corepresenting orientations of $E$ in the sense of \cite[Def.~7.2.7]{LurieEll2} (see \cite[\S 4.3, \S 7.2]{LurieEll2}).
The construction works in families, and for the universal elliptic curve $\mathcal{E}$ over Lurie's moduli stack of strict elliptic curves $\Mell$, $\K_{\mathcal{E}/\Mell}$ is a sheaf of $E_\infty$-ring spectra which, regarded as a sheaf on the étale site of the classical moduli of elliptic curves $\Mell^{\heartsuit}$, is the Goerss--Hopkins--Miller sheaf $\OO^{\mathrm{top}}$, with global sections $\mathrm{TMF}$ (Theorem~\ref{thm:Otop}).

Equivariantly, for $G$ a diagonalizable group scheme over $k$, $\HH_A^G(*) \in \Rot^G(A)$ similarly determines an equivariant refinement $\K_A^G$ of $\K_A$, by formation of Bott-periodized coinvariants ``in families''.
As we will see in Subsect.~\ref{ssec:HH on XAG}, the work of Lurie and Gepner--Meier \cite{LurieEllSurv,GepMei} shows that the $\K_A^G$ assemble into a cohomology theory $\Kcoh_A^G$ on $\mathsf{G}$-spaces, for $\mathsf{G}$ the compact abelian Lie group with character group $X^*(G)$, whose value on $\mathsf{G}/\mathsf{H}$ is (the pushforward of) $\K_A^H$.
In the elliptic case, $\Kcoh_E^G$ is precisely $\mathsf{G}$-equivariant elliptic cohomology in the sense of Lurie and Gepner--Meier \cite{LurieEllSurv,GepMei}.

\subsection{Loop spaces via categorified circles}\label{ssec:intro/circles}

Let $G$ be an affine fppf group scheme acting linearly on a scheme $X$ such that the quotient stack $\YY=X/G$ is \emph{perfect} in the sense of \cite{BFN}.
We now explain the main construction of this paper, the \emph{$A$-loop space} of $\YY$.

\subsubsection{Loop spaces and cohomology theories}

Our starting point is the observation that Hochschild homology $\operatorname{HH}(\YY)$ is encoded by the \emph{free loop space} $\MAP(S^1, \YY) \simeq \YY \times_{\YY \times \YY} \YY$ as the ring of functions,
\[
  \OO(\MAP(S^1, \YY)) \simeq
  \Gamma(\YY, \OO_\YY \otimes_{\OO_{\YY \times \YY}} \OO_\YY) \simeq
  \operatorname{HH}(\YY).
\]
Moreover, the $S^1$-equivariant structure on $\operatorname{HH}(\YY)$, used to produce negative and periodic cyclic homology, arises from the $S^1$-action by loop rotation.
See \cite{BNloop,TVChern}.

In the additive case, the \emph{$(-1)$-shifted tangent bundle} $\TT[-1]\YY$ plays an analogous role: its ring of functions is the derived Dolbeault complex,
\[ \OO(\TT[-1]\YY) \simeq \Gamma\bigl(\YY, \Sym_{\OO_\YY}(\LL_\YY[1])\bigr) \simeq \bigoplus_{n} L\Omega^n_\YY[n], \]
i.e. the derived de Rham complex, minus differential.
Moreover, $\TT[-1]\YY$ can also be written as the mapping stack $\MAP(B\GGhat_a, \YY)$, and the translation action of $B\GGhat_a$ encodes the de Rham differential, while the fibrewise scaling action of $\GGm$ on $\TT[-1]\YY$ encodes the Hodge filtration.
See \cite[Props.~4.4, 4.7, \S1.6.3]{BNloop}, Theorem~\ref{thm:LGa LGm}, and Corollary~\ref{cor:Ga comparison sheared}.

This leads us to the following formulation of the problem:

\begin{question}\label{que:intro/S^1}
	For any $1$-dimensional commutative group scheme $A$ over $k$, does there exist a canonical commutative group stack $S^1_A$ satisfying the following properties:
	\begin{enumerate}
		\item the structure sheaf of the mapping stack $\MAP(S^1_A, X/G)$, with its loop rotation action, encodes a well-behaved cohomology theory\footnote{%
		  Here $(-)^{tS^1_A}$ denotes a suitable Tate construction, which is also to be defined.
			In characteristic zero, we can ask for an algebra isomorphism $\OO(S^1_A)^\vee \simeq k[\epsilon]/(\epsilon^2)$, $|\epsilon| = 1$, so that $\OO(\MAP(S^1_A, X/G))$ can be regarded as a mixed complex and the usual Tate construction can be applied.
		} \[ X \mapsto \OO(\MAP(S^1_A, X/G))^{tS^1_A} \] on $k$-varieties with $G$-action (satisfying e.g. Mayer--Vietoris, K\"unneth, \dots);
		\item for $A=\GGa$ (over $k=\bQ$) and $A=\GGm$, $S^1_A$ recovers $B\GGhat_a$ and $B\bZ \simeq S^1$, respectively; and,
		\item for $A=E$ the resulting cohomology recovers Grojnowski's torus-equivariant elliptic cohomology \cite{Gr}.
	\end{enumerate}
\end{question}

\subsubsection{Categorified circles}

To address Question~\ref{que:intro/S^1}, let us point out that in the additive and multiplicative cases, the relationship between the commutative group stacks $A$ and $S^1_A$ becomes much more transparent after categorification.

To relate $\GGm$ and $B\bZ$, for example, consider the line bundle $L_{\GGm}$ on $B\bZ \times \GGm$ classified by the tautological character $\bZ \times \GGm \to \GGm$, $(n,t) \mapsto t^n$.
The Fourier--Mukai transform
\[
  \Phi_{\GGm} : \QCOH(\GGm) \to \QCOH(B\bZ),
  \quad F \mapsto \pr_{1,*}\bigl(\pr_2^*(F) \otimes L_{\GGm}\bigr),
\]
is an equivalence, and exchanges the tensor product $\otimes$ on the right with the convolution product
\[ F \star G := m_*(F \boxtimes G), \qquad m : \GGm \times \GGm \to \GGm \]
on the left.
On connective perfect complexes, it restricts to a symmetric monoidal equivalence
\[ \Phi_{\GGm} : (\Perf_0(\GGm), \star) \xrightarrow{\sim} (\Perf(B\bZ), \otimes) \]
where $\Perf_0(\GGm)$ is the category of connective perfect complexes with $0$-dimensional support (see Theorem~\ref{thm:HA vs QCoh checkA}).

For $\GGa$ in characteristic zero, the same holds with $L_{\GGa}$ the line bundle on $B\GGhat_a \times \GGa$ given by the character $\GGhat_a \times \GGa \to \GGm$, $(x,t) \mapsto \exp(xt)$: Fourier--Mukai restricts to an equivalence
\[ \Phi_{\GGa} : (\Perf_0(\GGa), \star) \xrightarrow{\sim} (\Perf(B\GGhat_a), \otimes), \]
where convolution is defined using the addition $m : \GGa \times \GGa \to \GGa$.

The category $\Perf_0(A)$, with the convolution product $F \star G := m_*(F \boxtimes G)$ categorifying the group law $m : A \times A \to A$, makes sense for every $A$.
This discussion suggests that the stack $S^1_A$ should have the property that $(\Perf(S^1_A), \otimes)$ is symmetric monoidally equivalent to $(\Perf_0(A), \star)$ for every $A$.
But for $A = E$ an elliptic curve, the classical Fourier--Mukai transform \cite{Mukai}, with kernel the Poincar\'e bundle on $\check{E} \times E$, restricts to a symmetric monoidal equivalence
\[ \Phi_{E} : (\Perf_0(E), \star) \xrightarrow{\sim} (\Perf^{\mathrm{hmg}}(\check{E}), \otimes) \]
where the right-hand side is the category of \emph{homogeneous} perfect complexes on the dual $\check{E}$ (see e.g. Theorem~\ref{thm:HA vs QCoh checkA}).
The requirement is thus
\[ (\Perf(S^1_E), \otimes) \simeq (\Perf^{\mathrm{hmg}}(\check{E}), \otimes), \]
which we do not know how to achieve.

\subsubsection{Mapping stacks via Tannaka duality}

In order to define the mapping stack $\MAP(S^1_A, \YY)$, it turns out to be sufficient to only have the data of $(\Perf(S^1_A), \otimes)$.
Indeed, since $\YY$ is perfect by assumption, Tannaka duality shows that a morphism of stacks $f : S^1_A \to \YY$ is the same data as that of a symmetric monoidal functor $\Perf(\YY) \to \Perf(S^1_A)$ (see \cite[Thm.~1.3]{BHL}).\footnote{%
  We write $\QCoh(\YY)$ and $\Perf(\YY)$ for the full subcategories of \emph{connective} objects of the stable $\infty$-categories $\QCOH(\YY)$, $\PERF(\YY)$.
  First, \cite[Thm.~1.3]{BHL} identifies $\Map(S^1_A, \YY)$ with symmetric monoidal colimit-preserving functors $\QCoh(\YY) \to \QCoh(S^1_A)$.
  Since $\YY$ is perfect, $\QCoh(\YY)$ is the ind-completion of $\Perf(\YY)$, hence this is equivalently the data of a symmetric monoidal functor $\Perf(\YY) \to \QCoh(S^1_A)$.
  Since symmetric monoidal functors preserve dualizable objects, any such lands in $\Perf(S^1_A)$.
}

For every $A$, taking $(\Perf(S^1_A), \otimes)$ to be the symmetric monoidal category $(\Perf_0(A), \star)$ by definition, we can thus define:

\begin{definition}
  The \emph{$A$-loop space} $\Loop_A(\YY)$ is the stack of symmetric monoidal functors $\Perf(\YY) \to \Perf(S^1_A)$.
  More precisely, $T$-points of $\Loop_A(\YY)$ are symmetric monoidal $\Perf(T)$-linear functors
  \[ \Perf(\YY_T) \to \Perf_0(A_T/T), \]
  where $\Perf_0(A_T/T)$ is the full subcategory of $\Perf(A_T)$ spanned by objects with finite support over $T$.
\end{definition}

We will see that $\Loop_A(\YY)$ is a substack of $\MAP(\check{A}, \YY)$ (\ref{sssec:loops as maps}), where $\check{A}$ is the \emph{$1$-shifted Cartier dual} (see \ref{sssec:HA trichotomy}); $\check{A}$ is $B\GGhat_a$, $B\bZ$ and the dual elliptic curve $\check{E}$ in the respective cases $A=\GGa$, $A=\GGm$, and $A=E$.
In particular, $\Loop_A(\YY)$ is equipped with a loop rotation action by $\check{A}$.

For $A=\GGa$ (in characteristic zero) and $\GGm$, $\Loop_A(\YY)$ recovers $\TT[-1]\YY$ and $\MAP(S^1, \YY)$, respectively.
The \emph{elliptic} loop space $\Loop_E(\YY)$ recovers a construction of Sibilla--Tomasini \cite{ST} in the torus-equivariant case $\YY=X/T$ (in characteristic zero), and has also been considered independently in the thesis of Devalapurkar \cite{Sanath}.

\subsubsection{Example: schemes}\label{sssec:intro/coh/HKR}

When $\YY=X$ is a \emph{scheme}, the $A$-loop space only depends on the formal group $\widehat{A}$ (see Theorem~\ref{thm:loop schemes}):
\begin{equation*}
  \Loop_A(X) \simeq \MAP(B(\mathbb{D}\widehat{A}), X),
\end{equation*}
where $\mathbb{D}\widehat{A}$ is the Cartier dual of $\widehat{A}$.
This description matches the one proposed by Moulinos--Robalo--To\"en~\cite[\S 6.3]{MRT} for $X$ affine (and $k$ discrete); see also \cite[Constr.~11.3]{LuriePrism}.

In characteristic zero, the formal logarithm identifies $\GGhat_a \simeq \GGhat_m \simeq \widehat{E}$, so the three loop spaces coincide:
\[ \Loop_{\GGa}(X) \simeq \Loop_{\GGm}(X) \simeq \Loop_E(X). \]

\subsubsection{Example: classifying stacks}\label{sssec:intro/coh/BG}
Equivariantly, $A$-loop spaces are not determined by just the formal group $\widehat{A}$.
For example, in characteristic zero, any symmetric monoidal functor $\Perf(B\GL_r) \to \Perf_0(A)$ is determined by the image of the standard representation (by Iwanari \cite{Iwa}), which is a sheaf of length $r$ on $A$ with $0$-dimensional support, so the $A$-loop space of $B\GL_r$ is
\[ \Loop_A(B\GL_r) \simeq \operatorname{Coh}_0^r(A), \]
the stack of $0$-dimensionally supported sheaves of length $r$ on $A$ (see Proposition~\ref{prop:LA(BGLr) char 0}).
For $r=1$ these are the skyscraper sheaves $\OO_a$, $a \in A$, whence $\Loop_A(B\GGm) \simeq A \times B\GGm$.
More generally, $\Loop_{\GGa}(BG) \simeq \mathfrak{g}/G$ for every affine algebraic group $G$ (Example~\ref{exam:equivariant loops Ga Gm}).

Over a general connective $E_\infty$-ring $k$, with $G$ linearly reductive, we have:
\begin{itemize}
\item For $A=\GGm$, $\Loop_{\GGm}(BG)$ is the adjoint quotient stack $G/G$.
\item For $A=E$ an elliptic curve, $\Loop_E(BG)$ is the stack $\Bun_G^\diamond(\check{E})$ (\ref{sssec:intro/coh/ell}); for $G=\GL_r$ in characteristic zero this agrees with the $\operatorname{Coh}_0^r(E)$ description above, via Fourier--Mukai (Remark~\ref{rem:HA vs QCoh checkA dzbl}).
\end{itemize}
We note that $\Loop_A(BG)$ is also well-defined and of interest without reductivity hypotheses.
For instance, for $B \subset \GL_r$ a Borel in characteristic zero, the map $\Loop_E(BB) \to \Loop_E(B\GL_r)$ recovers the elliptic Grothendieck--Springer resolution of Ben-Zvi--Nadler \cite[\S 2]{BNEll} (\ref{springer}).

\subsection{Loop cohomology}\label{ssec:intro/shv}

With $A$-loop spaces in hand, we now proceed to the construction of $A$-loop cohomology.
Let $G$ be an affine fppf group scheme acting on a scheme $X$.
At first approximation, the $A$-loop cohomology $\HH_A^G(X)$ is the structure sheaf of $\Loop_A(X/G)$, pushed forward along
\[ \pi : \Loop_A(X/G) \to \Loop_A(BG), \]
together with its $\check{A}$-equivariant structure.
The latter will be used to define its \emph{periodization}, $\Hcoh_A^G(X)$.

\subsubsection{Rotation modules}\label{sssec:intro/shv/def}

To make this precise, we need to be careful about the definition of $\check{A}$-equivariant structure on a quasi-coherent sheaf.
Because the classifying stack $*/\check{A} = B\check{A}$ is not perfect, the canonical functor
\[ \Ren : \INDPERF^{\check{A}}(*) = \INDPERF(B\check{A}) \to \QCOH(B\check{A}) = \QCOH^{\check{A}}(*) \]
is not an equivalence.
We define:

\begin{definition}
	The category of \emph{$G$-equivariant $A$-rotation modules} is
	\[ \Rot^G(A) := \INDPERF^{\check{A}}(\Loop_A(BG)). \]
	The \emph{$G$-equivariant $A$-loop cohomology} of $X$ is the $G$-equivariant $A$-rotation module
	\[ \HH_A^G(X) := \pi_*(\OO_{\Loop_A(X/G)}) \in \Rot^G(A). \]
\end{definition}

\subsubsection{Periodization}

Let us explain the effect of this ``renormalization'' in the nonequivariant case, i.e. for the category
\[ \Rot(A) := \INDPERF^{\check{A}}(\Loop_A(*)) \simeq \INDPERF(B\check{A}). \]
First, consider instead the category $\QCOH^{\check{A}}(*) \simeq \QCOH(B\check{A})$.
These are equivalent to modules over the \emph{rotation algebra} $\RA := \OO(\Omega A)$, the algebra of functions on the loop group $\Omega A := e \times_A e$ of $A$ at its identity section (\ref{sssec:RA}).
This is in fact a Hopf algebra, with coalgebra structure induced by the group structure of $A$, and there is a symmetric monoidal equivalence
\[ (\QCOH(B\check{A}), \otimes) \simeq (\MOD_{\RA}, \star), \]
where $\star$ is the convolution product defined by the comultiplication of $\RA$ (see Theorem~\ref{thm:QCoh checkA}).
In characteristic zero, we have $\RA \simeq k[\epsilon]/(\epsilon^2)$, $|\epsilon| = 1$, as algebras, so that $\RA$-modules are equivalent to mixed complexes.\footnote{%
  In general, $\RA$ is not a square-zero extension; for example, over the sphere spectrum, $\R_{\GGm} \simeq \mathbf{S}[S^1]$ is the spherical group ring (\ref{sssec:OO(Gm[-1])}, Remark~\ref{rem:Hopf map}).
}

The category $\QCOH(B\check{A})$ admits a canonical autoequivalence $(-)\Thom{1}$, an $A$-twisted form of the triangulated shift $F \mapsto F[2]$, together with the \emph{Euler map} $\alpha_A : F \to F\Thom{1}$, natural in $F \in \QCOH(B\check{A})$ (Subsect.~\ref{ssec:Euler}).
The \emph{periodization} of $F$ is the colimit
\[
  F \xrightarrow{\alpha_A}
  F\Thom{1} \xrightarrow{\alpha_A}
  F\Thom{2} \xrightarrow{\alpha_A}
  \cdots
  \to F[\alpha_A^{-1}].
\]
As written this colimit is always zero: the Euler map on $\OO_{B\check{A}}$ corresponds under $\QCOH(B\check{A}) \simeq \MOD_{\RA}$ to the boundary map $\partial$ of the augmentation triangle $\mathcal{I}_\eps \to \RA \xrightarrow{\eps} k$; this triangle is split on underlying $k$-modules, so on underlying $k$-modules the maps in the telescope above are all null-homotopic, hence the underlying $k$-module of $F[\alpha_A^{-1}] \simeq F \otimes \OO_{B\check{A}}[\alpha_A^{-1}]$ is zero.

Let us now renormalize: the above constructions restrict to $\PERF(B\check{A})$ and ind-extend to
\[ \INDPERF(B\check{A}) \simeq \Ind(\MOD_{\RA}^{\mathrm{perf}/k}), \]
the ind-completion of the category of $\RA$-modules that are perfect over $k$ (see Corollary~\ref{cor:IndPerf checkA}).
For $F \in \PERF(B\check{A})$, the colimit defining $F[\alpha_A^{-1}]$ is now formed in the formal sense, i.e. as an ind-object, and $\OO_{B\check{A}}[\alpha_A^{-1}]$ is no longer zero; in fact, the canonical functor $\Ren : \INDPERF(B\check{A}) \to \QCOH(B\check{A})$ annihilates precisely the Euler-periodic objects, i.e. those on which $\alpha_A$ acts invertibly.
For example, forming homotopy $\check{A}$-invariants gives the Tate construction,
\[ (F[\alpha_A^{-1}])^{h\check{A}} \simeq F^{t\check{A}}, \]
see Theorem~\ref{thm:Tate vs invariants of Bott}.
Here, the Tate construction is the cofibre of a norm map from coinvariants to invariants, defined in Subsect.~\ref{ssec:Tate}; for $A=\GGm$ it recovers (the renormalization of) the usual $S^1$-Tate construction (Example~\ref{ex:Tate multiplicative}).

The symmetric monoidal category $\Rot(A)$ acts on $\Rot^G(A)$.
We write $\Rot^{G,\mathrm{per}}(A)$ for the full subcategory of periodic objects, and define:

\begin{definition}
	The \emph{periodic $G$-equivariant $A$-loop cohomology} of $X$ is the rotation module
	\[ \Hcoh_A^G(X) := \HH_A^G(X)[\alpha_A^{-1}] \ \in \Rot^{G,\mathrm{per}}(A). \]
\end{definition}

Nonequivariantly, we can form $\check{A}$-invariants for the rotation module $\Hcoh_A(X)$ (since $\check{A}$ acts trivially on $\Loop_A(*) \simeq *$) and this recovers the Tate construction,
\[ \Hcoh_A(X)^{h\check{A}} \simeq \HH_A(X)^{t\check{A}}. \]
For the value at the point $X=*$, there is another ``dual'' type of periodization.
Upon taking $\check{A}$-\emph{co}invariants, the Euler map $\alpha_A$ splits and gives rise to a \emph{Bott class}
\[ \beta_A : \frav[2] \to (\OO_{B\check{A}}\Thom{1})_{h\check{A}} \to (\OO_{B\check{A}})_{h\check{A}} \simeq \OO(\Omega^2 A), \]
see Subsect.~\ref{ssec:Bott periodization}.
This turns out to be closely related to a construction of Lurie (see \cite[Constr.~4.3.7]{LurieEll2}, Lemma~\ref{lem:Bott vs Lurie}).
The $E_\infty$-$k$-algebra $\K_A$ is defined by Bott-periodizing the coinvariants\footnote{%
  Note that taking coinvariants of the Euler periodization directly, $(\OO_{B\check{A}}[\alpha_A^{-1}])_{h\check{A}}$, gives zero (Lemma~\ref{lem:coinv kill periodization}).
  This is why we need to periodize with respect to the Bott class instead.
}, i.e. by inverting $\beta_A$ in $\lam_A := \OO(\Omega^2 A)$:
\[\K_A := \lam_A[\beta_A^{-1}].\]
As discussed in Subsect.~\ref{ssec:intro/coh}, this recovers the complex K-theory spectrum for $A=\GGm$ over the sphere spectrum, and, for the universal elliptic curve over $\Mell$, recovers the sheaf of $E_\infty$-ring spectra $\OO^{\mathrm{top}}$.

Equivariantly, for $G$ diagonalizable, a generalization of this construction gives rise to the $\K_A^G$, which for $A=E$ assemble into $\mathsf{G}$-equivariant elliptic cohomology in the sense of Lurie and Gepner--Meier (see Subsect.~\ref{ssec:HH on XAG}).

\subsubsection{Loop sheaves}

A pleasant feature of $A$-loop cohomology is that it comes naturally equipped with a category of coefficients.
We define the category of \emph{$A$-loop sheaves} by
\[ \Shv_A^G(X) := \INDPERF^{\check{A}}(\Loop_A(X/G)), \]
see Subsect.~\ref{ssec:sheaves}, and denote by $\Shv_{A,\mathrm{c}}^{G}(X)$ the full subcategory of perfect complexes.
The category of $G$-equivariant rotation modules is then $\Rot^G(A) = \Shv_A^G(*)$.

Via the functorialities of $\Loop_A(-)$ and $\INDPERF(-)$, $\Shv_A^G(-)$ admits pullback and pushforward functors along $G$-equivariant morphisms, as well as induction and restriction functors for any group homomorphism $H \to G$.
Writing $\un$ for the unit of the monoidal structure, the definition of $A$-loop cohomology can be recast as
\[ \HH_A^G(X) \simeq a_* a^*(\un) \in \Rot^G(A) \]
for $a : X \to *$ the projection, and more generally we write
\[ \HH_A^G(X; \cF) := a_* (\cF) \in \Rot^G(A), \qquad \cF \in \Shv_A^G(X). \]
Similarly, periodic loop cohomology is calculated by the same constructions inside $\Shv_A^{G,\mathrm{per}}(X)$, the full subcategory of periodic objects.

For $A = \GGa$ in characteristic zero, $\Shv_{\GGa}^G(X)$ is closely related to the category of $D$-modules on $X/G$.
For example, taking into account the loop rescaling action of $\GGm$, Ben-Zvi and Nadler show that for $G=1$ and $X$ a smooth scheme, $\PERF^{\,\check{\GG}_a \rtimes \GGm}(\Loop_{\GGa}(X))$ identifies with a full subcategory of filtered $D$-modules on $X$; see \cite[Cor.~5.2]{BNloop}, as well as Subsect.~\ref{ssec:dR comparison} and \cite{Preygel,Chen}.
The periodic subcategory $\Shvper_{\GGa}(X)$ is closely related to $2$-periodic $D$-modules on $X$; see \cite[Cor.~1.15]{BNloop}.

Likewise, for $A = \GGm$ in characteristic zero, $\Shv_{\GGm}^{G,\mathrm{per}}(X)$ is closely related to the category of \emph{potent} $D$-modules on $X/G$ as introduced in \cite{BZNpotent}.

We expect that, in characteristic zero, the sheaf theory $\Shv_A^{G,\mathrm{per}}(-)$ is ``topological'':

\begin{conjecture}\label{conj:weave}
  Let $k$ be a $\bQ$-algebra.
  Then the assignment $X \mapsto \Shv_A^{G,\mathrm{per}}(X)$ defines a topological weave in the sense of \cite[\S 3]{KhanWeaves} on the category of $G$-quasi-projective schemes.
\end{conjecture}

While its construction is by nature ``coherent'' rather than topological, evidence for Conjecture~\ref{conj:weave} is that $\Shv_A^{G,\mathrm{per}}(-)$ satisfies homotopy invariance in the sense of \cite[Def.~3.14]{KhanWeaves} (see Theorem~\ref{thm:homotopy invariance shv}).
Moreover, recent work of Erlikhman--Schefers \cite{ES} establishes the localization property in the sense of \cite[Def.~3.4]{KhanWeaves} (a.k.a. the Kashiwara lemma) for potent D-modules.

\subsection{Related work}\label{ssec:intro/related}

The relationship between formal groups and cohomology theories has long been central to stable homotopy theory.
Every complex-oriented cohomology theory determines a formal group law over its coefficients, while the reverse passage from a formal group to a cohomology theory is much more delicate.
The Landweber exact functor theorem produces a homology theory from a formal group law satisfying a flatness condition, but not the $E_\infty$-ring structure of a multiplicative cohomology theory.
In the elliptic case, the modern perspective is due to Lurie, and its extension to an equivariant theory for finite groups and compact abelian Lie groups was carried out by Lurie and Gepner--Meier \cite{LurieEllSurv,LurieEll1,LurieEll2,LurieEll3,GepMei}.

In algebraic geometry, a general procedure to attach cohomology theories to formal groups is as follows.
For a formal group $\widehat{A}$ over a commutative ring $k$, Moulinos--Robalo--To\"en \cite{MRT} define a Hochschild-type theory $\operatorname{HH}_{\widehat{A}/k}(-)$ on $k$-schemes via functions on the stack of maps from the classifying stack of the Cartier dual $\mathbb{D}\widehat{A}$,
\[
  X \mapsto \operatorname{HH}_{\widehat{A}/k}(X) := \OO\bigl( \MAP(B(\mathbb{D}\widehat{A}), X) \bigr).
\]
When $A$ is an algebraization of $\widehat{A}$, our loop space $\Loop_{A/k}(X)$ collapses precisely to this mapping stack (see Theorem~\ref{thm:loop schemes}), so that $\operatorname{HH}_{\widehat{A}/k}(X)$ is our loop cohomology $\HH_{A/k}(X)$ over the discrete base $k$.\footnote{%
  To be precise, \cite{MRT} work in the setting of derived (rather than spectral) algebraic geometry.
}

One can thus think of loop cohomology as providing in particular a refinement $\HH_{A/\mathbf{S}}(-)$ of the above theory from discrete bases $k$ to the sphere spectrum (when $A$ is defined over $\mathbf{S}$).
A different spectral generalization of \cite{MRT} is proposed by Lurie in \cite[\S 11]{LuriePrism}: he predicts the existence of cohomology theories on $p$-adic formal schemes, which we denote here by $\widehat{\operatorname{THH}}^{\Delta}_{\widehat{A}}(-)$, which for $\widehat{A} = \widehat{\GG}_m$ recover $p$-complete THH.
One expects that these $p$-complete theories should glue with our rational $\HH_{A/\mathbf{Q}}(-)$ into an integral theory $\operatorname{THH}^{\Delta}_{\widehat{A}}(-)$; however, $p$-adically and outside the multiplicative case we expect this ``prismatic'' THH to differ from our ``vanilla'' loop cohomology in general.\footnote{%
  We thank Jacob Lurie for helpful explanations about this.
}

For the construction of \emph{equivariant} cohomology theories for algebraic varieties, new ideas are required.
Unlike in the nonequivariant case, the loop space $\Loop_A^G(X) = \Loop_A(X/G)$ depends on the entire group $A$, through the categorified circle $(\Perf(S^1_A),\otimes) = (\Perf_0(A), \star)$, and not only on its formal completion $\widehat{A}$.
For $A=E$ an elliptic curve over a field $k$ of characteristic zero, Sibilla--Tomasini \cite{ST} construct a torus-equivariant refinement of $\HH_{E/k}(-)$ using the substack of \emph{quasi-constant} maps in $\MAP(\check{E}, X/T)$ \cite[Def.~2.18]{ST}; in this situation, quasi-constancy coincides with the condition of being a \emph{loop} in the sense of (\ref{sssec:loops as maps}).
A primary motivation of our work was to extend this construction to general algebraic groups $G$, arbitrary $A$, and arbitrary characteristic bases.

Extending the theory outside characteristic zero raises several difficulties, the most substantial of which is finding the appropriate theory of periodization (\S\ref{sec:Euler}).
In characteristic zero, formality of the rotation algebra $\RA$ implies that our periodization reduces to Preygel's categorical Tate construction \cite[Def.~1.2.3]{Preygel} (see also \cite[Def.~2.7]{BNloop}, \cite[2.1.4]{BZNpotent}), inverting a degree $-2$ class $u$ (see Corollary~\ref{cor:IndPerf BcheckA Kosz}).
Over the sphere, it refines the $S^1$-Tate construction of homotopy theory \cite{GreenleesMay,NS}: we calculate the $\check{A}$-invariants of the periodization as an $\check{A}$-Tate construction $(-)^{t\check{A}}$, defined as the cofibre of a norm map from $\check{A}$-coinvariants to $\check{A}$-invariants (see Definition~\ref{defn:Tate}, Theorem~\ref{thm:Tate vs invariants of Bott}).

This paper supersedes our preprint \cite{ellhodge}, which introduced the loop spaces $\Loop_A(\YY)$ and the loop cohomology theories $\HH_A^G(-)$ over a field of characteristic zero.
The theory of periodization (Part~\ref{part:rotation}) and the extension to spectral bases such as the sphere spectrum $\mathbf{S}$ are new here.
We also correct here some serious errors in the treatment of stacks of tensor functors (specifically, it is necessary to work with \emph{sheaves} of tensor categories as in Subsect.~\ref{ssec:sheaves tensorcat}).

\subsection{Contents}\label{ssec:intro/contents}

Part~\ref{part:circles} is dedicated to the construction of categorified circles.
We begin with some preliminaries on spectral algebraic geometry in \S \ref{sec:perfect}, and develop the stacks of tensor functors used in the definition of loop spaces in \S \ref{sec:tensor}.
Finally, \S \ref{sec:Cartier} contains the construction of $\Perf(S^1_A)$.

Part~\ref{part:rotation} develops the theory of rotation-equivariant sheaves and periodization.
We introduce the rotation algebra $\RA$ in \S \ref{sec:RA} and study the structure of the category of ind-perfect complexes on $B\check{A}$ in \S \ref{sec:BcheckA}.
The Euler periodization and its relation to the Tate construction are studied in \S \ref{sec:Euler}, and the Bott periodization in \S \ref{sec:Bott per}.

The main constructions of the paper occupy Parts~\ref{part:loops} and~\ref{part:cohomology}.
Part~\ref{part:loops} is dedicated to $A$-loop spaces.
We give the definition in \S \ref{sec:loops} and specialize to the equivariant version in \S \ref{sec:loops BG}.
The case of $A$ an elliptic curve, or more generally an abelian variety, is treated in detail in \S \ref{sec:abvar}.
Finally, \S \ref{sec:loop properties} establishes general properties of $A$-loop spaces.

We come to the construction of $A$-loop cohomology $\HH_A^G$ in Part~\ref{part:cohomology}.
The definition is given in \S \ref{sec:cohomology}, its formal structure is described in \S \ref{sec:coh formal}, and its main properties are proven in \S \ref{sec:coh properties}.
Finally, \S \ref{sec:pushforwards} treats twists, pushforwards and Euler classes.

Appendix~\ref{app:strictification} records some facts about $E_\infty$-group stacks in spectral algebraic geometry, and their strictification.

\subsection{Conventions and notation}\label{ssec:intro/conv}

\subsubsection{}\label{sssec:intro/conv/base}

We work in the context of spectral algebraic geometry over a base $k$, a connective $E_\infty$-ring; see \cite{SAG}.
In particular, the reader may take $k$ to be an ordinary commutative ring.
When $k$ is a $\mathbf{Q}$-algebra, spectral algebraic geometry over $k$ collapses to derived algebraic geometry.

By default, all schemes and stacks are spectral in the sense of~\cite{SAG}.
Our terminology departs from Lurie's in the following respects:
\begin{enumerate}
	\item We will use the term ``smooth'' for the property Lurie calls \emph{fiber smooth} in \cite[Def.~11.2.5.5]{SAG}.
	The contrasting notion \emph{differentially smooth} of \cite[Def.~11.2.2.2]{SAG} will play only a minor role for us.

  \item We say a morphism is \emph{almost of finite presentation} if it is quasi-compact, quasi-separated, and locally almost of finite presentation in the sense of \cite[Def.~4.2.0.1]{SAG}.

  \item For group schemes, \emph{fppf} means flat and almost of finite presentation; see (\ref{sssec:fppf}).

  \item \emph{Commutative} group schemes are taken in the strict sense: that is, they are abelian group objects in the sense of \cite[Def.~1.2.4]{LurieEll1}.
  We also require them to be quasi-compact, separated, and (fiber) smooth; see Definition~\ref{def:comm grp scheme}.

  \item \emph{Abelian varieties} and \emph{elliptic curves} are taken in the strict sense of \cite[Defs.~1.5.1, 2.0.2]{LurieEll1}.
  For example, in our conventions, an elliptic curve over $k$ is a commutative group scheme which is proper and geometrically connected, of relative dimension $1$ over $k$.
  We write $\Mell$ for the moduli stack of elliptic curves, which is denoted $\mathcal{M}^{s}$ in \cite{LurieEll1}, and $\mathcal{M}_{\mathrm{Ell}}^{s}$ in \cite{LurieEll2}.

  \item An \emph{orientation} of a commutative group scheme $A$ over $S$ will mean for us a trivialization $\fra \simeq \OO_S$ (see Definition~\ref{def:orientation}).
  The preorientations, resp. orientations, of \cite{LurieEll2} will be called (resp. nondegenerate) \emph{periodicity data} (see Definition~\ref{def:periodicity datum}, Lemma~\ref{lem:Bott vs Lurie}).

  \item We never work with \emph{nonconnective} spectral Deligne--Mumford stacks in the sense of \cite[Def.~1.4.4.2]{SAG}, except in the proof of Theorem~\ref{thm:Otop}, Remark~\ref{rem:GM general base}, and Example~\ref{exam:KAG G=1}, where we explain that $\Specnc_{\Mell}(\K_{\mathcal{E}/\Mell})$, the (relative) nonconnective spectrum of our $\K_{\mathcal{E}/\Mell}$ (for $\mathcal{E}$ the universal elliptic curve over Lurie's $\mathcal{M}^{s}_{\mathrm{Ell}}$), reproduces precisely Lurie's nonconnective spectral stack $\mathcal{M}_{\mathrm{Ell}}^{\mathrm{or}}$ of oriented elliptic curves (see \cite[Prop.~7.2.10]{LurieEll2}).
\end{enumerate}

\subsubsection{}\label{sssec:intro/conv/io}

Throughout the paper, \emph{category} means $\infty$-category and \emph{groupoid} means $\infty$-groupoid.

Two kinds of categories appear in the construction of $A$-loop cohomology:
\begin{itemize}
\item The \emph{input}: $\QCoh(\YY)$, the symmetric monoidal category of \emph{connective} quasi-coherent sheaves on $\YY$.
\item The \emph{output}: $\Rot^G(A) := \INDPERF^{\check{A}}(\Loop_A(BG))$, the symmetric monoidal stable category of $\check{A}$-equivariant ind-perfect complexes on the loop space.
\end{itemize}
We write $\QCoh(\YY)$, $\Perf(\YY)$, $\Coh(\YY)$ for the connective (prestable) versions, and $\QCOH(\YY)$, $\PERF(\YY)$, $\COH(\YY)$ for their nonconnective (stable) counterparts.

\subsubsection{}

We use homological grading conventions.

\subsection{Acknowledgments}

The authors would like to thank Andrei Okounkov for the suggestion to look into these questions, and David Ben-Zvi, Sanath Devalapurkar, Jacob Lurie, Tasos Moulinos, David Nadler and Andrei Okounkov for comments on a draft.

E.B.\ is grateful to Wille Liu, Andy Jiang, Ian Grojnowski, Arkadij Bojko and particularly Chris Brav for helpful conversations related to this material.

A.A.K.\ thanks Kyoto University and RIMS for the hospitality and peaceful working conditions during his stay from April to August 2026, when much of his work on this paper was carried out, and acknowledges support from the grants AS-CDA-112-M01 and NSTC 112-2628-M-001-006.

\subsubsection*{Disclosure of software tools}

The authors used AI tools for literature search, for checking proofs and locating errors, and for typesetting assistance.

\addtocontents{toc}{\protect\setcounter{tocdepth}{2}}

\part{Categorified circles}\label{part:circles}

\section{Perfect and pluperfect stacks}\label{sec:perfect}

\subsection{Recollections on spectral algebraic geometry}
\label{ssec:dagreview}

For the reader's convenience we include a brief recollection on derived stacks and various categories of quasi-coherent sheaves on them.
We refer to \cite{GR,LurieHA,HAG} for textbook accounts.

\subsubsection{Stacks}

We let $\dAlgk$ denote the category of connective $E_\infty$-algebras over the base $k$.
When $k$ is a $\mathbf{Q}$-algebra, these are equivalently connective commutative dg-algebras.
A \emph{prestack} is by definition a functor: $$\YY: R \longmapsto \YY(R),$$ encoding for each $R \in \dAlgk$ a groupoid\footnote{%
  We recall our convention that \emph{groupoid} means $\infty$-groupoid.
} $\YY(R)$ of $R$-points.
A \emph{stack} is a prestack satisfying descent for the \'{e}tale topology. The category of such is denoted $\mathbf{St}_{k}$. Representable objects are called \emph{affine schemes}.
Recall that there is a notion of (higher) \emph{Artin} stacks (sometimes called \emph{geometric} stacks), see e.g. \cite[Chap.~2, \S 4]{GR}.

\subsubsection{Quasi-coherent sheaves and perfect complexes}

For a commutative $k$-algebra $R$ we write $\Mod_R$ for the category of $R$-modules.
These are connective module spectra over $R$ (or, when $k$ is a $\mathbf{Q}$-algebra, equivalently connective dg-modules). For a stack $\YY$, the category of quasi-coherent sheaves $\QCoh(\YY)$ is defined by $$\QCoh(\spc(R)) = \Mod_R$$ in the affine case, and in general by the limit \[\QCoh(\YY) = \lim_{R,y} \Mod_R \] over $R\in\dAlgk$ and $R$-points $y \in \YY(R)$. The category $\QCoh(\YY)$ is presentable, and there is a unique symmetric monoidal structure on $\QCoh(\YY)$ which is colimit-preserving in each argument and for which the pull-backs $\QCoh(\YY) \to \Mod(R)$ are symmetric monoidal for every $R\in\dAlgk$ and every morphism $\spc(R) \to \YY$.
The categories $\Mod_R$ and $\QCoh(\YY)$ are prestable in the sense of Lurie (see \cite[App.~C]{SAG}), and form the connective parts of $t$-structures on their stabilizations,
$$ \MOD_R \supset \Mod_R, \quad \text{resp.}\ \QCOH(\YY) \supset \QCoh(\YY). $$

We write $\Perf_R \subset \Mod_R$ for the full subcategory of \emph{perfect} $R$-modules, generated by $R$ under finite colimits and extensions. We write $\Perf(\YY) \subset \QCoh(\YY)$ for the full subcategory of \emph{perfect complexes} on $\YY$, i.e. the limit of $\Perf_R$ over $R\in\dAlgk$ and $y \in \YY(R)$. Its stabilization $$\PERF(\YY) \supset \Perf(\YY)$$ is the full subcategory of perfect complexes in $\QCOH(\YY)$.

When $\YY$ is Artin, the limits defining the various categories above can be taken over pairs $(R,y)$ where $R\in\dAlgk$ and $y : \spc(R) \to \YY$ is a \emph{smooth} morphism. It follows from \cite[Prop.~C.3.2.4]{SAG} that $\QCoh(\YY)$ is a Grothendieck prestable category in the sense of \cite[\S C.1.4]{SAG}.

\subsubsection{Affine stacks}

Any $E_\infty$-algebra $R$ over $k$, not necessarily connective, represents the stack
$$\Specnc(R) : B \mapsto \Map_{\dAlgk^{\mathrm{nc}}}(R, B),$$ where $\dAlgk^{\mathrm{nc}}$ is the category of nonconnective $E_\infty$-algebras over $k$.
The \emph{affine stack} associated to a stack $\YY$,
$$ \operatorname{Aff}(\YY) := \Specnc(\OO(\YY)), $$
has the universal property that any morphism $\YY\rightarrow S$ with $S$ an affine stack factors uniquely through $\YY \to \operatorname{Aff}(\YY)$. For example, when $k$ is an ordinary commutative ring, the affinization of an elliptic curve $E$ is non-canonically isomorphic to $B\GGa$, with an isomorphism corresponding to a trivialization of the one dimensional vector space $H^{1}(E,\mathcal{O})$. For a reference on affine stacks and affinization see \cite[\S 3]{BNloop}.

\subsubsection{Variant: derived algebraic geometry}\label{sssec:dagreview/animated}

Our choice to work in the context of spectral algebraic geometry in this paper is forced by our use of Tannaka duality (see Subsect.~\ref{ssec:stktenfun}), which relates spectral stacks to their symmetric monoidal stable $\infty$-categories of quasi-coherent sheaves.
However, one can indeed develop a parallel formalism in \emph{derived} algebraic geometry, the theory of stacks built out of animated commutative rings (as opposed to connective $E_\infty$-ring spectra), as follows.
The definition of \emph{tensor category} in Subsect.~\ref{ssec:tensorcat} should be modified by replacing symmetric monoidal $\infty$-categories by \emph{$\Theta$-categories} in the sense of Nuiten--To\"en \cite{NuitenToen} (which roughly speaking encode the extra structure of derived symmetric power operations).
We may then replace the invocations of \cite[Thm.~1.3]{BHL} by the Tannaka duality result of \cite[Thm.~6.1.5]{MundingerSahai}.

When $k$ is a $\bQ$-algebra, the two formalisms agree because derived symmetric powers then coincide with $E_\infty$ ones (see the introduction of \cite{NuitenToen}).
Over a general base, the animated formalism is better behaved in some respects; for example, the theory of twists and push-forwards developed in Sect.~\ref{sec:pushforwards} becomes much simpler there.
We do not develop animated loop spaces in detail here; they enter only in Theorem~\ref{thm:HKR fil}.

\subsection{Perfect and pluperfect stacks}

\subsubsection{Perfect stacks}

We recall the following definition from \cite{BFN}:

\begin{definition}\label{def:perfect}
  A stack $\YY$ is called \emph{perfect} if it has affine diagonal, $\OO_\YY$ is compact, and $\QCOH(\YY)$ is compactly generated by perfect complexes.
\end{definition}

For example, for a qcqs algebraic space $X$, $\QCOH(X)$ is compactly generated by a \emph{single} perfect complex by \cite[Thm.~9.6.3.2]{SAG}.
Perfect stacks are Tannakian by \cite{BHL} (see Subsect.~\ref{ssec:stktenfun}).

\subsubsection{Pluperfect stacks}

We introduce the following variant of perfectness, which requires a finite number of generators, but does not require perfect complexes to be compact.

\begin{definition} \label{def:pluperfect} We call $\YY$ \emph{pluperfect} if it is Tannakian (Definition~\ref{def:Tannakian}) and there exists a finite set $\{F_{i}\}$ of perfect complexes on $\YY$ which generate the derived category $\QCOH(\YY)$ under colimits, finite limits, direct summands, tensor products, and duals. \end{definition}

\begin{remark}
	We refer to a set $\{F_i\}$ as in Definition~\ref{def:pluperfect} as a \emph{set of perfect tensor generators} of $\QCOH(\YY)$.
\end{remark}

\begin{remark}\label{rem:perfect and pluperfect}
Suppose $\YY$ has affine diagonal and $\OO_{\YY}$ is compact.
Then $\YY$ is perfect and pluperfect if and only if $\QCOH(\YY)$ admits a finite set of perfect tensor generators.
\end{remark}

\begin{lemma}\label{lem:pluperfect base change}
  Let $k_0 \to k$ be a morphism of connective $E_\infty$-rings, $\YY_0$ a stack over $k_0$ admitting a faithfully flat qcqs morphism from an affine scheme, and $\YY := \YY_0 \times_{\spc(k_0)} \spc(k)$.
  If $\YY_0$ is perfect and pluperfect, then so is $\YY$.
\end{lemma}
\begin{proof}
  The projection $q : \YY \to \YY_0$ is affine, so perfectness holds by \cite[Prop.~3.21]{BFN}.
  As $\YY$ also admits a faithfully flat quasi-compact quasi-separated morphism from an affine scheme, it is Tannakian by \cite[Thm.~1.3]{BHL}.
  Let $\{F_i\}$ be a finite set of perfect complexes generating $\QCOH(\YY_0)$ under colimits, finite limits, direct summands, tensor products and duals.
  Then the full subcategory of $M \in \QCOH(\YY_0)$ with $q^*M$ in the closure of $\{q^*F_i\}$ under these operations is all of $\QCOH(\YY_0)$.
  Since $q$ is affine, $q_*$ is conservative and colimit-preserving, so every $N \in \QCOH(\YY)$ is the geometric realization of a simplicial object with terms $q^*(q_*q^*)^n q_*N$ (see \cite[Prop.~4.7.3.14]{LurieHA}).
  Thus $\{q^*F_i\}$ generates $\QCOH(\YY)$.
\end{proof}

\subsection{$G$-quasi-projective schemes}

\subsubsection{}\label{sssec:fppf}

Let $S$ be a stack over $k$.
By an \emph{fppf} group scheme $G$ over $S$, we mean a group stack that is schematic, flat, and almost of finite presentation.
We write $B_S(G) := S/G$ for the classifying stack, or simply $BG$ when there is no risk of ambiguity.

\subsubsection{}

An affine fppf group scheme $G$ over $S$ is \emph{linearly reductive} if for every affine $T$ over $S$, the pushforward $\QCOH(B_T(G_T)) \to \QCOH(T)$ is t-exact.

We say $G$ is \emph{reductive} if it is smooth, has reductive geometric fibres, and $G/G^\circ$ is finite over $S$ (where $G^\circ$ is the identity component).\footnote{%
  Unlike \cite[Exp.~XIX]{SGA3}, we do not assume any connectedness.
}
If $k$ is a $\bQ$-algebra, then every reductive group scheme is linearly reductive (see \cite[Thm.~9.7.6, Rem.~9.1.3]{AlperAdequate}).

\begin{definition}
  Let $G$ be an fppf group scheme over a stack $S$.
  A stack $X$ over $S$ with $G$-action is \emph{$G$-quasi-projective} if it is schematic and almost of finite presentation over $S$, and admits a $G$-linearized line bundle $L$ which is relatively ample over $S$ (i.e., $L|_{X_\cl}$ is relatively ample in the classical sense, after base change to any affine over $S_\cl$).
  It is \emph{$G$-projective} if it is $G$-quasi-projective and proper over $S$.
\end{definition}

\subsubsection{}

We record the following classical lemma for later use:

\begin{lemma}\label{lem:VarG projective embedding}
  Suppose $k$ is an ordinary commutative ring, $G$ is a linearly reductive group scheme over $k$, and $X$ is a classical $G$-quasi-projective scheme of finite presentation over $k$.
  Then, Nisnevich-locally on $\spc(k)$, there exists a finite projective representation $W$ of $G$ and a $G$-equivariant quasi-compact immersion $X \hookrightarrow \mathbf{P}(W)$.
\end{lemma}
\begin{proof}
  Let $L$ be a $G$-linearized ample line bundle on $X$.
  Then we may find sections $s_0, \ldots, s_n \in \Gamma(X, L^{\otimes d})$, for some $d > 0$, defining an immersion $i : X \hookrightarrow \mathbf{P}^n_{k}$.
  The proof of Lemma~\ref{lem:BG pluperfect lin red} shows that, after passing to a Nisnevich cover of $\spc(k)$, every object of $\QCOH(BG)^\heartsuit$ is a quotient of a direct sum of finite projective representations.
  Applied to the representation $\Gamma(X, L^{\otimes d})$, this gives a finite projective representation $W$ and a $G$-equivariant map $W \to \Gamma(X, L^{\otimes d})$ whose image contains $s_0, \ldots, s_n$.
  Since the $s_j$ generate $L^{\otimes d}$, the induced map $W \otimes \OO_{X} \to L^{\otimes d}$ is a $G$-equivariant surjection, and classifies a $G$-equivariant morphism $j : X \to \mathbf{P}(W)$.
  Lifting the $s_j$ to $W$ defines a linear projection $\mathbf{P}(W) \supset U \to \mathbf{P}^n_{k}$ on an open $U$ containing $j(X)$, through which $i$ factors as $X \xrightarrow{j} U \to \mathbf{P}^n_{k}$.
  The resulting morphism $X \to U \subset \mathbf{P}(W)$ is a quasi-compact immersion (since $U$ is quasi-compact and separated).
\end{proof}

\subsection{Perfectness and pluperfectness of global quotients}

For $1$-Artin stacks $\YY$, $\QCOH(\YY)$ is often perfect but rarely compactly generated by \emph{finitely} many perfect complexes. Nevertheless, we have\footnote{%
  See also \cite{DeligneMilne,Gross} for analogous statements for the abelian categories $\QCOH(\YY)^\heartsuit$.
}:

\begin{theorem}\label{thm:pluperfect quotient}
Let $G$ be an affine fppf group scheme over $k$ and $X$ a $G$-quasi-projective $k$-scheme.
\begin{enumerate}
\item If $k$ is a $\bQ$-algebra and $G$ admits a closed embedding $G \hookrightarrow \GL_{n,k}$ such that $\GL_{n,k}/G$ is $\GL_{n,k}$-quasi-projective\footnote{%
  e.g. $G$ is the base change of an affine algebraic group over a field of characteristic zero; cf. the proof of \cite[Cor.~3.22]{BFN}
}, then $X/G$ is perfect and pluperfect.
\item If $G$ is linearly reductive, then $X/G$ is perfect, and pluperfect Nisnevich-locally on $\spc(k)$.
\end{enumerate}
\end{theorem}

\begin{lemma}\label{lem:BG pluperf implies X/G pluperf}
  Let $G$ be an affine fppf group scheme and $X$ a $G$-quasi-projective scheme.
  If $BG$ is perfect, then so is $X/G$; if $BG$ is moreover pluperfect, then $X/G$ is perfect and pluperfect.
\end{lemma}
\begin{proof}
Since $X$ admits an ample line bundle, it has affine diagonal (this can be checked on classical truncations, where an argument is given in the proof of \cite[Cor.~2.3.3]{ThomasonTrobaugh}).
Then since $G$ is affine, $X/G$ also has affine diagonal.
The projection $p : X/G \to BG$ is qcqs representable, so $p_*$ preserves colimits.
Since $BG$ is perfect and hence $\OO_{BG}$ is compact, $\Gamma(X/G,-) \simeq \Gamma(BG, p_*(-))$ preserves colimits; that is, $\OO_{X/G}$ is compact.

Let $L$ be a $G$-linearized ample line bundle on $X$.
The perfect complexes $p^*M \otimes L^{\otimes a}$, for $M \in \PERF(BG)$ and $a \in \mathbf{Z}$, are compact because $\OO_{X/G}$ is.
To show that $X/G$ is perfect it will suffice to show they generate $\QCOH(X/G)$.
Let $N$ be an object right orthogonal to all their shifts.
Since the $M \in \PERF(BG)$ compactly generate $\QCOH(BG)$, adjunction gives
$p_*(N \otimes L^{\otimes -a}) = 0$ for all $a$, hence $\Gamma(X, N|_X \otimes L^{\otimes -a}) = 0$ for all $a$ by base change.
Since $L$ is ample, $X$ admits an affine open cover by the non-vanishing loci $X_s$ of $s \in \pi_0\Gamma(X, L^{\otimes n})$, $n\ge 1$.
For each,
\[\Gamma(X_s, N) \simeq \varinjlim \bigl(\Gamma(X, N|_X) \xrightarrow{s} \Gamma(X, N|_X \otimes L^{\otimes n}) \xrightarrow{s} \cdots\bigr)\]
vanishes, hence $N|_X = 0$ and hence $N = 0$.

Now suppose $\{F_i\}$ is a finite set of perfect tensor generators of $\QCOH(BG)$, and let us show that the objects $L$ and $p^*F_i$ form a set of tensor generators of $\QCOH(X/G)$.
Let $\mathcal{C} \subset \QCOH(X/G)$ be the full subcategory they generate under colimits, finite limits, direct summands, tensor products, and duals.
Since $p^*$ preserves these operations and $F_i$ tensor-generate, the full subcategory of $\QCOH(BG)$ spanned by $M$ with $p^*M$ lying in $\mathcal{C}$ is the entirety of $\QCOH(BG)$.
In particular, $\mathcal{C}$ contains the compact generators $p^*M \otimes L^{\otimes a}$, and hence $\mathcal{C} = \QCOH(X/G)$.
Finally, Remark~\ref{rem:perfect and pluperfect} now shows that $X/G$ is pluperfect.
\end{proof}

\begin{lemma}\label{lem:BG pluperfect lin red}
  Let $k$ be a connective $E_\infty$-ring and $G$ an affine fppf group scheme over $S=\spc(k)$.
  If $G$ is linearly reductive, then we have:
  \begin{enumerate}
    \item
    $\OO_{BG}$ is compact.

    \item\label{item:BG pluperfect lin red/compgen}
    $\QCOH(BG)$ is compactly generated by the finite projective representations.

    \item\label{item:BG pluperfect lin red/generators}
    Nisnevich-locally on $S$, there exists a vector bundle $V$ on $BG$ which is a perfect tensor generator of $\QCOH(BG)$.
  \end{enumerate}
  In particular, $BG$ is perfect, and pluperfect Nisnevich-locally on $S$.
\end{lemma}
\begin{proof}
  Since $G$ is affine, $BG$ has affine diagonal.
  Since $G$ is linearly reductive, i.e. $BG$ is of cohomological dimension zero over $S$, \cite[Prop.~9.1.5.3]{SAG} shows that $\OO_{BG}$ is compact, and the proof of \cite[Prop.~A.3.4]{BKRS} shows that $\QCOH(BG)$ is compactly generated by the finite projective representations.
  This shows that $BG$ is perfect.

  It remains to prove claim~(\ref{item:BG pluperfect lin red/generators}).
  Replacing $S$ by a Nisnevich cover if necessary, we may assume by \cite[Cor.~6.2(2)]{AHRLocalStructure} that $G_\cl$ admits a faithful finite projective representation $V_0$.
  Since the quotient $\GL(V_0)/G_\cl$ is affine by Matsushima (\cite[Rem.~2.5]{AHRLocalStructure}), \cite[Thm.~5.4, Prop.~5.5]{Gross} imply that every object of $\QCOH(BG_\cl)^\heartsuit$ is a quotient of a direct sum of objects $V_0^{\otimes m} \otimes V_0^{\vee, \otimes n}$, $m,n\ge 0$.

  We lift $V_0$, regarded as a vector bundle on $BG_\cl$, to a vector bundle $V$ on $\XX := BG$.
  By \cite[Cor.~19.2.1.4]{SAG}, there is an equivalence $\Perf(\XX) \simeq \lim_n \Perf(\XX_n)$, where $\XX_n := BG \times_S \spc(\tau_{\le n}k)$, and by \cite[Rem.~19.2.2.6]{SAG}, the obstruction to lifting from $\XX_n$ to $\XX_{n+1}$ lies in
  \[ H^{n+3}\bigl(BG_\cl, V_0^\vee \otimes V_0 \otimes_{\pi_0(k)} \pi_{n+1}(k) \bigr), \]
  which vanishes since $G$ is linearly reductive.\footnote{%
    Both \cite{SAG} statements globalize to the stack $\XX$ by fppf descent; see e.g. the proofs of \cite[Prop.~3.1.5]{HLP} and \cite[Lem.~5.1.11, Prop.~5.1.13]{HLP}, respectively.
  }
  We thus obtain a compatible collection of lifts $V_n$, and hence $V \in \Perf(\XX)$; by \cite[Prop.~2.7.3.2]{SAG}, $V$ is still a vector bundle.

  By the argument of \cite[Lem.~A.2.6(i)]{BKRS}, it follows that every finite projective representation of $G$ is a direct summand of a finite direct sum of objects $V^{\otimes m}\otimes V^{\vee,\otimes n}$, $m,n\ge0$ (by lifting a retraction from the level of classical truncations).
  By claim~(\ref{item:BG pluperfect lin red/compgen}) it follows that $V$ is a perfect tensor generator.
\end{proof}

\begin{lemma}\label{lem:BG pluperfect char 0}
  Suppose $k$ is a $\bQ$-algebra.
  Let $G$ be an fppf group scheme over $k$ admitting a closed embedding $G \hookrightarrow \GL_{n,k}$ such that $\GL_{n,k}/G$ is $\GL_{n,k}$-quasi-projective.
  Then $BG$ is perfect and pluperfect.
\end{lemma}
\begin{proof}
  Since $k$ is a $\bQ$-algebra, $\GL_{n,k}$ is linearly reductive and by Lemma~\ref{lem:BG pluperfect lin red}, $\OO_{B\GL_{n,k}}$ is compact.
  Moreover, taking $V_0$ in the proof of Lemma~\ref{lem:BG pluperfect lin red}(\ref{item:BG pluperfect lin red/generators}) to be the standard representation shows the existence of a perfect tensor generator $V \in \QCOH(B\GL_{n,k})$, globally on $\spc(k)$.
  Thus $B\GL_{n,k}$ is perfect and pluperfect by Remark~\ref{rem:perfect and pluperfect}.
  Then since $Q := \GL_{n,k}/G$ is $\GL_{n,k}$-quasi-projective, Lemma~\ref{lem:BG pluperf implies X/G pluperf} shows that $Q/\GL_{n,k} \simeq BG$ is also perfect and pluperfect.
\end{proof}

\begin{proof}[Proof of Theorem~\ref{thm:pluperfect quotient}]
  By Lemma~\ref{lem:BG pluperf implies X/G pluperf} we reduce to the case of $BG$, which is shown in Lemmas~\ref{lem:BG pluperfect char 0} and \ref{lem:BG pluperfect lin red}, respectively.
\end{proof}

\subsection{$G$-varieties}\label{ssec:Var_G}

\subsubsection{}

Fix a stack $S$ and an fppf group scheme $G$ over $S$.

\begin{definition}\label{def:VarG}
  A stack $X$ over $S$ with $G$-action is a \emph{relative $G$-variety} if it satisfies the following conditions:
  \begin{enumerate}
	  \item $X$ is schematic and almost of finite presentation over $S$;
	  \item $X_\cl$ is $G_\cl$-quasi-projective over $S_\cl$; and
		\item $(X/G) \times_S T$ is perfect, and pluperfect fppf-locally on $T$, for every affine $T$ over $S$.
  \end{enumerate}
  Note that these conditions are stable under base change in $S$.
  A \emph{$G$-variety} is a relative $G$-variety over $\spc(k)$.
\end{definition}

\begin{lemma}\label{lem:qproj implies variety}
	If $G$ is linearly reductive and $X$ is $G$-quasi-projective over $S$, then $X$ is a relative $G$-variety over $S$.
\end{lemma}
\begin{proof}
  Conditions (i) and (ii) hold by definition.
  For (iii), let $T$ be an affine over $S$.
  Since $X \to S$ is almost of finite presentation, hence \emph{a fortiori} qcqs, the base change $X_T$ is quasi-compact.
  If $L$ is a $G$-linearized $S$-ample line bundle on $X$, then $L_T$ is a $G_T$-linearized line bundle on $X_T$ which is relatively ample over the affine $T$, hence ample.
  Hence $X_T/G_T$ is perfect, and pluperfect Nisnevich-locally on $T$, by Theorem~\ref{thm:pluperfect quotient}(ii) applied with $k=\OO_T$.
\end{proof}

\subsubsection{}

Given an fppf group scheme $G$ over $k$, we write $\Stk^G$ for the category of stacks over $k$ with $G$-action, $\Sch^G$ for the full subcategory of qcqs schemes with $G$-action, and $\Var^G$ for the full subcategory of $G$-varieties.

\begin{example}\label{ex:VarG}
  Let $G$ be an affine fppf group scheme over $k$ and $X$ a $G$-quasi-projective $k$-scheme.
  If $G$ satisfies the hypotheses of Theorem~\ref{thm:pluperfect quotient}(i) or (ii), then $X \in \Var^G$.
\end{example}

\begin{example}
	Let $k_0$ be an ordinary commutative ring, $G_0$ a linearly reductive group scheme over $k_0$, and $X_0$ a classical $G_0$-quasi-projective scheme over $k_0$.
	Then for every connective $E_\infty$-algebra $k$ over $k_0$, the base change $X := X_0 \times_{\spc(k_0)} \spc(k)$ is $G$-quasi-projective over $k$, where $G := G_0 \times_{\spc(k_0)} \spc(k)$, and hence lies in $\Var^G$.
\end{example}

\begin{example}
	If $k$ is an ordinary commutative ring and $G$ is linearly reductive over $S=\spc(k)$, then $\Var^G$ includes all classical $G$-quasi-projective $k$-schemes.
	Moreover, if $k$ is an ordinary $\bQ$-algebra, then this holds without the linear reductivity hypothesis.
\end{example}

\section{Tensor functors}\label{sec:tensor}

\subsection{Tensor categories} \label{ssec:tensorcat}
We introduce now the categorical backdrop for our main constructions. We refer again to \cite[\S C.1.4]{SAG} for background on Lurie's theory of Grothendieck prestable $\infty$-categories.

\begin{definition} A \emph{tensor category} $\CC$ is a symmetric monoidal $k$-linear Grothendieck prestable category, where by convention the tensor product is assumed colimit-preserving in each argument. Given tensor categories $\CC$ and $\DD$, we denote by $\TensFun(\CC, \DD)$ the category of \emph{tensor functors}, i.e. symmetric monoidal $k$-linear colimit-preserving functors, from $\CC$ to $\DD$. \end{definition}

We denote by $\CCC$ the category of tensor categories and tensor functors. This is the category of commutative algebra objects in $\Cat$, Lurie's symmetric monoidal category of Grothendieck prestable $k$-linear categories (see \cite[\S C.4.2]{SAG}).

\begin{example} For an Artin stack $\YY$, $\QCoh(\YY)$ is a tensor category. Further, for a morphism $f:\YY\rightarrow\XX$ the induced functor $$f^{*}:\QCoh(\XX)\longrightarrow\QCoh(\YY)$$ is a tensor functor. \end{example}

\begin{example} Let us say a \emph{tensor $1$-category} is a Grothendieck abelian $1$-category equipped with a symmetric monoidal structure for which the tensor product commutes with colimits in each argument, and a \emph{tensor functor} between tensor $1$-categories is a colimit-preserving symmetric monoidal functor. Given a ``nice'' tensor $1$-category $A$\footnote{%
  More precisely, assume that $A$ admits a ``weakly flat descent structure'' in the sense of \cite[2.2, 3.1]{CisDeg}.
  For our purposes it is enough to consider the case of $A$ generated by compact projective objects, in which case $A$ is equivalent to the category of functors $(A_0)^{\mathrm{op}} \to \mathrm{Ab}$ that preserve finite products, where $A_0$ is the full subcategory of compact projectives (see e.g. \cite[\S 5.1]{CesScholze}).
}, the connective derived category $\DD_{\geq 0}(A)$ is a tensor category, and for any tensor functor $f : A \to B$, the left-derived functor $$Lf : \DD_{\geq 0}(A) \to \DD_{\geq 0}(B)$$ is a tensor functor. See \cite[Thm.~2.14, Prop.~3.2]{CisDeg} and \cite[Prop.~1.3.5.15, Cor.~1.3.4.26, Ex.~4.1.7.6]{LurieHA}. \end{example}

\begin{lemma}\label{lem:easyKunnethforQCoh}
  Let $\YY$ and $\YY'$ be Artin stacks over an affine scheme $S$, and consider the canonical symmetric monoidal functor
  \[ \QCoh(\YY) \otimes_{\QCoh(S)} \QCoh(\YY') \xrightarrow{\boxtimes} \QCoh(\YY \times_S \YY'). \]
  If $\QCoh(\YY')$ is dualizable in $\Cat_S$, then this is an equivalence.
\end{lemma}
\begin{proof}
  Note first that the claim holds if $\YY$ and $\YY'$ are both affine.
  Next suppose $\YY'=T$ is affine but $\YY$ is arbitrary.
  Then $\QCoh(T)$ is dualizable in $\Cat_S$, i.e., $\QCoh(T)\otimes_S -$ preserves limits.
  Writing $\YY$ as a colimit of affines $S_\alpha$, we have $\QCoh(\YY) \simeq \lim_\alpha \QCoh(S_\alpha)$ (by construction of $\QCoh(\YY)$).
  Thus the claim reduces in this case to the equivalence $$\QCoh(S_\alpha) \otimes_{\QCoh(S)} \QCoh(T) \simeq \QCoh(S_\alpha \times_S T).$$
  Finally, let $\YY$ and $\YY'$ be arbitrary and assume only that $\QCoh(\YY')$ is dualizable in $\Cat_S$.
  Then writing $\YY \simeq \varinjlim_\alpha S_\alpha$ again, the claim reduces to the equivalence $$\QCoh(S_\alpha) \otimes_{\QCoh(S)} \QCoh(\YY') \simeq \QCoh(S_\alpha \times_S \YY')$$ which follows from the previous step.
\end{proof}

\subsection{Sheaves of tensor categories}\label{ssec:sheaves tensorcat}

We will in fact have to pass to the following more general setup.

For an Artin stack $S$, the category $\Cat_S$ of \emph{$S$-linear Grothendieck prestable categories} is by definition the category of $\QCoh(S)$-modules in $\Cat$, where $\QCoh(S)$ is regarded as a commutative algebra object in $\Cat$. The category $\CCC_S$ of \emph{$S$-linear tensor categories} is by definition the category of commutative algebra objects in $\Cat_S$. For example, $\QCoh(\YY)$ is an $S$-linear tensor category for any Artin stack $\YY$ over $S$.

Given $S$-linear tensor categories $\CC$ and $\DD$, we denote the tensor product by $\CC \otimes_S \DD := \CC \otimes_{\QCoh(S)} \DD$ and by $$\TensFun_S(\CC, \DD)$$ the category of \emph{$S$-linear tensor functors} from $\CC$ to $\DD$, i.e., $\QCoh(S)$-linear symmetric monoidal colimit-preserving functors.

\begin{definition}
  We define the category $\stCCC$ of \emph{sheaves of tensor categories} as the lax limit of the diagram
  \[ T \mapsto \CCC_{T} \]
  indexed by affine schemes $T$.
  That is, a sheaf of tensor categories $\stCC$ associates:
  \begin{itemize}
      \item to every affine scheme $T$, a $T$-linear tensor category $\stCC(T)$;
      \item to every morphism of affine schemes $f : T' \to T$, a $T'$-linear tensor functor $$f^* : \stCC(T) \otimes_{\QCoh(T)} {\QCoh(T')} \to \stCC(T');$$
      \item a system of coherent homotopies up to which $f^*$ are compatible with composition in $f$.
  \end{itemize}
  We say that a sheaf of tensor categories is \emph{quasi-coherent} if for every morphism of affines $f$, the functor $f^*$ above is an equivalence.
\end{definition}

Note that $\stCCC$ is equivalently the category of commutative algebra objects in $\stCat$, the symmetric monoidal category defined as the lax limit of the diagram $T \mapsto \Cat_{T}$.

More generally still, for a base Artin stack $S$, the category $\stCCC_S$ of \emph{sheaves of $S$-linear tensor categories} is the lax limit of the diagram
$$ T \mapsto \CCC_T, $$
where $T$ ranges over affine schemes over $S$. This is the category of commutative algebra objects in the symmetric monoidal category $\stCat_S$\footnote{This should not be confused with $\operatorname{ShvCat}(S),$ as studied for example by Gaitsgory in \cite{1-aff}.}, the lax limit of $T \mapsto \Cat_T$.

\begin{example}
  For an $S$-linear tensor category $\CC$, we denote by $\CC \otimes_S \stQCoh_S \in \stCCC_S$ the quasi-coherent sheaf of $S$-linear tensor categories
  $$ T \mapsto (\CC \otimes \stQCoh_S)(T) := \CC \otimes_S \QCoh(T). $$
\end{example}

\begin{example}
  Let $S$ be an Artin stack.
  For an Artin stack $\YY$ over $S$ and an affine scheme $T$ over $S$, we denote by $\stQCoh_{\YY/S} \in \stCCC_S$ the sheaf of tensor categories
  $$ T \mapsto \stQCoh_{\YY/S}(T) := \QCoh(\YY \times_S T). $$
  By Lemma~\ref{lem:easyKunnethforQCoh} (applied to $\YY\times_S T$ and $T'$ over $T$), it is quasi-coherent; moreover, there is a canonical isomorphism $$\QCoh(\YY) \otimes_{S} \stQCoh_S \simeq \stQCoh_{\YY/S}$$ in $\stCCC$.

\end{example}

When $S=\spc(k)$, we drop the decorations ``$/S$'' from the notation.

\subsection{Stacks of tensor functors} \label{ssec:stktenfun}
Recall that for two stacks $\XX$ and $\YY$, there is a \emph{mapping stack}
\[\MAP(\XX,\YY) : T \mapsto \Map(\XX \times T, \YY).\]
We introduce an analogue of this construction which takes sheaves of tensor categories as input. Note that the output is taken to be a (pre-)stack, i.e., this is not an internal mapping object in sheaves of tensor categories.

\begin{definition} Let $\stCC$ and $\stDD$ be sheaves of $S$-linear tensor categories, where $S$ is an Artin stack. The \emph{$S$-stack of tensor functors} from $\stCC$ to $\stDD$ is defined by the assignment
\[ \TENSFUN_S(\stCC, \stDD) : T \mapsto \Map_{\stCCC_T}(\stCC|_T, \stDD|_T)^\simeq, \]
sending an affine scheme $T$ over $S$ to the groupoid of maps from $\stCC|_T$ to $\stDD|_T$ in the category of sheaves of $T$-linear tensor categories.
\end{definition}

\begin{remark}
  When $\stCC$ is quasi-coherent, a morphism of sheaves out of $\stCC|_T$ is uniquely determined by its value on the terminal object $T$, so we have:
  \[ \TENSFUN_S(\stCC, \stDD)(T) \simeq \TensFun_T(\stCC(T), \stDD(T))^\simeq. \]
\end{remark}

\begin{example} Let $\CC$ be a tensor category and $\stDD$ a sheaf of tensor categories. Then $\TENSFUN(\CC \otimes \stQCoh, \stDD)$ is the stack defined by
$$ T \mapsto \TensFun_T(\CC \otimes \QCoh(T), \stDD(T))^{\simeq} \simeq \TensFun(\CC, \stDD(T))^{\simeq}. $$ In particular, if $\stDD = \DD \otimes \stQCoh$ for some tensor category $\DD$, then $\TENSFUN(\CC \otimes \stQCoh, \DD \otimes \stQCoh)$ is given by $$ T \mapsto \TensFun(\CC, \DD \otimes \QCoh(T))^{\simeq}. $$ \end{example}

\begin{example} Let $\XX$ and $\YY$ be Artin stacks over $S=\spc(k)$. Then it follows from Lemma~\ref{lem:easyKunnethforQCoh} that $\TENSFUN(\stQCoh_{\YY}, \stQCoh_{\XX})$ is given by $$ T \mapsto \TensFun(\QCoh(\YY), \QCoh(\XX\times T)). $$ In particular, there is a canonical morphism of stacks
\begin{equation}\label{eq:MaptoTENSFUN}
  \MAP(\XX,\YY) \rightarrow \TENSFUN(\stQCoh_{\YY}, \stQCoh_{\XX}),
\end{equation}
sending $f \mapsto f^*$. \end{example}

\begin{definition}\label{def:Tannakian}
  A stack $\YY$ is \emph{Tannakian} if the morphism \eqref{eq:MaptoTENSFUN} is invertible for all stacks $\XX$ (or equivalently, for all affine schemes $\XX=S$).
\end{definition} \begin{remark} There are many natural stacks which are \emph{not} Tannakian, and the theory developed in this paper does not have much to say about their loop spaces. For example, for $n\geq 2$ the unit morphism $$\QCoh(*)\rightarrow\QCoh(B^n\GM)$$ is an equivalence and so $B^n\GM$ is not Tannakian. \end{remark}

\begin{example}
  If $\YY$ has quasi-affine diagonal and $\QCOH(\YY)$ is compactly generated, it is Tannakian by \cite[Thm.~1.3]{BHL}.\footnote{Note that the compact generation is known for a large class of algebraic stacks, see \cite[\S 1.8]{KhanKstack} for an overview.} \footnote{We remind the reader that $\QCoh(\XX)$ denotes the prestable category of \emph{connective} quasi-coherent complexes. In \cite{BHL}, the authors consider the full stable category $\QCOH(\XX)$ of (nonconnective) quasi-coherent complexes, and show that the map analogous to \eqref{eq:MaptoTENSFUN} is fully faithful with essential image spanned by the tensor functors that preserve connective objects.}
\end{example}

\begin{example} We include an example which is unrelated to our aims in this paper, just to illustrate the utility of the above definition. Let $D^*=\spc(k((z)))$ denote the formal punctured disc. Let $$\mathbf{D}^*\in\stCCC$$ denote the assignment defined by  $$\mathbf{D}^*:\bigl(S=\spc(A)\bigr)\mapsto\bigl(\QCoh(S\widehat{\times}D^*):=\Mod(A((z)))\bigr).$$ Then it can be checked that for $\YY$ Tannakian we have $$\mathcal{L}\YY\simeq\TENSFUN(\stQCoh_\YY,\mathbf{D}^*),$$ where $\mathcal{L}\YY$ denotes the meromorphic loop space. We note that the corresponding formula in stacks would be $\Map(D^*,\YY)=\mathcal{L}\YY$, which, despite being \emph{morally true}, is strictly false. Indeed it is precisely the failure of $\mathbf{D}^*$ to be quasi-coherent that obstructs the above. More prosaically, the natural morphism $A((z))\otimes_A A'\rightarrow A'((z))$ is an isomorphism iff $A'$ is perfect as an $A$-module. We might also consider $D^*$ as an element of $\stCCC$, via the obvious embedding. Then the family of morphisms $$\Mod(A)\otimes\Mod(k((z)))\rightarrow\Mod(A((z)))$$ defines an embedding $D^*\rightarrow\mathbf{D}^*$, which is informally a discretization of the topological object $\mathbf{D}^*$.\end{example}

\subsubsection{The coevaluation functor}\label{sssec:coev}

For stacks $\XX$ and $\YY$, recall that the mapping stack is equipped with a canonical evaluation morphism
\[ \ev : \MAP(\XX,\YY) \times \XX \to \YY. \]

Given sheaves of tensor categories $\stCC$ and $\stDD$, we have a natural isomorphism
\[ \Map(\YY, \TENSFUN(\stCC,\stDD)) \simeq \TensFun_\YY(\stCC(\YY), \stDD(\YY))^\simeq \]
for all stacks $\YY$, by descent.
Therefore, the identity morphism $$\operatorname{id}:\TENSFUN(\stCC,\stDD)\rightarrow\TENSFUN(\stCC,\stDD)$$ classifies a universal $\QCoh(\TENSFUN(\stCC,\stDD))$-linear tensor functor
$$\mathrm{coev}:\stCC\bigl(\TENSFUN(\stCC,\stDD)\bigr)\longrightarrow \stDD\bigl(\TENSFUN(\stCC,\stDD)\bigr).$$

\begin{example} Let $\CC$ and $\DD$ be tensor categories. Then the coevaluation functor for $\stCC := \CC \otimes \stQCoh$ and $\stDD := \DD \otimes \stQCoh$ is tensored up from the tensor functor
\[ \mathrm{coev} : \CC \to \DD \otimes \QCoh(\TENSFUN(\stCC,\stDD)). \] \end{example}

\subsubsection{Pull-backs of functor stacks}

\begin{lemma}\label{fibre} Given a diagram of sheaves of $S$-linear tensor categories $$\stCC_{0}\leftarrow \stCC_{01}\rightarrow \stCC_{1},$$ the induced square \[\begin{tikzcd}
\TENSFUN(\stCC_{0}\otimes_{\stCC_{01}}\stCC_{1},\stDD)\arrow{r} \arrow[swap]{d} & \TENSFUN(\stCC_{0},\stDD) \arrow{d} \\
\TENSFUN(\stCC_{1},\stDD)\arrow{r} & \TENSFUN(\stCC_{01},\stDD)
\end{tikzcd}\] is cartesian for every sheaf of $S$-linear tensor categories $\stDD$. \end{lemma}

\begin{proof} Since $\stCCC_S$ is the category of commutative algebra objects in $\stCat_S$, \cite[Prop.~3.2.4.7]{LurieHA} implies that the tensor product $\stCC_{0}\otimes_{\stCC_{01}} \stCC_{1}$ is the coproduct in $\stCCC_S$ of the diagram $\stCC_{0}\leftarrow \stCC_{01}\rightarrow \stCC_{1}$. The claim follows. \end{proof}

Taking $\stCC_{01} = \stQCoh_{S}$, we get:

\begin{corollary}\label{cor:tensfun of coproduct}
  Given quasi-coherent sheaves of $S$-linear tensor categories $\stCC_1$ and $\stCC_2$, and a sheaf of $S$-linear tensor categories $\stDD$, the canonical map
  \[
    \TENSFUN_S(\stCC_1 \otimes_{S} \stCC_2, \stDD)
    \to \TENSFUN_S(\stCC_1, \stDD) \times \TENSFUN_S(\stCC_2, \stDD)
  \]
  is invertible.
\end{corollary}

\subsubsection{Tensor $1$-categories}

Recall from Subsect.~\ref{ssec:tensorcat} that for tensor $1$-categories $A$ and $B$, we have the category $$\TensFun(A, B)$$ of tensor functors $A \to B$. For an ordinary commutative $k$-algebra $R$, we denote by $B \otimes R := B \otimes \Mod^\heartsuit_R$ the Lurie tensor product with the tensor $1$-category of ordinary $R$-modules; this is still a tensor $1$-category by \cite[Thm.~C.5.4.16]{SAG}. The classical stack of tensor functors is defined by
$$\TENSFUN(A,B) : R \mapsto \TensFun(A, B\otimes R)^\simeq.$$
We also have the classical stack $\TENSFUN(\DD_{\geq 0}(A),\DD_{\geq 0}(B))_{\cl}$, the classical truncation of the stack of tensor functors $\DD_{\geq 0}(A) \to \DD_{\geq 0}(B)$.
There is a canonical morphism $$\TENSFUN(A,B)\longrightarrow\TENSFUN(\DD_{\geq 0}(A),\DD_{\geq 0}(B))_{\cl},$$ sending a tensor functor $f : A \to B$ to the induced left-derived functor $Lf : \DD_{\geq 0}(A) \to \DD_{\geq 0}(B)$ on connective derived categories.
More precisely, on $R$-points it sends a tensor functor $f : A \to B \otimes \Mod^\heartsuit_R$ to $Lf : \DD_{\geq 0}(A) \to \DD_{\geq 0}(B) \otimes \Mod_R$, for every commutative $k$-algebra $R$.

\begin{lemma}\label{abelian} Let $A$ and $B$ be tensor $1$-categories. Suppose that the dualizable objects of $A$ are compact and projective, and generate $A$ under colimits. Then the canonical morphism $$\TENSFUN(A,B)\longrightarrow\TENSFUN(\DD_{\geq 0}(A),\DD_{\geq 0}(B))_{\cl}$$ is invertible. Moreover, both are equivalently the stack whose $R$-points are symmetric monoidal $k$-linear functors $A^{\mathrm{dzbl}} \to (B \otimes R)^{\mathrm{dzbl}}$, where $(-)^{\mathrm{dzbl}}$ denotes the full subcategory of dualizable objects. \end{lemma}

\begin{proof} The claim is that for every ordinary commutative $k$-algebra $R$, the functor
\[ \TensFun(A, B \otimes \Mod^\heartsuit_R)^\simeq \to \TensFun(\DD_{\geq 0}(A), \DD_{\geq 0}(B) \otimes \Mod_R)^\simeq \]
is an equivalence. We will show that, under the given assumptions, it is inverse to the functor sending a tensor functor $F$ to the tensor functor $X \mapsto \operatorname{H}_0(F(X))$.

Since $A_0 := A^{\mathrm{dzbl}} \subset A$ is stable under finite coproducts, the assumptions imply that $A$ is freely generated by $A_0$ under sifted colimits, in the $1$-categorical sense, while the same holds for $\DD_{\geq 0}(A)$ but in the $\infty$-categorical sense (see e.g. \cite[\S 5.1]{CesScholze}).
Let $F : \DD_{\geq 0}(A) \to \DD_{\geq 0}(B) \otimes \Mod_R$ be a tensor functor. We claim that its restriction $F|_{A_0}$ takes values in $B \otimes \Mod^\heartsuit_R$ (which is the full subcategory of discrete objects in $\DD_{\geq 0}(B) \otimes \Mod_R$, see the proof of \cite[Thm.~C.5.4.16]{SAG}). If $X$ is a dualizable object of $A$, then the functor $X \otimes (-)$ is right adjoint to $X^\vee \otimes (-)$ and is in particular left-exact, so that $X$ is flat. In particular, $X \otimes^L (-) \simeq X \otimes (-)$ so that $X$ remains dualizable in $\DD_{\geq 0}(A)$. Since $F : \DD_{\geq 0}(A) \to \DD_{\geq 0}(B) \otimes \Mod_R$ was symmetric monoidal, it follows that $F(X)$ is also dualizable and hence flat by the same argument as above. In particular, $F(X)$ must be discrete. Thus, $F$ restricts to the functor
\[ F : A_0 \to B \otimes \Mod^\heartsuit_R. \]
We let $f : A \to B \otimes \Mod^\heartsuit_R$ be its unique colimit-preserving extension to $A$. The assignment $F \mapsto f$ is clearly natural in $F$, inverse to $f \mapsto Lf$, and satisfies $\operatorname{H}_0 (F(X)) \simeq \operatorname{H}_0 (Lf(X)) \simeq f(X)$ for every $X$ in $A$. The claim follows.

For the last assertion, note that any tensor functor $f : A \to B \otimes \Mod^\heartsuit_R$ preserves colimits and dualizable objects, hence is determined by its restriction $A_0 \to (B \otimes \Mod^\heartsuit_R)^{\mathrm{dzbl}}$, and conversely every symmetric monoidal $k$-linear functor $A_0 \to (B \otimes \Mod^\heartsuit_R)^{\mathrm{dzbl}}$ extends uniquely as above. \end{proof}

\subsubsection{Tensor localizations} We are interested in computing $\TENSFUN(C,D)$ \emph{locally} on $C$.

\begin{definition} If $\CC$ and $\CC'$ are tensor categories and $$j^{*}:\CC\rightarrow \CC'$$ is a tensor functor admitting a fully faithful right adjoint $j_{*}$, we call $j^{*}$ a \emph{tensor localization}.\footnote{The adjoint $j_{*}$ is of course not assumed to be symmetric monoidal.}\end{definition}

\begin{definition}\label{cover} Let $\CC$ be a tensor category. A (finite) \emph{cover} of $\CC$ is a (finite) jointly conservative family of tensor localizations $$\{j^{*}_{\alpha}: \CC\rightarrow \CC_{\alpha}\}_{\alpha\in A}.$$ \end{definition}

\begin{lemma}\label{Cart} Let $j^{*}: \CC\rightarrow \CC'$ be a tensor localization. Let $\iota: \stDD\rightarrow \stDD'$ be a morphism of sheaves of tensor categories such that $\iota(T)$ is fully faithful for every affine scheme $T$. Then the commutative square of stacks \[\begin{tikzcd}
\TENSFUN(\CC' \otimes \stQCoh, \stDD)\arrow{r} \arrow[swap]{d} & \TENSFUN(\CC' \otimes \stQCoh, \stDD') \arrow{d} \\
\TENSFUN(\CC \otimes \stQCoh, \stDD)\arrow{r} & \TENSFUN(\CC \otimes \stQCoh, \stDD')
\end{tikzcd}\] is cartesian. \end{lemma}

\begin{proof} We may check the Cartesian condition after evaluation at any affine $T$. Taking $\DD := \stDD(T)$ and $\DD' := \stDD'(T)$ and replacing $\CC$ and $\CC'$ by $\CC \otimes \QCoh(T)$ and $\CC' \otimes \QCoh(T)$, respectively, it is thus enough to show that for every fully faithful tensor functor $\iota: \DD \to \DD'$, the square of groupoids
\[\begin{tikzcd}
\TensFun(\CC', \DD)^\simeq\arrow{r} \arrow[swap]{d} & \TensFun(\CC', \DD')^\simeq \arrow{d} \\
\TensFun(\CC, \DD)^\simeq\arrow{r} & \TensFun(\CC, \DD')^\simeq
\end{tikzcd}\]
is cartesian.
Since $\iota$ is fully faithful, the horizontal arrows are injective on $\pi_0$ and bijective on all higher homotopy groups. It will thus suffice to show that the functor \[ \TensFun(\CC', \DD)^\simeq \to \TensFun(\CC, \DD)^\simeq \times_{\TensFun(\CC, \DD')^\simeq} \TensFun(\CC', \DD')^\simeq \] is essentially surjective and bijective on connected components of mapping groupoids. This amounts to the assertion that for any commutative square of solid arrows \[\begin{tikzcd}
\CC \arrow{r}{j^{*}} \arrow[swap]{d}{\phi} & \CC' \arrow{d}{\psi}\arrow[dashed,swap]{ld}{\chi} \\
\DD\arrow{r}{\iota} & \DD',
\end{tikzcd}\]
the dashed arrow can be filled in by an essentially unique tensor functor $\chi: \CC'\rightarrow \DD$. In fact, any such $\chi$ will automatically be a tensor functor: since $\iota$ is fully faithful, colimit-preserving and symmetric monoidal, this follows from the fact that $\iota \chi \simeq \psi j^*$ is a tensor functor.

First note that if $\chi$ is such a filling, then the commutativity $\chi j^* \simeq \phi$ and the adjunction counit provide a canonical isomorphism $\phi j_* \simeq \chi j^* j_* \simeq \chi$. It will thus suffice to show that $\chi := \phi j_{*} : \CC' \to \DD$ indeed fills the above diagram. The counit of the adjunction provides a natural transformation $$\iota\chi = \iota \phi j_{*} \simeq \psi j^{*} j_{*} \xrightarrow{\sim} \psi,$$ which is invertible because $j_{*}$ is fully faithful. Similarly, we claim that the natural transformation
\[ \phi \to \phi j_{*} j^* = \chi j^* \]
provided by the unit is also invertible. This may be checked after applying $\iota$ on the left: since $\iota$ is fully faithful and colimit-preserving, it admits a right adjoint $\iota^R$ with $\iota^R \iota \simeq \id$. Under the identification $\iota \phi \simeq \psi j^*$, $\iota \phi \to \iota \phi j_{*} j^*$ is inverse to the counit isomorphism $\psi j^* j_{*} j^* \to \psi j^*$. \end{proof}

\subsection{Tensor categorical points}

We will now introduce a special class of tensor categories which behave like ``points'' in the world of tensor categories.

\begin{definition} We say that a tensor category $\DD$ is an \emph{integral domain} if for non-zero objects $d_{0},d_{1}$ it holds that $d_{0}\otimes d_{1}$ is also non-zero. \end{definition}

\begin{example} The tensor category $\Mod_{R}$ is an integral domain iff $R$ is a field. If $\YY$ is a reduced stack whose underlying topological space is a point, then $\QCoh(\YY)$ is an integral domain. Indeed, there exists a flat surjection $\pi:\spc(\kappa)\rightarrow\YY$ for some field $\kappa$ (see \cite[Tag~06MN]{StacksProject}), and $\pi^{*}$ is a conservative tensor functor to the integral domain $\Mod_{\kappa}$. \end{example}

\begin{lemma}\label{local} If $\DD$ is an integral domain, and we are given a finite cover $\{j^{*}_{\alpha}:\CC\rightarrow \CC_{\alpha}\}$ of $\CC$, then any tensor functor $F : \CC \to \DD$ necessarily factors through one of the $\CC_{\alpha}$. Equivalently, there exists an $\alpha$ such that $F_\alpha := F j_{\alpha,*} : \CC_\alpha \to \DD$ is a tensor functor and the natural transformation $F \to F j_{\alpha,*}j_{\alpha}^* =F_\alpha j_{\alpha}^*$ is invertible. \end{lemma}
\begin{proof} Since $j_\alpha^*$ is a tensor localization, such a factorization $$F_\alpha : \CC_\alpha \to \DD$$ exists (as a tensor functor) if and only if $F_{\alpha} \simeq F j_{\alpha,*}$. We first show the existence of $F_\alpha$ as a colimit-preserving functor, and then check that it is endowed with a symmetric monoidal structure. Let $K_{\alpha}$ denote the kernel of the localization $j_{\alpha}^*$. Assuming that $F$ did not factor through any $j_{\alpha}^*$ we could then pick $d_{\alpha}$ in $K_{\alpha}$ so that $F(d_{\alpha})$ is non-zero. However $\bigotimes_{\alpha\in A}d_{\alpha}\simeq 0,$ as it vanishes upon localization at each fixed $\alpha$ and the cover is by definition jointly conservative. It follows since $F$ is symmetric monoidal that $$\bigotimes_{\alpha\in A}F(d_{\alpha}) \simeq F\bigl(\bigotimes_{\alpha\in A}d_{\alpha}\bigr) \simeq 0.$$  We obtain thus a contradiction as $\DD$ is an integral domain.

It remains to show that the induced functor $F_{\alpha} : \CC_{\alpha}\rightarrow \DD$ is symmetric monoidal. As the right adjoint of a symmetric monoidal functor, $j_{\alpha,*}$ is automatically lax symmetric monoidal and so we have natural maps $\mathbf{1}_{\CC} \to j_{\alpha,*}(\mathbf{1}_{\CC_\alpha})$ and $$j_{\alpha,*}(x)\otimes j_{\alpha,*}(y)\rightarrow j_{\alpha,*}(x\otimes y)$$ for objects $x$ and $y$. Since $F$ is symmetric monoidal, it will suffice to show that these maps become invertible after applying $F$.  Note that the cone of any of these maps automatically lies in $K_{\alpha}=\operatorname{ker}(j_{\alpha}^{*})$, as $j_{\alpha}^{*}$ is monoidal and $j_{\alpha}^{*}j_{\alpha,*}\simeq\id$. In particular, it vanishes after applying $F \simeq F_\alpha j_\alpha^*$. The claim follows. \end{proof}

\begin{remark} Let $\CC$ be an integral domain. Applying Lemma~\ref{local} to the identity functor, we see that for any finite cover $\{j^{*}_{\alpha}: \CC \rightarrow \CC_{\alpha}\}$, there exists an $\alpha$ such that $\id \to j_{\alpha,*}j_\alpha^*$ is invertible. But this implies that $j_\alpha^* : \CC \to \CC_\alpha$ is an equivalence. In particular, integral domains admit no nontrivial finite covers. More strongly, they admit no non-trivial tensor localizations. Equivalently, if $\CC$ is an integral domain then $\CC$ contains no non-trivial $\otimes$-idempotent algebras. In particular the locale formed by idempotent algebras (cf.~\cite{BKS}) is trivial. \end{remark}

\subsection{Spectra of tensor categories}

\begin{definition}
  The \emph{spectrum} of a sheaf of $S$-linear tensor categories $\stCC$ is the stack over $S$
  \[ \SpecTens_S(\stCC) := \TENSFUN_S(\stCC, \stQCoh_S).\]
\end{definition}

\begin{remark}
  The identity $\id_{\SpecTens(C)}$ is classified by a tautological tensor functor
  \[ \tau_{C} : C \to \QCoh(\SpecTens(C)). \]
\end{remark}

\begin{example}
  For every stack $\YY$ there is a canonical map
  \[ \begin{multlined}
    \YY \simeq \MAP(\spc(k), \YY)\\
    \qquad\to \TENSFUN(\QCoh(\YY), \Mod_k) = \SpecTens(\QCoh(\YY)),
  \end{multlined} \]
  which is invertible when $\YY$ is Tannakian.
  In particular, $\tau_{\QCoh(\YY)}$ is an equivalence.
\end{example}

\begin{example}
  We have $\SpecTens(\Mod_R) \simeq \spc(R)$ for every commutative ring $R$.
\end{example}

By Corollary~\ref{cor:tensfun of coproduct}, we have:

\begin{proposition}\label{prop:SpecTens monoidal}
  For any sheaves of $S$-linear tensor categories $\stCC$ and $\stDD$, the canonical map
  \[
    \SpecTens_S(\stCC \otimes_S \stDD) \to \SpecTens_S(\stCC) \times_S \SpecTens_S(\stDD)
  \]
  is invertible.
  In particular, $\SpecTens_S(-) : \CCC_S \to \St_S^{\mathrm{op}}$ is symmetric monoidal.
\end{proposition}

\section{Categorified circles}\label{sec:Cartier}

In our approach to elliptic cohomology, the key player is the categorified circle associated to an elliptic curve $E$, a tensor category $\HE$ that should be thought of as (the ind-completion of) the category of perfect complexes on the putative circle $S^1_E$ (cf. Subsect.~\ref{ssec:intro/circles}).
More generally, we will define $\HA$ for any commutative group scheme $A$, and $\HAS$ for a relative commutative group scheme $A \to S$.

One way to motivate the construction is as follows.
For an elliptic curve $E$, the Fourier--Mukai transform \cite{Mukai} on $\QCOH(E)$ exchanges the tensor product with the convolution product, given by the formula
$$ \mathcal{F} \star \mathcal{G} := m_*( \mathcal{F} \boxtimes \mathcal{G} ) $$
where $m$ is the group law on $E$.
However, convolution does not restrict to the connective subcategory $\QCoh(E)$ because $E$ is not affine (and hence $m_*$ is not right t-exact).
We thus introduce the subcategory $\QCoh_0(E)$ of sheaves with $0$-dimensional (= finite) support, which is the largest subcategory of $\QCoh(E)$ that supports the convolution product.

\subsection{Sheaves with $0$-dimensional support}\label{ssec:QCoh0}

We begin with a more general setup. Let $S$ be a stack and $\pi : X \to S$ a quasi-compact separated schematic morphism of stacks.

\begin{definition}\label{perf0} We denote by
\[\PERF_{0}(X/S) \subseteq \PERF(X) \]
the full subcategory spanned by perfect complexes supported on a cocompact\footnote{i.e., $Z$ is a closed subset of the underlying topological space of $X$ such that $X\setminus Z$ is quasi-compact} closed subset $Z \subseteq X$ which is finite over $S$, and by
\[ \QCOH_{0}(X/S) = \Ind \PERF_{0}(X/S) \]
its ind-completion.
We also denote by $\Perf_{0}(X/S) \subseteq \PERF_{0}(X/S)$ the full subcategory spanned by connective perfect complexes, and by $\QCoh_{0}(X/S) = \Ind \Perf_{0}(X/S)$ its ind-completion. \end{definition}

\begin{remark}
  We regard $\QCOH_0(X/S)$ as a full subcategory of $\QCOH(X)$ via the fully faithful functor
  \[
    \QCOH_0(X/S) = \Ind\PERF_0(X/S) \hookrightarrow \QCOH(X),
  \]
  which restricts to $\QCoh_0(X/S) \to \QCoh(X)$.
\end{remark}

\begin{remark}
  If $X$ is a qcqs scheme, then the canonical functor
  \[ \varinjlim_Z \QCOH(X~\mathrm{on}~Z) \to \QCOH_0(X/S) \]
  is an equivalence, where $Z \subseteq X$ ranges over cocompact closed subsets which are finite over $S$, and $\QCOH(X~\mathrm{on}~Z)$ denotes the full subcategory of quasi-coherent complexes $F$ such that $F|_{X\setminus Z}\simeq 0$, and the colimit is taken in the category of stable presentable categories and colimit-preserving functors.
  Indeed, writing $\PERF_0(X/S)$ as the union (= filtered colimit) of $\PERF(X~\mathrm{on}~Z)$, this follows from \cite[Lem.~5.3.2.9]{LurieHA} and the fact that the canonical functors
  \[ \Ind\PERF(X~\mathrm{on}~Z) \to \QCOH(X~\mathrm{on}~Z) \]
  are equivalences.
  The latter is a well-known result of Thomason, generalized to the derived setting in \cite[\S 5.2]{ToenAzumaya}.
  More generally, it suffices that $X$ is absolutely perfect in the sense of \cite[Def.~1.47]{KhanKstack}.
\end{remark}

\begin{example} When $S=\spc(k)$, we drop the decoration ``$/S$'' from the notation. Note that if $X$ is of finite type over $\spc(k)$, then $F \in \PERF(X)$ belongs to $\PERF_0(X)$ if and only if $F|_{X\setminus Z} \simeq 0$ for some closed subset $Z \subseteq X$ of dimension $\le 0$ (i.e., $0$-dimensional or empty). Note also that if $X$ is an integral curve over $\spc(k)$, then $\QCoh_0(X)$ is identified with the full subcategory of $\QCoh(X)$ spanned by \emph{torsion} quasi-coherent sheaves, i.e., the kernel of $\eta^* : \QCoh(X) \to \QCoh(\eta)$ where $\eta \in X$ is the generic point. If moreover $k$ is an algebraically closed field, then the category decomposes in $\Cat$ as $$\QCoh_0(X) \simeq \bigoplus_{x\in X(k)} \QCoh(X~\mathrm{on}~x).$$ \end{example}

\begin{lemma}\label{lem:push finite support}
Let $f : X \to Y$ be a quasi-compact quasi-separated schematic morphism and $Z\subseteq X$ a cocompact closed subset finite over $Y$.
Then
\[
  f_*:\QCOH(X~\mathrm{on}~Z)\to\QCOH(Y)
\]
is conservative and t-exact.
\end{lemma}
\begin{proof}
By base change for $f_*$, we may assume $Y$ is affine and hence $X$ is a qcqs scheme.
Since $f_*$ is left t-exact and has finite cohomological dimension by \cite[Cor.~3.4.2.3]{SAG}, it suffices to show that $f_*$ sends $\QCoh(X~\mathrm{on}~Z)^\heartsuit$ to $\QCoh(Y)^\heartsuit$.
For this, we may replace $X$ and $Y$ by their classical truncations.
Since $Z \subseteq X$ is cocompact, we may choose some subscheme structure which makes the inclusion finitely presented; let $I$ denote its ideal sheaf, which is thus of finite type.
Given $F' \in \QCoh(X~\mathrm{on}~Z)^\heartsuit$, we may write $F'$ as the filtered colimit of $i_{n,*}(F'_n)$ where $F'_{n}$ is the subsheaf of $F'$ annihilated by $I^n$, and $i_n : Z_n \hookrightarrow X$ is the inclusion of the $n$th infinitesimal neighbourhood of $Z$.
Since $f_*$ preserves colimits, we have $f_*(F') \simeq \varinjlim_n (f|_{Z_n})_*(F'_n)$.
Since filtered colimits are t-exact, and each $(f|_{Z_n})_*$ is t-exact because $f|_{Z_n}$ is affine, we find that $f_*(F')$ is discrete.

The same argument shows conservativity on the heart: if $f_*(F')=0$, left t-exactness implies $(f|_{Z_n})_*(F'_n)=0$ for every $n$, hence $F'_n=0$ because $f|_{Z_n}$ is affine, and hence $F'=0$.
By the t-exactness just proven above, this then implies conservativity on the entire category.
\end{proof}

\begin{proposition}\label{prop:push on QCoh0}
  Let $\pi_X: X \to S$ and $\pi_Y: Y \to S$ be qc separated schematic morphisms and $f: X \to Y$ an $S$-morphism.
  If $f$ is almost of finite presentation and of finite Tor-amplitude, then $f_* : \QCOH(X) \to \QCOH(Y)$ restricts to colimit-preserving functors
  \[
    f_* : \QCOH_0(X/S) \to \QCOH_0(Y/S)
  \]
  and $f_* : \QCoh_0(X/S) \to \QCoh_0(Y/S)$.
\end{proposition}
\begin{proof}
  By the universal property of ind-completion it is enough to show that $f_* : \QCOH(X) \to \QCOH(Y)$ sends $\PERF_0(X/S)$ to $\PERF_0(Y/S)$ and $\Perf_0(X/S)$ to $\Perf_0(Y/S)$.

  Let $F \in \PERF(X)$ be supported on a closed subset $Z \subseteq X$ which is finite over $S$.
  Because $\pi_Y$ is separated and $\pi_Y \circ f|_Z : Z \to S$ is finite, it follows that $f|_Z : Z \to Y$ is finite and hence $f(Z) \subseteq Y$ is a closed subset finite over $S$.
  By base change, $f_*(F)|_{Y\setminus f(Z)} \simeq 0$ follows from $F|_{X\setminus f^{-1}(f(Z))} \simeq 0$.
  To conclude, it remains to show:
  \begin{enumerate}
    \item The quasi-coherent complex $f_*(F) \in \QCOH(Y)$ is perfect.
    \item If $F$ is connective, then $f_*(F)$ is connective.
  \end{enumerate}
  By base change for $f_*$, we may assume $Y$ is affine and hence $X$ is a quasi-compact separated scheme.
  Lemma~\ref{lem:push finite support} gives (ii).

  It remains to prove (i).
  Since $f$ is almost of finite presentation, $f|_Z$ is proper, and $F \in \PERF(X~\mathrm{on}~Z)$, we may apply \cite[Prop.~5.6.5.2]{SAG} to deduce that $f_*(F)$ is pseudo-coherent.
  It remains to show that $f_*(F)$ is of finite Tor-amplitude, i.e., that $f_*(F) \otimes G \simeq f_*(F \otimes f^*(G))$ is cohomologically bounded for any $G \in \QCoh(Y)^\heartsuit$.
  Since $f$ is of finite Tor-amplitude, $f^*(G)$ is cohomologically bounded.
  Since $F$ is perfect and lies in $\QCOH(X~\mathrm{on}~Z)$, $F \otimes f^*(G)$ is thus cohomologically bounded and lies in $\QCOH(X~\mathrm{on}~Z)$.
  Thus by Lemma~\ref{lem:push finite support}, $f_*(F \otimes f^*(G))$ is cohomologically bounded as desired.
\end{proof}

\begin{proposition}\label{prop:QCOH0Kunneth}
  Let $\pi_X: X \to S$ and $\pi_Y: Y \to S$ be separated schematic morphisms.
  Assume $S$ is perfect and $X$ and $Y$ are absolutely perfect (in the sense of \cite[Def.~1.47]{KhanKstack}).
  Then the canonical functor
  \[ \QCOH_0(X/S) \otimes_{\QCOH(S)} \QCOH_0(Y/S) \to \QCOH_{0}(X \times_S Y/S) \]
  is an equivalence.
  Moreover, it is natural with respect to push-forwards (Proposition~\ref{prop:push on QCoh0}) in $X$ and $Y$.
\end{proposition}
\begin{proof}
  This functor is the ind-completion of the canonical functor
  \[ \PERF_0(X/S) \otimes_{\PERF(S)} \PERF_0(Y/S) \to \PERF_{0}(X \times_S Y/S) \]
  given by $F \otimes G \mapsto F \boxtimes G$.
  This is well-defined because if $Z \subseteq X$ and $Z' \subseteq Y$ are closed subsets finite over $S$ such that $F|_{X\setminus Z} \simeq 0$ and $G|_{Y\setminus Z'} \simeq 0$, then $Z\times_S Z' \subseteq X\times_S Y$ is a closed subset finite over $S$, and $(F \boxtimes G)|_{X\times_S Y \setminus Z\times_S Z'} \simeq 0$.

  Since the tensor product commutes with colimits in each argument, the left-hand side is the colimit of the categories
  \[ \PERF(X~\mathrm{on}~Z) \otimes_{\PERF(S)} \PERF(Y~\mathrm{on}~Z') \]
  over closed subsets $Z \subseteq X$ and $Z' \subseteq Y$ which are finite over $S$.
  On the other hand, the right-hand side can be identified by a cofinality argument with the colimit of $\PERF(X \times_S Y~\mathrm{on}~Z\times_S Z')$.
  Indeed, for any closed subset $W \subseteq X \times_S Y$ which is finite over $S$, the closed subsets $Z := \mathrm{pr}_1(W) \subseteq X$ and $Z' := \mathrm{pr}_2(W) \subseteq Y$ are finite over $S$, and we have $W \subseteq Z \times_S Z'$; the separatedness assumptions on $\pi_X$ and $\pi_Y$ guarantee that $\mathrm{pr}_1|_W : W \to X$ and $\mathrm{pr}_2|_W : W \to Y$ are proper.

  It remains to show that the functor
  \[\begin{multlined}
    \QCOH(X~\mathrm{on}~Z) \otimes_{\QCOH(S)} \QCOH(Y~\mathrm{on}~Z')\\
    \to \QCOH(X \times_S Y~\mathrm{on}~Z \times_S Z')
  \end{multlined}\]
  is an equivalence for all $Z$ and $Z'$; since $S$ is perfect and $X$ and $Y$ are absolutely perfect, the equivalence on perfect objects will follow.
  Since the Lurie tensor product preserves Verdier sequences in each argument, this reduces to the version without supports proven in \cite[Thm.~1.2(i)]{BFN}.

  The naturality of the functor with respect to push-forwards follows from the invertibility of the canonical functors
  \[ f_*(-) \boxtimes g_*(-) \to (f \times g)_*(- \boxtimes -) \]
  on $\QCOH(X) \times \QCOH(Y)$, where $f : X \to X'$ and $g : Y \to Y'$ are $S$-morphisms, which is a standard application of the projection and base change formulas.
\end{proof}

\begin{corollary}\label{cor:QCoh0 fun} For any Artin stack $S$, the assignment $$X \mapsto \QCoh_0(X/S)$$ determines a symmetric monoidal functor from the category of Artin $S$-stacks $X \to S$ which are smooth, schematic and of finite presentation, to the category $\Cat_S$ of $S$-linear Grothendieck prestable categories (and $S$-linear colimit-preserving functors). \end{corollary}
\begin{proof}
  This follows from \cite[Ex.~C.1.4.4]{SAG}, Proposition~\ref{prop:push on QCoh0}, and Proposition~\ref{prop:QCOH0Kunneth}.
\end{proof}
\begin{corollary} For every smooth schematic morphism of finite presentation $X \to S$, the object $\QCoh_0(X/S) \in \Cat_S$ is equipped with a canonical cocommutative comonoid structure. \end{corollary}
\begin{proof} Follows from Corollary~\ref{cor:QCoh0 fun} in view of the tautological cocommutative comonoid structure on $X \to S$, with comultiplication defined by the diagonal $\Delta : X \to X \times_S X$ and counit defined by the projection $X \to S$. \end{proof}

\begin{remark} Informally speaking, $\QCoh_0(X/S)$ is a ``categorical cocommutative coalgebra'', with comultiplication
\[\begin{multlined}
\QCoh_0(X/S)\\
\xrightarrow{\Delta_*} \QCoh_0(X \times_S X/S)
\xrightarrow{\boxtimes^{-1}} \QCoh_0(X/S) \otimes_{\QCoh(S)} \QCoh_0(X/S),
\end{multlined}\]
and counit $\pi_{*} : \QCoh_0(X/S) \to \QCoh_0(S/S) \simeq \QCoh(S)$. \end{remark}

\begin{remark} The construction $\QCoh_0(X/S)$ does not satisfy base change along morphisms $T \to S$. That is, the canonical tensor functor
\[ \phi : \QCoh_0(X/S)\otimes_S \QCoh(T) \to \QCoh_0({X \times_S T/T}), \quad F \otimes G \mapsto F \boxtimes_S G\]
is typically not an equivalence. For example, take $S=\spc(k)$ and $T$ a qcqs scheme. Then $\phi$ is the ind-completion of a fully faithful functor
\[ \Perf_0(X) \otimes \Perf(T) \hookrightarrow \Perf_0(X\times T/T) \]
which can be identified with an inclusion of two full subcategories of $\Perf(X\times T)$:
\begin{itemize}
\item
The source is the full subcategory spanned by perfect complexes $F$ such that $F|_{X\times T\setminus Z\times T} \simeq 0$ for some closed subset $Z \subseteq X$ finite over $\spc(k)$.
This follows from the commutativity of tensor products with colimits in each argument, and the K\"unneth formula of \cite[Thm.~1.2(1)]{BFN}.

\item
The target is the full subcategory of $\Perf(X\times T)$ spanned by perfect complexes $F$ such that $F|_{X\times T\setminus W} \simeq 0$ for some $W \subseteq X \times T$ finite over $T$.
\end{itemize}
Take $T=X'$ to be another copy of $X$. Then the target of $\phi$ contains the structure sheaf of the diagonal, $\OO_{\Delta} \in \Perf(X \times X'~\mathrm{on}~\Delta)$. However, $\OO_{\Delta}$ does not lie in the source of $\phi$ since $\Delta$ is not contained in $Z \times X$ for any $Z \subseteq X$ finite over $\spc(k)$.

Note that by Corollary~\ref{cor:HA (x) QCoh} the inclusion $\QCOH_0(X) \hookrightarrow \QCOH(X)$ corresponds to an object $F_\Delta \in \QCOH_0(X) \otimes \QCOH(X)$.
Under the above description, it is given by the colimit
\[ F_\Delta \simeq \varinjlim_{Z\subseteq X} \Gamma_{Z\times X}(\OO_{\Delta}). \]
This object is neither compact nor invertible. \end{remark}

\begin{lemma}\label{lem:duality with proper supports}
  Let $p : X \to S$ be a smooth morphism, separated schematic of finite type, and let $Z\subseteq |X|$ be a cocompact closed subset whose reduced closed substack is proper over $S$.
  There are natural isomorphisms
  \[
    \Hom(p_*P, M) \simeq \Hom(P, p^*M \otimes \omega_{X/S})
  \]
  for $P\in\QCOH(X~\mathrm{on}~Z)$ and $M\in\QCOH(S)$.
\end{lemma}
\begin{proof}
  See \cite[Thm.~3.2.16, Rem.~3.2.17, Cor.~3.2.18]{JiangThesis} (for which the standing truncatedness assumption in \emph{op. cit.} is not required).
\end{proof}

\subsection{The tensor category \texorpdfstring{$\HA$}{H\_A}}
\label{ssec:HA}

We now define, for a commutative group scheme $A$, the categorified circle $\HA$.
This will be the formal definition of $\IndPerf(S^1_A)$ for the putative circle $S^1_A$.

\subsubsection{Commutative group schemes}\label{sssec:comm grps}

Fix a stack $S$ over $k$.

\begin{definition}\label{def:comm grp scheme}
  Let $\mathcal{C}_{S}$ be the category of stacks $X$ over $S$ such that $X \to S$ is smooth, separated, quasi-compact, and schematic.
	A \emph{commutative group scheme} over $S$ is an abelian group object of $\mathcal{C}_S$ in the sense of \cite[Def.~1.2.4]{LurieEll1}.
\end{definition}

Informally speaking, a commutative group scheme $\pi : A \to S$ is equipped with an operation $m : A \times_S A \to A$, an inversion map $\iota : A \to A$ over $S$, and an identity section $e : S \to A$, all satisfying the axioms of a commutative group up to coherent homotopy.
Moreover, the $T$-points $A(T)$ form a topological abelian group up to homotopy for every $T$ over $S$, rather than merely a grouplike $E_\infty$-space (see \cite[\S1, Rem.~1.2.10]{LurieEll1}).

We write $[n]:A\rightarrow A$ for scalar multiplication by $n\in\mathbf{Z}$, and $A[n] := \operatorname{fib}([n] : A \to A)$ for the subgroup of $n$-torsion points, which is finite locally free over $S$ whenever $[n]$ is.

\subsubsection{}\label{sssec:fra}

We write
\[ \fra := \omega_{S/A}[d], \qquad \frav := \fra^{\otimes -1} \simeq \omega_{S/A}^\vee[-d], \]
where $\omega_{S/A} := e^!(\OO_A)$ is the relative dualizing complex of the identity section $e : S \to A$, an invertible object of $\QCOH(S)$ sitting in degree $d$, where $d$ is the relative dimension of $A \to S$.
By \cite[Cor.~6.4.2.8]{SAG}, we have equivalently
\[
  \fra \simeq e^*\omega_{A/S}^\vee[d],
  \qquad
  \frav \simeq e^*\omega_{A/S}[-d],
\]
where $\omega_{A/S}$ is the relative dualizing complex of $A \to S$.
When $k$ is a $\bQ$-algebra, $\fra$ coincides with the determinant of the tangent at the identity, $\det(e^*T_{A/S})$. (In general, $e^*T_{A/S}$ is not perfect, as $A \to S$ is not \emph{differentially smooth} in the sense of \cite[Def.~11.2.2.2]{SAG}.)

\begin{definition}\label{def:orientation}
	An \emph{orientation} of $A$ over $S$ is an isomorphism $\fra \simeq \OO_S$.
\end{definition}

\begin{remark}
There is a canonical isomorphism $m^*\omega_{A/S}\simeq\pr_1^*\omega_{A/S}$.
(Indeed, factor $m = \pr_1\circ t : A \times_S A \to A$, where $t(a,b)=(ab,b)$, and use the identities $t^*\omega_{\pr_2}\simeq\omega_{\pr_2}$ and $\omega_{\pr_2}\simeq\pr_1^*\omega_{A/S}$.)
Pulling back along $e\times\operatorname{id}:A\to A\times_S A$, this gives in particular
\begin{equation}\label{eq:omegaAS}
  \omega_{A/S}\simeq\pi^*e^*\omega_{A/S}\simeq\pi^*\frav[d].
\end{equation}
See also \cite[Prop.~5.2.3]{LurieEll1}.
\end{remark}

\subsubsection{}\label{sssec:HA Hopf}

In the cartesian monoidal category of stacks over $S$, $A \to S$ is a commutative monoid object via $m : A\times_S A \to A$ and a cocommutative comonoid object via $\Delta : A \to A \times_S A$. These structures are compatible: $\Delta$ is a morphism of monoids because the operation $m$ is commutative. The inversion map $\iota : A \to A$ defines an antipode. In summary, $A$ is a commutative and cocommutative Hopf algebra (in the category of stacks over $S$).
By Corollary~\ref{cor:QCoh0 fun}, we have:

\begin{corollary}\label{cor:QCoh0 Hopf}
  The $S$-linear category $\QCoh_0(A/S)$ is equipped with the structure of an $S$-linear \emph{Hopf tensor category}, i.e., an $S$-linear category together with the structure of an abelian cogroup object in $\CCC_S$.
\end{corollary}

\begin{definition}
  The \emph{categorified circle} associated to the commutative group scheme $A$ over $S$ is the sheaf of $S$-linear tensor categories $\stHAS \in \stCCC_S$ defined by the assignment
  $$ T \mapsto \Hcat_{A\times_S T/T} := \QCoh_0(A \times_S T/T). $$
  By Corollary~\ref{cor:QCoh0 Hopf}, it carries the structure of a cogroup object in $\stCCC_S$.
\end{definition}

Note also that for any morphism of commutative group schemes $f : A \to A'$ over $S$, the pointwise pushforward $f_* : \stHAS \to \stHcat_{A'/S}$ is an $S$-linear tensor functor compatible with Hopf algebra structures.

\begin{remark}
  As we have seen in Subsect.~\ref{ssec:QCoh0}, the canonical map $\HAS \otimes_S \stQCoh_S \to \stHAS$ is typically not an equivalence even on underlying categories.
  It is this discrepancy which forces us to work in the generality of sheaves of tensor categories.
\end{remark}

We may spell out the Hopf algebra structure on $\HAS$ as follows.
To simplify notation, we take $S=\spc(k)$.
\begin{itemize}
  \item The underlying category is $\QCoh_0(A)$.
  \item The \emph{multiplication map} $\Hcat_{A\times A} \simeq \HA \otimes \HA \to \HA$ is (under the K\"unneth equivalence of Proposition~\ref{prop:QCOH0Kunneth}) the tensor functor $m_* : \QCoh_0(A \times A) \to \QCoh_0(A)$, where $m$ is a morphism of groups since $A$ is commutative.
  \item The \emph{unit map} $\eta : \Mod_k \to \HA$ is the tensor functor $$e_* : \QCoh_0(\spc(k)) \to \QCoh_0(A).$$
  \item The \emph{comultiplication map} $\nabla : \HA \to \HA \otimes \HA \simeq \Hcat_{A\times A}$ is the tensor functor $\Delta_{A,*} : \QCoh_0(A) \to \QCoh_0(A \times A)$, where $\Delta_A : A \to A \times A$ is the diagonal.
  \item The \emph{counit map} $\varepsilon : \HA \to \Mod_k$ is the tensor functor $\pi_* : \QCoh_0(A) \to \QCoh_0(\spc(k))$.
  \item The \emph{antipode map} $s : \HA \to \HA$ is the tensor functor $\iota_* : \QCoh_0(A) \to \QCoh_0(A)$.
\end{itemize}

In particular, the symmetric monoidal structure on $\HAS$ is defined by convolution:
  $$F\star G:=m_{*}(F \boxtimes_S G).$$
The monoidal unit is $\mathcal{O}_{e}$, the structure sheaf at the origin.

\begin{remark} For any two points $a,b \in A(S)$, we have $\OO_a \star \OO_b \simeq \OO_{m(a,b)}$. (By way of contrast, recall that $\mathcal{O}_{a}\otimes\mathcal{O}_{b}=0$ when $a\ne b$.) \end{remark}

\begin{example}
  Let $S=\spc(k)$ where $k$ is an algebraically closed field and let $A$ be a $1$-dimensional commutative group scheme over $S$.
  Denote by $\Hcat_{A,0}$ the tensor category $\QCoh(A~\mathrm{on}~e)$ equipped with the convolution product, and denote by $\grMod^{A(k)}_{k}$ the tensor category of $A(k)$-graded objects in $\Mod_k$.
  Then the equivalence $\QCoh_0(A) \simeq \bigoplus_{a\in A(k)} \QCoh(A~\mathrm{on}~a)$ in $\Cat$ enhances to an equivalence in $\CCC$
  $$\Hcat_{A,0} \otimes \grMod^{A(k)}_{k} \to \HA, \qquad F \otimes k\{a\} \mapsto F \star \OO_a$$
  where $k\{a\}$ denotes the $1$-dimensional $k$-vector space placed in degree $a$.
\end{example}

\subsubsection{The additive and multiplicative cases}

Let $A=\GGaS \to S$ be the additive group over $S$. Then $\Hcat_{\GGaS/S}$ is equivalent to the symmetric monoidal category of pairs $(V, T)$ where $V \in \IndPerf(S)$ and $T$ is an endomorphism of $V$. The tensor product is given by $$(V,T_V) \otimes (W,T_W) = (V\otimes W, T_V \otimes \id + \id \otimes T_W).$$

For $A=\GGmS \to S$ the multiplicative group over $S$, $\Hcat_{\GGmS/S}$ is equivalent to the symmetric monoidal category of pairs $(V, T)$ where $V \in \IndPerf(S)$ and $T$ is an automorphism of $V$. The tensor product is given by $$(V,T_V) \otimes (W,T_W) = (V\otimes W, T_V \otimes T_W).$$

\subsubsection{Global sections}

Regarding $\HAS$ as $\IndPerf(S^1_A)$ for the putative circle $S^1_A$, one expects a functor of global sections $\Gamma_{S^1_A}$.\footnote{%
  More precisely, this would be ``renormalized'' global sections: the ind-extension of global sections on $\Perf(S^1_A)$.
}
To define it, recall the unit of the tensor structure on $\HAS$ (\ref{sssec:HA Hopf}), i.e., the unique $\QCOH(S)$-linear symmetric monoidal colimit-preserving functor $$\eta:\QCOH(S)\rightarrow\QCOH_{0}(A/S), \quad V \mapsto e_*V.$$ Its right adjoint
$$\Gamma_{S^1_A} : \QCOH_{0}(A/S) \rightarrow \QCOH(S), \quad F \mapsto \Hom_{\QCOH(A)}(\OO_{e}, F)$$
restricts along $\HAS = \QCoh_0(A/S) \subseteq \QCOH_0(A/S)$ to
$$\Gamma_{S^1_A} : \HAS\rightarrow\QCOH_{0}(A/S)\xrightarrow{\Gamma_{S^1_A}}\QCOH(S).$$

\begin{remark} Note that $\eta$ restricts to a tensor functor $$\eta:\QCoh(S) \rightarrow\HAS.$$ In contrast, $\Gamma_{S^1_A}$ does not necessarily preserve connective objects. \end{remark}

\subsubsection{Integrality of $\HA$}

Let $A$ be a commutative group scheme over $S=\spc(K)$ for some field extension $K/k$.

\begin{lemma}\label{HE int} The tensor category $\Hcat_{A/K}$ is an integral domain.\end{lemma}
\begin{proof} Recall that non-zero objects of $\Hcat_{A/K}$ are quasi-coherent sheaves whose cohomology sheaves have $0$-dimensional support. Given two such $F$ and $G$, the claim is that $F \star G$ is non-zero. Since the convolution product is left-exact in each argument (by left-exactness of pull-back along the projections $A \times A \to A$ and push-forward along $m : A \times A \to A$), we have $H^i(F \star G) \simeq H^i(F) \star H^i(G)$. We may thus assume that $F$ and $G$ are discrete. Since they have $0$-dimensional support, we may find points $x,y \in A$ and injections $\OO_{x}\rightarrow F$ and $\OO_{y}\rightarrow G$. These give rise to an injection $\OO_{m(x,y)} \rightarrow F \star G$ (using left-exactness of convolution again), so $F \star G$ is nonzero as claimed. \end{proof}

\subsection{Dualizable objects of \texorpdfstring{$\HA$}{H\_A}}
We study duality internally to the tensor category $\HA$.
This analysis will also yield a proof of the following external self-duality:

\begin{theorem}\label{thm:Hcat self dual} The stable category $\QCOH_0(A/S)$ is self-dual in the symmetric monoidal $\infty$-category $\operatorname{Pr}^{\mathrm{L,st}}_S$ of stable presentable $S$-linear categories. \end{theorem}

\begin{remark}
The self-duality is a stable phenomenon: the convolution dual $\mathbb{D}^\star F \simeq \iota^*\mathbb{D} F\otimes\pi^*\frav[\dim(A)]$ preserves $\PERF_0(A/S)$ but does \emph{not} preserve the connective subcategory $\Perf_0(A/S)$ in general: for example, $$\mathbb{D}^\star(\OO_a[1]) \simeq \OO_{\iota(a)}[-1].$$
\end{remark}

By definition of the Lurie tensor product, this implies:

\begin{corollary}\label{cor:HA (x) QCoh}
  For every Artin stack $X$, there is a canonical equivalence in $\operatorname{Pr}^{\mathrm{L,st}}_S$
  \[
    \QCOH_0(A/S) \otimes_S \QCOH(X)
    \simeq \operatorname{Fun}^{\mathrm{L}}_S(\QCOH_0(A/S), \QCOH(X)),
  \]
  contravariantly functorial in $X$, where the right-hand side denotes the category of $S$-linear colimit-preserving functors.
\end{corollary}

The convolution product $\star$ commutes with colimits in each argument separately, and so we have internal hom objects $$\hhoms_{A}^{\star}(F,G) \in \HAS,$$ defined by natural isomorphisms $$\Hom(F\star G,H)\simeq\Hom(F,\hhoms_{A}^{\star}(G,H)).$$

Since the base $S$ is generally left implicit, we write $\Hom_{\CC}(-,-)$ for the $\OO_S$-linear mapping object of an $\OO_S$-linear category $\CC$, and reserve the underlined $\hhoms$ for internal homs relative to a nontrivial base.
We note that we have decorated internal hom objects in $\HAS$ with a superscript $\star$ to emphasize that the monoidal product is convolution. We will write $\hhoms_{A}$ for the internal hom object with respect to $\otimes.$ If $F$ is an object of $\QCoh(A),$ we denote its pre-dual object $$\mathbb{D}F:=\hhoms_{A}(F,\mathcal{O}).$$ If $F$ is an object of $\HAS$, then we denote its pre-dual object $$\mathbb{D}^{\star}F:=\hhoms_{A}^{\star}(F,\mathcal{O}_{e}).$$ We remind the reader that in a closed symmetric monoidal category, if an object $F$ is dualizable then its pre-dual is necessarily also a dual.

\begin{lemma} There are functorial isomorphisms $$\hhoms_{A}^{\star}(F,G)\simeq \pi_{1,*}\,\hhoms_{A^{2}}(\pi_{2}^{*}\,F,m^{\times}G),$$ where $m^{\times}$ denotes a right adjoint\footnote{Note that $m$ is qcqs schematic since $A \to S$ was assumed qcqs schematic.} to $m_*$. \end{lemma}

\begin{proof} We have \begin{align*}
  \Hom(H\star F,G) &= \Hom\bigl(m_{*}(\pi_{1}^{*}H\otimes\pi_{2}^{*}F),G\bigr)\\
  &\simeq \Hom(\pi_{1}^{*}H\otimes\pi_{2}^{*}F,m^{\times}G)\\
  &\simeq \Hom\bigl(\pi_{1}^{*}H,\hhoms_{A^{2}}(\pi_{2}^{*}F,m^{\times}G)\bigr)\\
  &\simeq \Hom\bigl(H,\pi_{1,*}\hhoms_{A^{2}}(\pi_{2}^{*}F,m^{\times}G)\bigr)
\end{align*} for all $H$. \end{proof}

\begin{corollary}\label{duals} Let $A$ be a commutative group scheme over $S$ as in (\ref{sssec:comm grps}). An object $F \in \QCOH_0(A/S)$ is dualizable with respect to convolution if and only if it belongs to $\PERF_{0}(A/S)$, i.e. if and only if it is dualizable as an object of $\QCOH(A)$. In that case the dual is given by $$\mathbb{D}^{\star}F\simeq \iota^{*}\,\mathbb{D}F \otimes \pi^{*}\frav[\operatorname{dim}(A)].$$\end{corollary}

\begin{proof}
Suppose $F \in \QCOH_0(A/S)$ is $\star$-dualizable.
To show that $F$ is perfect, we may assume $S$ is affine.
The unit $\OO_e$ belongs to $\PERF_0(A/S)$ and is therefore compact.
Since $\Hom(F,-) \simeq \Hom(\OO_e, \mathbb{D}^\star F\star(-))$ and convolution preserves colimits, $F$ is compact and hence perfect.

Conversely, let $F \in \PERF_0(A/S)$ with support contained in a cocompact closed subset $Z \subseteq A$ finite over $S$.
Using the shearing automorphism
\[
  \varsigma : A\times_S A \to A\times_S A,
  \qquad (a,b) \mapsto (a,ab),
\]
we can write
\[
  F\star G \simeq \pr_{2,*}(\pr_1^*F\otimes q^*G),
\]
where $q := \pr_2 \circ \varsigma^{-1} : A \times_S A \to A$.
Note that $\pr_2$ is finite on $Z\times_S A$, and $\pr_1^*F\otimes q^*G$ vanishes away from $Z\times_S A$.
We may thus apply Lemma~\ref{lem:duality with proper supports} to obtain, for all $G,H\in\QCOH(A)$,
\begin{align*}
  \Hom(F\star G,H)
  &\simeq \Hom(\pr_1^*F\otimes q^*G,
               \pr_2^*H\otimes\pr_1^*\omega_{A/S})\\
  &\simeq \Hom(q^*G,
               \pr_1^*(\mathbb{D}F\otimes\omega_{A/S})
               \otimes\pr_2^*H)\\
  &\simeq \Hom\bigl(G,
               q_*(\pr_1^*(\mathbb{D}F\otimes\omega_{A/S})
               \otimes\pr_2^*H)\bigr)\\
  &\simeq \Hom(G,E\star H),
\end{align*}
where $E := \iota^*(\mathbb{D}F\otimes\omega_{A/S})$ and the last isomorphism uses $q=m\circ(\iota\times\operatorname{id})$.
In other words, $L := F\star(-)$ is left adjoint to $R := E\star(-)$.
By \eqref{eq:omegaAS}, we have $E \simeq \iota^*\mathbb{D}F\otimes\pi^*\frav[d]$.

Let $\beta_{H,K} : R(H)\star K\to R(H\star K)$ be adjoint to $\varepsilon_H\star\operatorname{id}_K$, where $\varepsilon: LR \to \operatorname{id}$ is the counit.
To show that $E$ is a convolution dual of $F$, it remains to show that $\beta$ is invertible.
We may assume $S$ affine.
By K\"unneth, the adjunction $(L,R)$ gives rise to an adjunction $(L', R') := (L\otimes\operatorname{id}, R\otimes\operatorname{id})$ on $\QCOH(A\times_S A)$.
By construction, there is a canonical isomorphism $\pr_{1,*} \circ R' \simeq R \circ \pr_{1,*}$.
Writing $m = \pr_1\circ t : A \times_S A \to A$, where $t(a,b)=(ab,b)$, we also have $L' \circ t_* \simeq t_* \circ L'$, and hence by passing to right adjoints (since $t$ is an isomorphism), $t_* \circ R' \simeq R' \circ t_*$.
Their composite $m_* \circ R' \simeq R \circ m_*$ evaluates to $\beta_{H,K}$ on $H\boxtimes K$, whence the claim.
\end{proof}

\begin{remark} For a point $a\in A(S)$, it is not hard to show that the above formula reduces to $$\mathbb{D}^{\star}\mathcal{O}_{a}\simeq\mathcal{O}_{\iota(a)}.$$  \end{remark}

\begin{proof}[Proof of Theorem~\ref{thm:Hcat self dual}]
For an $S$-linear small stable $\infty$-category $C$, the dual of $\operatorname{Ind}(C)$ in $\operatorname{Pr}^{\mathrm{L,st}}_S$ is $\operatorname{Ind}(C^{\operatorname{op}})$.
By Corollary~\ref{duals}, $\mathbb{D}^{\star}$ determines an equivalence $\PERF_0(A/S)^{\operatorname{op}}\simeq\PERF_0(A/S)$.
Passing to ind-completions proves the claim.
\end{proof}

\subsection{Invertible objects of \texorpdfstring{$\HA$}{H\_A}}
\label{ssec:Pic HA}

We study the stack of $\star$-invertible objects $\PIC(\stHcat_{A/S})$, whose $T$-points form the $\infty$-groupoid $\Pic(\Hcat_{A_T/T})$.

\subsubsection{}

There is a canonical homomorphism of $E_\infty$-group stacks over $S$,
\begin{equation}\label{eq:A to PicH}
  A \to \PIC(\stHcat_{A/S}),
\end{equation}
which for every $T$ over $S$ sends $a : T \to A$ to the $\star$-invertible object $a_*\OO_T \in \Hcat_{A_T/T}$.
This is classified by $\Delta_{A,*} \OO_{A} \in \Hcat_{A\times_S A/A}$.
Over $A \times_S A$, we have
\[
  \pr_1^*(\Delta_\OO) \star \pr_2^*(\Delta_\OO)
  \simeq
  m_A^*(\Delta_\OO).
\]
This isomorphism and all its associativity, symmetry, and unit coherences are induced by the commutative group structure of $A$, so \eqref{eq:A to PicH} is a homomorphism of $E_\infty$-group stacks.

\subsubsection{}

The morphism of sheaves of tensor categories $e_* : \stQCoh_S \to \stHAS$ induces a homomorphism of $E_\infty$-group stacks over $S$,
\begin{equation}\label{eq:PIC to PIC H}
	B\GL_{1,S} \simeq \PIC_S(\stQCoh_S) \to \PIC_S(\stHcat_{A/S}).
\end{equation}

This gives rise to a canonical homomorphism  of $E_\infty$-group stacks,
\begin{equation}\label{eq:A x Pic to PicH}
	A\times_S B\GL_{1,S} \xrightarrow{\ \eqref{eq:A to PicH} \times \eqref{eq:PIC to PIC H}\ }
	\PIC(\stHcat_{A/S})\times_S\PIC(\stHcat_{A/S}) \xrightarrow{\ \star\ }
	\PIC(\stHcat_{A/S}),
\end{equation}
given informally by $(a,L) \mapsto a_*L$ where $a : T \to A_T$ and $L$ is a line bundle on $T$.

\begin{proposition}\label{prop:Pic HA classical}\leavevmode

(i) On classical truncations, \eqref{eq:A x Pic to PicH} induces an isomorphism of group stacks over $S_\cl$,
\[
  A_{\mathrm{cl}} \times_{S_\cl} B\GG_{m,S_{\cl}}
  \xrightarrow{\ \sim\ }
  \PIC(\stHcat_{A/S})_{\mathrm{cl}}.
\]

(ii) If $A \to S$ is of relative dimension $1$ and $k$ is a $\bQ$-algebra, then \eqref{eq:A x Pic to PicH} is an isomorphism.
\end{proposition}

\begin{lemma}\label{lem:rational etale rectification}
Let $R$ be a connective $E_\infty$-ring, $R[t] = R \otimes \Sigma^\infty_+(\mathbf{N})$ the polynomial $E_\infty$-$R$-algebra on a degree zero generator $t$, and $R[t] \to B$ an étale morphism of connective $E_\infty$-$R$-algebras.
If $R$ is a $\bQ$-algebra, then for every connective $E_\infty$-$R$-algebra $C$, the forgetful map
\[
  \Map_{\CAlg_R}(B,C)
  \to
  \Map_{\operatorname{Alg}_R}(B,C)
\]
is invertible.
\end{lemma}
\begin{proof}
Consider the commutative square
\[\begin{tikzcd}
  \Map_{\CAlg_R}(B,C) \ar[r] \ar[d] &
  \Map_{\operatorname{Alg}_R}(B,C) \ar[d] \\
  \Map_{\CAlg_R}(R[t],C) \ar[r] &
  \Map_{\operatorname{Alg}_R}(R[t],C).
\end{tikzcd}\]
Since $R$ is a $\bQ$-algebra, $R[t]$ is free as an $E_\infty$-algebra, and its underlying $E_1$-algebra is free on the same generator (see \cite[Props.~7.1.4.6, 7.1.4.11]{LurieHA}).
The respective universal properties thus identify the lower horizontal map above with the identity of $\Map_{\MOD_R}(R,C)$.

Taking fibres vertically at a fixed $E_\infty$-$R$-algebra map $R[t] \to C$ yields the map
\[
  \Map_{(\CAlg_R)_{R[t]/}}(B,C) \to
  \Map_{(\operatorname{Alg}_R)_{R[t]/}}(B,C).
\]
Since $R[t]\to B$ is \'{e}tale, Corollaries~7.5.4.6 and 7.5.3.13 of \cite{LurieHA} identify both source and target, respectively, with the set
\[
  \Hom_{\operatorname{CAlg}_{\pi_0R[t]}}(\pi_0B, \pi_0C),
\]
such that the forgetful map is identified with the identity.

Thus the lower horizontal map in the square above is invertible, as is the induced map on every vertical fibre.
It follows that the upper horizontal map is invertible.
\end{proof}

\begin{proof}[Proof of Proposition~\ref{prop:Pic HA classical}]
We first prove (i), for which we may take $S$ classical and $T$ a classical scheme over $S$.
Let $F \in \Hcat_{A_T/T}$ be a $\star$-invertible object.
Since $F$ is $\star$-invertible, it lies in $\Perf_0(A_T/T)$ by Corollary~\ref{duals}.
Since pushforward along $\pi : A_T \to T$ determines a symmetric monoidal functor $\Hcat_{A_T/T} \to \QCoh(T)$, $\pi_*F$ is an $\otimes$-invertible object of $\QCoh(T)$, i.e., a line bundle $L$ over $T$.
By Lemma~\ref{lem:push finite support}, $\pi_*$ is $t$-exact and conservative on $\QCoh_0$.
Since $L$ is discrete, $F$ is discrete and may therefore be regarded as an ordinary quasi-coherent sheaf on the classical scheme $A_T$.

Write $F \simeq i_*(F|_Z)$ where $i : Z \hookrightarrow A_T$ is the scheme-theoretic support.
Then $L \simeq f_*(F|_Z)$, where $f$ is the finite morphism $\pi \circ i: Z \hookrightarrow A_T \to T$, so that $L$ acquires the structure of a quasi-coherent $B$-module over the quasi-coherent $\OO_T$-algebra $B := f_*\OO_Z$, encoded by the associative $\OO_T$-algebra homomorphism
\[
  \rho: B \to
  \underline{\operatorname{End}}_{\OO_T}(L)
  \simeq \OO_T.
\]
Since $T$ is classical, $B$ and $\OO_T$ are discrete and $\rho$ is automatically a homomorphism of commutative $\OO_T$-algebras.
It thus determines a section $T \to Z$ of $f : Z \to T$, whence a section $a : T \to A_T$ of $A_T \to T$, with $F \simeq a_*L$ by construction.
The assignment $F \mapsto (a, L)$ determines an inverse to \eqref{eq:A x Pic to PicH} on $T$-points, and is natural in $T$.

We now prove (ii).
Let $T=\spc(R)$ be affine over $S$ and let $F \in \Pic(\Hcat_{A_T/T})$ be a $T$-point.
As above, $F$ lies in $\Perf_0(A_T/T)$ and $L := \pi_*F \in \QCoh(T)$ is a line bundle.
Applying (i), we obtain a section $a_\cl:T_\cl\to A_{T_\cl}$ and an isomorphism
\[
  F|_{T_\cl} \simeq (a_\cl)_*(L|_{T_\cl}).
\]

Given a Zariski cover $T' \twoheadrightarrow T$ with \v{C}ech nerve $T'_\bullet$, denote by $K_{F,T'_n}$ the fibre of
\[
  \bigl(A \times_S B\GL_{1,S}\bigr)(T'_n) \xrightarrow{\ \eqref{eq:A x Pic to PicH}\ }
  \PIC\bigl(\stHcat_{A/S}\bigr)(T'_n)
\]
at $F|_{T'_n}$.
Since \eqref{eq:A x Pic to PicH} is a map of stacks, the totalization of $K_{F,T'_\bullet}$ is the fibre $K_{F,T}$, so it will suffice to show that $K_{F,T'_n}$ is contractible for every $n$.

Since $A_T \to T$ is fiber-smooth of relative dimension $1$, we may take the Zariski cover $T' \twoheadrightarrow T$ to be such that, after replacing $T$ by any term of $T'_\bullet$, there exists an affine open $j:U=\spc(B)\hookrightarrow A_T$ whose classical truncation is an open neighbourhood of the image of $a_\cl$, along with an \'{e}tale coordinate $q:U\to\mathbf{A}^1_T=\spc(R[t])$.
Replacing $T$ by any $T'_n$, we may assume that such exist over $T$.
Since $F$ is supported on $U$, we have $F\simeq j_*j^*F$.
The fibre $K_{F,T}$ is the groupoid of pairs $(a',L')$ together with an isomorphism $\theta:a'_*L'\simeq F$.
Applying $\pi_*$ to $\theta$ identifies $L'$ with $L$, and the space of pairs consisting of a line bundle $L'$ and an isomorphism $L'\simeq L$ is contractible.
Part (i) shows that the classical truncation of $a'$ is $a_\cl$, so $a'$ factors through $U$.
As above, the restriction $F|_U$ corresponds to a $B$-module structure on $L$, encoded by $\rho : B \to \operatorname{End}_R(L) \simeq R$, where the isomorphism $R \to \operatorname{End}_R(L)$ sends $r$ to scalar multiplication $r : L \to L$.
After identifying $L'$ with $L$, the remaining datum is an $E_\infty$-$R$-algebra map $B \to R$, together with an identification of the underlying $E_1$-map with $\rho$.
That is to say, $K_{F,T}$ is identified with the fibre at $\rho$ of the forgetful map
\[
  \Map_{\CAlg_R}(B,R)
  \to
  \Map_{\operatorname{Alg}_R}(B,R),
\]
which is contractible by Lemma~\ref{lem:rational etale rectification} (since $R$ is a $\bQ$-algebra).
\end{proof}

\begin{remark}
	For $A = \GGaS$ in characteristic zero, \eqref{eq:A x Pic to PicH} is an isomorphism
	\[
    \GGaS \times_S B\GL_{1,S} \xrightarrow{\ \sim\ }
    \PIC(\stHcat_{\GGaS/S}).
	\]
	However, for $A = \GGm$, we instead have
	\[
    \PIC(\stHcat_{\GGmS/S}) \simeq
    \GL_{1,S} \times_S B\GL_{1,S},
	\]
	and \eqref{eq:A x Pic to PicH} is identified with the map $\GGmS \times_S B\GL_{1,S} \to \GL_{1,S} \times_S B\GL_{1,S}$, which is an isomorphism if and only if $S$ is defined over $\bQ$.
	Indeed, an object of $\Hcat_{\GGmS/S}(T)$ is a pair $(V,u)$ consisting of an object $V\in\IndPerf(T)$ and an automorphism $u$ of $V$.
  Such a pair is $\star$-invertible precisely
  when $V$ is a line bundle, and the automorphisms of a line bundle are parametrized by $\GL_{1,S}(T)$.
\end{remark}

\subsection{\texorpdfstring{$1$}{1}-shifted Cartier duals}

\subsubsection{The stack of $1$-shifted characters}\label{sssec:CharOne}

The following definition was first studied by Deligne in \cite{Deligne}.

\begin{definition}
  Let $A \to S$ be a relative commutative group stack.
  The \emph{stack of $1$-shifted characters} of $A$ is the group-theoretic mapping stack
  \[ \CharOne(A/S) := \HOM_S^{E_\infty}(A, B\GL_{1,S}). \]
\end{definition}

By definition, a $T$-point of $\CharOne(A/S)$ is a homomorphism of $E_\infty$-group stacks $A_T := A \times_S T \to B\GL_{1,S}$ over $S$.
Such a morphism amounts to a multiplicative line bundle on $A_T$, i.e., a line bundle $L \to A_T$ together with a specified isomorphism $$m^{*}L\rightarrow L\boxtimes_{T} L$$ which is associative and commutative up to coherent homotopy. Note that formation of $\CharOne(A/S)$ is stable under base change in $S$.

\begin{remark} The stack $\CharOne(A/S)$ is equipped with a canonical commutative group structure: its underlying $E_\infty$-group structure is given in terms of multiplicative line bundles by tensor product, and it is strict because $A$ is, since \eqref{eq:str adjunction} gives $\CharOne(A/S) \simeq \HOM_S(A, (B\GL_{1,S})^{\mathrm{str}})$. Note that the unit $\mathcal{O}_A$ is equipped with a tautological multiplicative structure. The stabilizer of $\CharOne(A/S)$ at the corresponding point is the usual Cartier dual,
$$\operatorname{Aut}_{\CharOne(A/S)}(\mathcal{O}_A) \simeq \chi(A/S) := \HOM_S^{E_\infty}(A, \GL_{1,S}).$$
In particular, when there are no nontrivial multiplicative line bundles on $A$, we have $\CharOne(A/S)\simeq B\chi(A/S).$ \end{remark}

\begin{example} For $A=\GGm$ we have canonical isomorphisms
\[
  \chi(\GGm) \simeq \mathbf{Z},
  \qquad \CharOne(\GGm) \simeq B\mathbf{Z}.
\]
\end{example}

\begin{example}
For $A=\GGa$ over a $\bQ$-algebra $k$, we have
$$ \chi(\GGa) \simeq \GGhat_{a}, \quad\textrm{but}\quad \CharOne(\GGa) \not\simeq B\GGhat_{a}$$ even over the base $k=\spc(\mathbf{C})$, due to the nontriviality of $\operatorname{Ext}^1(\GGa, \GM)$, see \cite{RibeiroRosengarten}. Over an ordinary commutative ring $k$, we have
\[ \chi(\GGa) \simeq \GGhat_{a}^\sharp, \]
where $\GGa^\sharp$ is the PD-hull of $\GGa$, i.e., $\spc(k[x, x^2/2!, x^3/3!, \ldots])$ where $x$ is the coordinate of $\GGa$ (see \cite[\S 2.4]{BhattFGauges}). There is always a canonical monomorphism of group stacks
$$ B \GGhat_a^\sharp \hookrightarrow \CharOne(\GGa), $$
the inclusion of the connected component of the identity element.
\end{example}

\begin{example}
If $A \to S$ is a relative abelian variety, we have
\[
  \chi_S(A) \simeq e,
  \qquad \CharOne(A/S) \simeq \check{A},
\]
where $e$ is the trivial group scheme over $S$, and $\check{A}$ is the dual abelian variety.
\end{example}

\subsubsection{The Fourier--Mukai transform}
\label{sssec:Phi}

Consider the diagram \[\begin{tikzcd}
  & \CharOne(A/S) \times_S {A} \ar[swap]{ld}{\pr_{1}}\ar{d}{\ev_A}\ar{rd}{\pr_{2}} &
  \\
  \CharOne(A/S) & B\GL_{1,S} & A
\end{tikzcd}\]
where $\pr_i$ are the projections, and $\ev_{A}$ is the evaluation map.
Recall that $B\GL_{1,S}$ is equipped with the universal line bundle $\OO(1)$, and $\ev_A^*(\OO(1))$ is the universal multiplicative line bundle.

The \emph{Fourier--Mukai transform} is the $S$-linear colimit-preserving functor
\[ \Phi_{A/S} : \QCOH(A) \to \QCOH(\CharOne(A/S)) \]
defined by $F \mapsto \pr_{1,*}(\pr_{2}^*(F) \otimes \ev_{A}^*(\OO(1)))$.

\begin{lemma}\label{lem:FM t-exact} The restriction $\Phi_{A/S}|_{\QCOH_0(A/S)}$ is t-exact. In particular, it sends $\HAS$ to $\QCoh(\CharOne(A/S))$.
\end{lemma}

\begin{definition}\label{def:FMA}
  We call $\Phi_{A/S} : \HAS \to \QCoh(\CharOne(A/S))$ the \emph{Fourier--Mukai transform} for $A$.
\end{definition}

\begin{proposition} The functor $\Phi_{A/S} : \HAS \to \QCoh(\CharOne(A/S))$ is a tensor functor. \end{proposition}
\begin{proof} The symmetric monoidality follows from the base change and projection formulas for $*$-pushforward, and the multiplicativity of $\ev_{A}^*(\OO(1))$. \end{proof}

\subsubsection{A Tannakian dual}
\label{sssec:Tannakian dual}

By Proposition~\ref{prop:SpecTens monoidal}, the Hopf algebra structure on $\HAS$ (Corollary~\ref{cor:QCoh0 Hopf}) endows its spectrum
$$ \SpecTens_S(\stHAS) $$
with the structure of a commutative group stack over $S$.
There is a canonical group homomorphism
\begin{equation}\label{eq:CharOne to SpecTens}
  \CharOne(A/S) \to \SpecTens_S(\stHAS)
\end{equation}
over $S$, classifying the morphism of sheaves of tensor categories
\[ \Phi_{A/S} : \stHAS \to \stQCoh_{\CharOne(A/S)} \]
defined by the pointwise Fourier--Mukai transform (Definition~\ref{def:FMA}).
On $T$-points, \eqref{eq:CharOne to SpecTens} sends a multiplicative line bundle $L \to A_T$ to the morphism of sheaves $\stHAS|_T \to \stQCoh_T$ whose value on an affine $T'$ over $T$ is the tensor functor
$$ \Hcat_{A_{T'}/T'} \rightarrow \QCoh(T'), \quad F \mapsto \pi_{T',*}\bigl(F \otimes L|_{A_{T'}}\bigr). $$

\begin{remark} The group stack $\SpecTens_S(\stHAS)$ seems to be rather pathological, and the morphism \eqref{eq:CharOne to SpecTens} is typically far from invertible. \end{remark}

\subsubsection{The $1$-shifted Cartier dual}\label{sssec:check A}

\begin{definition}\label{def:abstract 1-shifted Cartier dual}
  A \emph{$1$-shifted Cartier dual pair} over $S$ consists of a tuple $(A, \check{A}, \jmath)$, where:
  \begin{itemize}
    \item $A \to S$ is a relative commutative group scheme (smooth, separated, and of finite presentation);
    \item $\check{A} \to S$ is a relative commutative group stack;
    \item $\jmath : \check{A} \longrightarrow \CharOne(A/S)$ is a homomorphism of commutative group stacks over $S$,
  \end{itemize}
  subject to the following conditions:
  \begin{enumerate}[label={\upshape(\alph*)}]
    \item\label{item:abstract 1-shifted Cartier dual/FM}
    The tensor functor
    \begin{equation*}\label{eq:Phi check A st}
      \Phi_{(A,\check{A},\jmath)/S} : \stHAS
      \xrightarrow{\Phi_{A/S}} \stQCoh_{\CharOne(A/S)}
      \xrightarrow{\jmath^*} \stQCoh_{\check{A}},
    \end{equation*}
    defined by restricting the Fourier--Mukai transform (Definition~\ref{def:FMA}) along $\jmath$, is pointwise fully faithful.
    \item\label{item:abstract 1-shifted Cartier dual/sharp}
    Consider the projection $\check{p} : \check{A} \to S$.
    The functor $\check{p}_* : \QCOH(\check{A}) \to \QCOH(S)$ preserves perfect complexes and satisfies the projection formula and base change against $*$-pullback; moreover, $\check{p}^*$ admits a left adjoint $\check{p}_\sharp : \QCOH(\check{A}) \to \QCOH(S)$ also satisfying the projection formula and base change against $*$-pullback.
  \end{enumerate}
\end{definition}

\begin{remark}
When there is no risk of confusion, we will write $\Phi_{A/S}$ instead of $\Phi_{(A,\check{A},\jmath)/S}$.
\end{remark}

\begin{definition}\label{def:morphism shifted Cartier pairs}
	A \emph{morphism} $(A, \check{A}, \jmath) \to (A', \check{A}', \jmath')$ consists of group homomorphisms $f : A \to A'$ and $g : \check{A}' \to \check{A}$ over $S$, along with a commutative square
	\[\begin{tikzcd}
	  \check{A}' \ar{r}{g}\ar{d}{\jmath'}
		& \check{A} \ar{d}{\jmath}
		\\
		\CharOne(A'/S) \ar{r}{f^*}
		&  \CharOne(A/S)
	\end{tikzcd}\]
	in the category of commutative group stacks.
\end{definition}

\subsubsection{The additive, multiplicative, and elliptic cases}
\label{sssec:HA trichotomy}

In each of our cases of interest, there is a canonical choice for the $1$-shifted Cartier dual.
Namely, for $\GGaS \to S$, $\GGmS \to S$, and a relative elliptic curve $E \to S$, we define the $1$-shifted Cartier dual by:
\[
\check{\GG}_{a,S} := B\GGhat^\sharp_{a,S},
\quad \check{\GG}_{m,S} := B\mathbf{Z}_S = S \times S^1,
\quad \text{and }\check{E},
\]
respectively, where $\check{E} \to S$ denotes the dual elliptic curve.
More generally, for $A \to S$ a relative abelian variety, its $1$-shifted Cartier dual is the dual abelian variety $\check{A} \to S$ (cf. \cite[Prop.~6.8.2]{LurieEll1}).

As we have seen in (\ref{sssec:CharOne}), there is a canonical homomorphism $\jmath : \check{A} \hookrightarrow \CharOne(A/S)$, which is in fact invertible in all but the case $A=\GGaS$, for which it is the inclusion of the identity component.\footnote{%
Note conversely that $\CharOne(\check{A}/S) \simeq A$ does hold for each of $A \in \{ \GGaS, \GGmS, E \}$.
}
Let us check condition~\ref{item:abstract 1-shifted Cartier dual/FM}, on the Fourier--Mukai transform.
By allowing $S$ to vary and unravelling definitions, it suffices to show that the tensor functor
\begin{equation}\label{eq:Phi check A}
  \Phi_{(A,\check{A},\jmath)/S} : \HAS \to \QCoh(\check{A}),
  \qquad F \mapsto \pr_{1,*}(\pr_2^*(F) \otimes \ev^*_A(\OO(1)))
\end{equation}
is fully faithful, where $\pr_i$ denote the projections of $\check{A} \times_S A$ and $\ev_A$ denotes the evaluation map
\[ \ev_A : \check{A} \times_S A \hookrightarrow \CharOne(A/S) \times_S A \xrightarrow{\ev_A} B\GL_{1,S}. \]

\begin{theorem} \label{thm:HA vs QCoh checkA} Let $A \to S$ be either $\GGaS \to S$ (over a $\bQ$-algebra), $\GGmS \to S$, or a relative abelian variety. Then the functor \eqref{eq:Phi check A} is fully faithful. Moreover, in the additive and multiplicative cases it determines tensor equivalences
\[ \Hcat_{\GGaS/S} \simeq \IndPerf(B\GGhat_{a,S}^\sharp) \subseteq \QCoh(B\GGhat_{a,S}^\sharp), \]
\[ \Hcat_{\GGmS/S} \simeq \IndPerf(B\bZ_S) \subseteq \QCoh(B\bZ_S). \]
If $A \to S$ is proper, it determines a symmetric monoidal equivalence
\[ \Hcat_{A/S} \simeq \IndPerf^{\,\text{hmg}/S}(\check{A}) \subseteq \QCoh(\check{A}), \]
where $\Perf^{\,\text{hmg}/S}(\check{A})$ is the full subcategory spanned by perfect complexes $G\in\Perf(\check{A})$ such that the cohomology sheaves of the $S$-fibres of $G$ are homogeneous vector bundles.
In particular, if $A=E \to S$ is a relative elliptic curve, we have
\[ \Hcat_{E/S} \simeq \IndPerf^{\,\ss,0}(\check{E}) \subseteq \QCoh(\check{E}), \]
where $\Perf^{\,\ss,0}(\check{E})$ is the full subcategory spanned by perfect complexes $G\in\Perf(\check{E})$ such that the cohomology sheaves of the $S$-fibres of $G$ are semistable vector bundles of degree zero.
\end{theorem}

\begin{proof}

  For $A$ an abelian variety, $\Phi_{A/S}$ is the restriction of the classical Fourier--Mukai transform $\QCoh(A) \to \QCoh(\check{A})$ of \cite{Mukai}.
  Recall that the latter is an equivalence, and sends the skyscraper sheaf $\OO_x$ for an $S$-point $x : S \to A$ to the degree zero Poincar\'e line bundle $\mathcal{P}_x$.
  It then follows by Mukai's characterization of homogeneous vector bundles \cite[Exam.~3.2]{Mukai} that it sends discrete coherent sheaves on $A$ that are torsion, fibrewise over $S$, to discrete coherent sheaves on $\check{A}$ that are locally free and homogeneous, fibrewise over $S$.
  Since $\Phi_{A/S}$ is t-exact, we deduce that $\Phi_{A/S}$ determines an equivalence $\Perf_0(A/S) \simeq \Perf^{\,\text{hmg}/S}(\check{A})$, so the claim follows.
  In the elliptic curve case, recall that homogeneity of each geometric fibre is equivalent to semistability of degree zero (see \cite[p.~159]{Mukai} and \cite[Thm.~5]{At}).

  In the additive and multiplicative cases we use the descriptions of $\HAS$ given in Subsect.~\ref{ssec:HA}.
  We write $S=\spc(k)$ to simplify the notation.

  For $A=\GGm$, we identify $\Hcat_{\GM}$ with the ind-completion of the category of objects $V \in \Perf(S)\simeq \Mod_{k}^{\mathrm{perf}}$ equipped with an automorphism $T$ representing the action of the coordinate $t \in \OO(\GGm)$.
  For a stack $X$, $\IndPerf(X \times B\bZ)$ is the ind-completion of the category of objects of $\Perf(X)$ equipped with an equivariant structure for the discrete group $\bZ$, i.e., an automorphism representing the action of $1 \in \bZ$.
  To show that the Fourier--Mukai transform $\Phi_{\GGm}$ is identified with this tautological equivalence $\Hcat_{\GM} \simeq \IndPerf(B\bZ)$, observe that the multiplicative line bundle on $B\bZ \times \GGm$ classified by the isomorphism $B\bZ \to \CharOne(\GGm)$ is the trivial line bundle $\OO_{\GGm}$ equipped with the $\bZ$-equivariant structure where $1 \in \bZ$ acts by multiplication by the coordinate $t$ of $\GGm$.
  Thus for $F \in \Hcat_{\GGm}$ corresponding to $(V,T)$, where $T$ is the action of the coordinate $t$, $\pr_2^*(F) \otimes \ev_{\GGm}^*(\OO(1)) \in \QCoh(B\bZ \times \GGm)$ corresponds to the same object $F$ equipped with the $\bZ$-equivariant structure where $1\in\bZ$ acts by multiplication by $t$.
  Its pushforward along $\pr_{1,*}$ thus corresponds to $(V,T)$.

  For $A=\GGa$ (over a $\bQ$-algebra), we identify $\Hcat_{\GGa}$ with the ind-completion of the category of objects $V \in \Perf(S)\simeq \Mod_{k}^{\mathrm{perf}}$ equipped with an endomorphism $T$ representing the action of the coordinate $x \in \OO(\GGa)$.
  On the other hand, $\IndPerf(B\GGhat_{a}^\sharp)$ is the ind-completion of the category of objects of $\Perf(S)$ equipped with an equivariant structure for the formal group $\GGhat_{a}^\sharp$.
  Since the coordinate ring of $\GGhat_{a}^\sharp$ is the formal divided power algebra $\widehat{\Gamma}(t)$, such an equivariant structure amounts to a continuous coaction $\rho : V \to V \widehat{\otimes} \widehat{\Gamma}(t) \simeq \prod_{n \ge 0} V \otimes t^{[n]}$.
  Using Cartier duality to identify the continuous dual of $\widehat{\Gamma}(t)$ with the polynomial algebra $k[x]$, a continuous $\widehat{\Gamma}(t)$-comodule structure on $V$ is equivalent to a $k[x]$-module structure, specified by an endomorphism $T$ of $V$.
  We have described an equivalence $\Hcat_{\GGa} \simeq \IndPerf(B\GGhat_{a}^\sharp)$, and the proof that this coincides with the Fourier--Mukai transform $\Phi_{\GGa}$ follows similarly to the multiplicative case, using the observation that the canonical map $B\GGhat_{a}^\sharp \to \CharOne(\GGa)$ classifies the character $\exp(tx) = \sum_{n \ge 0} t^{[n]} x^n$.
\end{proof}

\begin{remark} The fact that, for an elliptic curve $E$, $\HE$ is not obviously of the form $\IndPerf(\YY)$ for any stack $\YY$ is precisely the reason we are forced to adopt the Tannakian perspective in this paper. \end{remark}

\begin{remark}\label{rem:HA vs QCoh checkA dzbl}
  In the situation of Theorem~\ref{thm:HA vs QCoh checkA} with $A=E$, $\Phi_{E/S}$ is $t$-exact (Lemma~\ref{lem:FM t-exact}) and restricts to an equivalence
  \[ \Hcat_{E/S}^{\heartsuit} \simeq \IndPerf^{\,\ss,0}(\check{E})^{\heartsuit} \]
  on hearts.
  Further restricting to dualizable objects yields the equivalence
  \[ \Hcat_{E/S}^{\heartsuit,\mathrm{dzbl}} \simeq \operatorname{\textsc{Vect}^{\,\ss,0}}(\check{E}), \]
  since a dualizable object of $\IndPerf^{\,\ss,0}(\check{E})^{\heartsuit}$ is precisely a vector bundle on $\check{E}$ which is fibrewise semistable of degree zero.
\end{remark}

To finish the verification that $\GGa$, $\GGm$, and $E$ have $1$-shifted Cartier duals $B\GGhat_{a}^\sharp$, $B\bZ$, and $\check{E}$, respectively, it remains to check condition~\ref{item:abstract 1-shifted Cartier dual/sharp} for the projection $\check{p} : \check{A} \to S$.
For $A=E$ (or more generally for a relative abelian variety $A \to S$), $\check{p}$ is smooth, proper and representable, so the claims are standard (see e.g. \cite[Exam.~A.7]{weilres}).

For the cases of $A=\GG_{m,S} \to S$ and $A=\GG_{a,S} \to S$, we use the algebraic descriptions of $\QCOH(\check{A})$ as in the proof of Theorem~\ref{thm:HA vs QCoh checkA}.
For example, $\QCOH(B\bZ_S)$ is identified with the category of $\OO_S[t,t^{-1}]$-modules in $\QCOH(S)$, and $\check{p}^* : \QCOH(S) \to \QCOH(B\bZ_S)$ is identified with restriction of scalars along $\OO_S[t,t^{-1}] \to \OO_S$, $t \mapsto 1$.
Its left adjoint $\check{p}_\sharp$ corresponds thus to extension of scalars along the same map, i.e., $(V,T) \mapsto \operatorname{cofib}(T-\id_V)$, while its right adjoint $\check{p}_*$ corresponds to coextension of scalars, i.e., $(V,T) \mapsto \operatorname{fib}(T-\id_V)$.

\subsubsection{Action of $\check{A}$ on $\stHA$}\label{sssec:actonHA}

The group stack $\SpecTens_S(\stHAS)$ is equipped with a natural action on $\stHAS$ via the pointwise comultiplication $\nabla : \stHAS \to \stHAS \otimes_S \stHAS$.
Along the homomorphism $\CharOne(A/S) \to \SpecTens_S(\stHAS)$, this tautological action restricts to the following action of $\CharOne(A/S)$: a $T$-point classifying a multiplicative line bundle $L \to A_T$ acts by
\[ \Hcat_{A_{T'}/T'} \to \Hcat_{A_{T'}/T'}, \qquad F \mapsto F \otimes L|_{A_{T'}} \]
for every affine $T'$ over $T$.

In particular, for any $1$-shifted Cartier dual pair $(A, \check{A}, \jmath)$, we obtain by restriction a canonical action of $\check{A}$ on the sheaf of tensor categories $\stHAS$.

\part{Rotation-equivariant sheaves}\label{part:rotation}

\section{The rotation algebra}\label{sec:RA}

We fix a base stack $S$ and a $1$-shifted Cartier dual pair $(A, \check{A}, \jmath)$ over $S$.
In this section we study several Hopf algebras associated to $A$ and $\check{A}$.
These include the algebras of functions on the groups $\Omega A$ and $\check{A}$, and others derived from these.
As we will see in \S \ref{sec:Euler}, these will provide algebraic models for the category of quasi-coherent sheaves on the classifying stack $B\check{A}$.

\subsection{Conventions}

\subsubsection{}\label{sssec:module conventions}

For an $E_\infty$-algebra $R$ in $\QCOH(S)$:

\begin{itemize}
	\item We write $\MOD_R$ for the symmetric monoidal category of $R$-modules in $\QCOH(S)$.
	\item We write $\MOD_R^{\mathrm{perf}} \subseteq \MOD_R$ for the full subcategory of \emph{dualizable} objects; these are equivalently those which, locally on $S$, lie in the thick subcategory generated by $R$.
	\item Since $\OO_S$ is the unit object of $\QCOH(S)$, we have $\MOD_{\OO_S} = \QCOH(S)$ and $\MOD_{\OO_S}^{\mathrm{perf}} = \PERF(S)$.
	\item We write $\MOD_R^{\mathrm{perf}/\OO_S} \subseteq \MOD_R$ for the full subcategory of $R$-modules whose underlying $\OO_S$-module is perfect.
\end{itemize}
When $R$ is moreover an $E_\infty$-Hopf algebra in $\QCOH(S)$:
\begin{itemize}
  \item We write $(\MOD_R, \star)$ for the category $\MOD_R$ equipped with the \emph{convolution} symmetric monoidal structure.
	Given $V, W \in \QCOH(S)$ with actions $\alpha_V : R \otimes V \to V$ and $\alpha_W$, the convolution product $V \star W$ is $V \otimes W \in \QCOH(S)$, equipped with the action
  \[\begin{multlined}
     R \otimes V \otimes W \xrightarrow{\Delta \otimes \operatorname{id}} R \otimes R \otimes V \otimes W
     \\\;\simeq\; R \otimes V \otimes R \otimes W \xrightarrow{\alpha_V \otimes \alpha_W} V \otimes W,
  \end{multlined}\]
  where $\Delta$ is the comultiplication map.
  The monoidal unit is $\OO_S$ with trivial $R$-action.
	In particular, the forgetful functor
	$$(\MOD_{R}, \star) \to (\MOD_{\OO_S}, \otimes)$$
	is symmetric monoidal.
\end{itemize}

\begin{lemma}\label{lem:dual mod conv}
  Let $R$ be an $E_\infty$-Hopf algebra in $\QCOH(S)$.
  An object of $(\MOD_R, \star)$ is dualizable if and only if it lies in $\MOD_{R}^{\mathrm{perf}/\OO_S}$.
\end{lemma}
\begin{proof}
  Consider the restriction of scalars functor $\varpi_* : (\MOD_R, \star) \to (\MOD_{\OO_S}, \otimes)$ along the unit $\varpi : \OO_S \to R$.
  Since $\varpi_*$ is symmetric monoidal, it sends dualizable objects to perfect $\OO_S$-modules; that is, every dualizable object $M \in (\MOD_R, \star)$ has perfect underlying $\OO_S$-module $\varpi_*(M)$.

  Conversely, let $M \in \MOD_R^{\mathrm{perf}/\OO_S}$, given by $V := \varpi_*(M) \in \MOD_{\OO_S}^{\mathrm{perf}}$ and action $\alpha : R \otimes V \to V$.
  The ``weak'' dual of $M$ in the closed symmetric monoidal category $(\MOD_R, \star)$ is $M^\vee := \hhoms_R^{\star}(M, \OO_S)$, where $\OO_S$ is regarded with trivial $R$-action (and $\hhoms_R^{\star}$ denotes the internal hom for the convolution monoidal structure).
  This is the $\OO_S$-linear dual $V^\vee$, equipped with the contragredient $R$-action
  $$ R \otimes V^\vee \xrightarrow{\chi \otimes \id} R \otimes V^\vee \xrightarrow{\alpha^\vee} V^\vee $$
  where $\chi : R \to R$ is the antipode of $R$.
  The antipode axioms imply that the evaluation and coevaluation maps for $V$ are $R$-linear, and hence lift to maps $M^\vee \star M \to \OO_S$ and $\OO_S \to M \star M^\vee$.
  The triangle identities showing that $M^\vee$ is dual to $M$ now follow from those exhibiting $\varpi_*(M^\vee) = V^\vee$ as the dual of $\varpi_*(M) = V$, since $\varpi_*$ is conservative, symmetric monoidal, and commutes with $(-)^\vee$.
\end{proof}

\subsubsection{}

We will write the relative algebra of functions of any $p : X \to S$ as
\[ \OO(X) := p_*(\OO_X) \in \QCOH(S), \]
when there is no risk of confusion.
When $X=G$ is a relative group stack and $p_*$ satisfies the base change and projection formulas, $\OO(G)$ carries a natural Hopf algebra structure.
Indeed, we then have the K\"unneth isomorphism
\[ \OO(G) \otimes_{\OO_S} \OO(G) \simeq \OO(G \times_S G), \]
and the Hopf algebra structure on $\OO(G)$ is given as follows:
\begin{itemize}
  \item the multiplication is the composite
  $$ \OO(G) \otimes_{\OO_S} \OO(G) \simeq \OO(G \times_S G) \xrightarrow{\Delta^*} \OO(G),$$
  where $\Delta : G \to G \times_S G$ is the diagonal;
  \item the unit is $p^* : \OO_S \to \OO(G)$;
  \item the comultiplication is the composite
  $$\OO(G) \xrightarrow{m^*} \OO(G \times_S G) \simeq \OO(G) \otimes_{\OO_S} \OO(G),$$
  where $m : G \times_S G \to G$ is the multiplication;
  \item the counit is $e^* : \OO(G) \to \OO_S$, where $e : S \to G$ is the unit section;
  \item the antipode is $\iota^* : \OO(G) \to \OO(G)$, where $\iota : G \to G$ is the inversion map.
\end{itemize}

\subsection{The rotation algebra}\label{ssec:RA}

\subsubsection{}\label{sssec:RA}

We write $\Omega A := e \times_A e$ for the loop space of $A$ based at the identity section $e$, with its induced group structure.
The \emph{rotation algebra} is the Hopf algebra $\RA := \OO(\Omega A)$.
We denote its $\OO_S$-linear dual by $\RA^\vee$.

\subsubsection{}\label{sssec:maps out of RA}

Since $A$ is separated over $S$, $\Omega A$ is affine over $S$.
For every $E_\infty$-algebra $R$ over $\OO_S$, we have
\begin{equation*}
	\Map_{\CAlg_{\OO_S}}(\RA, R) \;\simeq\; \Map_S(\spc(R), \Omega A) \;\simeq\; \Omega A(R),
\end{equation*}
where $\CAlg_{\OO_S}$ denotes the category of (quasi-coherent) $E_\infty$-algebras over $\OO_S$.

\subsubsection{}
Equivalently, we have
\[ \RA \simeq e^*e_*\OO_S, \qquad \RA^\vee \simeq e^!e_*\OO_S \simeq \RA \otimes \fra[-d], \]
where $d$ is the relative dimension of $A \to S$.
We have used the canonical isomorphism $e^!(-) \simeq e^*(-) \otimes \omega_{S/A}$ (see \cite[Cor.~6.4.2.7]{SAG}), where $\omega_{S/A} \simeq \fra[-d]$ is the relative dualizing complex of $e : S \to A$.

\subsubsection{Additive splitting}\label{sssec:funcheckA add split}

Suppose $A \to S$ has relative dimension $1$.
Then $e : S \to A$ is a virtual Cartier divisor with $e^*\OO_A(-e) \simeq \frav$ (recall $\frav \simeq e^*\omega_{A/S}[-1]$).
The exact triangle $\OO_A(-e) \to \OO_A \to e_*\OO_S$, after applying $e^*$ and rotating, becomes
$$ \OO_S \longrightarrow e^*e_*\OO_S \longrightarrow \frav[1]. $$
The unit $\OO_S \to e^*e_*\OO_S$ admits a retraction via the counit (by one of the triangle identities for the adjunction), so this triangle yields the \emph{additive} splitting
$$\RA \simeq e^*e_*\OO_S \simeq \OO_S \oplus \frav[1].$$
Dually, $\RA^\vee \simeq \OO_S \oplus \fra[-1]$ in $\QCOH(S)$.

\subsubsection{Functions on \texorpdfstring{$\check{A}$}{Ǎ}}\label{sssec:funcheckA}

We describe the Hopf algebra $\OO(\check{A})$ of functions on the $1$-shifted Cartier dual in terms of $\RA$.

\begin{proposition} \label{prop:fun check A}
  There are canonical isomorphisms of Hopf algebras
  \[
    \OO(\check{A})
    \simeq \operatorname{End}_{\HA}(\OO_{e})
    \simeq \RA^\vee.
  \]
\end{proposition}

\begin{proof}
Since $\OO_{e} := e_*(\OO_S) \in \HA$ is the unit, the Hopf algebra structure on $\operatorname{End}_{\HA}(\OO_{e})$ is induced by the Hopf algebra structure on the category $\HA$.
The first isomorphism then comes from the symmetric monoidal fully faithful functor $\Phi_{A/S} : \HA \to \QCoh(\check{A})$ of Theorem~\ref{thm:HA vs QCoh checkA}.
Then we observe that
$$ \Hom_{\HA}(e_*\OO_S, e_*\OO_S) \simeq \Hom_{\OO_S}(\OO_S,e^!e_*\OO_S)$$ $$\simeq e^!e_*\OO_S \simeq \RA^\vee$$ and we are done.
\end{proof}

\subsubsection{}\label{sssec:Omega A hat}

Write $\widehat{A}$ for the formal completion of $A$ along the identity section, given by the functor
$\widehat{A}(R) := \operatorname{fib}\big( A(R) \to A(\pi_0 (R)_{\mathrm{red}}) \big)$ on $\operatorname{CAlg}^{\mathrm{cn}}_{\OO_S}$, as in \cite[Constr.~7.1.1]{LurieEll2}.
Note that $\Omega A$ only depends on the formal group $\widehat{A}$: the morphism of group objects $\widehat{A} \rightarrow A$ induces an isomorphism
\[
  \Omega \widehat{A}
  \xrightarrow{\ \sim\ } \Omega A.
\]
In particular, the Hopf algebra structure on $\RA$ (or equivalently $\RA \simeq \OO(\Omega\widehat{A})$) is determined in characteristic zero by the formal group law of $A$.

\subsubsection{}

The previous computation leads to the following description of the affinization of the $1$-shifted Cartier dual:

\begin{corollary}\label{cor:Aff check}
  There is a canonical isomorphism
  \[ \operatorname{Aff}(\check{A}) \simeq \operatorname{Aff}(B(\mathbb{D}\widehat{A})). \]
\end{corollary}

\begin{remark}\label{rem:D hat A unipotence}
  If $S=\spc(k)$ where $k$ is a field, then $\mathbb{D}\widehat{A}$ is unipotent, so $B(\mathbb{D}\widehat{A})$ is an affine stack by \cite[Thm.~2.4.1]{Toen} and the isomorphism becomes
  \[ \operatorname{Aff}(\check{A}) \simeq B(\mathbb{D}\widehat{A}). \]
  For example, if $k$ is of characteristic zero then $\mathbb{D}\widehat{A} \simeq (\GGa)^{\dim A}$ is always unipotent.
  For $A=\GGm$ over a field $k$ of characteristic $p>0$, $\mathbb{D}\widehat{\GG}_m \simeq \mathbf{Z}_p$ is also unipotent: it is the unipotent completion of $\mathbf{Z}$ (see the proof of \cite[Thm.~2.5.1]{Toen}), so $B\mathbf{Z}_p$ is an affine stack (cf. \cite[Prop.~3.3.4]{MRT}); more generally, $\mathbb{D}\widehat{A}$ is unipotent for every $A$ (it is the limit of the Cartier duals of the Frobenius kernels $\widehat{A}[F^n]$, which are connected; see \cite[\S8.5]{Waterhouse} and \cite[IV, \S2, Prop.~2.3(b)]{DemazureGabriel}).\footnote{%
    We thank Tasos Moulinos for explaining this to us.
  }
\end{remark}

\begin{proof}[Proof of Corollary~\ref{cor:Aff check}]
  Since $\OO(\mathbb{D}\widehat{A}) \simeq \OO(\widehat{A})^\vee$, \v{C}ech descent gives
  \[
    \OO(B(\mathbb{D}\widehat{A}))
    \simeq \operatorname{Tot}(\OO((\mathbb{D}\widehat{A})^{\times \bullet}))
    \simeq \operatorname{Tot}(\OO(\widehat{A}^{\times \bullet})^\vee).
  \]
  Here $(-)^\vee$ on the right-hand side denotes the \emph{continuous}\footnote{The word \emph{continuous} is meant only to evoke the duality between, say, linearly compact complete topological vector spaces and discrete topological vector spaces.} linear dual
  \[
    (-)^\vee : \operatorname{Pro}(\PERF(S))^{\mathrm{op}}
    \simeq \Ind(\PERF(S)^{\mathrm{op}})
    \to \Ind(\PERF(S)) \simeq \QCOH(S),
  \]
  the ind-extension of the linear dual on perfect complexes.
  By construction, we have
  \[
    \operatorname{Tot}(\OO(\widehat{A}^{\times \bullet})^\vee)
    \simeq \big| \OO(\widehat{A}^{\times \bullet}) \big|^\vee
    \simeq \OO(\Omega\widehat{A})^\vee
    \simeq \RA^\vee.
  \]
  Thus by Proposition~\ref{prop:fun check A} we have the algebra isomorphism
  \[
    \OO(B(\mathbb{D}\widehat{A}))
    \simeq \RA^\vee
    \simeq \OO(\check{A}).
  \]
  Applying $\Specnc(-)$ to both sides yields the first isomorphism.
\end{proof}

\subsection{The double rotation algebra}\label{ssec:modules}

\subsubsection{}

Let $\varpi : \Omega A \to S$ denote the projection, $\eps : S \to \Omega A$ the identity section, and $\Omega^2(A) := \eps \times_{\Omega A} \eps$ the iterated loop space of $A$ based at $\eps$.
The \emph{double rotation algebra} $\lam_A := \OO(\Omega^2 A)$ is the Hopf algebra of functions on $\Omega^2 A$. We denote its $\OO_S$-linear dual by $\lamco_A$.

\subsubsection{}\label{sssec:maps out of lam}

Since $\Omega^2 A$ is affine over $S$, we have
\begin{equation*}
	\Map_{\CAlg_{\OO_S}}(\lam_A, R) \;\simeq\; \Map_S(\spc(R), \Omega^2 A) \;\simeq\; \Omega^2 A(R),
\end{equation*}
for every $E_\infty$-algebra $R$ over $\OO_S$, where $\CAlg_{\OO_S}$ denotes the category of (quasi-coherent) $E_\infty$-algebras over $\OO_S$.

\subsubsection{}

On functions, $\varpi : \Omega A \to S$ and $\eps : S \to \Omega A$ induce the unit and augmentation of $\RA$, respectively:
\[
  \varpi : \OO_S \to \RA,
  \qquad
  \eps : \RA \to \OO_S.
\]
We have
\[
  \lam_A \simeq \OO_S \otimes_{\RA} \OO_S,
  \qquad
  \lamco_A \simeq \operatorname{End}_{\RA}(\OO_S).
\]

\subsubsection{}

The restriction of scalars functor $\eps_* : \MOD_{\OO_S} \to \MOD_{\RA}$ admits left and right adjoints,
\[
  \eps^*(M) := M \otimes_{\RA} \OO_S,
  \qquad \eps^!(M) := \Hom_{\RA}(\OO_S, M),
\]
respectively.
In those terms we may write
\[ \lam_A \simeq \eps^*\eps_*\OO_S, \qquad \lamco_A \simeq \eps^!\eps_*\OO_S. \]

\begin{proposition}\label{prop:eps^!(R)}
  Suppose $A \to S$ has relative dimension one.
  Then there is a canonical isomorphism in $\QCOH(S)$,
  \[
    \eps^!(\RA) \;\simeq\; \frav[1],
  \]
  where $\frav \simeq e^*\omega_{A/S}[-1]$ (see \ref{sssec:comm grps}).
\end{proposition}
\begin{proof}
  We denote by $\varpi^*$ and $\varpi^!$ the extension and coextension of scalars functors $\MOD_{\OO_S} \to \MOD_{\RA}$, respectively.
  By base change, $\Omega A \to S$ is a closed immersion, almost of finite presentation and of finite Tor-amplitude, with relative dualizing complex $\omega_{\Omega A/S} \simeq \varpi^*\omega_{S/A} \simeq \varpi^*\fra[-1]$ by \cite[Prop.~6.4.2.1]{SAG}.
  The canonical morphism $\varpi^*(-) \otimes \omega_{\varpi} \to \varpi^!(-)$ is thus invertible (see \cite[Cor.~6.4.2.7]{SAG}).

  Since $\eps \circ \varpi = \id_{\OO_S}$, we have the chain of canonical isomorphisms:
  \begin{align*}
    \eps^!(\RA)
    &\simeq \eps^!(\varpi^*(\OO_S))\\
    &\simeq \eps^!(\varpi^!(\OO_S) \otimes \varpi^*(\frav[1]))\\
    &\simeq \eps^!(\varpi^!(\OO_S)) \otimes \frav[1]\\
    &\simeq \OO_S \otimes \frav[1] \simeq \frav[1],
  \end{align*}
  using the invertibility of $\fra$ (see \ref{sssec:comm grps}).
\end{proof}

\subsubsection{Norm map}\label{sssec:norm Mod_RA}

We define the \emph{norm map}
\begin{equation}\label{eq:norm}
  \Nm_A : \eps^*(-) \otimes \frav[1] \to \eps^!(-)
\end{equation}
by combining the isomorphism of Proposition~\ref{prop:eps^!(R)} with the composite
\[\begin{aligned}
  \eps^*(M) \otimes \eps^!(\RA)
    &\xrightarrow{\mathrm{unit}}
      \eps^!\eps_*\!\left(\eps^*(M) \otimes \eps^!(\RA)\right)\\
    &\xrightarrow{\,\ \sim\ \,} \eps^!\!\left(M \otimes \eps_*\eps^!(\RA)\right) \\
    &\xrightarrow{\mathrm{counit}} \eps^!(M).
\end{aligned}\]
where the isomorphism is the projection formula.

\begin{remark}\label{norm iso} $\Nm_A$ is an equivalence evaluated on $\RA$, by definition and the triangle identities for the adjunction $(\varepsilon_*,\varepsilon^!)$. Note that it is not an equivalence in general, despite $\RA$ being a generator, as $\eps^!$ is discontinuous.\end{remark}

\subsubsection{Augmentation triangle}

We denote by $\mathcal{I}_\eps \in \MOD_{\RA}^{\mathrm{perf}/\OO_S}$ the fibre of the augmentation map $\eps : \RA \to \OO_S$.
We call the exact triangle
\begin{equation}\label{eq:aug triangle}
	\mathcal{I}_\eps
	\to \RA
	\to \eps_*(\OO_S)
\end{equation}
the \emph{augmentation triangle}.
The underlying $\OO_S$-module $\varpi_*\mathcal{I}_\eps$ is $\frav[1]$, by (\ref{sssec:funcheckA add split}).

\begin{lemma}\label{lem:Ieps inv}
  Suppose $A \to S$ has relative dimension one.
  Then the object $\mathcal{I}_\eps$ is invertible in the symmetric monoidal category $(\MOD^{\mathrm{perf}/\OO_S}_{\RA}, \star)$.
\end{lemma}
\begin{proof}
  The functor $\varpi_* : (\MOD_{\RA}, \star) \to (\QCOH(S), \otimes)$ is conservative and symmetric monoidal.
  Since $\mathcal{I}_{\eps}$ lies in $\MOD_{\RA}^{\mathrm{perf}/\OO_S}$ and hence is dualizable (Lemma~\ref{lem:dual mod conv}), there is a canonical evaluation map $\mathcal{I}_{\eps} \star \mathcal{I}_{\eps}^\vee \to \eps_*(\OO_S)$ which $\varpi_*$ sends to the evaluation map for $\frav[1]$.
  The claim follows.
\end{proof}

Applying $\eps^*$ and rotating \eqref{eq:aug triangle} yields the exact triangle
\begin{equation}\label{eq:eps^* aug triangle}
\OO_S
\to \eps^*\eps_*(\OO_S)
\simeq \lam_A
\xrightarrow{\eps^* (\partial)} \eps^* \mathcal{I}_{\eps}[1],
\end{equation}
which is split:

\begin{lemma}\label{lem:lam splitting}
  There is a canonical splitting of $\OO_S$-modules
  \[ \lam_A \simeq \OO_S \oplus \eps^* \mathcal{I}_{\eps}[1]. \]
\end{lemma}
\begin{proof}
  By the triangle identities, the counit map $\eps^*\eps_*(\OO_S) \to \OO_S$ provides a retraction of the left arrow $\eps^*(\mathrm{unit})$ in \eqref{eq:eps^* aug triangle}.
\end{proof}

\begin{remark}
  More generally, tensoring with $V \in \QCOH(S)$ yields $\eps^*\eps_*(V) \simeq V \oplus V \otimes_{\OO_S} \eps^* \mathcal{I}_{\eps}[1]$ by the projection formula.
\end{remark}

\subsection{Examples}

Let $(A, \check{A}, \jmath)$ be a $1$-shifted Cartier dual pair over $S$, with $A \to S$ of relative dimension one.
The additive splitting $\RA \simeq \OO_S \oplus \frav[1]$ (\ref{sssec:funcheckA add split}) provides a canonical map of $E_\infty$-algebras
\begin{equation}\label{eq:Sym to RA}
  \Sym^*_{\OO_S}(\frav[1]) \to \RA
\end{equation}
which is \emph{not} invertible in general (see e.g. Remark~\ref{rem:Hopf map}).

\subsubsection{Characteristic zero}\label{sssec:funcheckA char 0}

When $k$ is a $\mathbf{Q}$-algebra, the free $E_\infty$-algebra $\Sym_{\OO_S}^*(\frav[1])$ collapses to the free (strictly) commutative dg-algebra.
In particular, \eqref{eq:Sym to RA} is invertible on underlying $\OO_S$-modules, and hence is an isomorphism of $E_\infty$-algebras.
Moreover, $\RA \simeq \Sym^*(\frav[1])$ is an isomorphism of Hopf algebras, where:
\begin{itemize}
  \item the underlying algebra is the exterior algebra $\Lambda^*(\frav)$;
  \item the comultiplication map is the unique algebra map extending $\Delta(\varpi) = \varpi \otimes 1 + 1 \otimes \varpi$ for $\varpi \in \frav[1]$.
\end{itemize}
Dually, $\RA^\vee \simeq \Gamma^*(\fra[-1])$ as Hopf algebras, where $\Gamma^*$ is the divided-power Hopf algebra on $\fra[-1]$ (cf.~\cite[Sect.~3]{BCN}):
\begin{itemize}
  \item the multiplication on divided powers $\gamma_n(\nu)$ ($\nu \in \fra[-1]$, $n \ge 0$) is given by $\gamma_n(\nu) \cdot \gamma_m(\nu) = \binom{n+m}{n}\,\gamma_{n+m}(\nu)$;
  \item the comultiplication is $\Delta(\gamma_n(\nu)) = \sum_{i+j=n} \gamma_i(\nu) \otimes \gamma_j(\nu)$.
\end{itemize}
Summarizing:
\begin{align*}
  \RA &\simeq \Lambda^*(\frav), & \RA^\vee &\simeq \Gamma^*(\fra[-1]), \\
  \lam_A &\simeq \Sym^*(\frav[2]), & \lamco_A &\simeq \Sym^*(\fra[-2]).
\end{align*}
A choice of orientation $\fra \simeq \OO_S$ further identifies:
\begin{align*}
  \RA &\simeq \OO_S[\epsilon]/(\epsilon^2), & \RA^\vee &\simeq \OO_S[\nu]/(\nu^2), \\
  \lam_A &\simeq \OO_S[\beta_A], & \lamco_A &\simeq \OO_S[u].
\end{align*}
where the degrees are $|\epsilon|=1$, $|\nu|=-1$, $|\beta_A|=2$, and $|u|=-2$.

\subsubsection{Classical bases}\label{sssec:funcheckA for k discrete}

Suppose now that $S$ is a classical Artin stack.
In this case, \eqref{eq:Sym to RA} factors through the $E_\infty$-algebra maps
\[
  \Sym^*(\frav[1])
  \to \Sym^{\le 1}(\frav[1])
  \xrightarrow{\sim} \RA,
\]
inducing an isomorphism $\RA \simeq \Sym^{\le 1}(\frav[1])$.
In other words, as an augmented $E_\infty$-algebra over $\OO_S$,
$$\RA \simeq \OO_S \oplus \frav[1]$$
is the trivial square-zero extension in the sense of \cite[Rem.~7.3.4.16]{LurieHA}.

Indeed, for $n \ge 2$ the $n$-ary operations are all null-homotopic on $\frav[1]^{\otimes n}$ for degree reasons: when $\OO_S$ is discrete, $\frav[1]$ is concentrated in degree $1$, so that $(\frav[1]^{\otimes n})_{h\Sigma_n}$ is at least $2$-connective, and the space of maps $\Map_{\OO_S}((\frav[1]^{\otimes n})_{h\Sigma_n}, \frav[1])$ is contractible.
In other words, the non-unital $E_\infty$-structure on $\frav[1]$ is trivial, and the claim follows from \cite[Prop.~5.4.4.10]{LurieHA}.
We thus have:
\begin{align*}
  \RA &\simeq \Lambda^*(\frav), & \RA^\vee &\simeq \Gamma^*(\fra[-1]), \\
  \lam_A &\simeq \Gamma^*(\frav[2]), & \lamco_A &\simeq \Sym^*(\fra[-2]).
\end{align*}
Here $\Lambda^*$, $\Gamma^*$, and $\Sym^*$ denote the classical exterior, divided-power, and symmetric Hopf algebras over the discrete ring $\OO_S$.

\subsubsection{Multiplicative case}\label{sssec:OO(Gm[-1])}

For $A = \GGm$ over the sphere spectrum $k=\mathbf{S}$\footnote{%
By $\GGm$ over a connective $E_\infty$-ring $k$, we always mean the ``flat'' variant with $\OO(\GGm) = k[\bZ] = k \otimes \Sigma^\infty_+(\bZ)$.
},
$$\R_{\GGm} := \OO(\Omega\GGm) \simeq \mathbf{S} \otimes_{\mathbf{S}[\bZ]} \mathbf{S} \simeq \mathbf{S}[\mathsf{B}\bZ]$$
is the spherical group ring of the classifying space $\mathsf{B}\bZ$.
Identifying $\mathsf{B}\mathbf{Z} \simeq S^1$ and $\mathsf{B}^2\mathbf{Z} \simeq \mathbf{C}\mathbf{P}^\infty$,\footnote{%
  We write $\mathsf{B}$ to distinguish classifying \emph{space}s from our notation for classifying stacks $BG$.
} we have:
\begin{align*}
  \RA &\simeq \mathbf{S}[S^1], & \RA^\vee &\simeq \mathbf{S}^{S^1}, \\
  \lam_A &\simeq \mathbf{S}[\mathbf{C}\mathbf{P}^\infty], & \lamco_A &\simeq \mathbf{S}^{\mathbf{C}\mathbf{P}^\infty},
\end{align*}
where the left-hand and right-hand sides are spherical chains and cochains, respectively, on $S^1$ and $\mathbf{C}\mathbf{P}^\infty$.

\begin{remark} \label{rem:Hopf map} In this case, the non-invertibility of \eqref{eq:Sym to RA} corresponds to the fact that the spherical group ring $\mathbf{S}[S^1]=\Sigma^\infty_+ S^1$ is not a square zero extension of the base $\mathbf{S}$, essentially because the multiplication $S^1\times S^1\to S^1$ records the Hopf map $\eta$. Indeed after suspending and identifying $\Sigma (S^1\times S^1)\simeq S^2\vee S^2\vee S^3$ we see $\eta$ explicitly as the component $S^3\to\Sigma S^1=S^2.$\end{remark}

\subsection{Koszul duality}\label{ssec:indperf_RA}

\subsubsection{$\lamco_A$-modules}\label{sssec:eps_lam}

Since the functor $\eps_* : (\MOD_{\OO_S}, \otimes) \to (\MOD_{\RA}, \star)$ is symmetric monoidal, its right adjoint $\eps^! : \MOD_{\RA} \to \MOD_{\OO_S}$ is lax symmetric monoidal.
As such, it lifts tautologically to the category of modules over $\eps^!(\eps_*\OO_S) \simeq \hends_{\RA}(\eps_*\OO_S) = \lamco_A$:
\[
  \eps^!_{\Lambda} : \MOD_{\RA} \to \MOD_{\lamco_A}.
\]
This has a left adjoint, symmetric monoidal and $\INDPERF(S)$-linear, given by
\[
  \eps_*^{\Lambda} : \MOD_{\lamco_A} \to \MOD_{\RA},
  \quad N \mapsto N \otimes_{\lamco_A} \OO_S.
\]

The adjunction $(\eps_*, \eps^!)$ factors as
\[\begin{tikzcd}[column sep=large]
  \eps_*: \MOD_{\OO_S}
  \arrow[r, shift left=3pt, "\mathrm{free}"]
  &
  \MOD_{\lamco_A}
  \arrow[l, shift left=3pt, "\mathrm{oblv}"]
  \arrow[r, shift left=3pt, "\eps_*^{\Lambda}"]
  &
  \MOD_{\RA} : \eps^!.
  \arrow[l, shift left=3pt, "\eps^!_{\Lambda}"]
\end{tikzcd}\]

\subsubsection{Ind-perfect \texorpdfstring{$\RA$}{R\_A}-modules}
\label{sssec:renorm}

We now renormalize the category $\MOD_{\RA}$, replacing it with the ind-completion of the category $\MOD_{\RA}^{\mathrm{perf}/\OO_S}$ of $\RA$-modules whose underlying $\OO_S$-module is perfect.

The functor $\eps_* : \QCOH(S) \to \MOD_{\RA}$ sends $\PERF(S)$ to $\MOD_{\RA}^{\mathrm{perf}/\OO_S}$.
If $S$ is perfect, it thus extends uniquely to a colimit-preserving functor
$$\epsb_* : \QCOH(S) \to \Ind(\MOD_{\RA}^{\mathrm{perf}/\OO_S}).$$
Let $\epsb^* : \Ind(\MOD_{\RA}^{\mathrm{perf}/\OO_S}) \to \QCOH(S)$ denote the unique colimit-preserving functor extending $\eps^* : \MOD_{\RA}^{\mathrm{perf}/\OO_S} \to \MOD_{\OO_S}$, and similarly for $\epsb^!$.\footnote{Note that neither $\eps^*$ nor $\eps^!$ sends $\MOD_{\RA}^{\mathrm{perf}/\OO_S}$ to $\MOD_{\OO_S}^{\mathrm{perf}}$.}
The functors $\epsb^*$ and $\epsb^!$ are then left and right adjoint to $\epsb_*$, respectively.

Since $\epsb^!$ agrees with $\eps^!$ on $\OO_S$-perfect modules, Proposition~\ref{prop:eps^!(R)} also gives $\epsb^!(\RA) \simeq \frav[1]$ in $\QCOH(S)$.

\subsubsection{Koszul duality I}

\begin{proposition}\label{prop:Kosz I}\leavevmode
Suppose that the stack $S$ is perfect.
The adjoint functors $(\eps_*^{\Lambda}, \eps^!_{\Lambda})$ introduced in (\ref{sssec:eps_lam}) restrict to a symmetric monoidal $\INDPERF(S)$-linear equivalence of categories
\[ \MOD_{\lamco_A}^{\mathrm{perf}} \simeq \langle \eps_*(\MOD_{\OO_S}^{\mathrm{perf}}) \rangle_{\MOD_{\RA}} \]
where the target is the thick subcategory of $\MOD_{\RA}$ generated by $\eps_*(\MOD_{\OO_S}^{\mathrm{perf}})$.
\end{proposition}
\begin{proof}
The unit map $\id \to \eps^!_{\Lambda} \eps_*^{\Lambda}$ is tautologically an isomorphism on the unit $\lamco_A$.
By $\OO_S$-linearity it follows that it is an isomorphism on $\lamco_A \otimes V$ for every $V \in \MOD_{\OO_S}^{\mathrm{perf}}$.
Since $S$ is perfect, $\MOD_{\lamco_A}^{\mathrm{perf}} \subset \MOD_{\lamco_A}$ is generated as a thick subcategory by $\lamco_A \otimes V$ for $V \in \MOD_{\OO_S}^{\mathrm{perf}}$.
This shows that $\eps_*^{\Lambda} : \MOD_{\lamco_A}^{\mathrm{perf}} \to \MOD_{\RA}$ is fully faithful, with essential image the thick subcategory generated by $\eps_*^{\Lambda}(\lamco_A \otimes V) = \eps_*(V)$ over $V \in \MOD_{\OO_S}^{\mathrm{perf}}$.
\end{proof}

\subsubsection{Koszul duality II: the formal case}

We now assume that $\RA$ is \emph{formal} in the following sense; recall that this is the case when $k$ is an $E_\infty$-$\bQ$-algebra, or any ordinary commutative ring (\ref{sssec:funcheckA for k discrete} and \ref{sssec:funcheckA char 0}).

\begin{definition}\label{def:formal}
  Suppose $A \to S$ has relative dimension one, so that $\RA \simeq \OO_S \oplus \frav[1]$ as $\OO_S$-modules (\ref{sssec:funcheckA add split}).
  We say the rotation algebra $\RA$ is \emph{formal} if this splitting is one of $E_\infty$-algebras, i.e.\ $\RA$ is the trivial square-zero extension of $\OO_S$ by $\frav[1]$; more precisely, if the canonical map \eqref{eq:Sym to RA} factors through an isomorphism
  \[ \Sym^{\le 1}(\frav[1]) \xrightarrow{\ \sim\ } \RA. \]
\end{definition}

\begin{lemma}\label{lem:eps triangle}
  If $\RA$ is formal, then there is a canonical exact triangle
  \[
    \eps_*(\eps^*(-) \otimes \frav[1])
    \to \id
    \xrightarrow{\mathrm{unit}} \eps_*\eps^*(-)
  \]
  of endofunctors of $\MOD_{\RA}$.
  In particular, there is a canonical isomorphism $\mathcal{I}_{\eps} \simeq \eps_*(\frav[1])$ in $\MOD_{\RA}^{\mathrm{perf}/\OO_S}$.
\end{lemma}
\begin{proof}
  Since $\RA$ is formal, the multiplication map $\RA \otimes_{\OO_S} \RA \to \RA$ is null-homotopic on the summand $\frav[1] \otimes_{\OO_S} \frav[1]$ (see \cite[Rem.~7.3.4.16(iii)]{LurieHA}).
  It follows that on the summand $\frav[1]$ of $\RA$, the $\RA$-action factors through $\eps : \RA \to \OO_S$, so that the fibre of $\eps : \RA \to \OO_S$ identifies with $\eps_*(\frav[1])$ in $\MOD_{\RA}$.
  We thus have the exact triangle in $\MOD_{\RA}$,
  \begin{equation}\label{eq:OO(A[-1]) triangle}
    \eps_*(\frav)[1] \to \RA \to \eps_*(\OO_S).
  \end{equation}
  Tensoring with $M \in \MOD_{\RA}$ and applying the projection formula $\eps_*(V) \otimes_{\RA} M \simeq \eps_*(V \otimes_{\OO_S} \eps^*(M))$ yields the claim.
\end{proof}

\begin{remark}
  The exact triangle $\eps_*(\frav[1]) \to \RA \to \eps_*(\OO_S)$ is split at the level of underlying $\OO_S$-modules, but not in $\MOD_{\RA}$.
  Indeed, the obstruction to splitting is precisely the Euler map which will be introduced in \S \ref{sec:Euler}.
\end{remark}

\begin{corollary}\label{cor:sqzero thick=perf/OS}
  Suppose $S$ is a classical affine scheme which is regular noetherian (e.g. the spectrum of a field).
  Then $\MOD_{\RA}^{\mathrm{perf}/\OO_S}$ is generated as a thick subcategory of $\MOD_{\RA}$ by $\eps_*(\OO_S)$.
\end{corollary}
\begin{proof}
  On $\pi_0$, $\eps$ induces an isomorphism $\pi_0(\RA) \simeq \OO_S$, hence $\eps_* : \MOD_{\OO_S} \to \MOD_{\RA}$ induces an equivalence on hearts (for the standard t-structures).
  Consider the restricted functor $F := \eps_* : \MOD_{\OO_S}^{\mathrm{perf}} \to \MOD_{\RA}^{\mathrm{perf}/\OO_S}$.
  Every $M \in \MOD_{\RA}^{\mathrm{perf}/\OO_S}$ has finitely many nonzero homotopy groups, each of which is again a perfect $\OO_S$-module since $S$ is regular noetherian, so $F$ is a t-exact functor between categories with bounded t-structures that induces an equivalence on hearts.
  The claim now follows by induction on the number of nonzero homotopy groups.
\end{proof}

\begin{corollary}[Koszul duality II]\label{cor:Kosz II}
If $S$ is a classical Artin stack which is regular and locally noetherian, then the adjoint functors $(\eps_*^{\Lambda}, \eps^!_{\Lambda})$ restrict to a symmetric monoidal $\INDPERF(S)$-linear equivalence of categories
\[ \MOD_{\lamco_A}^{\mathrm{perf}} \simeq \MOD_{\RA}^{\mathrm{perf}/\OO_S}. \]
\end{corollary}
\begin{proof}
Since both sides satisfy descent in $S$, and $\eps_*^{\Lambda} = (-)\otimes_{\lamco_A}\eps_*(\OO_S)$ commutes with base change in $S$, we may assume $S$ is affine.
Then the claim follows from Proposition~\ref{prop:Kosz I} and Corollary~\ref{cor:sqzero thick=perf/OS}.
\end{proof}

\subsubsection{Ind-perfect \texorpdfstring{$\lamco_A$}{Lambda_A}-modules}\label{sssec:indperf Lambda modules}

We renormalize the preceding discussion.
Since $\epsb_* : \QCOH(S) \to \Ind(\MOD_{\RA}^{\mathrm{perf}/\OO_S})$ is symmetric monoidal, its right adjoint $\epsb^!$ is lax symmetric monoidal.
As such, it lifts tautologically to the category of modules over $\epsb^!\epsb_*(\OO_S) \simeq \lamco_A$,
\[
  \epsb^!_{\Lambda} : \Ind(\MOD_{\RA}^{\mathrm{perf}/\OO_S}) \to \MOD_{\lamco_A},
\]
the colimit-preserving extension of $\eps^!_{\Lambda}|_{\MOD_{\RA}^{\mathrm{perf}/\OO_S}}$ in the sense of (\ref{sssec:renorm}).
By the projection formula, it sends $\epsb_*(V)$ to the free module $\lamco_A \otimes_{\OO_S} V$, for $V \in \QCOH(S)$.
Its left adjoint, symmetric monoidal and $\INDPERF(S)$-linear, is
\[
  \epsb_*^{\Lambda} : \MOD_{\lamco_A} \to \Ind(\MOD_{\RA}^{\mathrm{perf}/\OO_S}),
  \quad N \mapsto N \otimes_{\lamco_A} \eps_*(\OO_S).
\]
The renormalized adjunction $(\epsb_*, \epsb^!)$ thus factors as
\[\begin{tikzcd}[column sep=large]
  \epsb_*: \QCOH(S)
  \arrow[r, shift left=3pt, "\mathrm{free}"]
  &
  \MOD_{\lamco_A}
  \arrow[l, shift left=3pt, "\mathrm{oblv}"]
  \arrow[r, shift left=3pt, "\epsb_*^{\Lambda}"]
  &
  \Ind(\MOD_{\RA}^{\mathrm{perf}/\OO_S}) : \epsb^!.
  \arrow[l, shift left=3pt, "\epsb^!_{\Lambda}"]
\end{tikzcd}\]
When $S$ is a classical Artin stack which is regular and locally noetherian, the functors $\epsb_*^{\Lambda}$ and $\epsb^!_{\Lambda}$ are mutually inverse equivalences, by Corollary~\ref{cor:Kosz II}.

\section{Ind-perfect complexes on $B\check{A}$}\label{sec:BcheckA}

In this section we set up the sheaf categories where $A$-loop cohomology will reside.
We first develop some general theory of ind-perfect complexes (Subsect.~\ref{ssec:indperf}).
We then apply the algebraic models of \S\ref{sec:RA} to describe ind-perfect complexes on the classifying stack $B\check{A}$ (Subsects.~\ref{ssec:indperf_RA}--\ref{ssec:rotation}), which underlie the Euler periodization developed in \S \ref{sec:Euler}.

\subsection{Ind-perfect complexes}
\label{ssec:indperf}

For a stack $\XX$, we consider the ind-completion
$$\INDPERF(\XX) \supset \PERF(\XX).$$
There is always a colimit-preserving symmetric monoidal functor
$$\Ren_\XX : \INDPERF(\XX) \to \QCOH(\XX),$$
which is an equivalence when $\XX$ is perfect. In order to form the Euler periodization for $A$-loop cohomology, we will need the case $\XX = B\check{A}$, the classifying stack of the $1$-shifted Cartier dual, for which $\Ren_\XX$ is typically not an equivalence (see Remark~\ref{rem:pi not conservative}). The present subsection records the basic functoriality of $\INDPERF$ needed in the sequel.

We denote by $\Ren_\XX^R$ the right adjoint\footnote{We note that $\Ren_\XX^R$ is colimit-preserving if and only if $\OO_\XX$ is compact.} of $\Ren_\XX$.
Note that $\Ren_\XX$ and $\Ren_\XX^R$ restrict to mutually inverse equivalences on the respective embedded copies of $\PERF(\XX)$. We will identify the embedded copies of $\PERF(\XX)$ along these equivalences and so we will not attempt to distinguish them notationally. Note that $\Ren_\XX$ is symmetric monoidal, though $\Ren_\XX^R$ is not.

Both categories $\INDPERF$ and $\QCOH$ are endowed with the usual adjunctions $(f^*, f_*)$ for morphisms $f:\XX\rightarrow\YY$. We will write $f^*$ and $ f_*$ for the functors on $\INDPERF$,
$$ f^* : \INDPERF(\YY) \to \INDPERF(\XX), \qquad f_* : \INDPERF(\XX) \to \INDPERF(\YY), $$
and $f^*_{\QCOH}$ and $f_*^{\QCOH}$ for their counterparts on $\QCOH$.

\begin{lemma}\label{lem:indperf push easy}
For every morphism $f:\XX\rightarrow\YY$, the functor $f_*$ preserves colimits and satisfies the projection formula.
\end{lemma}
\begin{proof}
  Since $f^*$ preserves compact objects (= perfect complexes), its right adjoint preserves colimits.
  For the second claim, we consider the canonical map
  \[
    f_*(F) \otimes G \longrightarrow f_*(F \otimes f^*G)
  \]
  defined as the composite
  \[
    f_*(F) \otimes G
    \xrightarrow{\mathrm{unit}} f_*f^*\bigl(f_*(F) \otimes G\bigr)
    \simeq f_*\bigl(f^*f_*(F) \otimes f^*G\bigr)
    \xrightarrow{\mathrm{counit}} f_*\bigl(F \otimes f^*G\bigr).
  \]
  We show that this map is invertible for all $F \in \INDPERF(\XX)$ and $G \in \INDPERF(\YY)$.
  Since the tensor product is colimit-preserving in each variable, as are $f^*$ and $f_*$, this map is compatible with colimits in $G$.
  We may thus assume that $G$ is perfect, and hence dualizable in $\INDPERF(\YY)$.
  It is enough to check that for every $G' \in \INDPERF(\YY)$, the map becomes invertible after applying $\Map_{\INDPERF(\YY)}(G', -)$.
  Indeed, it is then identified with the composite isomorphism
  \begin{align*}
    \Map_{\INDPERF(\YY)}(G', f_*(F) \otimes G)
    &\simeq \Map_{\INDPERF(\YY)}(G' \otimes G^\vee, f_*(F)) \\
    &\simeq \Map_{\INDPERF(\XX)}(f^*G' \otimes f^*(G)^\vee, F) \\
    &\simeq \Map_{\INDPERF(\XX)}(f^*G', F \otimes f^*(G)) \\
    &\simeq \Map_{\INDPERF(\YY)}(G', f_*(F \otimes f^*(G))),
  \end{align*}
  so the claim follows.
\end{proof}

\subsection{Base change formula}

Next we will establish base change for $f_*$, under suitable assumptions on $f$.
More precisely, given a fibre square of stacks,
\begin{equation}\label{eq:indperf fibsq}
  \begin{tikzcd}
    \tilde{\XX} \arrow{d}{\tilde{g}}\arrow{r}{\tilde{f}} & \tilde{\YY}\arrow{d}{g} \\
    \XX\arrow{r}{f} & \YY,
  \end{tikzcd}
\end{equation}
there is a canonical natural transformation on $\INDPERF$,
\[
  \operatorname{Ex}^*_* : g^* f_*
  \xrightarrow{\mathrm{unit}} \tilde{f}_* \tilde{f}^* g^* f_*
  \simeq \tilde{f}_* \tilde{g}^* f^* f_*
  \xrightarrow{\mathrm{counit}} \tilde{f}_* \tilde{g}^*.
\]

\begin{theorem}\label{thm:bcindperf}
  For a fibre square of stacks as in \eqref{eq:indperf fibsq}, suppose either:
  \begin{enumerate}
    \item $f$ is qcqs and representable, and $g$ is proper,
    representable, and quasi-smooth; or
    \item $f$ is strongly cohomologically proper (see Definition~\ref{def:str coh proper}) and of
    finite Tor-amplitude.
  \end{enumerate}
  Then the square satisfies base change for $\INDPERF$.
\end{theorem}

\begin{remark}\label{rem:indperf open base change}
For pushforward along a non-proper morphism, e.g. a quasi-compact open immersion or affine morphism, the base change formula in $\INDPERF$ can fail even over $k=\mathbf{Q}$.
For example, consider the stacks
\[
  \XX := BE\times\spc(k[t,z]),
  \qquad \YY := BE\times\spc(k[t]),
\]
where $E$ is an oriented elliptic curve over $\bQ$, and $t$ and $z$ have homological degree $0$ and $2$, respectively.
Then for the cartesian square
\[
\begin{tikzcd}
  \XX[t^{-1}] \arrow{d}{j'}\arrow{r}{\tilde{f}}
  & \YY[t^{-1}]\arrow{d}{j}\\
  \XX\arrow{r}{f} & \YY,
\end{tikzcd}
\]
both base change maps
\[
  f^*j_*\OO_{\YY[t^{-1}]}\longrightarrow j'_*\OO_{\XX[t^{-1}]},
  \qquad
  j^*f_*\OO_\XX\longrightarrow \tilde{f}_*\OO_{\XX[t^{-1}]},
\]
are invertible in $\QCOH$, but not in $\INDPERF$.
Indeed, by Corollary~\ref{cor:IndPerf BcheckA Kosz} we may identify
\[
  \INDPERF(\YY)\simeq\MOD_{k[t,u]},\quad
  \INDPERF(\YY[t^{-1}])\simeq\MOD_{k[t,t^{-1},u]},
\]
with $u$ in degree $-2$.
Then $j^*$ is localization at $t$, and $j_*\OO_{\YY[t^{-1}]}\simeq\OO_\YY[t^{-1}]$.
The first map thus becomes
\[
  \OO_\XX[t^{-1}] \longrightarrow j'_* \OO_{\XX[t^{-1}]}.
\]
Applying $\Map_{\INDPERF(\XX)}(\OO_\XX, -)$ to both sides and using the compactness of $\OO_\XX$ in $\INDPERF(\XX)$, this becomes the canonical map
\[
  \OO(\XX)[t^{-1}] \longrightarrow \OO(\XX[t^{-1}])
\]
where we have identified $\Map_{\INDPERF(\XX)}(\OO_\XX, \OO_\XX) \simeq \OO(\XX)$ since $\OO_\XX$ is perfect, and similarly for $\XX[t^{-1}]$.
Now using Theorem~\ref{thm:QCoh checkA}, one computes
\[
  \OO(\XX) \simeq \prod_{n\ge 0} k[t,z] \cdot u^n,
  \qquad \OO(\XX[t^{-1}]) \simeq \prod_{n\ge 0} k[t,t^{-1},z] \cdot u^n.
\]
Passing to $\pi_0$, the map in question becomes
\[
  \left(\prod_{n\geq0}k[t]\right)[t^{-1}]
  \longrightarrow \prod_{n\geq0}k[t,t^{-1}],
\]
which is not surjective: e.g. the sequence $(t^{-n})_{n\geq0}$ is not in the image (since its denominators are unbounded).
The argument for the other base change map is very similar.
\end{remark}

The proof of Theorem~\ref{thm:bcindperf} will occupy the rest of this subsection.

\begin{lemma}\label{lefts} There is a canonical natural isomorphism of functors $$\operatorname{Ex}^*_{\Ren} : \Ren_\XX f^* \simeq f_{\QCOH}^*\Ren_\YY:\INDPERF(\YY)\rightarrow\QCOH(\XX).$$\end{lemma}

\begin{proof} All functors involved are left adjoints and hence colimit-preserving. Since $\PERF$ generates $\INDPERF$ by definition, it suffices to check the claim on $\PERF$, where it is clear.\end{proof}

Passing to right adjoints in Lemma~\ref{lefts}, we obtain:
\begin{corollary} \label{rights} There is a canonical natural isomorphism of functors $$\operatorname{Ex}_*^{\Ren^R} : f_*\Ren_\XX^R\simeq \Ren_\YY^R f^{\QCOH}_*.$$\end{corollary}

\begin{construction}\label{constr:Ex f^* PhiR}
For every morphism $f : \XX \to \YY$ the natural transformation
$\operatorname{Ex}^*_{\Ren^R} : f^*\Ren_\YY^R \to \Ren_\XX^R f^*_{\QCOH}$ is defined as the composite
\[
  f^*\Ren_\YY^R
  \xrightarrow{\mathrm{unit}} \Ren_\XX^R\Ren_\XX f^*\Ren_\YY^R
  \simeq \Ren_\XX^R f^*_{\QCOH}\Ren_\YY\Ren_\YY^R
  \xrightarrow{\mathrm{counit}} \Ren_\XX^R f^*_{\QCOH},
\]
where the middle isomorphism is Lemma~\ref{lefts}.
\end{construction}

\begin{definition}\label{def:str coh proper}
We call a map $f : \XX \to \YY$ \emph{strongly cohomologically proper} if it is representable by Artin stacks, almost of finite presentation, universally of finite cohomological dimension (see \cite[Def.~1.18]{KhanKstack}), and universally cohomologically proper in the sense of \cite[Def.~1.26]{KhanKstack}.
\end{definition}

Note that this class is stable under base change and composition. Any proper representable morphism is strongly cohomologically proper.

\begin{lemma}\label{lem:pushindperf} Assume that $f$ is strongly cohomologically proper and of finite Tor-amplitude. Then the canonical natural transformation
\[
  \operatorname{Ex}_*^{\Ren} : \Ren_\YY f_*
  \xrightarrow{\mathrm{unit}} f^{\QCOH}_* f^*_{\QCOH} \Ren_\YY f_*
  \simeq f^{\QCOH}_* \Ren_\XX f^* f_*
  \xrightarrow{\mathrm{counit}} f^{\QCOH}_* \Ren_\XX,
\]
where the middle isomorphism is Lemma~\ref{lefts}, is invertible.
\end{lemma}\begin{proof} Under our assumptions all functors are colimit-preserving and we need only check this on perfect objects. Under the assumptions, $f^{\QCOH}_*$ preserves perfect complexes by \cite[Prop.~1.31]{KhanKstack}. The same is true of $f_*$, as $f^*$ has a right adjoint at the level of $\PERF$ and $f_*$ is the ind-extension of this to $\INDPERF$. Agreement on perfect objects is then clear.\end{proof}

Although $\Ren_\XX^R$ is not symmetric monoidal, it is naturally a morphism of $\PERF$-modules:

\begin{lemma}\label{perftens} For every $b \in \QCOH(\XX)$ and $c\in\PERF(\XX)$, there are natural isomorphisms $$\Ren_\XX^{R}(b)\otimes c\simeq \Ren_\XX^R(b\otimes c).$$  \end{lemma}

\begin{proof} If $c$ is a dualizable object of a symmetric monoidal category we write $c^\vee$ for its dual. Then we compute $$\Hom_{\INDPERF}(a,\Ren_\XX^R(b)\otimes c)\simeq\Hom_{\INDPERF}(a\otimes c^\vee,\Ren_{\XX}^R(b))$$ $$\simeq \Hom_{\QCOH}(\Ren_\XX(a\otimes c^\vee),b)\simeq\Hom_{\QCOH}(\Ren_\XX(a)\otimes c^\vee,b)$$ $$\simeq\Hom_{\QCOH}(\Ren_{\XX}(a),b\otimes c)\simeq\Hom_{\INDPERF}(a,\Ren_\XX^R(b\otimes c)).\qedhere$$\end{proof}

\begin{corollary}\label{cor:proper pullback PhiR}
Let $f : \XX \to \YY$ be a proper, representable, quasi-smooth morphism.
Then the exchange transformation of Construction~\ref{constr:Ex f^* PhiR} is an isomorphism
$$\operatorname{Ex}^*_{\Ren^R} : f^*\Ren_\YY^R\simeq \Ren_\XX^R f^*_{\QCOH}.$$
\end{corollary}

\begin{proof} Under our assumptions $f^{\QCOH}_*$ and $f_*$ have \emph{right} adjoints given respectively by $$f_{\QCOH}^!:=f_{\QCOH}^*\otimes\omega_f,$$ $$f^!:=f^*\otimes\omega_f,$$ where, as $\omega_f$ is invertible, we have not notationally distinguished its images in $\INDPERF$ and $\QCOH$. Passing to right adjoints from the isomorphism $\Ren_\YY f_*\simeq f^{\QCOH}_*\Ren_{\XX}$ we obtain an isomorphism $$f^!\Ren_\YY^R\simeq \Ren_\XX^R f_{\QCOH}^!.$$ Now we have by definition $$(f^*\Ren^R_\YY)\otimes\omega_f\simeq\Ren_\XX^R(f_{\QCOH}^*\otimes\omega_f).$$ By Lemma~\ref{perftens} we identify the right hand side with $(\Ren^R_{\XX} f_{\QCOH}^*)\otimes\omega_f$ and cancelling off the $\omega_f$ factors we obtain the desired claim.  \end{proof}

\begin{theorem}\label{AppMain} For a fibre square as in \eqref{eq:indperf fibsq}, suppose either
\begin{enumerate}
\item $g$ is proper, representable, and quasi-smooth, or
\item $f$ is strongly cohomologically proper (Definition~\ref{def:str coh proper}) and of finite Tor-amplitude.
\end{enumerate}
If the square satisfies base change for $\QCOH$, then it satisfies base change for $\INDPERF$.
\end{theorem}

\begin{proof}
It suffices to check that $g^*f_*(M) \to \tilde{f}_*\tilde{g}^*(M)$ is invertible for every $M \in \PERF(\XX)$, since all functors involved are colimit-preserving.

In case (i), we have the isomorphisms
\begin{align*}
  g^* f_* \Ren_\XX^R
  \simeq g^* \Ren_\YY^R f^{\QCOH}_*
  \simeq \Ren_{\tilde{\YY}}^R g^*_{\QCOH} f^{\QCOH}_*
  &\simeq \Ren_{\tilde{\YY}}^R \tilde{f}^{\QCOH}_* \tilde{g}^*_{\QCOH}\\
  &\simeq \tilde{f}_* \Ren_{\tilde{\XX}}^R \tilde{g}^*_{\QCOH}
  \simeq \tilde{f}_* \tilde{g}^* \Ren_\XX^R
\end{align*}
by Corollary~\ref{rights}, Corollary~\ref{cor:proper pullback PhiR} for $g$, base change for $\QCOH$, Corollary~\ref{rights} for $\tilde{f}$, and Corollary~\ref{cor:proper pullback PhiR} for $\tilde{g}$, in order.
The claim follows by writing $M \simeq \Ren_\XX^R \Ren_\XX(M)$ (since $M$ is perfect).

Under assumption (ii), the proof of Lemma~\ref{lem:pushindperf} shows that $f^{\QCOH}_*$ and $\tilde{f}^{\QCOH}_*$ preserve $\PERF$ and agree there with $f_*$ and $\tilde{f}_*$. By Lemma~\ref{lefts}, $g^*$ and $\tilde{g}^*$ also agree with $g^*_{\QCOH}$ and $\tilde{g}^*_{\QCOH}$ on $\PERF$. For $M \in \PERF(\XX)$,
\[
g^* f_*(M) \simeq g^*_{\QCOH} f^{\QCOH}_*(M) \simeq \tilde{f}^{\QCOH}_* \tilde{g}^*_{\QCOH}(M) \simeq \tilde{f}_* \tilde{g}^*(M),
\]
the middle isomorphism by base change for $\QCOH$.
\end{proof}

\begin{proof}[Proof of Theorem~\ref{thm:bcindperf}] By Theorem~\ref{AppMain} we reduce to the analogous statements for $\QCOH$. In case (i) this is \cite[Prop.~6.3.4.1(b)]{SAG}; in case (ii) it is \cite[Prop.~A.1.5]{HLP}. \end{proof}

\subsection{The functor \texorpdfstring{$f^!$}{f\textasciicircum !}}\label{ssec:indperf dualizing complex}

\subsubsection{$f^!$, $f^\dag$, and dualizing complexes}

Let $f : \XX \to \YY$ be a proper representable morphism, or more generally a strongly cohomologically proper morphism, so that $f_*^{\QCOH}$ preserves colimits.
Then $f_*^{\QCOH}$ admits a right adjoint $f^!_{\QCOH}$, and we have the \emph{relative dualizing complex}
\[
  \omega_f^{\QCOH} := f^!_{\QCOH}(\OO_\YY) \quad \in \QCOH(\XX),
\]
see \cite[Cor.~6.3.4.3, \S 6.4.2]{SAG} and \cite[Prop.~A.1.5]{HLP}.

\begin{definition}
Let $f : \XX \to \YY$ be a strongly cohomologically proper morphism of stacks.
\begin{enumerate}
\item We let $f^! : \INDPERF(\YY) \to \INDPERF(\XX)$ denote the right adjoint of the colimit-preserving functor $f_*$ (Lemma~\ref{lem:indperf push easy}).

\item We let $f^\dag : \INDPERF(\YY) \to \INDPERF(\XX)$ denote the composite $f^\dag(-) := \Ren^R \circ f^!_{\QCOH} \circ \Ren (-)$.
\end{enumerate}
\end{definition}

\begin{lemma}\label{lem:f^dag vs f^!}
Let $f : \XX \to \YY$ be a strongly cohomologically proper morphism of finite Tor-amplitude.
Then there is a canonical isomorphism $f^\dag(-) \simeq f^! \bigl( \Ren^R \Ren(-) \bigr)$.
In particular, $f^\dag(F) \simeq f^!(F)$ for all $F \in \PERF(\YY)$.
\end{lemma}
\begin{proof}
Passing to right adjoints from Lemma~\ref{lem:pushindperf}, we have $f^! \circ \Ren^R \simeq \Ren^R \circ f^!_{\QCOH}$, whence the first claim.
The second follows since $\Ren^R \Ren$ is the identity on $\PERF(\YY)$.
\end{proof}

\begin{definition}
For a strongly cohomologically proper morphism $f : \XX \to \YY$, we define $\omega_f \in \INDPERF(\XX)$ as $f^\dag(\OO) \simeq \Ren^R \bigl( \omega_f^{\QCOH} \bigr)$.
\end{definition}

When $f$ is of finite Tor-amplitude, $\omega_f \simeq f^!(\OO)$.
When $\omega_f^{\QCOH}$ is perfect\footnote{%
For example, this is the case when $k$ is a $\bQ$-algebra and $f$ is quasi-smooth.
}, we have $\omega_f \simeq \omega_f^{\QCOH}$.

\begin{proposition}\label{prop:f^! formula indperf}
	Let $f : \XX \to \YY$ be a strongly cohomologically proper morphism.
	Then for every $F \in \PERF(\YY)$ there is a canonical natural isomorphism
	\[
    f^\dag(F) \simeq \omega_f \otimes f^*(F) \quad \in \INDPERF(\XX).
	\]
	In particular, if $f$ is of finite Tor-amplitude, then the same formula holds for $f^!(F)$.
\end{proposition}
\begin{proof}
  Passage to right adjoints from \cite[Cor.~6.3.4.3(2)]{SAG} yields $f^!_{\QCOH}(F) \simeq \omega_f^{\QCOH} \otimes f^*_{\QCOH}F$.
  Then the first claim follows from Lemmas~\ref{lefts} and~\ref{perftens}, and the second from Lemma~\ref{lem:f^dag vs f^!}.
\end{proof}

\begin{remark}\label{rem:f^dag lax functoriality}
	Let $f : \XX \to \YY$ and $g : \YY \to \ZZ$ be strongly cohomologically proper morphisms.
	Then we have $(gf)^\dag \simeq \Ren^R f^!_{\QCOH} g^!_{\QCOH} \Ren$, and the counit $\Ren\Ren^R \to \id$ induces a canonical map $f^\dag g^\dag \to (gf)^\dag$.
  Similarly, there is a canonical map $\id \to \id^\dag = \Ren^R\Ren$.
  These maps are both invertible if $\Ren$ is an equivalence.
\end{remark}

\subsubsection{Exchange transformations}\label{sssec:Ex f^* f^dag indperf}

Consider a cartesian square of stacks
\[\begin{tikzcd}
  \XX' \ar{r}{g}\ar[swap]{d}{p} & \YY' \ar{d}{q} \\
  \XX \ar{r}{f} & \YY
\end{tikzcd}\]
with $f$ proper representable.
In $\QCOH$ we have the base change formula $\operatorname{Ex}^*_* : q^*f_* \simeq g_*p^*$ \cite[Prop.~6.3.4.1(b)]{SAG} and the projection formula $\operatorname{Ex}^{\otimes,*}_* : f_*(-) \otimes (-) \simeq f_*(- \otimes f^*(-))$ \cite[Cor.~6.3.4.3(2)]{SAG}.
By adjunction, these give rise to the exchange transformations
\begin{align*}
  \operatorname{Ex}^{*,!}_\otimes &: f^*(-) \otimes f^!(-) \to f^!(- \otimes -),\\
  \operatorname{Ex}^{*,!} &: p^*f^! \to g^!q^*
\end{align*}
of functors on $\QCOH$, see e.g. \cite[A.1.2, A.2.7]{weavelisse}.

\begin{construction}\label{constr:Ex f^* f^dag}
On $\INDPERF$, we define the exchange transformations
\begin{align*}
  \operatorname{Ex}^{*,\dag}_\otimes &: f^*(-) \otimes f^\dag(-) \to f^\dag(- \otimes -),\\
  \operatorname{Ex}^{*,\dag} &: p^*f^\dag \to g^\dag q^*
\end{align*}
as the composites
\begin{align*}
  f^*F \otimes \Ren^R f^!_{\QCOH}\Ren G
  &\xrightarrow{\mathrm{unit}} f^*\Ren^R\Ren F \otimes \Ren^R f^!_{\QCOH}\Ren G\\
  &\xrightarrow{\operatorname{Ex}^*_{\Ren^R}} \Ren^R f^*_{\QCOH}\Ren F \otimes \Ren^R f^!_{\QCOH}\Ren G\\
  &\xrightarrow{\mathrm{lax}} \Ren^R\big(f^*_{\QCOH}\Ren F \otimes f^!_{\QCOH}\Ren G\big)\\
  &\xrightarrow{\operatorname{Ex}^{*,!}_\otimes} \Ren^R f^!_{\QCOH}(\Ren F \otimes \Ren G)
  \simeq \Ren^R f^!_{\QCOH}\Ren(F \otimes G)
\end{align*}
and
\[
  p^*\Ren^R f^!_{\QCOH}\Ren
  \xrightarrow{\operatorname{Ex}^*_{\Ren^R}} \Ren^R p^*_{\QCOH} f^!_{\QCOH}\Ren
  \xrightarrow{\operatorname{Ex}^{*,!}} \Ren^R g^!_{\QCOH} q^*_{\QCOH}\Ren
  \xsimeq{\operatorname{Ex}^*_\Ren} \Ren^R g^!_{\QCOH}\Ren q^*,
\]
where the arrow labelled $\mathrm{lax}$ comes from the lax symmetric monoidal structure of $\Ren^R$, $\operatorname{Ex}^*_{\Ren^R}$ is Construction~\ref{constr:Ex f^* PhiR}, and $\operatorname{Ex}^*_\Ren$ is Lemma~\ref{lefts}.
\end{construction}

\begin{construction}\label{constr:omega composite indperf}
Let $f : \XX \to \YY$ and $g : \ZZ \to \XX$ be proper representable morphisms.
We have the canonical map in $\INDPERF(\ZZ)$,
\[
  g^*\omega_f \otimes \omega_g
  \xrightarrow{\operatorname{Ex}^{*,\dag}_\otimes} g^\dag f^\dag(\OO)
  \xrightarrow{\mathrm{lax}} (f\circ g)^\dag(\OO) = \omega_{f\circ g},
\]
where the map labelled $\mathrm{lax}$ is the lax functoriality of Remark~\ref{rem:f^dag lax functoriality}.
\end{construction}

\subsubsection{Trace maps}\label{sssec:trace indperf}

When $f: \XX \to \YY$ is proper representable, or more generally strongly cohomologically proper, we have the counit maps
\[
  f_*^{\QCOH} f^!_{\QCOH}(F) \longrightarrow F,
  \qquad F \in \QCOH(\YY).
\]
For $F = \OO$ this becomes $f_*^{\QCOH} (\omega_f^{\QCOH}) \longrightarrow \OO$.

\begin{construction}\label{constr:trace indperf}
Applying the functors $\Ren$ and $\Ren^R$ on the right and left, respectively, we get for every $F \in \INDPERF(\YY)$ the map
\[
  f_* f^\dag (F) :=
  f_*\Ren^R f^!_{\QCOH} \Ren(F) \simeq \Ren^R f_*^{\QCOH} f^!_{\QCOH} \Ren(F) \xrightarrow{\mathrm{counit}} \Ren^R \Ren(F),
\]
where $\Ren^R$ commutes with $f_*$ by Corollary~\ref{rights}.
If $F \in \PERF(\YY)$, then $\Ren^R\Ren(F) \simeq F$ and we get the canonical natural map in $\INDPERF(\YY)$
\[
  \operatorname{tr}_f : f_* f^\dag (F) \longrightarrow F.
\]
The source can moreover be identified with $f_*(\omega_f \otimes f^*F)$ (Proposition~\ref{prop:f^! formula indperf}).
\end{construction}

\subsection{\texorpdfstring{$\check{A}$}{Ǎ}-equivariant ind-perfect complexes}\label{ssec:rotation}

We fix a $1$-shifted Cartier dual pair $(A, \check{A}, \jmath)$ over a stack $S$.
In this subsection we describe the category of ind-perfect complexes on the classifying stack $B\check{A} := S/\check{A}$.
We will denote by
$$\sfp : B\check{A} \to S, \qquad \sfe : S \twoheadrightarrow B\check{A}$$
the projection and quotient maps, respectively.

\subsubsection{Quasi-coherent sheaves}

We first describe the entire category of quasi-coherent sheaves on $B\check{A}$.

\begin{theorem}\label{thm:QCoh checkA}
  There are canonical equivalences
  \begin{equation*}
    \QCOH(B\check{A})
    \;\simeq\; \MOD_{\RA}
    \;\simeq\; \COMOD_{\RA^\vee},
  \end{equation*}
  where $\RA$ is the Hopf algebra introduced in (\ref{sssec:RA}).
\end{theorem}

\begin{remark}
The symmetric monoidal structure on $\QCOH(B\check{A})$ corresponds to the convolution monoidal structure on $\MOD_{\RA}$, whose product is $\star$ and unit is $\eps_*(\OO_S)$; see (\ref{sssec:module conventions}).
(Somewhat dually, the standard tensor product on $\MOD_{\RA}$ corresponds to the convolution product on $\QCOH(B\check{A})$, defined by $F \star G := \mathsf{m}_\sharp(F \boxtimes G)$ where $\mathsf{m}$ is the multiplication on $B\check{A}$.)

For an $E_\infty$-Hopf algebra $R$ in $\QCOH(S)$, the category $\COMOD_R$ of $R$-comodules in $\QCOH(S)$ is equipped with a symmetric monoidal structure where the tensor product $\otimes$ is defined using the multiplication of $R$, dually to the convolution product $\star$.
When $R$ is perfect as an $\OO_S$-module, the $\OO_S$-linear dual $R^\vee$ is also a Hopf algebra, and there is a canonical symmetric monoidal equivalence $(\MOD_R, \star) \simeq (\COMOD_{R^\vee}, \otimes)$.
\end{remark}

\begin{proof}[Proof of Theorem~\ref{thm:QCoh checkA}]
  Consider the adjunction
  \[
    \sfe^* : \QCOH(B\check{A}) \to \QCOH(S),
    \quad \sfe_* : \QCOH(S) \to \QCOH(B\check{A})
  \]
  induced by the atlas $\sfe : S \twoheadrightarrow B\check{A}$.
  It follows from Definition~\ref{def:abstract 1-shifted Cartier dual}\ref{item:abstract 1-shifted Cartier dual/sharp} and \cite[Vol.~I, Chap.~1, Lem.~2.6.4]{GR} that $\sfe_*$ satisfies base change against $*$-pullback, and that $\sfe^*$ admits a left adjoint $\sfe_\sharp$ also satisfying base change against $*$-pullback.
  In view of the equivalence $\QCOH(B\check{A}) \simeq \operatorname{Tot}(\QCOH(\check{A}^\bullet))$, $\sfe^*$ is clearly conservative.
  We may thus apply the Barr--Beck--Lurie theorem (see \cite[Thm.~4.7.3.5]{LurieHA}) to conclude that the adjunction $(\sfe^*, \sfe_*)$ is comonadic, i.e., $\QCOH(B\check{A})$ is identified with the category of comodules for the comonad $\sfe^*\sfe_*$.
  By the base change and projection formulas, we have $\sfe^*\sfe_* \simeq \check{p}_*\check{p}^* \simeq \OO(\check{A}) \otimes (-)$.

  Since $\RA$ is perfect as an $\OO_S$-module, the dual description follows.
  Alternatively, the dual form of Barr--Beck--Lurie shows that the adjunction $(\sfe_\sharp, \sfe^*)$ is monadic, so that $\QCOH(B\check{A})$ is also identified with the category of modules for the monad $\sfe^*\sfe_\sharp \simeq \check{p}_\sharp \check{p}^* \simeq \check{p}_\sharp(\OO_{\check{A}}) \otimes (-)$.
  Since $\check{p}_*$ preserves perfect complexes (Definition~\ref{def:abstract 1-shifted Cartier dual}\ref{item:abstract 1-shifted Cartier dual/sharp}), \cite[Rem.~A.4]{weilres} identifies $\check{p}_\sharp(\OO_{\check{A}}) \simeq \check{p}_*(\OO_{\check{A}})^\vee \simeq \OO(\check{A})^\vee$, which by Proposition~\ref{prop:fun check A} is $\RA$.
\end{proof}

Passing to endomorphisms of the unit on either side of the equivalence $(\QCOH(B\check{A}), \otimes) \simeq (\MOD_{\RA}, \star)$ yields:

\begin{corollary}\label{cor:O B checkA via End}
  There is a canonical isomorphism of $E_\infty$-algebras in $\QCOH(S)$,
  \[
    \lamco_A := \operatorname{End}_{\RA}(\OO_S) \;\simeq\; \OO(B\check{A}).
  \]
\end{corollary}

\subsubsection{Ind-perfect complexes}

In view of Lemma~\ref{lem:dual mod conv}, taking dualizable objects on both sides of Theorem~\ref{thm:QCoh checkA} and ind-completing gives:

\begin{corollary}\label{cor:IndPerf checkA}
  There are canonical symmetric monoidal equivalences
  \begin{equation*}
    (\INDPERF(B\check{A}), \otimes)
      \;\simeq\; (\Ind(\MOD^{\operatorname{perf}/\OO_S}_{\RA}), \star)
    \;\simeq\; (\Ind(\COMOD^{\operatorname{perf}/\OO_S}_{\RA^\vee}), \otimes)
  \end{equation*}
  where $\COMOD^{\operatorname{perf}/\OO_S}_{\RA^\vee}$ denotes the full subcategory of $\RA^\vee$-comodules whose underlying $\OO_S$-module is perfect.
\end{corollary}

\subsubsection{Koszul dual side}

Assume $S$ is perfect.
Since $\sfp^* : \INDPERF(S) \to \INDPERF(B\check{A})$ is symmetric monoidal, its right adjoint $\sfp_*$ is lax symmetric monoidal, and lifts to modules over the algebra $\sfp_*(\OO_{B\check{A}}) \simeq \lamco_A$.
By the projection formula, we have $\sfp_*^{\Lambda}(\sfp^* V) \simeq \lamco_A \otimes_{\OO_S} V$ for $V \in \QCOH(S)$.

By Corollary~\ref{cor:Kosz II}, when $S$ is moreover regular, this lift is an equivalence:

\begin{corollary}\label{cor:IndPerf BcheckA Kosz}
  Suppose $S$ is a classical Artin stack which is regular, locally noetherian, and perfect.
  Then $\sfp_*$ lifts to a symmetric monoidal, $\INDPERF(S)$-linear equivalence
  \[
    \sfp_*^{\Lambda} : (\INDPERF(B\check{A}), \otimes) \xrightarrow{\ \sim\ } (\MOD_{\lamco_A}, \otimes).
  \]
\end{corollary}

\begin{remark}
  Recall that $\lamco_A \simeq \Sym^*_{\OO_S}(\fra[-2]) \simeq \OO_S[u]$ modulo the choice of an orientation, where $|u| = -2$.
\end{remark}

\begin{remark}\label{rem:pi not conservative}
  Unlike for $\QCOH$, the forgetful functor $\sfe^*: \INDPERF(B\check{A}) \rightarrow \INDPERF(S)$ is \emph{not conservative} on $\INDPERF$.
  This fails even for $S=\spc(k)$ with $k$ a field of characteristic zero.
  Under the identification $\lamco_A \simeq k[u]$, $|u|=-2$ (modulo a choice of orientation), $\sfe^*$ corresponds to $(-) \otimes_{k[u]} k$ along the augmentation $u \mapsto 0$, which kills the non-zero object $k[u, 1/u]$. This motivates Euler periodization below.
\end{remark}

\subsection{\texorpdfstring{$\check{A}$}{Ǎ}-invariants and coinvariants}\label{ssec:checkA invariants}

\subsubsection{}

Consider the functor $$\sfp^*_{\QCOH} : \QCOH(S) \to \QCOH(B\check{A}),$$ induced by the projection $\sfp : B\check{A} \to S$.
Under the equivalence $(\QCOH(B\check{A}), \otimes) \simeq (\MOD_{\RA}, \star)$ of Theorem~\ref{thm:QCoh checkA}, it corresponds to $\eps_* : \MOD_{\OO_S} \to \MOD_{\RA}$, restriction of scalars along the augmentation $\eps : \RA \to \OO_S$.
Indeed this follows from the fact that $\sfp_*^{\QCOH}$ is identified with $\eps^!$ as both compute homs out of the respective monoidal units. It follows then that $\sfp^*$ moreover admits a left adjoint $\sfp_{\sharp}^{\QCOH}$ corresponding to $\eps^*$, extension of scalars.

\subsubsection{}\label{sssec:sfp push}

On ind-perfect complexes the functor
$$\sfp^* : \INDPERF(S) \to \INDPERF(B\check{A})$$
admits a right adjoint $\sfp_*$, which is easily checked to be the renormalization of the functor $\sfp_*$ on $\QCOH$. There is also a functor $\sfp_\sharp$ obtained as the renormalization of the identically named functor on $\QCOH$. \begin{remark} We caution the reader that $\sfp_\sharp$ is no longer the left adjoint to $\sfp^*$, which in general will not even have a left adjoint. For example let us work over $S=\text{Spec}(k)$, with $k$ a $\mathbf{Q}$-algebra. Then choosing an orientation we can identify $\INDPERF(B\check{A})=\MOD(k[u])$, so that $\sfp^*=\otimes_k k[u].$ Now if $x$ has homological degree $2$ then we observe $$\sfp^*k[[x]]=k[[x]]\otimes k[u]\neq k[u][[x]],$$ as for example the right hand side contains the element $\sum u^n x^n$ in $\pi_0$, whereas the left hand side has no such element. In particular $\sfp^*$ does not preserve limits and so has no left adjoint.  \end{remark}
Under the equivalence $(\INDPERF(B\check{A}), \otimes) \simeq (\Ind(\MOD_{\RA}^{\mathrm{perf}/\OO_S}), \star)$ of Corollary~\ref{cor:IndPerf checkA}, the functors $\sfp_\sharp$, $\sfp^*$, and $\sfp_*$ correspond to the renormalized functors of (\ref{sssec:renorm}):
\[
  \sfp_\sharp \leftrightarrow \epsb^*, \qquad
  \sfp^* \leftrightarrow \epsb_*, \qquad
  \sfp_* \leftrightarrow \epsb^!.
\]
By construction, they all agree with their counterparts on $\QCOH$ when applied to perfect complexes.
In particular we have
\begin{align*}
  \sfp_\sharp(\OO_{B\check{A}}) &\simeq \epsb^*\epsb_*(\OO_S) \simeq \OO_S \otimes_{\RA} \OO_S \simeq \lam_A,\\
  \sfp_*(\OO_{B\check{A}}) &\simeq \epsb^!\epsb_*(\OO_S) \simeq \operatorname{End}_{\RA}(\OO_S) \simeq \lamco_A.
\end{align*}

\subsubsection{Co/invariants}

We define the functors of \emph{$\check{A}$-coinvariants} and \emph{$\check{A}$-invariants} on $\INDPERF(B\check{A})$ by
\[
  (-)_{h\check{A}} := \sfp_\sharp(-),
  \qquad (-)^{h\check{A}} := \sfp_*(-),
\]
respectively.

\begin{example}
	For $A=\GGm$ over the sphere spectrum $k=\mathbf{S}$, the functor $(-)_{hB\bZ}$ corresponds under the equivalence $$\INDPERF(B^2\bZ) \simeq \Ind\bigl(\MOD_{\mathbf{S}[S^1]}^{\mathrm{perf}/\mathbf{S}}\bigr)$$ to homotopy $S^1$-coinvariants, and $(-)^{hB\bZ}$ to the colimit-preserving extension of homotopy $S^1$-invariants on $\mathbf{S}[S^1]$-modules which are perfect over $\mathbf{S}$.
\end{example}

\subsubsection{}\label{sssec:map to coinvariants}

The atlas $\sfe : S \twoheadrightarrow B\check{A}$ induces the functor
\[ \sfe^* : \INDPERF(B\check{A}) \to \INDPERF(S), \]
forgetting\footnote{We use the word \emph{forget} here, despite the fact that the functor in question is not conservative.} the $\check{A}$-equivariant structure on an object $F \in \INDPERF(B\check{A})$.
Under the equivalence $(\INDPERF(B\check{A}), \otimes) \simeq (\Ind(\MOD_{\RA}^{\mathrm{perf}/\OO_S}), \star)$ of Corollary~\ref{cor:IndPerf checkA}, $\sfe^*$ corresponds to the renormalized restriction of scalars $\varpib_*$ along the unit $\varpi : \OO_S \to \RA$ (in the sense of (\ref{sssec:renorm})).
Its left and right adjoints then correspond to the renormalized extension and coextension $\varpib^*$ and $\varpib^!$, respectively:
\[
  \sfe_\sharp \leftrightarrow \varpib^*, \qquad
  \sfe^* \leftrightarrow \varpib_*, \qquad
  \sfe_* \leftrightarrow \varpib^!.
\]
In particular, $\sfe_\sharp(\OO_S) \simeq \RA$ and $\sfe_*(\OO_S) \simeq \RA^\vee$.

There is a natural map
\[
  \sfe^*(F)
  \simeq \sfp_\sharp \sfe_\sharp \sfe^*(F)
  \xrightarrow{\mathrm{counit}} \sfp_\sharp(F)
  = F_{h\check{A}}
\]
from the underlying object to its $\check{A}$-coinvariants.
Dually,
\[
  F^{h\check{A}} = \sfp_*(F)
  \xrightarrow{\mathrm{unit}} \sfp_* \sfe_* \sfe^*(F)
  \simeq \sfe^*(F).
\]

\subsubsection{$\sharp$-projection formula}
The projection formula
$$ \sfp_\sharp(F \otimes \sfp^*V) \simeq \sfp_\sharp(F)\otimes V $$
corresponds to the identity
$$ \eps^*(M \star \eps_*V) \simeq \eps^*(M) \otimes V. $$
This holds for any $\RA$-module $M$ and any $\OO_S$-module $V$, since the convolution product $M \star \eps_*V$ is equivalent to the $\OO_S$-module $M \otimes V$ with $\RA$ acting only on $M$.
Therefore, $\sfp_\sharp$ satisfies the projection formula on $\QCOH$.
Extending by colimits, so does the renormalization on $\INDPERF$.

\subsubsection{$\sharp$-base change formula}
The base change formula
$$ f^* \sfp_\sharp (F) \simeq \sfp'_{\sharp} f'^* (F), $$
where $f : S' \to S$ is a morphism and $f'$ its base change, corresponds to the identity
$$ f^* \eps^*(M) \simeq \eps'^* f'^*(M), $$
where we write $\eps' : \R_{A_{S'}/S'} = \RA \otimes_{\OO_S} \OO_{S'} \to \OO_{S'}$ for the base change.
This holds for any $\RA$-module $M$, so $\sfp_\sharp$ satisfies the base change formula on $\QCOH$ and so does the renormalization on $\INDPERF$.

\subsubsection{$*$-projection formula}
The projection formula
$$ \sfp_*(F \otimes \sfp^*V) \simeq \sfp_*(F)\otimes V $$
corresponds to the identity
$$ \eps^!(M \star \eps_*V) \simeq \eps^!(M) \otimes V. $$
When $V \in \QCOH(S)$ is \emph{dualizable}, i.e. $V \in \PERF(S)$, then symmetric monoidality of $\eps_*$ means $\eps_*V$ is $\star$-dualizable in $\MOD_{\RA}$. We have:
\begin{align*}
	\eps^!(M \star \eps_*V)
	&\simeq \Hom_{\RA}(\OO_S, M \star \eps_*V)\\
	&\simeq \Hom_{\RA}(\eps_*(V)^\vee, M)\\
	&\simeq \Hom_{\RA}(\eps_*(V^\vee), M)\\
	&\simeq \Hom_{\OO_S}(V^\vee, \eps^!(M))\\
	&\simeq \eps^!(M) \otimes V,
\end{align*}
using the symmetric monoidality of $\eps_*$ with respect to the convolution product on $\MOD_{\RA}$.
Thus $\sfp_*$ satisfies the projection formula for all $F \in \QCOH(B\check{A})$ and $V \in \PERF(S)$, and its renormalization satisfies the projection formula for all $F \in \INDPERF(B\check{A})$ and $V \in \INDPERF(S)$.
(Alternatively, the $\INDPERF$ version is a special case of Lemma~\ref{lem:indperf push easy}.)

\subsubsection{$*$-base change formula}
The base change formula
$$ f^* \sfp_* (F) \simeq \sfp'_{*} f'^* (F), $$
where $f : S' \to S$ is a morphism and $f'$ its base change, corresponds to the identity
$$ f^* \eps^!(M) \simeq \eps'^! f'^*(M), $$
where we write $\eps' : \R_{A_{S'}/S'} = \RA \otimes_{\OO_S} \OO_{S'} \to \OO_{S'}$ for the base change, and we regard $f'$ as the map of algebras $f' : \RA \to \R_{A_{S'}/S'}$.
Assuming that $\OO_{S'} \in \PERF(S)$, let us show the formula in $\QCOH$.
Taking limits, we may assume that $S$ and $S'$ are affine.
We have
\[
  f^* \eps^!(M)
  \simeq \Hom_{\RA}(\eps_*\OO_S, M) \otimes_{\OO_S} \OO_{S'}
  \simeq \Hom_{\RA}(\eps_*\OO_S, M) \otimes_{\RA} \R_{A_{S'}/S'}.
\]
It follows from the assumption that $\R_{A_{S'}/S'}$ is perfect, and hence dualizable, as an $\RA$-module.
Thus the tensor distributes over the hom, and the base change formula $f'^*(\eps_*\OO_S) \simeq \eps_*'\OO_{S'}$ identifies the right-hand side above with $\eps'^! (M \otimes_{\RA} \R_{A_{S'}/S'}) \simeq \eps'^! f'^*(M)$.
Thus $\sfp_*$ satisfies the base change formula on $\MOD^{\mathrm{perf}/\OO_S}_{\RA}$, hence so does its renormalization on $\INDPERF$.

\subsubsection{Base change for $\sfe^*$}\label{sssec:forget rotation equivariance}

Recall the quotient map $\sfe : S \twoheadrightarrow B\check{A}$.
For a stack $\YY$ with $\check{A}$-action, we also write $\sfe : \YY \twoheadrightarrow \YY/\check{A}$ for the quotient map, and still refer to
\[
  \sfe^* : \INDPERF^{\check{A}}(\YY) \to \INDPERF(\YY)
\]
as the ``forgetting rotation-equivariance'' functor (even though it is not conservative, see Remark~\ref{rem:pi not conservative}).
We now show that it is compatible with pushforward on $\INDPERF$.

\begin{lemma}\label{lem:sfe base change}
  Let $f : \XX \to \YY$ be an $\check{A}$-equivariant morphism of stacks with $\check{A}$-action, and $\bar{f} : \XX/\check{A} \to \YY/\check{A}$ the induced morphism.
  \begin{enumerate}
    \item On $\INDPERF$, $\sfe^*$ admits a left adjoint $\sfe_\sharp$, and $\bar{f}^*\sfe_\sharp \simeq \sfe_\sharp f^*$.
    \item There is a canonical isomorphism $\sfe^*\bar{f}_* \simeq f_*\sfe^*$ of functors $\INDPERF^{\check{A}}(\XX) \to \INDPERF(\YY)$.
  \end{enumerate}
\end{lemma}

\begin{proof}
  On $\QCOH$, $\sfe^*$ admits a left adjoint $\sfe_\sharp$ satisfying base change against $*$-pullback (by the argument in the proof of Theorem~\ref{thm:QCoh checkA}).
  We claim that $\sfe_\sharp$ preserves perfect complexes.
  For this it suffices to show $\sfe^*\sfe_\sharp$ preserves perfect complexes, since perfect complexes descend along $\sfe$.
  Base change in $\QCOH$ gives $\sfe^*\sfe_\sharp \simeq \mathrm{pr}_\sharp\,\mathrm{act}^*$, where $\mathrm{pr}, \mathrm{act} : \check{A} \times \YY \to \YY$ are the projection and action morphisms.
  Since $\mathrm{pr}_\sharp$ preserves perfect complexes by \cite[Rem.~A.4]{weilres}, the claim follows.

  Thus the restriction of $\sfe_\sharp$ to $\PERF$ ind-extends to a left adjoint of $\sfe^*$ on $\INDPERF$, and the base change isomorphism also extends by colimits from $\PERF$.
  Claim~(ii) follows from (i) by passing to right adjoints.
\end{proof}

\section{Euler periodization}\label{sec:Euler}

Throughout this and the following section we assume that $A \to S$ has relative dimension one.

\subsection{Euler periodization}\label{ssec:Euler}

\subsubsection{Thom twist}\label{sssec:Thom}

Recall that $\sfe : S \twoheadrightarrow B\check{A}$ denotes the quotient map.
We consider on $\INDPERF(B\check{A})$ the endofunctor $(-)\Thom{1}$ defined as the cofibre
\begin{equation}\label{eq:Euler triangle}
  \sfe_\sharp \sfe^*(-)
  \xrightarrow{\mathrm{counit}} (-)
  \to (-)\Thom{1}.
\end{equation}
It follows from Definition~\ref{def:abstract 1-shifted Cartier dual}\ref{item:abstract 1-shifted Cartier dual/sharp} that $\sfe_\sharp$ satisfies the projection formula, so this exact triangle is obtained by tensoring with its value on the structure sheaf,
\begin{equation}\label{eq:Euler O triangle}
	\sfe_\sharp (\OO_S) \xrightarrow{\mathrm{counit}} \OO_{B\check{A}} \to \OO_{B\check{A}}\Thom{1}.
\end{equation}
We refer to \eqref{eq:Euler triangle}, or \eqref{eq:Euler O triangle}, as the \emph{Euler triangle}.

\begin{remark}\label{rem:Euler = rotated aug triangle}
  Under the $\INDPERF(S)$-linear, symmetric monoidal equivalence $(\INDPERF(B\check{A}), \otimes) \simeq (\Ind(\MOD_{\RA}^{\mathrm{perf}/\OO_S}), \star)$ of Corollary~\ref{cor:IndPerf checkA}, the counit $\sfe_\sharp(\OO_S) \to \OO_{B\check{A}}$ is identified with the augmentation
	$$ \eps : \RA \simeq \varpi^*(\OO_S) \simeq \varpi^*\varpi_*\eps_*(\OO_S) \to \eps_*(\OO_S). $$
	Hence the Euler triangle \eqref{eq:Euler O triangle} corresponds to
	\[
	  \RA
		\to \eps_*(\OO_S)
		\xrightarrow{\partial} \mathcal{I}_\eps[1],
	\]
	the right-rotated augmentation triangle \eqref{eq:aug triangle}.
\end{remark}

Since $\mathcal{I}_{\eps}$ is an invertible object of $(\MOD_{\RA}^{\mathrm{perf}/\OO_S}, \star)$ (see Lemma~\ref{lem:Ieps inv}), we have:

\begin{lemma}\label{lem:Thom inv}
	The object $\OO_{B\check{A}}\Thom{1} \in \INDPERF(B\check{A})$ is $\otimes$-invertible.
\end{lemma}

It follows that the twist $(-)\Thom{1}$ is an autoequivalence, and we have the family of autoequivalences of $\INDPERF(B\check{A})$,
\[
  F \mapsto F\Thom{n} \simeq F \otimes (\OO_{B\check{A}}\Thom{1})^{\otimes n},
  \qquad n \in \bZ.
\]

\begin{remark}\label{rem:sfe^*O<1>}
  Applying $\sfe^*$ to the Euler triangle gives:
  \[
    \sfe^*\sfe_\sharp (\OO_S)
    \xrightarrow{\mathrm{counit}} \sfe^*\OO_{B\check{A}}
    \to \sfe^*\OO_{B\check{A}}\Thom{1}.
  \]
  Translating, this is the split exact triangle of $\OO_S$-modules,
  \[
    \RA \xrightarrow{\eps}
    \OO_S \to
    \frav[2],
  \]
  corresponding to the $\OO_S$-module splitting $\RA \simeq \OO_S \oplus \frav[1]$ (\ref{sssec:funcheckA add split}).
  In particular, we have $\sfe^*\OO_{B\check{A}}\Thom{1} \simeq \frav[2]$ in $\PERF(S)$.
\end{remark}

\begin{remark}\label{rem:Thom=2-shift}
	If $\RA$ is formal, then the isomorphism $\mathcal{I}_{\eps}[1] \simeq \eps_* \frav[2]$ in $\MOD_{\RA}^{\mathrm{perf}/\OO_S}$ (Lemma~\ref{lem:eps triangle}) corresponds to an identification $\OO_{B\check{A}}\Thom{1} \simeq \sfp^* \frav [2]$ in $\PERF(B\check{A})$, hence
	\[
		F\Thom{1} \simeq F \otimes \sfp^* \frav[2],
		\qquad \text{for all $F \in \INDPERF(B\check{A})$.}
	\]
	Up to a choice of orientation $\fra \simeq \OO_S$, the twist $(-)\Thom{1}$ is thus identified with the triangulated shift $(-)[2]$ in this case.
\end{remark}

\subsubsection{The Euler map}\label{sssec:Euler map}

By construction, there is a canonical map in $\PERF(B\check{A})$,
\[
  \alpha_A : \OO_{B\check{A}} \to \OO_{B\check{A}}\Thom{1},
\]
which we call the \emph{Euler map}.
More generally, we have
\[
  \alpha : F \to F\Thom{1}, \qquad F \in \INDPERF(B\check{A}).
\]

\begin{remark}\label{rem:sfe^*beta}
  In $\MOD_{\RA}^{\mathrm{perf}/\OO_S}$, the Euler map corresponds to the boundary map $\eps_*(\OO_S) \to \mathcal{I}_{\eps}[1]$ of the exact triangle
  \[
    \mathcal{I}_{\eps} \to \RA \to \eps_*(\OO_S),
  \]
  see Lemma~\ref{lem:Ieps inv}.
  The fact that this is naturally split on underlying $\OO_S$-modules corresponds to the fact that $\sfe^*(\alpha_A) : \OO_S \to \frav[2]$ is naturally null-homotopic (see Remark~\ref{rem:sfe^*O<1>}).
\end{remark}

\subsubsection{Euler periodization}\label{sssec:Euler periodization}

The \emph{Euler periodization} of $\INDPERF(B\check{A})$ is the Bousfield localization at the Euler map $\alpha_{A} : \OO_{B\check{A}} \to \OO_{B\check{A}}\Thom{1}$, in the sense of \cite[\S 5.5.4]{LurieHTT}:
\[
  \INDPERF(B\check{A})[\alpha_A^{-1}].
\]
By \cite[Prop.~4.8.2.7]{LurieHA}, $\INDPERF(B\check{A})[\alpha_A^{-1}]$ inherits a symmetric monoidal structure, and we denote the symmetric monoidal localization functor by $F \mapsto F[\alpha_A^{-1}]$.

More generally, for any $\INDPERF(B\check{A})$-linear category $\mathcal{C}$, we define its \emph{Euler periodization} as the tensor product over $\INDPERF(B\check{A})$:
\[
  \mathcal{C}[\alpha_A^{-1}] := \mathcal{C} \otimes_{\INDPERF(B\check{A})} \INDPERF(B\check{A})[\alpha_A^{-1}].
\]
For instance, this applies to $\INDPERF^{\check{A}}(\YY) := \INDPERF(\YY/\check{A})$ where $\YY$ is a stack with $\check{A}$-action, with $\INDPERF(B\check{A})$-linear structure obtained by pulling back the map $\YY/\check{A}\to B\check{A}$.

We may regard $\mathcal{C}[\alpha_A^{-1}]$ as the full subcategory of \emph{$\alpha_A$-periodic} objects in $\mathcal{C}$, in the sense of \cite[\S 3]{HoyoisCdh}.
By \cite[Lem.~3.3]{HoyoisCdh}, the Euler periodization of an object $F \in \mathcal{C}$ can be described by the following telescope construction:
\[
  F[\alpha_A^{-1}]
  \;=\; \varinjlim_n \left( F \xrightarrow{\alpha_A} F\Thom{1} \xrightarrow{\alpha_A} F\Thom{2} \xrightarrow{\alpha_A} \cdots \right).
\]

\subsubsection{}

The telescope formula shows that there are canonical isomorphisms
\[
  F[\alpha_A^{-1}] \;\simeq\; F \otimes \OO_{B\check{A}}[\alpha_A^{-1}],
  \qquad \text{for all } F \in \INDPERF(B\check{A}),
\]
and $\OO_{B\check{A}}[\alpha_A^{-1}]$ is idempotent as an $E_\infty$-algebra in $\INDPERF(B\check{A})$.
In particular,
\[ \INDPERF(B\check{A})[\alpha_A^{-1}] \;\simeq\; \MOD_{\OO_{B\check{A}}[\alpha_A^{-1}]}\bigl( \INDPERF(B\check{A}) \bigr); \]
see \cite[\S 4.8.2]{LurieHA}.

\subsection{Euler periodization vs. $u$-inversion}

Suppose that $S$ is perfect and $\RA$ is formal (Definition~\ref{def:formal}).
Using the lax symmetric monoidal lift $\sfp_*^{\Lambda} : \INDPERF(B\check{A}) \to \MOD_{\lamco_A}$ constructed before Corollary~\ref{cor:IndPerf BcheckA Kosz}, Euler periodization can be described as an algebraic localization at a (twisted) degree $-2$ class $u$ in $\lamco_A$.

Observe that we have
\[ \sfp_*^{\Lambda}(\OO_{B\check{A}}) \simeq \lamco_A, \qquad \sfp_*^{\Lambda}(\OO_{B\check{A}}\Thom{1}) \simeq \lamco_A \otimes_{\OO_S} \frav[2]. \]
The right-hand isomorphism follows from the projection formula in view of $\OO_{B\check{A}}\Thom{1} \simeq \sfp^*\frav[2]$ (which uses that $\RA$ is formal, see Remark~\ref{rem:Thom=2-shift}).
The Euler map $\alpha_A : \OO_{B\check{A}} \to \OO_{B\check{A}}\Thom{1}$ thus gives rise to a $\lamco_A$-linear map
$$\sfp_*^{\Lambda}(\alpha_A) : \lamco_A \to \lamco_A \otimes_{\OO_S} \frav[2],$$
or equivalently by adjunction an $\OO_S$-linear map
\[ u : \fra[-2] \to \lamco_A, \]
which we call the \emph{cohomology Euler class}.

\begin{proposition}\label{prop:Bott as localization in sqzero case}
Suppose $S$ is perfect and $\RA$ is formal.
For every $F \in \INDPERF(B\check{A})$ there is a canonical isomorphism of $\lamco_A$-modules
\[ \sfp_*^{\Lambda}(F[\alpha_A^{-1}]) \simeq \sfp_*^{\Lambda}(F)[u^{-1}], \]
the algebraic localization of $\sfp_*^{\Lambda}(F)$ at $u$ in the sense of \cite[Prop.~4.3.17]{LurieEll2}.
\end{proposition}
\begin{proof}
For any $F \in \INDPERF(B\check{A})$, we have $$\sfp_*^{\Lambda}(F\Thom{1}) \simeq \sfp_*^{\Lambda}(F \otimes \sfp^*\frav[2]) \simeq \sfp_*^{\Lambda}(F) \otimes_{\OO_S} \frav[2]$$ by Remark~\ref{rem:Thom=2-shift} and the projection formula (Lemma~\ref{lem:indperf push easy}).
The Euler map $\alpha_A : F \to F\Thom{1}$ is sent to $\sfp_*^{\Lambda}(F) \to \sfp_*^{\Lambda}(F) \otimes_{\OO_S} \frav[2]$, multiplication by $u$, by naturality of the lax symmetric monoidal structure of $\sfp_*$ applied to $\id_F \otimes \alpha_A$.
Since $\sfp_*^{\Lambda}$ is colimit-preserving,
\[
  \sfp_*^{\Lambda}(F[\alpha_A^{-1}])
  \simeq \varinjlim \bigl( \sfp_*^{\Lambda}(F) \xrightarrow{\id \otimes u} \sfp_*^{\Lambda}(F) \otimes \frav[2] \xrightarrow{\id \otimes u} \cdots \bigr).
\]
This is the localization of $\sfp_*^{\Lambda}(F)$ at $u$ by construction (see \cite[Prop.~4.3.17]{LurieEll2}).
\end{proof}

\subsection{The Tate construction}\label{ssec:Tate}

For any $F \in \INDPERF(B\check{A})$, we define the \emph{norm map} in $\INDPERF(S)$,
\begin{equation}
	\Nm_A : \sfp_\sharp(F) \otimes \frav[1] \to \sfp_*(F).
\end{equation}
Under the equivalence with $\Ind(\MOD_{\RA}^{\mathrm{perf}/\OO_S})$, this is the renormalization
\[
  \Nm_A^{\mathrm{ren}} : \epsb^*(-) \otimes \frav[1] \to \epsb^!(-)
\]
of the norm map defined in (\ref{sssec:norm Mod_RA}).\footnote{%
  That is, it agrees with \eqref{eq:norm} on $\MOD_{\RA}^{\mathrm{perf}/\OO_S}$ and is extended to the ind-completion by colimits.
}
For $M \in \Ind(\MOD_{\RA}^{\mathrm{perf}/\OO_S})$, we write $(-)^{t\check{A}}_{\mathrm{ren}}$ for the cofibre of the renormalized norm map.

\begin{remark}\label{renorm norm iso} Again, the renormalized norm map is an equivalence evaluated on $\RA$. Note that both target and source of the renormalized norm map now \emph{are} continuous, but now $\RA$ is no longer a generator, and so again $\Nm_A^{\mathrm{ren}}$ need not be an equivalence generally.\end{remark}

\begin{definition}\label{defn:Tate}
	The \emph{Tate construction} $(-)^{t\check{A}} : \INDPERF(B\check{A}) \to \INDPERF(S)$ is the cofibre
	\[
	  \sfp_\sharp(F) \otimes \frav[1]
		\xrightarrow{\Nm_A} \sfp_*(F)
		\longrightarrow F^{t\check{A}}.
	\]
\end{definition}

\begin{example}\label{ex:Tate multiplicative}
For $A=\GGm$, the identification $\R_{\GGm}\simeq k[S^1]$ of (\ref{sssec:OO(Gm[-1])}) gives an equivalence $\MOD_{\R_{\GGm}}\simeq\operatorname{Fun}(BS^1,\MOD_k)$ under which \eqref{eq:norm} is the norm map $M_{hS^1}[1]\to M^{hS^1}$ constructed in \cite{GreenleesMay,NS}.
Thus its cofibre is the usual $S^1$-Tate construction.
In particular, on $k[S^1]$-modules perfect over $k$, the renormalized Tate construction of Definition~\ref{defn:Tate} also agrees with the usual $S^1$-Tate construction.

More generally, it agrees with $M^{tS^1}$ on $\Ren^R_{B\check{\GG}_m}(M)$ for every $M \in \operatorname{Fun}(BS^1, \MOD_k)$: write $\Ren^R_{B\check{\GG}_m}(M) \simeq \varinjlim P$ over the $k[S^1]$-modules $P$ perfect over $k$ with a map $P \to M$; the source of the renormalized norm is $\varinjlim P_{hS^1}[1] \simeq M_{hS^1}[1]$, its target is $\sfp_*\Ren^R_{B\check{\GG}_m}(M) \simeq M^{hS^1}$ by adjunction, and by naturality the map is the norm of $M$.
\end{example}

\subsection{\texorpdfstring{$\check{A}$}{Ǎ}-co/invariants and the Tate construction}\label{ssec:invariants vs Tate}

In this subsection we will study the close relationship between the Tate construction and $\check{A}$-invariants of the Euler periodization: see Theorem~\ref{thm:Tate vs invariants of Bott}.

\subsubsection{Coinvariants of the Euler periodization}\label{sssec:coinv Bott per}

We begin by showing that the Euler periodization is invisible to $\check{A}$-\emph{co}invariants.

\begin{lemma}\label{lem:coinv kill periodization}
  For every $F \in \INDPERF(B\check{A})$, we have
  \[
    \bigl( F[\alpha_A^{-1}] \bigr)_{h\check{A}} \;\simeq\; 0.
  \]
\end{lemma}
\begin{proof}
  Since the expression is colimit-preserving in $F$, we may assume $F \in \PERF(B\check{A})$.
  The claim being local on $S$ by the $\sharp$-base change formula (Subsect.~\ref{ssec:checkA invariants}), we may also assume $S$ affine.
  Let $M \in \MOD_{\RA}^{\mathrm{perf}/\OO_S}$ be the object corresponding to $F$.
  Under the equivalence of Corollary~\ref{cor:IndPerf checkA}, the twist $F \mapsto F\Thom{1}$ corresponds to $M \mapsto M \star \mathcal{I}_{\eps}[1]$ and $\sfp_\sharp$ corresponds to $\epsb^*$.
  The claim is thus that the telescope $\varinjlim_n \epsb^*(M \star \mathcal{I}_\eps[1]^{\star n})$ vanishes.

  The perfect $\OO_S$-module $\varpi_*(M)$ is bounded below, say $c$-connective for some $c \in \bZ$.
  Then each $M \star \mathcal{I}_\eps[1]^{\star n}$ is $(c+2n)$-connective.
  Indeed, this can be detected on underlying $\OO_S$-modules; since $$\varpi_* : (\MOD_{\RA}, \star) \to (\QCOH(S), \otimes)$$ is symmetric monoidal (\ref{sssec:module conventions}) and $\varpi_*(\mathcal{I}_\eps[1]) \simeq \frav[2]$, we have
  \[
    \varpi_*\bigl( M \star \mathcal{I}_\eps[1]^{\star n} \bigr)
    \;\simeq\;
    \varpi_*(M) \otimes \fra^{\vee, \otimes n}[2n],
  \]
  which is $(c+2n)$-connective since $\frav$ is of Tor-amplitude $[0,0]$.

  Now since $\OO_S$ and $\RA$ are connective, $M \star \mathcal{I}_\eps[1]^{\star n}$ is still $(c+2n)$-connective after applying $\eps^*(-) = \OO_S \otimes_{\RA} (-)$.
  The telescope $\varinjlim_n \epsb^*(M \star \mathcal{I}_\eps[1]^{\star n})$ is thus a sequential colimit of objects of strictly increasing connectivity.
  Each homotopy group is then a sequential colimit of groups which eventually vanish, hence vanishes.
\end{proof}

\begin{remark} Work over a base $S=\text{Spec}(k)$ with $k$ a regular noetherian $\mathbf{Q}$-algebra, and choose an orientation. Then we have $\INDPERF(B\check{A})\simeq \MOD(k[u])$, and the functor $\sfp_\sharp$ is determined by $\sfp_\sharp k[u]=k[\beta]$, where $\beta$ is in degree $2$, and the induced $k[u]$-module structure has $u$ acting as $\partial_\beta$. In this case the above lemma amounts to the easy observation $$\sfp_\sharp k[u,1/u]\simeq \operatorname{colim}\Big( k[\beta]\xrightarrow{\partial_\beta} k[\beta]\xrightarrow{\partial_\beta}k[\beta]\xrightarrow{\partial_\beta}\cdots\Big)\simeq 0.$$\end{remark}

\subsubsection{Invariants of the Euler periodization}

In the situation of Proposition~\ref{prop:Bott as localization in sqzero case}, recall that $\sfp_*^{\Lambda}$ is a lift of the functor $\sfp_* = (-)^{h\check{A}}$ along the forgetful functor $\MOD_{\lamco_A} \to \MOD_{\OO_S}$.
Therefore, we have:

\begin{corollary}\label{cor:invariants of Bott vs beta inversion}
  Suppose $S$ is perfect and $\RA$ is formal.
  There is a canonical isomorphism in $\MOD_{\OO_S}$,
  \[
    (F[\alpha_A^{-1}])^{h\check{A}} \simeq \sfp_*^{\Lambda}(F)[u^{-1}].
  \]
\end{corollary}

\subsubsection{Tate construction vs. invariants of Euler periodization}

We will prove the following result about the Tate construction:

\begin{theorem}\label{thm:Tate vs invariants of Bott}
  For every $F \in \INDPERF(B\check{A})$, we have:

  {\rm(i)} There is a canonical isomorphism
  \[
    F^{t\check{A}}
    \;\simeq\;
    \sfp_*\bigl( F[\alpha_A^{-1}] \bigr)
    \;=\;
    \bigl( F[\alpha_A^{-1}] \bigr)^{h\check{A}}.
  \]

  {\rm(ii)} The localization map $F \to F[\alpha_A^{-1}]$ induces a canonical isomorphism
  \[
    F^{t\check{A}}
    \;\simeq\;
    \bigl( F[\alpha_A^{-1}] \bigr)^{t\check{A}}.
  \]
\end{theorem}

Theorem~\ref{thm:Tate vs invariants of Bott} expresses the Tate construction as the telescope
\[
  \sfp_*(F)
  \xrightarrow{\alpha} \sfp_*(F\Thom{1})
  \xrightarrow{\alpha} \sfp_*(F\Thom{2})
  \xrightarrow{\alpha} \cdots
  \longrightarrow F^{t\check{A}}.
\]
When $\RA$ is formal, the identification $\OO_{B\check{A}}\Thom{1} \simeq \sfp^*(\frav[2])$ (Remark~\ref{rem:Thom=2-shift}) allows us to rewrite this as the telescope
\[
  \sfp_*(F)
  \xrightarrow{u} \sfp_*(F)\Thom{1}
  \xrightarrow{u} \sfp_*(F)\Thom{2}
  \xrightarrow{u} \cdots
  \longrightarrow F^{t\check{A}},
\]
i.e., as the algebraic localization $\sfp_*(F)[u^{-1}]$, recovering Proposition~\ref{prop:Bott as localization in sqzero case}.

\subsubsection{Proof of Theorem~\ref{thm:Tate vs invariants of Bott}}

We begin with several preparatory lemmas.

\begin{lemma}\label{lem:induced beta null}
  We have $\sfe^*\bigl( \OO_{B\check{A}}[\alpha_A^{-1}] \bigr) \simeq 0$.
  Moreover,
  \[
    \sfe_\sharp(W)\bigl[\alpha_A^{-1}\bigr] \;\simeq\; 0
    \qquad \text{for every } W \in \INDPERF(S).
  \]
\end{lemma}
\begin{proof}
  Recall that $\sfe^*(\alpha_A) : \OO_S \to \frav[2]$ is null-homotopic (Remark~\ref{rem:sfe^*beta}).
  Since $\sfe^*$ is symmetric monoidal, the same holds for the induced map $\sfe^*(\alpha_A \otimes \id_{\OO\Thom{n}}) \simeq \sfe^*(\alpha_A) \otimes \sfe^*(\OO\Thom{n})$ for every $n$.
  As a sequential colimit along null-homotopic maps, $\sfe^*(\OO_{B\check{A}}[\alpha_A^{-1}])$ thus vanishes.
  Finally, by the projection formula we deduce
  \[
    \sfe_\sharp(W)[\alpha_A^{-1}]
    \simeq \sfe_\sharp(W) \otimes \OO_{B\check{A}}[\alpha_A^{-1}]
    \simeq \sfe_\sharp\bigl( W \otimes \sfe^*(\OO_{B\check{A}}[\alpha_A^{-1}]) \bigr)
    \simeq 0. \qedhere
  \]
\end{proof}

Let $P_A \in \INDPERF(B\check{A})$ denote the fibre of the localization map,
\begin{equation}\label{eq:isotropy triangle}
  P_A \longrightarrow \OO_{B\check{A}} \longrightarrow \OO_{B\check{A}}[\alpha_A^{-1}].
\end{equation}

\begin{lemma}\label{lem:P induced}
For every $F \in \INDPERF(B\check{A})$, the object $F \otimes P_A$ lies in the localizing subcategory of $\INDPERF(B\check{A})$ generated by the objects $\sfe_\sharp(W)$, $W \in \INDPERF(S)$.
\end{lemma}
\begin{proof}
  By exactness of colimits, $P_A$ is the colimit over $n$ of the objects $P_n := \operatorname{fib}(\alpha^{\circ n} : \OO_{B\check{A}} \to \OO_{B\check{A}}\Thom{n})$.
  Using the fact that $F \otimes (-)$ is colimit-preserving, and the projection formula $F \otimes \sfe_\sharp(W) \simeq \sfe_\sharp(\sfe^*F \otimes W)$, it is enough to show that each $P_n$ lies in the thick subcategory generated by the objects $\sfe_\sharp(\fra^{\vee,\otimes j}[2j])$, $0 \le j < n$.

  By the octahedral axiom applied to the composite $\OO \to \OO\Thom{n} \to \OO\Thom{n+1}$, there are exact triangles
  \[
    P_n \longrightarrow P_{n+1} \longrightarrow \operatorname{fib}\bigl( \alpha_A : \OO\Thom{n} \to \OO\Thom{n+1} \bigr).
  \]
  Tensoring the triangle $\sfe_\sharp (\OO_S) \xrightarrow{\mathrm{counit}} \OO_{B\check{A}} \to \OO_{B\check{A}}\Thom{1}$ \eqref{eq:Euler O triangle} with $\OO\Thom{n}$ and using the projection formula for $\sfe_\sharp$ and $\sfe^*(\OO\Thom{n}) \simeq \fra^{\vee, \otimes n}[2n]$ (Remark~\ref{rem:sfe^*O<1>}) gives
  \[
    \sfe_\sharp\bigl( \fra^{\vee, \otimes n}[2n] \bigr)
    \longrightarrow \OO\Thom{n}
    \xrightarrow{\alpha_A} \OO\Thom{n+1}.
  \]
  Since $P_1 \simeq \sfe_\sharp(\OO_S)$ (see \eqref{eq:Euler O triangle}), the claim follows by induction.
\end{proof}

\begin{remark} Work over the base $S=\text{Spec}(k)$, with $k$ a regular noetherian ring and choose an orientation identifying $\INDPERF(B\check{A})=\MOD(k[u])$. Then the above amounts to the observation that $k[u,1/u]/k[u]\otimes_{k[u]} M$ is locally $u$-torsion for any $M$. \end{remark}

\begin{lemma}\label{lem:norm invertible induced}
For every $G$ in the localizing subcategory of $\INDPERF(B\check{A})$ generated by the objects $\sfe_\sharp(W)$, $W \in \INDPERF(S)$, the norm map
\[
  \Nm_A(G) : \sfp_\sharp(G) \otimes \frav[1] \longrightarrow \sfp_*(G)
\]
is invertible.
In particular, by Lemma~\ref{lem:P induced}, this holds for $G = F \otimes P_A$ for every $F \in \INDPERF(B\check{A})$.
\end{lemma}
\begin{proof}
Since both sides are colimit-preserving in $G$, we may assume $G = \sfe_\sharp(W)$ with $W \in \PERF(S)$.
The assertion is moreover local on $S$ (by the base change formulas of Subsect.~\ref{ssec:checkA invariants}), so we may assume $S = \spc(R)$ affine, in which case the claim reduces by continuity to $W=\OO_S.$
Under the equivalence $\INDPERF(B\check{A})\simeq\Ind(\MOD_{\RA}^{\mathrm{perf}/\OO_S})$ the object $\sfe_\sharp\OO_S$ corresponds to $\RA$, and the claim reduces to the fact that the renormalized norm map is the identity for $\RA$, as per Remark~\ref{renorm norm iso}.
\end{proof}

\begin{remark} Work over the base $S=\text{Spec}(k)$ with $k$ a regular noetherian ring and choose an orientation, identifying $\INDPERF(B\check{A})=\MOD(k[u])$. Then the above just comes down to shifted self-duality of the $k[u]$ module $k=k[u]/u$, which is clear from the evident resolution. \end{remark}

\begin{proof}[Proof of Theorem~\ref{thm:Tate vs invariants of Bott}]
  We begin with the exact triangle $F \otimes P_A \to F \to F[\alpha_A^{-1}]$ \eqref{eq:isotropy triangle}.
  The naturality of the norm map yields the commutative diagram with exact rows:
  \[\begin{tikzcd}
    \sfp_\sharp(F \otimes P_A) \otimes \frav[1] \ar{r}\ar{d}
    & \sfp_\sharp(F) \otimes \frav[1] \ar{r}\ar{d}
    & \sfp_\sharp(F[\alpha_A^{-1}]) \otimes \frav[1] \ar{d}
    \\
    \sfp_*(F \otimes P_A) \ar{r}
    & \sfp_*(F) \ar{r}
    & \sfp_*(F[\alpha_A^{-1}]).
  \end{tikzcd}\]
  By Lemma~\ref{lem:norm invertible induced}, the left-hand vertical map is invertible.
  By Lemma~\ref{lem:coinv kill periodization}, the top right term vanishes, so the top left horizontal map is invertible.
  We can now read off:

  {\rm(i)}
  The cofibre of the middle vertical map is identified with the bottom right term, i.e., $F^{t\check{A}} \simeq \sfp_*(F[\alpha_A^{-1}])$.

  {\rm(ii)}
  The cofibre of the left-hand vertical map is zero, so the induced map on cofibres in the middle and right columns $F^{t\check{A}} \to (F[\alpha_A^{-1}])^{t\check{A}}$ is invertible.
\end{proof}

\section{Bott periodization}\label{sec:Bott per}

The Euler periodization of \S \ref{sec:Euler} inverts the Euler map at the level of $\check{A}$-equivariant ind-perfect complexes, yielding an operator
\[
  \INDPERF(B\check{A}) \to \INDPERF(B\check{A}), \qquad F \mapsto F[\alpha_A^{-1}].
\]
One can then form $\check{A}$-invariants to obtain a functor
\[
  \INDPERF(B\check{A}) \to \INDPERF(S), \qquad F \mapsto \bigl( F[\alpha_A^{-1}] \bigr)^{h\check{A}},
\]
which by Theorem~\ref{thm:Tate vs invariants of Bott} computes the Tate construction $F^{t\check{A}}$.
On the other hand, by Lemma~\ref{lem:coinv kill periodization} it is invisible to $\check{A}$-\emph{co}invariants as $(F[\alpha_A^{-1}])_{h\check{A}}$ vanishes.
In this section we introduce a second, complementary construction that lives on the coinvariants side.

By (\ref{sssec:sfp push}), the $\check{A}$-coinvariants of the unit are computed by the double rotation algebra of Subsect.~\ref{ssec:modules}, via the equivalence in $\INDPERF(S)$,
\[
  (\OO_{B\check{A}})_{h\check{A}} = \sfp_{\sharp}(\OO_{B\check{A}}) \simeq \lam_A.
\]
While the Euler map does not formally induce an operator on $\lam_A$, we will see that its canonical splitting on coinvariants produces a (twisted) degree $2$ class $\beta_A$ in $\lam_A$, the \emph{Bott class}.
The \emph{Bott periodization} is the algebraic localization $\lam_A[\beta_A^{-1}]$ at this class (Subsect.~\ref{ssec:Bott periodization}).

In characteristic zero, a choice of orientation $\fra \simeq \OO_S$ identifies $\lam_A \simeq \Sym^*(\frav[2]) \simeq \OO_S[\beta_A]$, $|\beta_A| = 2$ (\ref{sssec:funcheckA char 0}), so the Bott periodization is $\lam_A[\beta_A^{-1}] \simeq \OO_S[\beta_A^{\pm 1}]$.
On the other hand, the $\check{A}$-invariants of the Euler periodization of $\OO_{B\check{A}}$ are $\lamco_A[u^{-1}] \simeq \OO_S[u^{\pm 1}]$ (Corollary~\ref{cor:invariants of Bott vs beta inversion}).
Thus in this case, the Bott-periodized coinvariants are the $\OO_S$-linear dual of the invariants of the Euler periodization.

Over the sphere spectrum, Bott periodization recovers topological invariants of interest: in the multiplicative case, periodic complex K-theory (Subsect.~\ref{ssec:KU}), and in the elliptic case, the Goerss--Hopkins--Miller sheaf $\OO^{\mathrm{top}}$ of the moduli stack of elliptic curves, whose global sections are $\mathrm{TMF}$ (Subsect.~\ref{ssec:TMF}).

We continue to assume that $(A, \check{A}, \jmath)$ is a $1$-shifted Cartier dual pair over $S$, with $A \to S$ of relative dimension $1$.

\subsection{Bott periodization}\label{ssec:Bott periodization}

\subsubsection{}

Recall that the Euler map $\alpha_A$ sits, after taking $\check{A}$-coinvariants, in the exact triangle
\[
  \OO_S
  \to \lam_A
  \simeq (\OO_{B\check{A}})_{h\check{A}}
  \xrightarrow{(\alpha_A)_{h\check{A}}} (\OO_{B\check{A}}\Thom{1})_{h\check{A}}.
\]
This is split by the augmentation $\lam_A \to \OO_S$ (see Lemma~\ref{lem:lam splitting}), so we may identify $(\OO_{B\check{A}}\Thom{1})_{h\check{A}}$ with the fibre of the augmentation map $\lam_A \to \OO_S$.
Then $(\alpha_A)_{h\check{A}}$ admits the section
\begin{equation}\label{eq:sigma section}
  \sigma : \bigl(\OO_{B\check{A}}\Thom{1}\bigr)_{h\check{A}}
  \simeq \operatorname{fib}\bigl(\lam_A \to \OO_S\bigr)
  \longrightarrow \lam_A.
\end{equation}

\subsubsection{}

By (\ref{sssec:map to coinvariants}), there is a canonical map from the underlying module of $\OO_{B\check{A}}\Thom{1}$ to its coinvariants:
\begin{equation}\label{eq:c map}
	c : \frav[2] \simeq \sfe^* (\OO_{B\check{A}}\Thom{1}) \to \sfp_\sharp(\OO_{B\check{A}}\Thom{1}) = (\OO_{B\check{A}}\Thom{1})_{h\check{A}}.
\end{equation}

\subsubsection{Bott class}

Composing \eqref{eq:sigma section} and \eqref{eq:c map}, we obtain the canonical map
\begin{equation}\label{eq:Bott class}
	\beta_A : \frav[2] \xrightarrow{c} \bigl(\OO_{B\check{A}}\Thom{1}\bigr)_{h\check{A}} \xrightarrow{\sigma} \lam_A
\end{equation}
in $\INDPERF(S)$.
We call this the \emph{Bott class} for $A$.

\subsubsection{}

Under the equivalence $\INDPERF(B\check{A}) \simeq \Ind(\MOD_{\RA}^{\mathrm{perf}/\OO_S})$, $\beta_A$ reads:
\[
  \beta_A : \frav[2] \simeq \varpi_* \mathcal{I}_\eps[1] \xrightarrow{c} \eps^*  \mathcal{I}_\eps[1] \xrightarrow{\sigma} \lam_A.
\]

\subsubsection{Bott periodization}\label{sssec:Bott periodization}

The \emph{Bott-periodized double rotation algebra} is defined as the $E_\infty$-algebra
\begin{equation}
  \K_A := \bigl(\OO_{B\check{A}}\bigr)_{h\check{A}}[\beta_A^{-1}] \simeq \lam_A[\beta_A^{-1}]
\end{equation}
obtained by localizing the double rotation algebra $\lam_A$ at $\beta_A$.

More generally, for $V \in \QCOH(S)$, the projection formula gives
\[
  (\sfp^* V)_{h\check{A}} \;=\; \sfp_\sharp \sfp^* (V) \;\simeq\; V \otimes_{\OO_S} \lam_A
\]
and we define the \emph{Bott-periodized coinvariants} of $\sfp^* V$ by
\[
  (\sfp^* V)_{h\check{A}}[\beta_A^{-1}] \;:=\; V \otimes_{\OO_S} \lam_A[\beta_A^{-1}].
\]

More precisely, multiplication by $\beta_A : \frav[2] \to \lam_A$ determines a $\lam_A$-module map $\beta_A^{\sharp} : \lam_A \to \lam_A \otimes_{\OO_S} \fra[-2]$.
The Bott-periodized double rotation algebra is defined by the telescope
\[
  \lam_A
  \xrightarrow{\beta_A^\sharp} \lam_A \otimes \fra[-2]
  \xrightarrow{\beta_A^\sharp} \lam_A \otimes \fra^{\otimes 2}[-4]
  \xrightarrow{\beta_A^\sharp} \cdots
  \longrightarrow  \lam_A[\beta_A^{-1}].
\]
See \cite[Prop.~4.3.17]{LurieEll2}.

\subsubsection{}\label{sssec:Bott RA}

Under the equivalence $\INDPERF(B\check{A}) \simeq \Ind(\MOD_{\RA}^{\mathrm{perf}/\OO_S})$, the Bott class
\[
  \beta_A : \varpi_* \mathcal{I}_\eps[1] \xrightarrow{c} \eps^*  \mathcal{I}_\eps[1] \xrightarrow{\sigma} \lam_A,
\]
can be described as follows.

The Euler map $\alpha_A$ is the boundary map $\partial : \eps_* \OO_S \to \mathcal{I}_\eps [1]$, and after applying $\eps^*$ it sits in the split exact triangle (Lemma~\ref{lem:lam splitting})
\[
  \OO_S
  \to \lam_A
  \xrightarrow{\eps^*(\partial)} \eps^* \mathcal{I}_{\eps}[1].
\]
The map $\eps^*(\partial)$ admits the section
\begin{equation}\label{eq:sigma RA}
  \sigma : \eps^*\mathcal{I}_{\eps}[1]
  \simeq \operatorname{fib}\bigl(\lam_A \to \OO_S\bigr)
  \longrightarrow \lam_A.
\end{equation}

Now $\beta_A$ is the composite of $\sigma$ with the natural transformation
\begin{equation}\label{eq:c RA}
  c : \varpi_{*}(-)
  \simeq \eps^* \varpi^* \varpi_{*}(-) \xrightarrow{\ \eps^*(\mathrm{counit})\ } \eps^*(-).
\end{equation}

\subsubsection{Universal property}

Recall (\ref{sssec:maps out of lam}) that for every $R \in \CAlg_{\OO_S}$, we have natural isomorphisms
\begin{equation*}
	\Map_{\CAlg_{\OO_S}}(\lam_A, R) \;\simeq\; \Omega^2 A(R).
\end{equation*}

\begin{definition}\label{def:periodicity datum}
  A \emph{periodicity datum} for $A$ over $R$ is a point $x \in \Omega^2 A(R)$.
\end{definition}

Given a periodicity datum $x$, corresponding to $x : \lam_A \to R$, consider the $\OO_S$-linear map
\[
  \beta_{A,x} : \frav[2] \xrightarrow{\beta_A} \lam_A \xrightarrow{x} R.
\]
We say that $x$ is \emph{nondegenerate} if the induced $R$-linear map $\beta_{A,x}^{\sharp} : R \otimes \frav[2] \to R$ is invertible.
We write $\Omega^2 A(R)^{\mathrm{nd}} \subset \Omega^2 A(R)$ for the subgroupoid spanned by nondegenerate points.

\begin{proposition}\label{prop:bott per univ prop}
  The localization map $\lam_A \to \lam_A[\beta_A^{-1}]$ induces natural isomorphisms
  \[
    \Map_{\CAlg_{\OO_S}}\bigl(\lam_A[\beta_A^{-1}], R\bigr)
    \;\simeq\;
    \Omega^2 A(R)^{\mathrm{nd}}
    \;\subseteq\; \Omega^2 A(R),
  \]
	for every $R \in \CAlg_{\OO_S}$.
\end{proposition}

\subsubsection{Classical bases}\label{sssec:Bott classical}

For a classical base $S$, Bott periodization can be described as follows.
For simplicity, assume $S = \spc(k)$ and choose an orientation $\fra \simeq k$.
By (\ref{sssec:funcheckA for k discrete}), $\lam_A \simeq \Gamma^*(k[2])$ is the divided power algebra.
This is the graded $k$-algebra with generators $\gamma_n$ in degree $2n$ for $n \ge 0$, satisfying $\gamma_n \gamma_m = \binom{n+m}{n} \gamma_{n+m}$.
The Bott class is the first divided power $\beta_A = \gamma_1$, which acts as $\beta_A \cdot \gamma_n = (n+1) \gamma_{n+1}$.
The localization $\lam_A[\beta_A^{-1}]$ is thus concentrated in even homogeneous degrees, with degree $2n$ component given by the telescope
\[
  k \xrightarrow{\cdot (n+1)} k \xrightarrow{\cdot (n+2)} k \to \dots,
\]
i.e., the rationalization $k \otimes \bQ$.
Thus $\lam_A[\beta_A^{-1}] \simeq (k \otimes \bQ)[\beta_A^{\pm 1}]$.
In particular, this vanishes in characteristic $p$.

Below we will see that over non-discrete bases, like the sphere spectrum, the construction is much more interesting.

\subsection{Alternative description of the Bott class}

Inspired by \cite[\S 4.3]{LurieEll2}, we give another description of the Bott class $\beta_A$.

\subsubsection{}\label{sssec:linearization map}

Denote by $\CAlg_{\OO_S}$ the category of (quasi-coherent) $E_\infty$-$\OO_S$-algebras.
The universal property of $\lam_A \simeq \OO_S \otimes_{\RA} \OO_S$ as the pushout in $\CAlg_{\OO_S}$ provides the canonical map\footnote{%
  As per our usual convention, when regarding $\RA$ as an $\OO_S$-module we usually write it as $\varpi_* \RA$.
}
\[\begin{aligned}
  \Map_{\CAlg_{\OO_S}}(\lam_A, \lam_A)
  &\xrightarrow{~\sim~} \Omega \Map_{\CAlg_{\OO_S}}(\RA, \lam_A)\\
  &\xrightarrow{\mathrm{oblv}} \Omega \Map_{\OO_S}(\varpi_* \RA, \lam_A)\\
  &\xrightarrow{~\sim~} \Map_{\OO_S}(\varpi_* \RA[1], \lam_A),
\end{aligned}\]
where the second map is passing to underlying $\OO_S$-modules and the last isomorphism is translation of the base point.\footnote{%
  The base point in the first and second lines is the augmentation $\eps : \RA \to \OO_S$ post-composed with the unit $\OO_S \to \lam_A$.
}
The identity of $\lam_A$ thus gives rise to a canonical map
\begin{equation}\label{eq:linearization map RA}
  \mathrm{lin} : \varpi_* \RA[1] \to \lam_A
  \qquad \text{in } \MOD_{\OO_S}.
\end{equation}
This is equivalently an additive map $\OO(\Omega A)[1] \to \OO(\Omega^2 A)$, or again
\begin{equation}\label{eq:linearization map}
  \mathrm{lin} : \sfe^* \sfe_\sharp (\OO_S)[1] \to \eps^*\eps_*\OO_S
  \qquad \text{in } \INDPERF(S).
\end{equation}

\begin{lemma}\label{lem:Bott vs Lurie}
  The Bott class $\beta_A : \frav[2] \to \lam_A$ of \eqref{eq:Bott class} is canonically homotopic to the map in $\INDPERF(S)$,
  \begin{equation}\label{eq:Bott Lurie}
    \beta_A^{\mathrm{Lur}} : \frav[2]
    \simeq \sfe^* \OO_{B\check{A}}\Thom{1}
    \xrightarrow{\ \sfe^*(\partial)\ } \sfe^* \sfe_\sharp \OO_{S} [1]
    \xrightarrow{\ \mathrm{lin}\ } \lam_A
  \end{equation}
   where $\partial : \OO_{B\check{A}}\Thom{1} \to \sfe_\sharp \OO_{S} [1]$ denotes the boundary map of the Euler triangle \eqref{eq:Euler O triangle}.
\end{lemma}

Equivalently, \eqref{eq:Bott Lurie} is the map
\begin{equation}\label{eq:Bott Lurie 2}
	\beta_A^{\mathrm{Lur}} : \frav[2] \simeq \varpi_* \mathcal{I}_\eps [1] \xrightarrow{\text{can}[1]} \varpi_* \RA [1] \xrightarrow{\ \eqref{eq:linearization map RA}\ } \lam_A
\end{equation}
in $\MOD_{\OO_S}$.

\subsubsection{}

To prove Lemma~\ref{lem:Bott vs Lurie} it will be useful to consider the more general setup of any connective $E_\infty$-$\OO_S$-algebra $R$ with an augmentation $\eps_R : R \to \OO_S$; we write $\varpi_R : \OO_S \to R$ for the unit.
We set $\lam_R := \OO_S \otimes_R \OO_S$ and write $\phi_R : \lam_R \to \OO_S$ for its augmentation.

The same construction in (\ref{sssec:linearization map}) gives for every $R' \in \CAlg_{\OO_S}$ a natural map
\[
  \mathrm{lin}_{R,R'} : \Map_{\CAlg_{\OO_S}}(\lam_R, R')
  \to \Map_{\OO_S}(\varpi_{R,*} R[1], R').
\]
In particular, the identity of $\lam_R$ gives rise to a canonical $\OO_S$-linear map $\mathrm{lin}_R : \varpi_{R,*} R[1] \to \lam_R$.

\begin{lemma}\label{lem:phi lin null}
  The composite
  \[
    \varpi_{R,*} R [1]
    \xrightarrow{\mathrm{lin}_R} \lam_R
    \xrightarrow{\phi_R} \OO_S
  \]
  is null-homotopic.
\end{lemma}
\begin{proof}
  By the Yoneda lemma, this composite is the image of $\phi_R \in \Map_{\CAlg_{\OO_S}}(\lam_R, \OO_S)$ under the map $\mathrm{lin}_{R,\OO_S}$.
  Tracing through the construction, we see that $\phi_R$ corresponds in $\Omega \Map_{\CAlg_{\OO_S}}(R, \OO_S)$ to the constant self-homotopy of $\eps : R \to \OO_S$, which after translation at the base point becomes the zero map.
\end{proof}

Let $\mathcal{I}_R$ denote the fibre of $\eps_R : R \to \OO_S$, and consider the rotated triangle
\[
  R \to \eps_{R,*} \OO_S \xrightarrow{\partial} \mathcal{I}_R[1].
\]
Applying $\eps_R^*$ yields the exact triangle
\[
  \OO_S \to \lam_R \xrightarrow{\eps^*(\partial)} \eps^*\mathcal{I}_R[1],
\]
where the second map admits the section
\begin{equation*}
  \sigma_R : \eps^*\mathcal{I}_R[1]
  \simeq \operatorname{fib}\bigl(\phi_R : \lam_R \to \OO_S\bigr)
  \longrightarrow \lam_R.
\end{equation*}

\begin{lemma}\label{lem:mozambique}
  The map $\mathrm{lin}_{R} : \varpi_{R,*} R[1] \to \lam_R$ is canonically homotopic to the composite
  \[
    \varpi_{R,*} R[1]
    \xrightarrow{\mathrm{lin}_{R}} \lam_R
    \xrightarrow{\eps^*(\partial)} \eps_R^* \mathcal{I}_R[1]
    \xrightarrow{\sigma_R} \lam_R.
  \]
\end{lemma}
\begin{proof}
In view of Lemma~\ref{lem:phi lin null} it suffices to observe that, by construction of $\sigma_R$, any $\OO_S$-linear map $\psi : V \to \lam_R$ equipped with a null-homotopy of $\phi \circ \psi : V \to \lam_R \to \OO_S$ factors as follows:
\[
  \psi : V
  \xrightarrow{\psi} \lam_R
  \xrightarrow{\partial} \eps^* \mathcal{I}_R [1]
  \xrightarrow{\sigma_R} \lam_R.
\]
Indeed, writing $\tau_R : \OO_S \to \lam_R$ for the unit, the identity $\id_{\lam_R} \simeq \tau_R \circ \phi_R + \sigma_R \circ \partial$ shows
\[
  \psi \simeq \tau_R \circ (\phi_R \circ \psi) + \sigma_R \circ (\partial \circ \psi) \simeq \sigma_R \circ \partial \circ \psi.
  \qedhere
\]
\end{proof}

We let $c_R$ denote the natural transformation
\[
  c_R : \varpi_{R,*}(-)
  \simeq \eps_R^* \varpi_R^* \varpi_{R,*}(-) \xrightarrow{\ \eps_R^*(\mathrm{counit})\ } \eps_R^*(-).
\]

\begin{lemma}\label{lem:zimbabwe}
The map $c_R : \varpi_{R,*} \mathcal{I}_R [1] \to \eps_R^* \mathcal{I}_R [1]$ is canonically homotopic to the composite
\begin{equation}\label{eq:zimbabwe}
  \varpi_{R,*} \mathcal{I}_R [1]
  \xrightarrow{\mathrm{can}[1]} \varpi_{R,*} R [1]
  \xrightarrow{\mathrm{lin}_R} \lam_R
  \xrightarrow{\eps^*(\partial)} \eps_R^* \mathcal{I}_R [1]
\end{equation}
\end{lemma}
\begin{proof}
Set $I_R := \varpi_{R,*} \mathcal{I}_R \in \MOD_{\OO_S}$ and consider the augmented $E_\infty$-algebra
\[ R_0 := \Sym^*_{\OO_S}(I_R) \xrightarrow{\ \eps_{R_0}\ } \OO_S, \]
with augmentation induced by the zero map $I_R \to 0$, and $I_{R_0} := \varpi_{R_0,*} \mathcal{I}_{R_0} \simeq \Sym^{\ge 1}(I_R)$.

With values in augmented $E_\infty$-algebras, $V \mapsto \Sym^*_{\OO_S}(V)$ is left adjoint to the functor $(R, \eps_R) \mapsto \operatorname{fib}(\eps_R : R \to \OO_S)$.
We write $\delta : R_0 \to R$ for the unique map such that
\begin{equation}\label{eq:zimbabwe 2}
	I_R \hookrightarrow \Sym^{\ge 1}(I_R) \simeq I_{R_0} \xrightarrow{\delta'} I_R
\end{equation}
is the identity, where $\delta' := I_\delta$.
Since all constructions involved in the statement are natural in $R$ as an  augmented $E_\infty$-algebra, we have the commutative square
\[\begin{tikzcd}[column sep=large]
  I_{R_0}[1] = \varpi_{{R_0},*} \mathcal{I}_{R_0}[1] \arrow{r}{c_{R_0}} \arrow{d}{\delta'}
  & \eps_{R_0}^* \mathcal{I}_{R_0}[1] \arrow{d}{\delta''} \\
  I_{R}[1] = \varpi_{R,*} \mathcal{I}_R[1] \arrow{r}{c_R}
  & \eps_R^* \mathcal{I}_R[1],
\end{tikzcd}\]
where the vertical arrows are induced by $\delta$, and similarly for \eqref{eq:zimbabwe} in place of $c$.
Restricting along $I_R[1] \hookrightarrow I_{R_0}[1]$ and using the fact that \eqref{eq:zimbabwe 2} is the identity, we obtain in particular
\[
  c_R \;\simeq\; \delta'' \circ c_{R_0}\big|_{I_R[1]},
  \qquad
  c'_R \;\simeq\; \delta'' \circ c'_{R_0}\big|_{I_R[1]},
\]
where $c'$ is \eqref{eq:zimbabwe}.
It therefore suffices to show that $c_{R_0}$ and $c'_{R_0}$ agree after restriction along $I_R[1] \hookrightarrow I_{R_0}[1]$.

Thus let us take $R := \Sym^*(V)$ for $V \in \MOD_{\OO_S}$ connective, and show that $c_R$ and $c'_R$ agree on $V[1]$.
We have $\varpi_{R,*} \mathcal{I}_R \simeq \Sym^{\ge 1}_{\OO_S}(V)$ and
\[
  \lam_{R} := \OO_S \otimes_{R} \OO_S \simeq \Sym_{\OO_S}^*(V[1]),
\]
and
\[
  \eps_{R}^* \mathcal{I}_R [1]
  \simeq \operatorname{cofib}(\OO_S \to \lam_R)
  \simeq \Sym_{\OO_S}^{\ge 1}(V[1]).
\]
The map $c'_R$ \eqref{eq:zimbabwe} takes the form
\[
  \Sym^{\ge 1}(V)[1]
  \hookrightarrow \Sym^*(V)[1]
  \xrightarrow{\mathrm{lin}_R} \Sym^*(V[1])
  \twoheadrightarrow \Sym^{\ge 1}(V[1]).
\]
On the summand $V[1] \simeq \Sym^1(V)[1]$, it is easy to check by tracing through the construction of $\mathrm{lin}_R$ that this is just the inclusion $V[1] \hookrightarrow \Sym^{\ge 1}(V[1])$.

It remains to check the same for $c_R$.
By the construction of relative tensor products  (see \cite[Constr.~4.4.2.7]{LurieHA}) there is a colimit diagram
\[
  \cdots
  \rightrightrightrightarrows \OO_S \otimes R^{\otimes 2} \otimes \mathcal{I}_R
  \rightrightrightarrows \OO_S \otimes R \otimes \mathcal{I}_R
  \rightrightarrows \OO_S \otimes \mathcal{I}_R
  \rightarrow \OO_S \otimes_R \mathcal{I}_R,
\]
where the tensor products are implicitly over $\OO_S$.
The last map is precisely $c_R : \varpi_{R,*}\mathcal{I}_R[1] \to \eps_R^*\mathcal{I}_R[1]$ (up to a harmless shift).
Note that the map $c$ preserves weights if $\Sym^n$ is regarded with weight $n$ as usual, so the restriction $c|_{V[1]}$ only depends on its weight one component.
In weight one, the simplicial diagram $\OO_S \otimes R^{\otimes \bullet} \otimes \mathcal{I}_R$, $n \ge 0$, is constant with value $V$: since $\mathcal{I}_R$ sits in weights $\ge 1$, any positive weight contribution from $R^{\otimes \bullet}$ would result in total weight $\ge 2$.
In particular, we conclude that the last map is the identity in weight $1$; that is, its restriction to its weight $1$ summand $V$, which is $V$, is the inclusion as claimed.
\end{proof}

\begin{proof}[Proof of Lemma~\ref{lem:Bott vs Lurie}]
  Combining Lemmas~\ref{lem:mozambique} and \ref{lem:zimbabwe} shows that the composite
  \[
    \varpi_* \mathcal{I}_R [1] \xrightarrow{\text{can}[1]} \varpi_* R [1] \xrightarrow{\ \mathrm{lin}_R\ } \lam_R
  \]
  is homotopic to
  \[
    \varpi_* \mathcal{I}_R [1] \xrightarrow{c} \eps_R^* \mathcal{I}_R [1] \xrightarrow{\ \sigma_R\ } \lam_R.
  \]
  Specializing to the case of $R = \RA$, the latter is the Bott map $\beta_A$ while the former is \eqref{eq:Bott Lurie 2}.
\end{proof}

\subsection{Multiplicative case: periodic complex K-theory}\label{ssec:KU}

Let $k = \mathbf{S}$ be the sphere spectrum and let $A = \GGm$ be the strict multiplicative group over $\mathbf{S}$.
We show that in this case the Bott-periodized double rotation algebra is precisely the periodic complex K-theory spectrum $\mathrm{KU}$.

\subsubsection{}

By definition, $\OO(\GGm) \simeq \mathbf{S}[\mathbf{Z}]$ and by (\ref{sssec:OO(Gm[-1])}),
$$
\R_{\GGm} \simeq \mathbf S[\mathsf{B}\bZ] \simeq \mathbf S[S^1], \qquad
\lam_{\GGm} \simeq \mathbf S[\mathsf{B}^2\bZ] \simeq \mathbf S[\mathbf{CP}^\infty].
$$
Fix the coordinate $t$ on $\GGm$, corresponding to the tautological character in $\bZ$.
The $E_\infty$-space identifications $\mathsf{B}\bZ \simeq S^1$ and $\mathsf{B}^2\bZ \simeq \mathbf{CP}^\infty$ are normalized to send $t$ to the canonical generators
$$
  [\,\mathrm{id}_{S^1}\,] \in \pi_1(S^1) = \bZ,
  \qquad
  [\,S^2 \simeq \mathbf{C}\mathbf{P}^1 \hookrightarrow \mathbf{C}\mathbf{P}^\infty\,] \in \pi_2(\mathbf{C}\mathbf{P}^\infty) = \bZ,
$$
where $\mathbf{C}\mathbf{P}^1$ is oriented by its complex structure.

\subsubsection{}\label{sssec:tau t}

The unit map $\varpi : \OO_S \to \R_{\GGm}$ is
$$\Sigma^\infty_+\bigl( * \hookrightarrow S^1 \bigr) \;\simeq\; \bigl( \,\mathbf{S} \to \mathbf{S}[S^1] \,\bigr),$$
so the additive splitting $\R_{\GGm} \simeq \OO_S \oplus \frav[1]$ (\ref{sssec:funcheckA add split}) corresponds to $\mathbf{S}[S^1] \simeq \mathbf{S} \oplus \mathbf{S}[1]$.
This gives a canonical identification $\frav[1] \simeq \mathbf{S}[1]$, or equivalently an orientation of $A=\GGm$:
\begin{equation}\label{eq:tau t}
  \tau_t : \fra \xrightarrow{\ \sim\ } \OO_S = \mathbf{S}.
\end{equation}
The inversion $t \mapsto t^{-1}$ of $\GGm$ exchanges $\tau_t$ with $-\tau_t$.

\subsubsection{}

The augmentation $\eps : \R_{\GGm} \to \OO_S$ is
\[
  \Sigma^\infty_+\bigl( S^1 \to * \bigr) \;\simeq\; \bigl( \,\mathbf{S}[S^1] \to \mathbf{S} \,\bigr).
\]
Consider the rotated augmentation triangle $\RA \to \eps_*\OO_S \to \mathcal{I}_{\eps}[1]$ \eqref{eq:aug triangle}.
On underlying spectra, this is the exact triangle
\begin{equation}\label{eq:Gm aug}
  \R_{\GGm} \simeq \Sigma^\infty_+ S^1
  \to \mathbf{S}
  \to \mathbf{S}[2].
\end{equation}
Applying instead $\eps^*$, we get the exact triangle
\begin{equation}\label{eq:Gm eps^* aug}
  \mathbf{S}
  \to \lam_{\GGm} \simeq \Sigma^\infty_+ \mathbf{C}\mathbf{P}^\infty
  \to \Sigma^\infty \mathbf{C}\mathbf{P}^\infty.
\end{equation}
The left-hand arrow is the unit of $\mathbf{S}[\mathbf{C}\mathbf{P}^\infty]$, which is $\Sigma^\infty_+$ of $\mathrm{pt} \to \mathbf{C}\mathbf{P}^\infty$; it is split by the augmentation, which is $\Sigma^\infty_+$ of $\mathbf{C}\mathbf{P}^\infty \to \mathrm{pt}$.
The section $\sigma : \eps^* \mathcal{I}_{\eps}[1] \to \lam_{\GGm}$ \eqref{eq:sigma RA} is thus
\begin{equation*}
  \sigma : \Sigma^\infty \mathbf{C}\mathbf{P}^\infty
  \simeq \operatorname{fib}\bigl(\Sigma^\infty_+ \mathbf{C}\mathbf{P}^\infty \to \mathbf{S}\bigr)
  \longrightarrow \Sigma^\infty_+ \mathbf{C}\mathbf{P}^\infty,
\end{equation*}
i.e., the inclusion of the reduced summand.

\subsubsection{}

The natural transformation $c : \varpi_*(-) \to \eps^*(-)$ \eqref{eq:c RA} provides a commutative diagram
\[\begin{tikzcd}
  \Sigma^\infty_+ S^1 \ar{r}\ar{d}{c}
  & \mathbf{S} \ar{r}\ar{d}{c}
  & \mathbf{S}[2] \ar{d}{c}\\
  \mathbf{S} \ar{r}
  & \Sigma^\infty_+ \mathbf{C}\mathbf{P}^\infty \ar{r}
  & \Sigma^\infty \mathbf{C}\mathbf{P}^\infty,
\end{tikzcd}\]
whose rows are the triangles \eqref{eq:Gm aug} and \eqref{eq:Gm eps^* aug}.
The left-hand square is $\Sigma^\infty_+$ of the cartesian square exhibiting the isomorphism $S^1 \simeq \Omega \mathbf{C}\mathbf{P}^\infty$.
It follows that $c : \mathbf{S}[2] \to \Sigma^\infty \mathbf{C}\mathbf{P}^\infty$ is $\Sigma^\infty$ of the map $f : S^2 \simeq \Sigma(S^1) \to \mathbf{C}\mathbf{P}^\infty$ adjoint to this isomorphism.
This $f$ is the inclusion of the bottom cell (2-skeleton) $S^2 \simeq \mathbf{C}\mathbf{P}^1 \hookrightarrow \mathbf{C}\mathbf{P}^\infty$.

\subsubsection{}

The Bott class $\beta_{\GGm} : \frav[2] \to \lam_{\GGm}$ \eqref{eq:Bott class} is by (\ref{sssec:Bott RA}) the composite
\begin{equation*}
	\beta_{\GGm} : \frav[2] \xrightarrow{c} \eps^*(\mathcal{I}_{\eps}[1]) \xrightarrow{\sigma} \lam_{\GGm}.
\end{equation*}
By the discussion above, this is the composite
\[
  \beta_{\GGm} : \mathbf{S}[2] \xrightarrow{c} \Sigma^\infty \mathbf{C}\mathbf{P}^\infty \xrightarrow{\sigma} \Sigma^\infty_+ \mathbf{C}\mathbf{P}^\infty
\]
where $c$ and $\sigma$ are as above.
This is precisely the usual Bott class in $\pi_2(\Sigma^\infty_+(\mathbf{C}\mathbf{P}^\infty))$, so $\K_{\GGm} := \lam_{\GGm}[\beta_{\GGm}^{-1}]$ identifies with $\Sigma^\infty_+(\mathbf{C}\mathbf{P}^\infty)[\beta_{\GGm}^{-1}]$.
By Snaith's theorem (see \cite{Snaith}, \cite[Thm.~6.5.1]{LurieEll2}), this is the periodic complex K-theory spectrum $\mathrm{KU}$:

\begin{theorem}\label{thm:Snaith}
There is a canonical isomorphism of $E_\infty$-rings $\K_{\GGm} \simeq \mathrm{KU}$.
\end{theorem}

\subsection{Elliptic case: periodic topological modular forms}\label{ssec:TMF}

We now specialize to the case of elliptic curves.
We will see that Bott periodization is intimately related to Lurie's work on \emph{orientations} of elliptic curves in \cite{LurieEll2}.
In particular, we will see that the double rotation algebra of the universal elliptic curve, and its Bott periodization, describe precisely Lurie's moduli stacks of preoriented and oriented elliptic curves, respectively.
Therefore, Bott periodization in this case recovers the Goerss--Hopkins--Miller sheaf $\OO^{\mathrm{top}}$, as well as $\mathrm{TMF}$, the periodic spectrum of topological modular forms.

\subsubsection{}

Let $S$ be a stack and take $A = E \to S$ to be a relative elliptic curve.%
\footnote{By our implicit conventions, this means $S$ is a (connective) spectral stack in the sense of Lurie, and $E \to S$ is a relative strict elliptic curve in the sense of \cite{LurieEll2}, i.e., an abelian group object in the full subcategory of spectral stacks over $S$ spanned by flat proper morphisms $X \to S$ which are almost of finite presentation, geometrically reduced, and geometrically connected.}
We denote the double rotation algebra and its Bott periodization by
\[
  \lam_{E/S} \quad \text{and} \quad \lam_{E/S}[\beta_{E/S}^{-1}],
\]
to emphasize the role of the base.
For example, we may take the universal elliptic curve $\mathcal{E} \to \Mell$ as in \cite[Thm.~2.0.3]{LurieEll1}\footnote{%
  In \cite{LurieEll1}, $\Mell$ is denoted $\mathcal{M}^{s}$.
  In \cite{LurieEll2}, it is denoted $\mathcal{M}_{\mathrm{Ell}}^{s}$.
}, for which we have
\[
  \lam_{E/S} \simeq \lam_{\mathcal{E}/\Mell} \otimes_{\OO_{\Mell}} \OO_S,
\]
and similarly for $\lam_{\mathcal{E}/\Mell}[\beta_{\mathcal{E}/\Mell}^{-1}]$.

\subsubsection{}

By construction, $\K_{\mathcal{E}/\Mell} := \lam_{\mathcal{E}/\Mell}[\beta_{\mathcal{E}/\Mell}^{-1}]$ is a quasi-coherent sheaf of $E_\infty$-algebras on $\Mell$.
Regarded as a sheaf on the étale site of $(\Mell)_\cl \simeq \Mell^{\heartsuit}$, the classical moduli stack of elliptic curves, it recovers the sheaf of $E_\infty$-algebras $\OO^{\mathrm{top}}$ constructed by Goerss--Hopkins--Miller (see \cite[Chap.~12]{DFHH}, \cite[\S4.2]{HopkinsICM}, \cite[Thm.~7.0.1]{LurieEll2}):

\begin{theorem}\label{thm:Otop}
  There is a canonical isomorphism of sheaves of $E_\infty$-algebras on the étale site of $\Mell^{\heartsuit}$,
  \[
    \K_{\mathcal{E}/\Mell}
    \simeq \OO^{\mathrm{top}}.
  \]
  In particular, its global sections recover the spectrum of periodic topological modular forms:
  \[
    \Gamma\bigl( \Mell^{\heartsuit}, \K_{\mathcal{E}/\Mell} \bigr)
    \simeq \mathrm{TMF}.
  \]
\end{theorem}

\subsubsection{Proof of Theorem~\ref{thm:Otop}}\label{sssec:proof of thm:Otop}

We define the stack of \emph{elliptic curves with periodicity data} as the relative spectrum of the universal double rotation algebra,
\[
  \Mell^{\mathrm{per}} := \spc_{\Mell}(\lam_{\mathcal{E}/\Mell}) \to \Mell,
\]
whose fibre over a given elliptic curve $E \to S$ is
\[
  \Mell^{\mathrm{per}} \times_{\Mell} S \simeq \spc_S(\lam_{E/S}).
\]
In other words, $\Mell^{\mathrm{per}}$ classifies elliptic curves $E \to S$ together with a periodicity datum
\[
  x \in \Map_{\CAlg_{\OO_S}}(\lam_{E/S}, \OO_S) \;\simeq\; \Omega^2 E(\OO_S).
\]
The right-hand description, which follows from (\ref{sssec:maps out of lam}), shows that $\Mell^{\mathrm{per}}$ is an explicit construction of the moduli stack whose functor of points is described by Lurie in \cite[Prop.~7.2.5]{LurieEll2}.\footnote{%
  That is, it is Lurie's \emph{moduli stack of preoriented elliptic curves}, denoted $\mathcal{M}_{\mathrm{Ell}}^{\mathrm{pre}}$ in \cite[Prop.~7.2.5]{LurieEll2}.
  We avoid Lurie's term \emph{preorientations} in order not to clash with our unrelated notion of orientations of $A \to S$.
}

The claim about $\OO^{\mathrm{top}}$ now follows from the proofs of Proposition~7.2.10 and Theorem~7.0.1 in \cite{LurieEll2}.
In fact, we have the following stronger assertion.
Consider the nonconnective relative spectrum
\[
  \Mell^{\mathrm{per},\mathrm{nd}} := \Specnc_{\Mell}(\K_{\mathcal{E}/\Mell}) \to \Mell,
\]
regarded as a nonconnective spectral stack; this is the moduli stack of \emph{elliptic curves with nondegenerate periodicity data}, classifying elliptic curves $E \to S$ together with a periodicity datum $x$ as above which is \emph{nondegenerate}: the composite
\[
  \beta_x : \frav[2] \xrightarrow{\ \beta_{E/S}\ } \lam_{E/S} \xrightarrow{\ x\ } \OO_S
\]
is invertible.\footnote{%
  Note that this condition requires $\OO_S$ to be nonconnective, which is why $\Mell^{\mathrm{per},\mathrm{nd}}$ must be regarded as a genuinely nonconnective spectral stack.
}
By Lemma~\ref{lem:Bott vs Lurie}, $\beta_x$ agrees with the construction of \cite[Constr.~4.3.7]{LurieEll2}.
Now the proof of \cite[Prop.~7.2.10]{LurieEll2} shows that this is precisely Lurie's \emph{moduli stack of oriented elliptic curves}, denoted $\mathcal{M}_{\mathrm{Ell}}^{\mathrm{or}}$ in \cite{LurieEll2}.
\qed

\part{Loop spaces}\label{part:loops}

\section{Loop spaces}\label{sec:loops}

Throughout this section, $\YY$ is always a \emph{Tannakian} stack (Definition~\ref{def:Tannakian}) unless otherwise stated.

\subsection{The definition}\label{ssec:loopsdef}

Let $(A, \check{A}, \jmath)$ be a $1$-shifted Cartier dual pair over $k$.
In Subsect.~\ref{ssec:HA} we have introduced the categorified circle associated to $A$, i.e. the sheaf of tensor categories $\stHA$, together with an action by the commutative group stack $\check{A}$ (\ref{sssec:actonHA}).

\begin{definition} Let $\stCC$ be a sheaf of tensor categories. The \emph{$A$-loop space} of $\stCC$ is the stack $$\Loop_A(\stCC) := \TENSFUN(\stCC, \stHA),$$ equipped with the induced $\check{A}$-action. For a stack $\YY$, its $A$-loop space is $$\Loop_A(\YY) := \Loop_A(\stQCoh_{\YY}),$$ the $A$-loop space of the sheaf of tensor categories $\stQCoh_{\YY}$. \end{definition}

\subsubsection{Loops as maps}
\label{sssec:loops as maps}

Given a sheaf of tensor categories $\stCC$, consider the stack of tensor functors $$\TENSFUN(\stCC, \stQCoh_{\check{A}})$$ as in Subsect.~\ref{ssec:stktenfun}. Composition with the Fourier--Mukai transform $\Phi_{A}: \stHA\rightarrow \stQCoh_{\check{A}}$ (\ref{sssec:check A}), which is pointwise fully faithful, induces a natural monomorphism of stacks \[ \Loop_A(\stCC) \hookrightarrow \TENSFUN(\stCC, \stQCoh_{\check{A}}). \]
For $\YY$ Tannakian, this becomes a natural monomorphism
\begin{equation}\label{eq:loop hook map}
  \Loop_A(\YY) \hookrightarrow \TENSFUN(\stQCoh_{\YY}, \stQCoh_{\check{A}}) \simeq \MAP(\check{A}, \YY).
\end{equation}

Thus, $\Loop_A(\YY)$ can equivalently be defined as the substack of maps $\gamma : \check{A} \to \YY$ for which $$\gamma^{*}:\QCoh(\YY)\rightarrow\QCoh(\check{A})$$ factors through the essential image of the Fourier--Mukai transform $\Phi_{A} : \HA \hookrightarrow \QCoh(\check{A})$.

\subsubsection{Evaluation at the origin}\label{sssec:ev_0}
There is a canonical map, natural in $\YY$,
\[ \ev_0 : \Loop_A(\YY) \to \YY, \qquad \gamma \mapsto \gamma(e_{\check{A}}) \]
induced by applying $\MAP(-, \YY)$ to the identity section $e : S \to \check{A}$.

\subsubsection{Functoriality}
Any morphism of sheaves of tensor categories $f: \stCC\rightarrow \stDD$ induces a morphism of $A$-loop spaces $\Loop_A(f):\Loop_A(\stDD) \rightarrow \Loop_A(\stCC).$ This assignment naturally assembles into a functor $$\Loop_A: \stCCC^{\mathrm{op}} \longrightarrow \check{A}\circlearrowright \Stk,$$ valued in the category of stacks with an action of $\check{A}$.

If $f:\YY\rightarrow\XX$ is a morphism of stacks, the tensor functor $f^*$ induces a morphism $\Loop_A(f) : \Loop_A(\YY) \to \Loop_A(\XX)$, and there is a functor $$\Loop_A: \Stk \longrightarrow \check{A} \circlearrowright \Stk.$$

Given a morphism of $1$-shifted Cartier dual pairs $(A, \check{A}, \jmath) \to (A', \check{A}', \jmath')$, the pushforward $\stHA \to \stHcat_{A'}$ is a tensor functor that induces a morphism
\[ \Loop_{A}(\stCC) \to \Loop_{A'}(\stCC) \]
which is natural in $\stCC$ and $\check{A}'$-equivariant (where the action on the source is through $\check{A}' \to \check{A}$).

\subsubsection{Relative loop spaces}\label{sssec:relative}

If $A \to S$ is a relative commutative group scheme over an Artin stack $S$, one can define similarly $\Loop_{A/S}(\stCC)$ and $\Loop_{A/S}(\YY)$ for any $S$-linear tensor category $\stCC$ or any Artin stack $\YY$ over $S$.
For example, we may take $A$ to be the universal elliptic curve $\mathcal{E}$ over the moduli stack $\Mell$ classifying smooth elliptic curves; cf. Subsect.~\ref{ssec:Weierstrass} below.

\subsubsection{The functor $\Theta$}\label{sssec:Theta}

Consider the correspondence
\[ \begin{tikzcd}
    & \check{A} \times \MAP(\check{A}, \YY) \ar{rd}{\eval}\ar[swap]{ld}{\mathrm{pr}} & \\
    \MAP(\check{A}, \YY) & & \YY,
\end{tikzcd} \]
where $\eval$ and $\pr$ are the evaluation and projection morphisms, respectively.
Using the restricted correspondence
\[ \begin{tikzcd}
    & \check{A} \times \Loop_A(\YY) \ar{rd}{\eval}\ar[swap]{ld}{\mathrm{pr}} & \\
    \Loop_A(\YY) & & \YY,
\end{tikzcd} \]
we define the functor
$$\Theta^A_{\YY} := \mathrm{pr}_* \eval^* : \PERF(\YY) \to \PERF(\Loop_A(\YY)).$$
Here $\mathrm{pr}_*$ preserves perfect complexes by condition~\ref{item:abstract 1-shifted Cartier dual/sharp} of Definition~\ref{def:abstract 1-shifted Cartier dual}.

Recall that $\check{A}$ acts by translation on $\Loop_A(\YY)$; if we let $\check{A}$ act by the inverse translation on the factor $\check{A}$ of $\check{A} \times \Loop_A(\YY)$, then $\eval$ is equivariant for the trivial $\check{A}$-action on $\YY$, so $\Theta^A_{\YY}$ naturally lifts to a functor
$$\Theta^{A}_{\YY} := \mathrm{pr}_* \eval^* : \PERF(\YY) \to \PERF^{\check{A}}(\Loop_A(\YY)).$$

\begin{remark} For any morphism $f:\XX\rightarrow\YY$, the square \[\begin{tikzcd}
  \PERF(\YY) \ar{r}{\Theta^A_\YY}\ar{d}{f^*} & \PERF^{\check{A}}(\Loop_A(\YY)) \ar{d}{\Loop_A(f)^*}\\
  \PERF(\XX) \ar{r}{\Theta^A_\XX} & \PERF^{\check{A}}(\Loop_A(\XX))
\end{tikzcd}\]
commutes. \end{remark}

\subsubsection{Artinness}\label{sssec:ArtA}

We make the following definition for convenience.

\begin{definition}
  We say that $\YY$ \emph{satisfies $(\textup{Map}\textsubscript{A})$} if the canonical monomorphism $\Loop_A(\YY) \hookrightarrow \MAP(\check{A}, \YY)$ (\ref{sssec:loops as maps}) is a quasi-compact open immersion.
\end{definition}

\begin{remark}
	If $A$ is a relative abelian variety over a stack $S$ and $\YY$ is a stack over $S$, we similarly say that $\YY$ \emph{satisfies $(\textup{Map}\textsubscript{A/S})$} if $\Loop_{A/S}(\YY) \hookrightarrow \MAP_S(\check{A}, \YY)$ (\ref{sssec:relative}) is a quasi-compact open immersion.
\end{remark}

In the additive and multiplicative cases, every Artin stack $\YY$ satisfies $(\textup{Map}\textsubscript{A})$ by Theorem~\ref{thm:LGa LGm}.
For the elliptic case, see Theorem~\ref{open in map pluperfect}.

\subsection{The additive and multiplicative cases}

We specialize to the additive and multiplicative cases, i.e., where $A$ is $\GGa$ and $\GGm$, respectively. (The additive case assumes $k$ is a $\bQ$-algebra.)

\begin{theorem}\label{thm:LGa LGm}
  For every stack $\YY$, we have:
  \begin{enumerate}
    \item
    If $\YY$ is Artin, then there are canonical isomorphisms
    \[
      \Loop_{\GGa}(\YY)
      \simeq \MAP(B\GGhat_a^\sharp, \YY)
      \simeq \TT[-1]\YY,
    \]
    where $\TT[-1]\YY$ is the $(-1)$-shifted tangent bundle of $\YY$.
    \item
    There are canonical isomorphisms
    \[\Loop_{\GGm}(\YY) \simeq \MAP(B\mathbf{Z}, \YY)
    \simeq L(\YY),\]
    where $L(\YY) := \MAP(S^1, \YY)$ is the free loop space of $\YY$.
  \end{enumerate}
\end{theorem}
\begin{proof}
  Recall that $\GGa$ and $\GGm$ have $1$-shifted Cartier duals $B\GGhat_a^\sharp$ and $B\bZ \simeq S^1$, respectively.
  Thus we have monomorphisms
  \[ \Loop_{\GGa}(\YY) \hookrightarrow \MAP(B\GGhat_a^\sharp, \YY), \qquad \Loop_{\GGm}(\YY) \hookrightarrow \MAP(S^1, \YY) \]
  by (\ref{sssec:loops as maps}), which we identify with post-composition with the Fourier--Mukai transform $\Phi_{A}$,
  \[
    \TENSFUN(\stQCoh_{\YY}, \stHA) \hookrightarrow \TENSFUN(\stQCoh_{\YY}, \stQCoh_{\check{A}}),
  \]
  for each $A \in \{\GGa, \GGm\}$.
  By Theorem~\ref{thm:HA vs QCoh checkA}, this is fully faithful with essential image spanned by tensor functors that factor through $$\stIndPerf_{\check{A}}.$$
  But this condition is automatic, since any such tensor functor preserves dualizable objects (= perfect complexes).
  This shows the isomorphisms $\Loop_{A}(\YY) \simeq \MAP(\check{A}, \YY)$ for both $A$.
  Finally, the isomorphism $\MAP(B\GGhat_a^\sharp, \YY) \simeq \TT[-1]\YY$ is well-known; in the affine case, it follows by comparing the functors of points of both sides, and in general one shows that the left-hand side satisfies smooth codescent in $\YY$.\footnote{%
    In the characteristic zero case one can alternatively argue as follows: a pointed morphism $(B\widehat{G},*)\rightarrow(\YY,y)$ is equivalent to a morphism of Lie algebras $\operatorname{Lie}(G)\rightarrow\TT_{y,\YY}[-1],$ where $\TT_{y,\YY}[-1]$ has the Lie algebra structure constructed in \cite{Hennion}; taking $G=\GGa^\sharp$ yields the claim.
    See also \cite{NaefSafronov}, who show the claim even without Artinness hypotheses on $\YY$.
  }
\end{proof}

\subsection{Loop spaces of schemes}

For a scheme, the $A$-loop space only depends on the formal group $\widehat{A}$:

\begin{theorem}\label{thm:loop schemes}
  Let $(A, \check{A}, \jmath)$ be a $1$-shifted Cartier dual pair.
  For every qcqs scheme $X$, there is a canonical isomorphism of stacks
  \[ \Loop_A(X) \simeq \MAP(B(\mathbb{D}\widehat{A}), X). \]
  Moreover, the $\check{A}$-action on $\Loop_A(X)$ factors through the affinization $\operatorname{Aff}(\check{A})$, and the isomorphism is $\operatorname{Aff}(\check{A})$-equivariant (hence in particular $B(\mathbb{D}\widehat{A})$-equivariant).
\end{theorem}
\begin{proof}
  Using Theorem~\ref{Zarcosh} we reduce to the case of $X$ affine. For an affine $T$, a $T$-point of $\Loop_A(X)$ is a $T$-linear tensor functor
  $$\QCoh(X \times T) \rightarrow \stHA(T) \simeq \Hcat_{A_T/T}, $$
  where $A_T := A\times T$.
  Since $X$ and $T$ are affine this is specified by an $\OO_T$-algebra homomorphism from $\OO_X \otimes \OO_T$ to the endomorphism algebra of the convolution unit $\OO_{e} \in \Hcat_{A_T/T}$, which by Proposition~\ref{prop:fun check A} is
  \[ \operatorname{End}_{\QCoh(A_T)}(\OO_{e}, \OO_{e}) \simeq \OO(\check{A}_T). \]
  Thus $\Loop_A(X)(T) \simeq \Map_{\dAlgk^{\mathrm{nc}}}(\OO_X \otimes \OO_T, \OO(\check{A}_T))$, i.e. the space of maps $\operatorname{Aff}(\check{A}) \times T \to X$, so $\Loop_A(X) \simeq \MAP(\operatorname{Aff}(\check{A}), X)$.
  Since $X$ is an affine scheme, the $\check{A}$-action on $\Loop_A(X)$ factors through $\operatorname{Aff}(\check{A})$.
  Moreover, the isomorphism $\operatorname{Aff}(\check{A}) \simeq \operatorname{Aff}(B(\mathbb{D}\widehat{A}))$ of Corollary~\ref{cor:Aff check}, and the universal property of affinization, give:
  $$\MAP(\operatorname{Aff}(\check{A}), X) \simeq \MAP(\operatorname{Aff}(B(\mathbb{D}\widehat{A})), X) \simeq \MAP(B(\mathbb{D}\widehat{A}), X),$$
  whence $\Loop_A(X) \simeq \MAP(B(\mathbb{D}\widehat{A}), X)$.
\end{proof}

\begin{corollary}\label{cor:hkr0}
  Suppose the base ring $k$ is a $\mathbf{Q}$-algebra and $A$ a $1$-dimensional commutative group scheme over $k$.
  Then we have:
  \begin{enumerate}
	\item For every qcqs scheme $X$, there is a canonical isomorphism of schemes
  \[ \Loop_A(X) \simeq \mathbf{V}_X(\LL_X \otimes \frav[1]), \]
  where $\mathbf{V}_X(-)$ denotes the stack of cosections (see \cite[App.~B]{weilres}) and $\fra$ is as in (\ref{sssec:comm grps}).

  \item The $\check{A}$-action on $\Loop_A(X)$ factors canonically through the affinization $\operatorname{Aff}(\check{A}) \simeq B\mathbf{V}_S(\fra)$, and the above upgrades to an isomorphism of $B\mathbf{V}_S(\fra)$-equivariant schemes.
  \end{enumerate}
  In particular, a choice of orientation $\fra \simeq k$ induces a $B\GGa$-equivariant isomorphism $\Loop_A(X) \simeq \Loop_{\GGa}(X) \simeq \TT[-1]X$.
\end{corollary}
\begin{proof}
  Since $A$ is a $1$-dimensional commutative group scheme over a $\mathbf{Q}$-algebra, the formal logarithm identifies $\widehat{A}$ with the formal completion of $\mathbf{V}_S(\frav)$, whence $\operatorname{Aff}(\check{A}) \simeq B(\mathbb{D}\widehat{A}) \simeq B\mathbf{V}_S(\fra)$.
  By Theorem~\ref{thm:loop schemes}, $\Loop_A(X) \simeq \MAP(B\mathbf{V}_S(\fra), X) \simeq \mathbf{V}_X(\LL_X \otimes \frav[1])$, the case $\fra \simeq \OO_S$ being that of $\GGa$.
\end{proof}

\begin{remark}
  For $A=\GGa$, Corollary~\ref{cor:hkr0} encodes the classical Hochschild--Kostant--Rosenberg theorem; see e.g. \cite{BNloop}.
\end{remark}

\subsection{HKR filtration}\label{ssec:HKR fil}

The assumption in Theorem~\ref{thm:loop schemes} that $X$ is a scheme, and the characteristic zero assumption in Corollary~\ref{cor:hkr0}, are both necessary.

In general we have the following statement, for which we adopt the setting of (\ref{sssec:dagreview/animated}).

\begin{theorem}\label{thm:HKR fil}
  Let $k$ be an animated commutative ring and $A$ a $1$-dimensional commutative derived group scheme.
  Denote by $A_{\text{dgn}} := D_{e/A}/\GGm \to S_{\text{dgn}} := \aff/\GGm$ the deformation to the normal bundle of the identity element $e$ (e.g. as in \cite{weilres}).
  Then $A_{\text{dgn}}$ is a relative commutative derived group scheme over $S_{\text{dgn}}$ degenerating $A$ into $\mathbf{V}_S(\frav)$ ($\simeq \GGa$ up to a choice of orientation), and for any stack $\YY$ the animated loop space
  \[ \Loop_{A_{\text{dgn}}/S_{\text{dgn}}}(\YY) \to S_{\text{dgn}} \]
  is a family degenerating $\Loop_{A}(\YY)$ into $\Loop_{\mathbf{V}_S(\frav)}(\YY)$ ($\simeq \Loop_{\GGa}(\YY)$ up to a choice of orientation).
\end{theorem}

\begin{remark}\label{rem:check A dgn qcqs}
  Let $A$ be $\GGm$ or an elliptic curve $E$.
  Then there is a $1$-shifted Cartier dual pair $(A_{\mathrm{dgn}}, \check{A}_{\mathrm{dgn}})$ if we take $\check{A}_{\mathrm{dgn}}$ to be the stack of relative $1$-shifted characters $\CharOne(A_{\mathrm{dgn}}/S_{\text{dgn}})$.
  The degeneration $\Loop_{A_{\mathrm{dgn}}/S_{\text{dgn}}} \to S_{\text{dgn}}$ is then equipped with an action of $\check{A}_{\mathrm{dgn}}$ over $S_{\text{dgn}}$.
\end{remark}

\begin{example}
  This is a generalization of the filtered circle $S^1_{\mathrm{Fil}}$ of \cite{MRT}: for any $A$ we may apply the above construction to the affinization $\operatorname{Aff}(\check{A})$ in place of $A$ to obtain the relative group stack
  \[
    \check{A}^{\mathrm{aff}}_{\mathrm{dgn}} \;:=\; \operatorname{Aff}(\check{A})_{\mathrm{dgn}}
  \]
  over $S_{\text{dgn}}$, degenerating
  $$ \check{A}^{\mathrm{aff}}_{\mathrm{dgn}}|_1 \simeq \operatorname{Aff}(\check{A}) \qquad \rightsquigarrow \qquad \check{A}^{\mathrm{aff}}_{\mathrm{dgn}}|_0 \simeq \operatorname{Aff}(\check{\GG}_a) = B\GGhat^\sharp_a = \check{\GG}_a \quad \footnote{over an ordinary commutative ring.}$$
  For $A=\GGm$ over $\mathbf{Z}_{(p)}$, this $\check{\GG}_{m,\mathrm{dgn}}^{\mathrm{aff}}$ recovers $S^1_{\mathrm{Fil}}$.
\end{example}

\begin{remark}
  For $\YY=X$ a qcqs scheme, Theorem~\ref{thm:loop schemes} identifies $\Loop_A(X) \simeq \MAP(B(\mathbb{D}\widehat{A}), X)$.
  In that guise, and over $k = \mathbf{Z}_{(p)}$, a construction equivalent to that of Theorem~\ref{thm:HKR fil} has been studied in \cite[\S 7]{MoulinosFiltered}.
\end{remark}

\subsection{Universal loop spaces}
\label{ssec:Weierstrass}

The moduli stack $\Mell$ of smooth elliptic curves is a spectral Deligne--Mumford stack (Goerss--Hopkins--Miller; see \cite{LurieEllSurv}).
The universal elliptic curve $\mathcal{E} \to \Mell$ defines a \emph{universal elliptic loop space}
$$ \Loop_{\mathcal{E}/\Mell}(\YY \times \Mell), $$
whose fibres are the elliptic loop spaces $\Loop_{E}(\YY)$ for every elliptic curve $E$.
The Deligne--Mumford compactification $\overline{\mathcal{M}}_{\mathrm{ell}}$, parametrizing elliptic curves with at worst nodal singularities, also admits a spectral structure.

Over an ordinary commutative ring $k$ with $2$ and $3$ invertible, these admit the following presentation.
Let $W \subset \mathbf{P}^2 \times \mathbf{A}^2$ be the universal Weierstrass cubic curve $Y^2Z = X^3 - 27 c_4 X Z^2 - 54 c_6 Z^3$ over $\mathbf{A}^2 = \spc k[c_4, c_6]$.
Let $W^{\mathrm{sm}} \to \mathbf{A}^2$ be the relative smooth locus.
Since the point $[0:1:0] \in \mathbf{P}^2$ is always smooth in the total space, $W^{\mathrm{sm}} \to \mathbf{A}^2$ admits a canonical section $e$.
By \cite[II, Prop.~2.7]{DeligneRapoport}, $W^{\mathrm{sm}} \to \mathbf{A}^2$ is a smooth $1$-dimensional relative commutative group scheme.

The morphism $W \to \mathbf{A}^2$ is $\GGm$-equivariant, where $\GGm$ acts on $\mathbf{A}^2$ by $\lambda \cdot (c_4, c_6) = (\lambda^4 c_4,\, \lambda^6 c_6)$ and on the total space by $[X:Y:Z]\mapsto[\lambda^{2}X:\lambda^{3}Y:Z]$.
Since $\lambda \in \GGm$ acts on $\Delta = c_4^3 - c_6^2$ by $\lambda \cdot \Delta = \lambda^{12}\Delta$, $W^{\mathrm{sm}}$ is $\GGm$-invariant and descends to a relative group scheme
\[ \mathcal{A}_{\mathrm{Weier}} := W^{\mathrm{sm}}/\GGm \;\longrightarrow\; \mathcal{S}_{\mathrm{Weier}} := \mathbf{A}^2/\GGm. \]
The quotient stacks
\[
  \Mell \simeq \{\Delta\neq 0\}/\GGm
  \;\subset\;
  \overline{\mathcal{M}}_{\mathrm{ell}} \simeq (\mathbf{A}^2 \setminus \{0\})/\GGm
  \;\subset\;
  \mathbf{A}^2/\GGm
\]
are, respectively, the moduli stack of smooth elliptic curves, the Deligne--Mumford compactification adding the nodal ($\GGm$) fiber, and the Weierstrass compactification adding the cuspidal ($\GGa$) fiber.

The \emph{Weierstrass loop space} of a stack $\YY$ is the relative loop space
$$\Loop_{\mathcal{A}_{\mathrm{Weier}}/\mathcal{S}_{\mathrm{Weier}}}(\YY \times \mathcal{S}_{\mathrm{Weier}}).$$
Over $\Mell$, this recovers the universal elliptic loop space.
The fibre over a geometric point $(c_4, c_6) \in \mathbf{A}^2$ where $\Delta=0$ and $c_4 \neq 0$ is isomorphic to $\Loop_{\GGm}(\YY)$, and the fibre over the origin $(0,0)$ is isomorphic to $\Loop_{\GGa}(\YY)$.

\subsection{K\"unneth for loop spaces}

Recall the notion of perfect stack (Definition~\ref{def:perfect}).

\begin{theorem}\label{Kun} Let $\XX$ be a stack and $\{\YY_i\}_{i=1,\ldots,n}$ a finite collection of stacks over $\XX$.  If $\XX$ and the $\YY_i$ are perfect fppf-locally on $S=\spc(k)$, then the canonical morphism
$$ \Loop_A(\YY_1 \times_{\XX} \cdots \times_{\XX} \YY_n) \longrightarrow \Loop_A(\YY_1) \times_{\Loop_A(\XX)} \cdots \times_{\Loop_A(\XX)} \Loop_A(\YY_n) $$
is an isomorphism in $\check{A}\circlearrowright \Stk$. \end{theorem}

\begin{proof} The claim is fppf-local on $S$, so we may assume that $\XX$ and the $\YY_i$ are perfect. By induction we may assume $n=2$. The perfectness assumptions guarantee that we have the equivalences \[\QCoh(\YY_1 \times_{\XX} \YY_2) \simeq \QCoh(\YY_1) \otimes_{\QCoh(\XX)} \QCoh(\YY_2), \] by passing to connective subcategories from \cite[Thm.~1.2(i)]{BFN}. Thus the claim follows from Lemma~\ref{fibre}. \end{proof}

We also have the following variant without perfectness (cf.\ \cite[Vol.~I, Chap.~3, Prop.~3.3.5]{GR}):

\begin{theorem}\label{thm:Kun affine}
  Let $\XX$ be a stack and $\YY_1$, $\YY_2$ stacks over $\XX$.
  If $\YY_2$ is affine over $\XX$, then the canonical morphism
  $$ \Loop_A(\YY_1 \times_{\XX} \YY_2) \longrightarrow \Loop_A(\YY_1) \times_{\Loop_A(\XX)} \Loop_A(\YY_2) $$
  is an isomorphism in $\check{A}\circlearrowright \Stk$.
\end{theorem}
\begin{proof}
  Write $p : \YY_1 \to \XX$ and $q : \YY_2 \to \XX$.
  Since $q$ is affine, $\QCoh(\YY_2) \simeq \Mod_{q_*\OO_{\YY_2}}(\QCoh(\XX))$ and, by base change, $\QCoh(\YY_1 \times_{\XX} \YY_2) \simeq \Mod_{p^*q_*\OO_{\YY_2}}(\QCoh(\YY_1))$.
  Thus \cite[Thm.~4.8.4.6]{LurieHA} gives
  \[ \QCoh(\YY_1 \times_{\XX} \YY_2) \simeq \QCoh(\YY_1) \otimes_{\QCoh(\XX)} \QCoh(\YY_2). \]
  Now the claim follows from Lemma~\ref{fibre}.
\end{proof}

\section{Equivariant loop spaces}\label{sec:loops BG}

\subsection{Equivariant loop spaces}
\label{ssec:LA(BG)}

Let $G$ be an fppf group scheme over $k$ acting on a stack $X \in \Stk$.
The $G$-equivariant geometry of $X$ is encoded by the morphism
\[
  a_X^G : X/G \to BG\footnote{Where we consider $BG$ as a pointed object.},
\]
since the fibre over $* \twoheadrightarrow BG$ recovers $X$ together with its $G$-action.

\begin{definition}\label{def:Loop_A^G}
The \emph{$G$-equivariant $A$-loop space of $X$} is $\Loop_A^G(X) := \Loop_A(X/G)$.
We regard it as a stack with $\check{A}$-action, equipped with an $\check{A}$-equivariant structure morphism
\[
  \pi_X^G := \pi_X^{A,G} := \Loop_A(a_X^G) : \Loop_A^G(X) \to \Loop_A^G(*).
\]
We use the same notation for the induced map on $\check{A}$-quotients.
\end{definition}

\begin{remark}
If $f : Y \to X$ is a $G$-equivariant morphism, there is an induced map $\Loop_A^G(Y) \to \Loop_A^G(X)$.
If $(A, \check{A}) \to (A', \check{A}')$ is a morphism of $1$-shifted Cartier dual pairs, there is an induced map $\Loop_{A}^G(X) \to \Loop_{A'}^G(X)$.
If $K \to G$ is a morphism of algebraic groups, there is an induced map $\Loop_A^K(*) \to \Loop_A^G(*)$.
\end{remark}

\begin{example}\label{exam:equivariant loops Ga Gm}
By Theorem~\ref{thm:LGa LGm} applied to $\YY=BG$, we have
\[
  \Loop_{\GGa}^G(*) \simeq \mathfrak{g}/G,
  \qquad
  \Loop_{\GGm}^G(*) \simeq G/G,
\]
where the actions are the adjoint and conjugation actions, respectively.
For the left-hand isomorphism we assume that $k$ is a $\bQ$-algebra and $G$ is smooth with Lie algebra $\mathfrak{g} := \operatorname{Lie}(G)$.
\end{example}

\subsection{Loop spaces for \texorpdfstring{$\GL_r$}{GLr} in characteristic zero}
\label{ssec:LA(BGLr) char 0}

In this subsection we take $k$ to be a field of characteristic zero, $S=\spc(k)$, and $G=\GL_r$ with $r \geq 1$.
Using a remarkable theorem of Iwanari \cite{Iwa}, showing that the category $\QCoh(B\GL_r)^\heartsuit \simeq \Rep(G)^\heartsuit$ is in this case tensor-generated by the standard representation $V_{\mathrm{std}}$, we give a description of the $\GL_r$-equivariant loop space of a point.

\subsubsection{}

We fix a $1$-dimensional commutative algebraic group $A$ over $k$.
We denote by $\operatorname{Coh}_0^r(A)$ the stack of $0$-dimensionally supported sheaves of length $r$ over $A$.
More precisely, its $T$-points are perfect complexes on $A\times_k T$ which are connective and flat over $T$, with finite support over $T$ and geometric fibres of length $r$.

\begin{proposition}\label{prop:LA(BGLr) char 0}
There is a canonical isomorphism of stacks
\[
  \Loop_A^{\GL_r}(*) \simeq \operatorname{Coh}_0^r(A).
\]
\end{proposition}

\begin{lemma}\label{lem:Coh_0^r wedge finite}
Let $T$ be an affine $k$-scheme and $F \in \PERF_0(A_T/T)$.
Then the following conditions are equivalent:
\begin{enumerate}
	\item $F$ and its $\star$-dual $\mathbb{D}^\star(F)$ are connective.
	\item $F$ is connective and flat over $T$.
	\item The pushforward $V := p_*(F)$ along $p : A_T \to T$ is a vector bundle.
\end{enumerate}
Moreover, under these equivalent conditions, $F$ has geometric fibres of length $r$, for some $r \geq 1$, if and only if $\wedge_\star^r F$ is invertible and $\wedge_\star^{r+1}F\simeq 0$, where $\wedge_\star^n$ denotes the $n$th exterior power taken with respect to convolution.
\end{lemma}
\begin{proof}
For every $F \in \PERF_0(A_T/T)$, $V := p_*(F)$ is perfect by Proposition~\ref{prop:push on QCoh0}.
By Lemma~\ref{lem:push finite support}, $p_*$ is t-exact and conservative on complexes with finite support over $T$.
Thus $F$ is connective if and only if $V$ is connective.
Note that $p_*(\mathbb{D}^\star F)\simeq V^\vee$, since $p_*$ is symmetric monoidal with respect to the convolution product on the source.
Hence also $\mathbb{D}^\star F$ is connective if and only if $V^\vee$ is connective.
This shows (i) $\Leftrightarrow$ (iii).

For every $N\in\QCoh(T)^\heartsuit$, the projection formula gives
\[
  p_*(F\otimes p^*N)\simeq V\otimes N.
\]
Since $F\otimes p^*N$ has finite support over $T$, Lemma~\ref{lem:push finite support} shows that $F\otimes p^*N$ is discrete if and only if $V\otimes N$ is discrete.
Allowing $N$ to vary, this means that $F$ is connective and flat over $T$ if and only if $V$ is connective and flat, i.e., a vector bundle (since $V$ is perfect).
This shows (ii) $\Leftrightarrow$ (iii).

For the last claim, suppose $F$ satisfies the equivalent conditions.
Then for any geometric point $t$ of $T$, the fibre $F_t$ is discrete with finite support over the algebraically closed field $k(t)$, hence has length
\[
  \operatorname{length}_{k(t)}(F_t)
  = \dim_{k(t)}\Gamma(A_t, F_t)
  = \dim_{k(t)}(V_t).
\]
Thus $F$ has geometric fibres of length $r$ if and only
if $\wedge^r V$ is invertible and $\wedge^{r+1}V\simeq0$.
It remains to show these conditions are equivalent to $\wedge_\star^r F$ being $\star$-invertible and $\wedge^{r+1}_\star F\simeq 0$.
Note that $P := \wedge_\star^r F$ is dualizable, since it is a direct summand of $F^{\star r}$, so its $\star$-invertibility is equivalent to invertibility of the evaluation map $P\star\mathbb{D}^\star P \to \OO_e$.
The claim now follows from the fact that $p_*$ is symmetric monoidal, and conservative on complexes with finite support over $T$.
\end{proof}

\begin{proof}[Proof of Proposition~\ref{prop:LA(BGLr) char 0}]
For every affine $T$ over $k$, Iwanari's theorem \cite[Thm.~1.7]{Iwa} identifies the groupoid of colimit-preserving $k$-linear symmetric monoidal functors
\[
  \phi : (\QCOH(B\GL_r), \otimes) \to (\QCOH_0(A_T/T), \star)
\]
with the groupoid of objects $F \in \PERF_0(A_T/T)$ such that $\wedge_\star^r F$ is invertible and $\wedge_\star^{r+1}F\simeq 0$, via the assignment $\phi \mapsto F := \phi(V_{\mathrm{std}})$.
Note here that $F$ is $\star$-dualizable, since $\phi$ is symmetric monoidal, hence perfect by Corollary~\ref{duals}.

This equivalence is natural in $T$.
Since $\Loop_A(B\GL_r)(T)$ is the groupoid of such functors $\phi$ that preserve connective objects, it remains to observe:
\begin{itemize}
\item Any $\phi$ as above preserves connective objects if and only if both $F$ and $\mathbb{D}^\star(F)$ are connective.
Indeed, the connective subcategory $\QCoh(B\GL_r)$ is generated under colimits, extensions, and direct summands by the objects $V_{\mathrm{std}}^{\otimes m} \otimes V_{\mathrm{std}}^{\vee,\otimes n}$ for $m,n \geq 0$ (see the proof of Lemma~\ref{lem:BG pluperfect lin red}).

\item By Lemma~\ref{lem:Coh_0^r wedge finite}, $F$ and $\mathbb{D}^\star(F)$ are connective if and only if $F$ is connective and flat over $T$; and in this case, since $\wedge_\star^r F$ is invertible and $\wedge_\star^{r+1}F\simeq 0$, $F$ has geometric fibres of length $r$.
\end{itemize}
That is, $\phi$ preserves connective objects if and only if $F\in \operatorname{Coh}_0^r(A)(T)$.
\end{proof}

\subsubsection{}

Proposition~\ref{prop:LA(BGLr) char 0} shows in particular that $\Loop_A^{\GL_r}(*)$ is a smooth classical Artin stack.
Moreover, recall that the map taking support cycles of length $r$ coherent sheaves,
\[
  \Loop_A^{\GL_r}(*) \simeq \operatorname{Coh}_0^r(A) \to A^{(r)},
  \qquad
  F \mapsto \sum_{a\in\operatorname{Supp}(F)} \operatorname{length}_{\OO_{A,a}}(F_a)[a],
\]
exhibits the symmetric power $A^{(r)}$ as a good moduli space (see \cite[Ex.~4.3.6]{HuybrechtsLehn} and \cite[Ex.~8.7]{Alper}).

This morphism is $\check{A}$-equivariant for the trivial action on $A^{(r)}$, since the action described in (\ref{sssec:actonHA}) is given by tensoring with multiplicative line bundles, under which support cycles are invariant.

\subsubsection{}

As another example, consider the case $G=SL_r$.
The $SL_r$-equivariant loop space $\Loop_A^{SL_r}(*)$ can be identified with the smooth classical Artin stack of pairs $(F,\epsilon)$, where $F\in\operatorname{Coh}_0^r(A)$ and $\epsilon : \wedge_\star^r F \simeq \OO_e$.

Indeed, K\"unneth (Theorem~\ref{thm:Kun affine}) applied to $BSL_r\simeq B\GL_r\times_{B\GGm}*$ identifies $\Loop_A^{SL_r}(*)$ with the fibre of the map
\[
  \Loop_A^{\GL_r}(*) \to \Loop_A^{\GGm}(*),
  \qquad F \mapsto \wedge_\star^r F,
\]
over the point corresponding to the convolution unit $\OO_e$.

A good moduli space is provided in this case by the support cycles map
\[
  \Loop_A^{SL_r}(*)\longrightarrow A^{(r)}_0,
  \qquad A^{(r)}_0 := A^{(r)}\times_A\{e\},
\]
where $A^{(r)}\to A$ is the summation map.

\subsection{Loop spaces for diagonalizable groups}
\label{ssec:LA(BG) diag}

Let $G$ be a strict diagonalizable group scheme over $S$.
That is,
\[
  G = D_S(M) := \spc_S(\OO_S[M])
  \simeq \HOM_S(M_S, \GG_{m,S}),
\]
where $M$ is a finitely generated abelian group; see \cite[Constr.~1.6.10 and Rem.~1.6.11]{LurieEll2}.
We also write $X^*(G) = M$ for the group of characters.

In this subsection we will describe the loop space $\Loop_A^G(*)$:

\begin{proposition}\label{prop:loop BG diag}
  There is a canonical $\check{A}$-equivariant isomorphism
  \begin{equation*}
    \Loop_{A/S}(BG)
    \xrightarrow{\ \sim\ }
    \HOM_S^{E_\infty}\bigl(M_S,\PIC(\stHcat_{A/S})\bigr)
  \end{equation*}
  of stacks over $S$.
  Moreover, this isomorphism is natural in the finitely generated abelian group $M$.
\end{proposition}

\begin{example}\label{exam:loop BGm vs Pic}
  For $G=\GGm$, we have
  \[
    \Loop_{A/S}(B\GGm) \simeq
    \HOM_S^{E_\infty}\bigl(\bZ_S, \PIC(\stHcat_{A/S})\bigr) \simeq
    \PIC(\stHcat_{A/S})^{\mathrm{str}}.
  \]
  Similarly, $\Loop_{A/S}(B\mu_n)$ is the strictification of the stack of $\star$-invertible objects $F$ of $\stHcat_{A/S}$ equipped with a trivialization $F^{\star n} \simeq \OO_e$.
\end{example}

\subsubsection{}

For any commutative group stack $H$ over $S$, there is a canonical homomorphism of commutative group stacks over $S$,
\begin{equation}\label{eq:B vs HOM}
  B\HOM_S(M_S, H) \to
  \HOM_S(M_S, BH).
\end{equation}

\begin{proposition}\label{prop:Homstr to BGm}
  For $H=\GG_{m,S}$, the map \eqref{eq:B vs HOM} induces the isomorphism
  \begin{equation}\label{eq:Homstr to BGm}
    BG = BD_S(M)
    \xrightarrow{\sim}
    \HOM_S(M_S, B\GG_{m,S})
  \end{equation}
   of commutative group stacks over $S$.
\end{proposition}

\begin{proof}
Choose a finite free resolution
\[
  0\longrightarrow\bZ^a\xrightarrow{\phi}\bZ^b
  \longrightarrow M\longrightarrow0
\]
and let $f : \GG_{m,S}^b \to \GG_{m,S}^a$ be the homomorphism induced by $\phi$.
Regarding this as a cofibre sequence in commutative group stacks over $S$, applying $\HOM_S(-, B\GG_{m,S})$ and $\HOM_S(-, \GG_{m,S})$ shows that the source and target of \eqref{eq:B vs HOM} identify respectively with
\[
  B\operatorname{fib}(f)
  \qquad\text{and}\qquad
  \operatorname{fib}(Bf).
\]
By \cite[Thm.~3.48]{Bunk} it will suffice to show that $f$ is an effective epimorphism.
Since $f$ is affine, it will suffice by \cite[Rem.~6.3.3.8]{SAG} to show that it is faithfully flat.
For this we may assume that $f$ is of the form
\[
  f: (\GG_{m,S})^b \to (\GG_{m,S})^a, \qquad
  (t_1,\ldots,t_b) \mapsto
  (t_1^{d_1},\ldots,t_a^{d_a}).
\]
Write $f=([d_1]\times\cdots\times[d_a])\circ p$, where $p:(\GG_{m,S})^b\to(\GG_{m,S})^a$ forgets the last $b-a$ coordinates.
Since the morphism $[d] : \GGm \to \GGm$ is faithfully flat for every $d>0$, as is $p$, the same holds for $f$.
\end{proof}

\subsubsection{}\label{sssec:QCohBG diag}

By \cite[\S 4]{MoulinosGeometry}, there is an equivalence of sheaves of categories over $S$,
\[
  \stQCoh_{BG/S} \simeq \stQCoh_{S}^{M\dashgr},
\]
where the right-hand side is the sheaf of tensor categories
$$ T \mapsto \stQCoh_{S}^{M\dashgr}(T) := \operatorname{Fun}(M^{\mathrm{disc}}, \QCoh(T)), $$
where $M^{\mathrm{disc}}$ denotes $M$ regarded as a discrete category; objects of the right-hand side may be identified with tuples $V = (V_m)_{m\in M}$, $V_m \in \QCoh(T)$.
This equivalence respects the natural functoriality of both sides with respect to group homomorphisms $M \to M'$: $\stQCoh_{BG/S}$ by pullback along $BG' = BD(M') \to BG = BD(M)$, and $\stQCoh_{S}^{M\dashgr}$ by left Kan extension.
Moreover, by Lemma~\ref{lem:easyKunnethforQCoh} both assignments send products in $M$ to tensor products of sheaves of categories, compatibly with the above equivalence.
In particular, the Hopf monoid structure on $M$, where the multiplication, comultiplication, unit, counit, and antipode are encoded respectively by
\[
  + : M^2\to M,\qquad
  \Delta : M\to M^2,\qquad
  0,\qquad
  M\to0,\qquad
-\mathrm{id}_M,
\]
gives rise to Hopf tensor structures on $\stQCoh_{BG/S}$ and $\stQCoh_{S}^{M\dashgr}$, which are identified under the above equivalence.
Identifying objects of $\stQCoh_{S}^{M\dashgr}(T)$ with tuples $V = (V_m)_{m\in M}$, $V_m \in \QCoh(T)$, the tensor product is Day convolution, with $(V \star W)_k = \bigoplus_{m+n=k} V_m \otimes W_n$.

\subsubsection{}\label{sssec:MdisctoPic}

For every tensor category $\CC$ and every $T$ over $S$, we have the canonical equivalences
\begin{equation*}
  \begin{aligned}
    \operatorname{Fun}_{\QCoh(T)}^{L,\otimes}
      \bigl(\stQCoh_{S}^{M\dashgr}(T), \CC\bigr)
    &\simeq \operatorname{Fun}^\otimes
      \bigl(M^{\mathrm{disc}}, \CC\bigr) \\
    &\simeq \operatorname{Fun}^\otimes
      \bigl(M^{\mathrm{disc}}, \operatorname{Pic}(\CC)\bigr),
  \end{aligned}
\end{equation*}
since every object $m\in M^{\mathrm{disc}}$ is $\otimes$-invertible (with inverse $-m$).
By the universal property of Day convolution, $\operatorname{Fun}(M^{\mathrm{disc}}, \CC)$ is the free $\CC$-linear cocomplete tensor category on $M^{\mathrm{disc}}$.

\subsubsection{Proof of Proposition~\ref{prop:loop BG diag}}

By (\ref{sssec:MdisctoPic}) we have, naturally in stacks $T$ over $S$,
\begin{align*}
  \TENSFUN_S\bigl(\stQCoh_S^{M\dashgr}, \stHAS\bigr)(T)
  &\simeq \TensFun_T\bigl(\operatorname{Fun}(M^{\mathrm{disc}},\QCoh(T)),\, \Hcat_{A_T/T}\bigr)^\simeq\\
  &\simeq \TensFun\bigl(M^{\mathrm{disc}},\, \Hcat_{A_T/T}\bigr)^\simeq\\
  &\simeq \TensFun\bigl(M^{\mathrm{disc}},\, \Pic(\Hcat_{A_T/T})\bigr)\\
  &\simeq \Map_{E_\infty}\bigl(M,\, \Pic(\Hcat_{A_T/T})\bigr).
\end{align*}
Combining this with (\ref{sssec:QCohBG diag}), we have:
\begin{align*}
  \Loop_{A/S}(BG)
  &= \TENSFUN_S\bigl(\stQCoh_{BG/S},\, \stHAS\bigr)\\
  &\simeq \TENSFUN_S\bigl(\stQCoh_{S}^{M\dashgr},\, \stHAS\bigr)\\
  &\simeq \HOM_S^{E_\infty}\bigl(M_S,\, \PIC(\stHAS)\bigr).
\end{align*}
This isomorphism is tautologically $\check{A}$-equivariant, since the $\check{A}$-action on both sides arises from the action on $\stHAS$.
\qed

\subsubsection{}

We now suppose given a morphism of stacks over $S$, strictly multiplicative in each variable,
\begin{equation}\label{eq:evstr}
  \ev_A^{\mathrm{str}}:
  \check A\times_S A
  \to
  B\GGmS
\end{equation}
lifting the evaluation map $\ev_A : \check{A} \times_S A \to B\GL_{1,S}$.
Equivalently, this is the datum of a $\GGm$-torsor $\mathcal{P}^{\mathrm{str}}_A$, strictly multiplicative in each variable, over $\check{A} \times_S A$, lifting the multiplicative line bundle $\mathcal{P}_A$ classified by $\ev_A$.

\subsubsection{}\label{sssec:tot P_A^str}

The total space of the $\GGm$-torsor $\mathcal{P}_A^{\mathrm{str}}$ is the fibre of the homomorphism of group stacks over $A$,
\[
  (\pr_2, \ev_A^{\mathrm{str}}) : \check{A} \times_S A \to A \times_S B\GGmS,
  \qquad
  (x,a) \mapsto (a, \mathcal{P}_A^{\mathrm{str}}|_{\{(x,a)\}}).
\]
It is thus a commutative group stack whose $T$-points are triples $(a,x,s)$ with $a$ a $T$-point of $A$, $x$ a $T$-point of $\check{A}$, and $s$ a section of the $\GGm$-torsor $(x,a)^*\mathcal{P}_A^{\mathrm{str}}$ over $T$; its group law is
\[
  (x,s) \cdot (y,t) = (x+y,\, s \otimes t).
\]

\subsubsection{}\label{sssec:P_A^str act on A x BG_m}

Note that $\mathcal{P}_A^{\mathrm{str}}$ acts on $A \times_S B\GGmS$ by translation along the homomorphism of group stacks over $A$,
\begin{equation}\label{eq:P_A^str act on A x BG_m}
  \mathcal{P}_A^{\mathrm{str}} \to \check{A} \times_S A \xrightarrow{(\pr_2, \ev_A^{\mathrm{str}})} A \times_S B\GGmS,
\end{equation}
i.e. by $(x,s) \cdot (a,Q) := (a,\, Q \otimes \mathcal{P}_A^{\mathrm{str}}|_{\{(x,a)\}})$.
Since $\mathcal{P}_A^{\mathrm{str}}$ is the fibre of $(\pr_2, \ev_A^{\mathrm{str}})$, \eqref{eq:P_A^str act on A x BG_m} is canonically null-homotopic as a homomorphism (of group stacks over $A$).
This null-homotopy makes the projection $\pr_2 : A \times_S B\GGmS \to B\GGmS$ equivariant for the trivial action on $B\GGmS$.

\subsubsection{}\label{sssec:O_A boxtimes O(1) equiv}

Write also $\OO(1)$ for the pullback along $B\GGmS \to B\GL_{1,S}$ of the tautological line bundle $\OO(1)$ (\ref{sssec:Phi}).

The fact that $\pr_2 : A \times_S B\GGmS \to B\GGmS$ is $\mathcal{P}_A^{\mathrm{str}}$-equivariant with respect to the trivial action on the target (\ref{sssec:P_A^str act on A x BG_m}) implies that the line bundle $\pr_2^*\OO(1) = \OO_A \boxtimes_S \OO(1)$ on $A \times_S B\GGmS$ tautologically acquires a $\mathcal{P}_A^{\mathrm{str}}$-equivariant structure.

\subsubsection{}

Consider the canonical homomorphism of $E_\infty$-group stacks over $S$,
\begin{equation}\label{eq:A x BGm to Pic}
	A \times_S B\GGmS \to
	A \times_S B\GL_{1,S}
	\xrightarrow{\ \eqref{eq:A x Pic to PicH}\ }
	\PIC(\stHAS),
\end{equation}
given informally by $(a,Q)\mapsto a_*(Q^{\mathrm{und}})$, where $Q^{\mathrm{und}}$ is the line bundle underlying the $\GGm$-torsor $Q$.
Recall that the $\check{A}$-action on $\PIC(\stHAS)$ is given by $x\cdot F := F\otimes \mathcal{P}_A|_{\{x\} \times A}$.
If $\check A$ acts on $A\times_S B\GGmS$ by
\[
  x \cdot (a,Q) :=
  \bigl(a, Q \otimes \mathcal{P}_A^{\mathrm{str}}|_{\{(x,a)\}} \bigr),
\]
then the projection formula for $a_*$ shows that \eqref{eq:A x BGm to Pic} is $\check A$-equivariant.

\subsubsection{}

Let $f_A:A\times_S B\GGmS\to\PIC(\stHAS)$ denote the homomorphism \eqref{eq:A x BGm to Pic}.
Post-composition with $f_A$ and Proposition~\ref{prop:loop BG diag} give the homomorphism
\begin{equation}\label{eq:yowdvbf}
  \HOM_S(M_S, A \times_S B\GGmS) \to
  \HOM_S^{E_\infty}(M_S, \PIC(\stHAS)) \simeq
  \Loop_A(BG).
\end{equation}

\subsubsection{}\label{sssec:XAG}

Set $\frX_A(G) := \HOM_S(M_S, A)$.
This is a commutative group scheme over $S$ which, when $G=T$ is a split torus of rank $r$, can be written
\[
  \frX_A(T) \simeq A^r \simeq A \otimes_{\bZ} \mathbf{X}_*(T).
\]
More generally, writing $M \simeq \bZ^r \oplus \bZ/n_1 \oplus \cdots \oplus \bZ/n_s$ with $n_1 \mid \cdots \mid n_s$, we have
\[
  \frX_A(G) \simeq A^r \times_S A[n_1] \times_S \cdots \times_S A[n_s],
\]
which is fppf over $S$ whenever each $[n_j] : A \to A$ is faithfully flat.\footnote{%
  This clearly holds for $A=\GGa$ (over $k$ a $\bQ$-algebra), $A=\GGm$, or $A$ an abelian scheme.
}
Since $\HOM_S(-, A)$ takes direct sums to products, this follows from the identifications
\[
  \HOM_S(\bZ_S, A) \simeq A, \qquad
  \HOM_S((\bZ/n)_S, A) \simeq A[n].
\]

\subsubsection{}

By Proposition~\ref{prop:Homstr to BGm}, the source of \eqref{eq:yowdvbf} is identified with $\frX_A(G) \times_S BG$.
Thus we obtain the canonical $\check{A}$-equivariant homomorphism of commutative group stacks over $S$,
\begin{equation}\label{eq:pAG}
	p_{A,G} : \frX_A(G) \times_S BG \to
	\Loop_A(BG).
\end{equation}
By construction, $p_{A,G}$ is determined by the data of \eqref{eq:A x BGm to Pic}.

\begin{proposition}\label{prop:pAG isomorphism}
The morphism $p_{A,G}$ is an isomorphism on classical truncations.
Moreover, if $A\to S$ has relative dimension $1$ and $k$ is a $\bQ$-algebra, then $p_{A,G}$ is an isomorphism.
\end{proposition}

\begin{proof}
Write $Y:=A\times_S B\GGmS$, $X:=\PIC(\stHAS)$ and let $f_A : Y \to X$ be the homomorphism \eqref{eq:A x BGm to Pic}.
Let $g_A : Y \to X^{\mathrm{str}}$ denote the composite of $f_A^{\mathrm{str}} : Y^{\mathrm{str}} \to X^{\mathrm{str}}$ with the unit map $\eta_Y : Y \to Y^{\mathrm{str}}$; then $p_{A,G}$ is identified with the composite
\[
  \HOM_S(M_S,Y)
  \xrightarrow{\ \HOM_S(-, g_A)\ }
  \HOM_S(M_S,X^{\mathrm{str}})
  \simeq
  \HOM_S^{E_\infty}(M_S,X),
\]
cf. Appendix~\ref{app:strictification}.
It will suffice to prove that both $f_A$ and $\eta_Y$ are isomorphisms on classical truncations, and isomorphisms when $A \to S$ is of relative dimension $1$ and $k$ is a $\bQ$-algebra.

We first consider the claim about classical truncations.
For $f_A$ this follows from Example~\ref{ex:strict units comparison} and Proposition~\ref{prop:Pic HA classical}(i).
For every classical affine $T$ over $S$, the space $A(T)$ is $0$-truncated because $A\to S$ is a group scheme, while $(B\GGmS)(T)$ is $1$-truncated because $B\GGmS$ is a $1$-Artin stack.
Thus $Y(T)$ is $1$-truncated, and Corollary~\ref{cor:strictification classical strict source} shows that $\eta_Y$ is an isomorphism on classical truncations.

Suppose next that $A\to S$ has relative dimension $1$ and that $k$ is a $\bQ$-algebra.
By Example~\ref{ex:strict units comparison} and Proposition~\ref{prop:Pic HA classical}(ii), $f_A$ is an isomorphism.
For $\eta_Y$ we will apply Corollary~\ref{cor:strictification strict source}, for which we need to check that $\pi_iY(T)$ is a $\bQ$-vector space for every $i\geq2$ and every affine $T$ over $S$.
Let $T=\spc(R)$ be an affine over $S$.
We check that $\pi_i A(T)$ and $\pi_i B\GGmS(T)$ are $\bQ$-vector spaces for $i\ge 2$:
\begin{itemize}
\item
For every positive integer $n$, we claim that multiplication by $n$ on $\pi_i A(T)$ is invertible.
Indeed, this is realized by the morphism $[n]:A\to A$, which is \'etale: its differential is multiplication by $n$, which is invertible because $\OO_A$ is a $\bQ$-algebra by assumption.
Since \'etale morphisms are $0$-truncated, this means that $[n]:A(T)\to A(T)$ has discrete fibres and hence induces an isomorphism on $\pi_i$ for every $i\geq2$.

\item
Since $k$ is a $\bQ$-algebra, Example~\ref{ex:strict units comparison} identifies $\GGm$ with the $E_\infty$-group stack $\GL_{1}$ representing units.
The inclusion $\GL_{1,S}(R)\subseteq\Omega^\infty R$ is a union of connected components, so for every $i\geq2$ we have:
\[
  \pi_i\bigl((B\GGmS)(T)\bigr)
  \simeq
  \pi_{i-1}\bigl(\GGmS(T)\bigr)
  \simeq
  \pi_{i-1}\bigl(\GL_{1,S}(R)\bigr)
  \simeq
  \pi_{i-1}(R),
\]
which is a module over the $\bQ$-algebra $\pi_0(R)$.
\end{itemize}
The claim follows.
\end{proof}

\begin{remark}
	For $A=\GGm$, Theorem~\ref{thm:LGa LGm} shows that $p_{\GGm,G}$ is an isomorphism unconditionally.
\end{remark}

\subsubsection{}
\label{sssec:LAG X_AG correspondence}

The commutative group stack $\frX_A(G)$ is related to the loop stack $\Loop_A(BG)$ via a canonical correspondence of $\check{A}$-equivariant stacks
\begin{equation}\label{eq:L_AG X_AG correspondence}
  \begin{tikzcd}
    & \frX_A(G) \times_S BG \ar[swap]{ld}{p_{A,G}}\ar{rd}{q_{A,G}} &
    \\
    \Loop_A(BG) & & \frX_A(G),
  \end{tikzcd}
\end{equation}
where $q_{A,G}$ is the projection, equivariant with respect to the trivial action on the target.
Since $G$ is diagonalizable, we may regard $q_{A,G}$ as the projection of a good moduli space: more precisely, under the equivalence of (\ref{sssec:QCohBG diag}), the functor
\[ q_{A,G,*} : \QCOH(\frX_A(G) \times_S BG) \to \QCOH(\frX_A(G)) \]
sends an $M$-graded quasi-coherent sheaf to its weight-zero piece (cf. \ref{sssec:Shv wt=0}).
Since $p_{A,G}$ is an isomorphism on classical truncations (Proposition~\ref{prop:pAG isomorphism}), it follows that the classical truncation of $\Loop_A(BG)$ has good moduli space $\frX_A(G)$.

\subsubsection{}\label{sssec:O(chi)}

For every character $\chi \in X^*(G)$, precomposition with $\chi : \bZ \to X^*(G)$ gives an $\check{A}$-equivariant morphism (see Proposition~\ref{prop:Homstr to BGm})
\begin{multline*}
  \bar{\chi} :
  \frX_A(G) \times_S BG
  \simeq \HOM_S(X^*(G)_S, A \times_S B\GGmS)
  \\
  \to \HOM_S(\bZ_S, A \times_S B\GGmS)
  \simeq A \times_S B\GGmS,
\end{multline*}
sending $(\xi, Q) \mapsto (\xi(\chi), \chi_*Q)$.

Let $\mathcal{P}_\chi$ be the pullback of $\mathcal{P}_A^{\mathrm{str}}$ along $\frX_A(G) \to A$, $\xi \mapsto \xi(\chi)$, whose points are $(\xi,x,s)$ with $s$ a section of $(x,\xi(\chi))^*\mathcal{P}_A^{\mathrm{str}}$.
Then $\mathcal{P}_\chi$ acts on $\frX_A(G) \times_S BG$ by $(\xi,x,s)\cdot(\xi,Q) := x\cdot(\xi,Q)$, and on $A \times_S B\GGmS$ through $\mathcal{P}_\chi \to \mathcal{P}_A^{\mathrm{str}}$.
Since $\mathcal{P}_\chi$ and $\mathcal{P}_A^{\mathrm{str}}$ act through their projections to $\check{A}$, $\bar{\chi}$ is $\mathcal{P}_\chi$-equivariant.

In particular, for every $\chi \in X^*(G)$ we have a line bundle
\[ \OO(\chi) := \bar{\chi}^*(\OO_A \boxtimes_S \OO(1)) \]
on $\frX_A(G) \times_S BG$, which is $\mathcal{P}_\chi$-equivariant by (\ref{sssec:O_A boxtimes O(1) equiv}).

\subsection{Good moduli spaces}\label{ssec:gms Ga Gm}

Let us take $S=\spc(k)$ and $G$ a linearly reductive affine fppf group scheme over $k$.
For simplicity, we will assume that $A=\GGa$ (over a $\bQ$-algebra), $A=\GGm$, or $A=E$ is an elliptic curve over $k$.

Write $\frL_A(G) := \Loop_A^G(*)_\cl$ for the classical truncation of the $G$-equivariant $A$-loop space of the point.
When it is $1$-Artin and admits a good moduli space in the sense of \cite{Alper}, we denote the latter by
\[
  \pi_{A,G} : \frL_A(G) \to \frM_A(G).
\]

\subsubsection{}

When $G$ is \emph{diagonalizable}, Proposition~\ref{prop:pAG isomorphism} shows that the correspondence \eqref{eq:L_AG X_AG correspondence} identifies $\pi_{A,G}$ with the projection
\[ q_{A,G} : \frX_A(G) \times BG \to \frX_A(G). \]

\subsubsection{}

For $A=\GGa$ (over a $\bQ$-algebra) and $A=\GGm$, Example~\ref{exam:equivariant loops Ga Gm} identifies the $G$-equivariant loop space as
\[
  \Loop_{\GGa}^G(*) \simeq \frL_{\GGa}(G) \simeq \mathfrak{g}/G, \qquad
  \Loop_{\GGm}^G(*) \simeq \frL_{\GGm}(G) \simeq G/G.
\]
In particular, they are classical $1$-Artin and, since $G$ is linearly reductive and $\mathfrak{g}$ and $G$ are affine, they admit good moduli spaces
\[
  \pi_{\GGa,G} : \mathfrak{g}/G \to \mathfrak{g} \modmod G,\qquad
  \pi_{\GGm,G} : G/G \to G \modmod G,
\]
which are the affines $\frM_{\GGa}(G) \simeq \spc k[\mathfrak{g}]^G$ and $\frM_{\GGm}(G) \simeq \spc k[G]^G$, respectively (see \cite[Thm.~13.2]{Alper}).
The elliptic case $A=E$ will be studied in Subsect.~\ref{ssec:gms elliptic} below.

\subsubsection{}

Let $k$ be a field, $G$ a split reductive algebraic group over $k$ with identity component $G^\circ$, maximal torus $T$, and Weyl group $W := N/T$ (where $N$ is the normalizer of $T$ in $G$); we assume that $|W|$ is coprime to the characteristic of $k$.
Under the hypothesis that $\frM_A(G)$ exists, we will show that its \emph{trivial} component (i.e., the component of the trivial loop) can be recovered from the case $G=T$ as $\frM_A(T)\modmod W$ (where by the diagonalizable case discussed above, the latter is given by $\frX_A(T)$).

The conjugation action of $N$ on $T$ induces a $W$-action on $\frL_A(T)$ for which the map $\frL_A(T) \to \frL_A(G) \to \frM_A(G)$ is $W$-invariant, hence factors through the quotient stack $\frL_A(T)/W$.
Since $T$ is a split torus, $\frM_A(T) \simeq \frX_A(T) \simeq A^{\operatorname{rank}(T)}$ (\ref{sssec:XAG}) is separated over $k$.
The following lemma then shows that $\frL_A(T)/W$ admits the good moduli space $\frM_A(T)\modmod W$.

\begin{lemma}\label{lem:gms descent}
  Let $H$ be a finite group of order coprime to the characteristic of $k$, regarded as a linearly reductive algebraic group over $k$.
  Let $\YY$ be a classical $1$-Artin stack of finite type over $k$, equipped with an action of $H$.
  If $\YY$ admits a good moduli space $p : \YY \to M$ with $M$ separated, then the $H$-action descends to $M$, the (separated) quotient stack $M/H$ admits a good moduli space $M \modmod H$, and the induced map $q : \YY/H \to M/H \to M \modmod H$ is a good moduli space, with $M \modmod H$ separated of finite type over $k$.
\end{lemma}
\begin{proof}
  The action descends by universality of good moduli spaces.
  Since $M$ is separated, and of finite type by \cite[Thm.~4.16(xi)]{Alper}, the quotient stack $M/H$ is of finite type.
  Since it has finite diagonal and linearly reductive stabilizers, it is separated and tame in the sense of \cite[Thm.~3.2]{AOV}; by Keel--Mori it thus admits a separated coarse moduli space $\rho : M/H \to M \modmod H$, which by \cite[Rem.~4.2, Ex.~8.1]{Alper} is also a good moduli space.

  In particular, $\rho$ is cohomologically affine and the unit map $\OO_{M \modmod H}\to \rho_*\OO_{M/H}$ is invertible.
  The analogous properties hold for the induced map $\bar{p} : \YY/H \to M/H$, by faithfully flat descent from $p$, and hence for the composite $q = \rho \circ \bar{p} : \YY/H \to M/H \to M\modmod H$.
  Since the target is an algebraic space, we deduce that $q$ is a good moduli space map, of finite type by \cite[Thm.~4.16(xi)]{Alper}.
\end{proof}

By the universal property of good moduli spaces we thus obtain a commutative diagram
\[\begin{tikzcd}
  \frL_A(T) \ar[twoheadrightarrow]{r}\ar[swap]{d}{\pi_{A,T}}
  & \frL_A(T)/W \ar{d}\ar{rd}
  &
  \\
  \frM_A(T) \ar[twoheadrightarrow]{r}
  & \frM_A(T)\modmod W \ar{r}
  & \frM_A(G).
\end{tikzcd}\]

Since $\frM_A(T)/W \simeq \frX_A(T)/W$ is connected, so is $\frM_A(T)\modmod W$ (see \cite[Thm.~4.16(viii)]{Alper}).
Assuming that $\frM_A(G)$ is locally noetherian, so that its connected components are open and closed, the map $\frM_A(T)\modmod W \to \frM_A(G)$ factors uniquely through a canonical map
\begin{equation}\label{eq:MAG vs MAT}
	\frM_A(T)\modmod W \longrightarrow \frM_A(G)^\circ,
\end{equation}
where $\frM_A(G)^\circ \subset \frM_A(G)$ is the component of the trivial loop.

In the cases $A=\GGa$ and $A=\GGm$, \eqref{eq:MAG vs MAT} recovers the well-known isomorphisms
\[
  \mathfrak{t} \modmod W \simeq \mathfrak{g} \modmod G, \qquad
  T \modmod W \simeq G^\circ \modmod G.
\]
That \eqref{eq:MAG vs MAT} is also invertible in the elliptic case $A=E$ will be shown in Subsect.~\ref{ssec:gms elliptic} below.

\section{The case of abelian varieties}\label{sec:abvar}

\subsection{Artinness}\label{ssec:abvar Artin}

In the additive and multiplicative cases, we have seen in Theorem~\ref{thm:LGa LGm} that the loop space $\Loop_A(\YY)$ is Artin for any Artin stack $\YY$.
We now show that the same holds in the elliptic case, or more generally for $A$ any abelian variety, at least when $\YY$ is pluperfect.

\begin{proposition}\label{prop:open in map implies Artin}
  Let $A$ be a relative abelian variety over a stack $S$.
  If $\YY$ is Artin over $S$ and satisfies $(\textup{Map}\textsubscript{A/S})$, then $\Loop_{A/S}(\YY)$ is Artin over $S$.
\end{proposition}
\begin{proof}
  It suffices to show that $\MAP_S(\check{A}, \YY)$ is Artin over $S$.
  Since this mapping stack is stable under base change in $S$, we may assume $S$ is affine.
  Since $\check{A}$ is a smooth proper scheme over $S$, the claim can then be deduced from Lurie's version of Artin's representability theorem (see e.g.\ the discussion following the proof of Theorem~2.5 in \cite{PTVV}).
\end{proof}

We now show that ($\textup{Map}\textsubscript{A/S}$) holds for pluperfect stacks $\YY$ (Definition~\ref{def:pluperfect}).
We will need the following two lemmas.

\begin{lemma}\label{lem:std suff for hmg}
  Let $S$ be an affine scheme, $A \to S$ a relative abelian scheme, and $\pi : \YY \to BGL_n$ a quasi-affine morphism of stacks.
  Denote by $V_{\mathrm{std}} \in \Perf(BGL_n)$ the standard representation.
  For any morphism $f : \check{A} \to \YY$, the following are equivalent:
  \begin{enumerate}
    \item $f^*(\pi^* V_{\mathrm{std}}) \in \Perf^{\,\mathrm{hmg}/S}(\check{A})$.
    \item $f^*(\pi^* V) \in \INDPERF^{\mathrm{hmg}/S}(\check{A})$ for all $V \in \INDPERF(BGL_n)$.
    \item $f^*(F) \in \INDPERF^{\mathrm{hmg}/S}(\check{A})$ for all $F \in \INDPERF(\YY)$.
  \end{enumerate}
\end{lemma}
\begin{proof}
  Let $f : \check{A} \to \YY$ be a morphism such that $E := f^*(\pi^* V_{\mathrm{std}}) \in \Perf^{\,\mathrm{hmg}/S}(\check{A})$.
  Let us show that $g^*(V) \in \INDPERF^{\mathrm{hmg}/S}(\check{A})$ for all $V \in \INDPERF(BGL_n)$, where $g = \pi \circ f : \check{A} \to BGL_n$.
  It is enough to show that $g^*(V) \in \Perf^{\,\mathrm{hmg}/S}(\check{A})$ for all $V \in \Perf(BGL_n)$.
  Replacing $S$ by a geometric point $s$, and $f$ by $f_s$, this amounts to showing that the cohomology sheaves of $g^*(V)$ are homogeneous vector bundles on $\check{A}$.

  Since $E$ is a homogeneous vector bundle on $\check{A}$, it is invariant under translations $t_a : \check{A} \to \check{A}$, i.e., $t_a^* E \simeq E$ for all geometric points $a$ of $\check{A}$.
  Under the equivalence of groupoids $\Map(\check{A}, BGL_n) \simeq \operatorname{Vect}_n(\check{A})$, this implies that the classifying morphism $g : \check{A}_s \to BGL_n$ is translation-invariant, i.e., $g \circ t_a \simeq g$ for all $a$.
  Thus $g \circ t_a \simeq g$ for all translation maps $t_a$ on $\check{A}$, so for \emph{any} $E \in \PERF(BGL_n)$ we have $t_a^* g^*(E) \simeq (g \circ t_a)^*(E) \simeq g^*(E)$, i.e., $g^*(E)$ is translation-invariant in $\PERF(\check{A})$.
  In particular, its cohomology sheaves are translation-invariant coherent sheaves, and hence by \cite[Exam.~3.2]{Mukai} are homogeneous vector bundles.
  Thus $g^*(E) \in \INDPERF^{\mathrm{hmg}/S}(\check{A})$ for all $E \in \INDPERF(BGL_n)$.

  Next let $F \in \INDPERF(\YY)$ and let us show that $F$ lies in the full subcategory
  \[ \mathcal{C} := \{ G \in \INDPERF(\YY) : f^*(G) \in \INDPERF^{\mathrm{hmg}/S}(\check{A})\}. \]
  Since $f^*$ is colimit-preserving and $\INDPERF^{\mathrm{hmg}/S}(\check{A})$ is closed under colimits, $\mathcal{C}$ is closed under colimits in $\INDPERF(\YY)$.
  By the first part, the essential image of $\pi^*$ lies in $\mathcal{C}$.
  Since $\pi$ is quasi-affine, the pushforward $\pi_* : \INDPERF(\YY) \to \INDPERF(BGL_n)$ is conservative and colimit-preserving.
  Therefore, $F$ can be written as the geometric realization of a simplicial object whose terms lie in the essential image of $\pi^*$ (see \cite[Prop.~4.7.3.14]{LurieHA}), and hence in $\mathcal{C}$.
  Since $\mathcal{C}$ is closed under colimits, we conclude that $F \in \mathcal{C}$.
\end{proof}

\begin{lemma}\label{0-open}
  Let $S$ be an affine scheme, $Y$ a proper and almost finitely presented scheme over $k$, and $F\in\Perf(Y\times S)$ a perfect complex.
  Then the set $S' \subseteq S$ of points $s \in S$ such that $s^*F = F|_{Y_{k(s)}}$ lies in $\Perf_{0}(Y_{k(s)})$ is a quasi-compact open subset of $S$, and $F|_{Y \times S'}$ lies in $\Perf_{0}(Y\times S'/S')$.
\end{lemma}

\begin{proof}
Let $Z$ denote the support of $F \in \Perf(Y\times S)$. Since $F$ is perfect, $Z$ is a closed subset of $Y \times S$ with quasi-compact complement (see \cite[Lem.~3.3]{ThomasonClassification}). Since $Y$ is proper, the projection $\pi|_Z: Z \subset Y \times S \to S$ is proper.

Note that $S'$ is the locus of points $s \in S$ where the fibre $Z_s = (\pi|_Z)^{-1}(s) = \operatorname{Supp}(s^* F)$ has dimension $\le 0$.
Since $\pi|_Z$ is proper, Chevalley's upper semi-continuity theorem implies that $S'$ is open in $S$.
It is also quasi-compact: the set of $z \in Z$ with $\dim_z Z_{\pi(z)} \ge 1$ is locally constructible by \cite[IV$_3$, Prop.~9.9.1(i)]{EGA}, so its image $S \setminus S'$ is constructible by Chevalley's theorem (see \cite[IV$_1$, Thm.~1.8.4]{EGA}), and any open constructible subset of an affine scheme is quasi-compact.
Then the restriction $\pi|_{Z \cap (Y \times S')} : Z \cap (Y \times S') \to S'$ is proper and quasi-finite, hence finite.
Thus $F|_{Y \times S'}$ is supported on a closed subset that is finite over $S'$, i.e., $F|_{Y \times S'} \in \Perf_{0}(Y\times S'/S')$ as desired.
\end{proof}

\begin{theorem}\label{open in map pluperfect} If $\YY$ is pluperfect fppf-locally over $S$, then $\YY$ satisfies $(\textup{Map}\textsubscript{A/S})$. In particular, $\Loop_{A/S}(\YY)$ is Artin over $S$. \end{theorem}

\begin{proof}
  The claim is fppf-local on $S$, so we may assume that $\YY$ is pluperfect.
  Let $\{F_i\}_i$ be a finite set of perfect tensor generators of $\QCOH(\YY)$.
  Then $\Loop_A(\YY)$ is the substack of $\MAP(\check{A}, \YY)$ of maps $f$ such that $f^*(F_i) \in \INDPERF^{\mathrm{hmg}}(\check{A})$ for all $i$, and the claim is that this condition is open.
  Let $f : S \times \check{A} \to \YY$ be a morphism with $S$ affine.
  Let $G_i := \Psi(f^* F_i)$ for $\Psi : \PERF(\check{A}_S) \to \PERF(A_S)$ the inverse Fourier--Mukai transform.
  For $s \in S$, $f_s^* F_i \in \INDPERF^{\mathrm{hmg}}(\check{A}_{s})$ if and only if $G_i|_{A_s} \in \PERF_0(A_{s})$.
  By Lemma~\ref{0-open}, for each $i$ the set $S'_i$ of $s \in S$ with $G_i|_{A_s} \in \PERF_0(A_s)$ is a quasi-compact open subset of $S$; hence $\Loop_A(\YY) \times_{\MAP(\check{A}, \YY)} S = \bigcap_i S'_i$ is a quasi-compact open subscheme of $S$.
\end{proof}

\subsection{Finiteness}

We also record the following finiteness property:

\begin{proposition}\label{prop:fin pres abvar}
  Let $A$ be an abelian variety over $k$, and let  $f : \XX \to \YY$ be a morphism of stacks satisfying $(\textup{Map}\textsubscript{A})$.
  If $f$ is locally of finite presentation, then so is $\Loop_A(f) : \Loop_A(\XX) \to \Loop_A(\YY)$.
\end{proposition}
\begin{proof}
  Factor $\Loop_A(f)$ as
  \[ \Loop_A(\XX) \hookrightarrow \ZZ \to \Loop_A(\YY), \qquad \ZZ := \MAP(\check{A}, \XX) \times_{\MAP(\check{A}, \YY)} \Loop_A(\YY). \]
  The second map is a base change of $\MAP(\check{A},f)$ which, since $\check{A}$ is quasi-compact and separated, is locally of finite presentation when $f$ is (one argues as in the proof of \cite[Prop.~19.1.3.1(4)]{SAG}).
  For the first map, note that $\ZZ \hookrightarrow \MAP(\check{A}, \XX)$ is an open immersion (base change of $\Loop_A(\YY) \hookrightarrow \MAP(\check{A}, \YY)$), and $\Loop_A(\XX) \hookrightarrow \MAP(\check{A}, \XX)$ is an open immersion by hypothesis.
  Since $\Loop_A(\XX) \subset \ZZ$ (by functoriality of $\Loop_A$: if $g \in \Loop_A(\XX)(T)$ then $f \circ g \in \Loop_A(\YY)(T)$), the inclusion $\Loop_A(\XX) \hookrightarrow \ZZ$ is an open immersion.
  Since open immersions are locally of finite presentation, the composite $\Loop_A(f)$ is locally of finite presentation.
\end{proof}

\subsection{Classifying stacks}\label{ssec:classifying stacks}

We now describe the elliptic loop space for $\YY = BG$, where $G$ is linearly reductive. We write $\HA^{\omega,\heartsuit} \subset \HA^{\heartsuit}$ for the abelian categories of coherent and quasi-coherent sheaves on $A$ with $0$-dimensional support, respectively.

\subsubsection{}

We take $S=\spc(k)$ where $k$ is a connective $E_\infty$-ring and $G$ is an fppf affine group scheme over $k$.

\begin{definition}
  Let $C$ be a smooth proper curve over $k$.
  The \emph{spectral moduli stack of $G$-bundles on $C$}, denoted $\operatorname{Bun}_G(C)$, is the mapping stack $\MAP(C, BG)$.
\end{definition}

The Lurie--Artin representability theorem implies that $\operatorname{Bun}_G(C)$ is an Artin stack (see \cite[Ch.~27]{SAG}).
Its classical truncation is the classical Artin stack of $G_\cl$-bundles on the classical curve $C_\cl$,
\[ \Bun_G(C)_\cl \simeq \frBun_{G_\cl}(C_\cl), \]
studied in Subsect.~\ref{ssec:BunG} below.
Note that even if $k$ is an ordinary commutative ring (so that $G$ and $C$ are classical), $\Bun_G(C)$ is typically not classical unless $k$ is a $\bQ$-algebra or $G$ is finite \'etale.

\subsubsection{}

Take $C=E$ to be an elliptic curve and assume that $G$ is linearly reductive.
Recall from (\ref{sssec:BunG diamond}) the classical substack
\[
  \frBun_{G_\cl}^\diamond(E_\cl)
  \subset \frBun_{G_\cl}(E_\cl)
\]
consisting of $G$-bundles $F$ with $F \times^G V$ semistable of degree zero for all finite-dimensional $V \in \Rep^{\heartsuit}(G_\cl)$.
This is open by Proposition~\ref{prop:frBunGdiamond open}.

\begin{definition}\label{def:BunG diamond}
  Denote by $\Bun_G^\diamond(E)$ the unique open substack of $
  \Bun_G(E)$ whose classical truncation is the open $\frBun_{G_\cl}^\diamond(E_\cl) \subset \frBun_{G_\cl}(E_\cl)$.
\end{definition}

\begin{theorem}\label{thm:loop BG bun spectral}
  Let $k$ be a connective $E_\infty$-ring, $G$ a linearly reductive group scheme over $k$, and $E$ an elliptic curve over $k$.
  Then there is a canonical isomorphism
  $$\Loop_E^G(*) \simeq \Bun_G^\diamond(\check{E}).$$
\end{theorem}
\begin{proof}
	By Theorem~\ref{open in map pluperfect}, we may identify $\Loop_E^G(*)$ with an open substack of $\operatorname{Bun}_G(\check{E})$.
  Since open substacks are determined by their classical truncations, it will suffice to identify its classical truncation with $\frBun_{G_\cl}^\diamond(\check{E}_\cl)$.

  Assuming without loss of generality that $k$ is an ordinary commutative ring, we show $\Loop_E^G(*)_\cl \simeq \frBun_G^\diamond(\check{E})$.
  Since the question is Nisnevich-local on $S=\spc(k)$, we may assume by \cite[Cor.~6.2(2), Rem.~2.5]{AHRLocalStructure} that $BG$ has the resolution property.
  Since $G$ is also linearly reductive, the dualizable objects of $\Rep^\heartsuit(G)$ are compact projective and generate under colimits (see e.g. \cite[Prop.~A.3.4]{BKRS}).
  Thus by Lemma~\ref{abelian} we have
  \[
    \Loop_{E}^G(*)_{\cl}
    \simeq \TENSFUN\bigl( \stRep^{\heartsuit}(G), \stHcat_{E}^{\heartsuit} \bigr),
  \]
  and this is equivalently the stack of symmetric monoidal $k$-linear functors $\stRep^{\heartsuit,\mathrm{dzbl}}(G) \to \stHcat_{E}^{\heartsuit,\mathrm{dzbl}}$ from the $1$-category of finite projective representations of $G$ to the $1$-category of $\star$-dualizable objects of $\stHcat_{E}^{\heartsuit}$.
  The Fourier--Mukai transform gives an equivalence $\stHcat_{E}^{\heartsuit,\mathrm{dzbl}} \simeq \textbf{\textsc{Vect}}^{\ss,0}_{\check{E}}$ with the $1$-category of vector bundles on $\check{E}$ which are fibrewise semistable of degree zero (Remark~\ref{rem:HA vs QCoh checkA dzbl}), so we have:
  \[
    \Loop_{E}^G(*)_{\cl}
    \simeq \TENSFUN\bigl( \stRep^{\heartsuit,\mathrm{dzbl}}(G), \textbf{\textsc{Vect}}^{\ss,0}_{\check{E}} \bigr).
  \]
  Applying Lemma~\ref{abelian} again, now with target $\QCoh(\check{E})^\heartsuit$, and Tannaka duality for the perfect stack $BG$, the right-hand side is the substack of $\frBun_{G}(\check{E})$ consisting of $G$-bundles $P$ such that $P \times^G V$ is fibrewise semistable of degree zero for every finite projective $V \in \Rep(G)^\heartsuit$.
  By Lemma~\ref{lem:diamond one rep}, this is precisely the condition that $P \in \frBun_G^\diamond(\check{E})$.
  Thus $\Loop_E^G(*)_\cl \simeq \frBun_{G_\cl}^\diamond(\check{E}_\cl)$, as required.
  (See also \cite[Thm.~5.5]{FGL} for the case of $k$ an algebraically closed field.)
\end{proof}

\subsubsection{}\label{sssec:LG P component}

Suppose that $k$ is algebraically closed and $G$ is smooth affine and linearly reductive.
By Theorem~\ref{thm:loop BG bun spectral}, we have
\[
  \frL_E(G):=\Loop_E^G(*)_\cl
  \simeq\frBun_G^\diamond(\check{E}).
\]
By (\ref{sssec:BunG decomposition}) we have a finite open
and closed decomposition
\[
  \frL_E(G)\simeq\coprod_{[P]}\frL_E(G)_P,
  \qquad
  \frL_E(G)_P\simeq\UU_P/H_P,
\]
where $[P]$ ranges over isomorphism classes of $\Gamma$-torsors on $
\check{E}$, for $\Gamma := G/G^\circ$, and $\UU_P$ and $H_P$ are as defined there.

\subsubsection{}\label{sssec:LG circ G^o}

Let $P_0$ be the trivial $\Gamma$-torsor on $\check{E}$.
By (\ref{sssec:U_P0 = BunGcirc}), we have
\[
  \frL_E(G)_{P_0}\simeq\frL_E(G^\circ)/\Gamma.
\]

Write $\frL_E(G)^\circ$ for the connected component of the trivial loop.
The $\Gamma$-action on $\frL_E(G^\circ)$ preserves the isomorphism class of
the trivial bundle, hence restricts to $\frL_E(G^\circ)^\circ$.
This component is open and closed (since $\frL_E(G^\circ)$ is locally
noetherian); its quotient by $\Gamma$ is therefore a connected open and closed substack of $\frL_E(G)_{P_0}$ containing the trivial loop.
Since $\frL_E(G)_{P_0}$ is open and closed in $\frL_E(G)$, we obtain
\[
  \frL_E(G)^\circ
  \simeq \frL_E(G^\circ)^\circ/\Gamma.
\]

\begin{remark}
  When $G$ is connected and $\pi_1(G)$ is torsion-free (e.g. $GL_r$, $SL_r$, or a torus), Proposition~\ref{prop:frBunGdiamond vs ss tors top type} below identifies $\frBun_G^\diamond(\check{E})$ with the locus $\frBun_G^{\ss,0}(\check{E})$ of semistable bundles of trivial topological type.
\end{remark}

\subsubsection{Elliptic Springer resolution}\label{springer}

Let $k$ be a field of characteristic zero.
For $B \subset GL_r$ the Borel subgroup of upper triangular matrices, the classifying stack $BB$ is perfect, with $\QCOH(BB)$ compactly generated by the irreducible representations of $B$ (see \cite[Thm.~A(c)]{HallRydh}).
These are the characters of $T = B/U$, where $U$ is the unipotent radical, as $V^U \neq 0$ for every representation $V \neq 0$ (see \cite[Exp.~XVII, Prop.~3.2]{SGA3}).
Thus $BB$ is pluperfect, with the $r$ coordinate characters as perfect tensor generators, and the proof of Theorem~\ref{open in map pluperfect} identifies $\Loop_E^B(*)$ with the substack of $\operatorname{Bun}_B(\check E)$ on which the associated character line bundles have fibrewise degree zero.
In other words,
\[
  \Loop_E^B(*)\simeq\operatorname{Bun}_B^0(\check E),
\]
the stack of flags whose graded line bundles have fibrewise degree zero.
Under the Fourier--Mukai transform $\FM^*_{E}$, these correspond to flags $0=F_0\subset\cdots\subset F_r$ with $F_i/F_{i-1}\in\operatorname{Coh}_0^1(E)$.
The map $\Loop_E^B(*)\to\Loop_E^{GL_r}(*)$ forgets the flag, and is precisely the elliptic Grothendieck--Springer resolution of \cite[\S 2]{BNEll}.

\subsection{Good moduli spaces}\label{ssec:gms elliptic}

Take $S=\spc(k)$ with $k$ an algebraically closed field, let $G$ be a linearly reductive smooth affine group scheme over $k$, and $A=E$ an elliptic curve over $k$.
As in (\ref{sssec:LG P component}), we denote the classical truncation of the elliptic loop space by
\[
  \frL_E(G) := \Loop_E^G(*)_{\cl} \simeq
  \frBun_G^\diamond(\check{E}).
\]
As in Subsect.~\ref{ssec:gms Ga Gm}, we denote by $\frM_{E}(G)$ the good moduli space of $\frL_E(G)$.
In this subsection we will show that $\frM_{E}(G)$ exists, and moreover that \eqref{eq:MAG vs MAT} is also an isomorphism in the elliptic case:

\begin{theorem}\label{thm:gms elliptic}\leavevmode

  (i) The classical $1$-Artin stack $\frL_E(G)$ admits a good moduli space $\pi_{E,G} : \frL_E(G) \to \frM_E(G)$, with $\frM_E(G)$ a separated algebraic space of finite type over $k$.

  (ii)
  If $T$ denotes a maximal torus of $G$, and $W$ the Weyl group, then the canonical map \eqref{eq:MAG vs MAT} is an isomorphism
  \[
    \frM_E(G)^{\circ} \simeq
    \frM_E(T)\modmod W \simeq
    \bigl(E \otimes_{\bZ} \mathbf{X}_*(T)\bigr) \modmod W.
  \]
\end{theorem}

\subsubsection{Proof of Theorem~\ref{thm:gms elliptic}(i)}

Claim~(i) follows from Theorem~\ref{thm:frBunGdiamond} using the
identification $\frL_E(G)\simeq\frBun_G^\diamond(\check{E})$ of
Theorem~\ref{thm:loop BG bun spectral}.
\qed

\subsubsection{Proof of Theorem~\ref{thm:gms elliptic}(ii)}

Suppose first that $G$ is connected.
In characteristic zero, Proposition~\ref{prop:frBunGdiamond vs ss tors top
type} identifies the locus of trivial topological type in $\frL_E(G)$ with $\frBun_G^{\ss,0}(\check{E})$.
By (\ref{sssec:Z^ss gms char 0}), the latter admits a projective good moduli space $M$.
By its universal property (see \cite[Thm.~4.16(vi)]{Alper}) and \cite[Thm.~1.2, Rem.~1.6]{Fratila}, extension of structure groups along $T\subset G$ induces an isomorphism $\frM_E(T)\modmod W \simeq M$.
This space is connected, as a quotient of $\frM_E(T)\simeq E\otimes_{\bZ}\mathbf{X}_*(T)$, and is open and closed in $\frM_E(G)$.
Since it contains the trivial loop, it is $\frM_E(G)^\circ$.

In positive characteristic, $G$ is a torus by Nagata (cf. \cite[Thm.~1.2]{HallRydh}).
Thus $G=T$ and $W=1$, while $\frM_E(G) \simeq E\otimes_{\bZ}\mathbf{X}_*(T)$ is connected, whence the claim.

For general $G$, write $\Gamma := G/G^\circ$ and $W^\circ := N_{G^\circ}(T)/T$.
By (i) and \cite[Lem.~4.14, Thm.~4.16(i--iii)]{Alper}, the good moduli spaces $\frL_E(G) \to \frM_E(G)$ and $\frL_E(G^\circ) \to \frM_E(G^\circ)$ restrict to separated good moduli spaces for the trivial loop components $\frL_E(G)^\circ$ and $\frL_E(G^\circ)^\circ$.
Thus (\ref{sssec:LG circ G^o}) and Lemma~\ref{lem:gms descent} yield
\[
  \frM_E(G)^\circ
  \simeq \frM_E(G^\circ)^\circ\modmod\Gamma.
\]
The isomorphism constructed above in the connected case,
\[
  \frM_E(T)\modmod W^\circ
  \simeq \frM_E(G^\circ)^\circ
\]
is equivariant with respect to $\Gamma\simeq W/W^\circ$, since extension of structure group along $T\subset G^\circ$ commutes with conjugation by $N_G(T)$.

Applying Lemma~\ref{lem:gms descent} again, we obtain
\[
  \frM_E(G)^\circ
  \simeq
  \bigl(\frM_E(T)\modmod W^\circ\bigr)\modmod\Gamma
  \simeq \frM_E(T)\modmod W.
\]
By construction, this isomorphism is inverse to the canonical map
\eqref{eq:MAG vs MAT}.
The last isomorphism follows from $\frM_E(T)\simeq E\otimes_{\bZ}\mathbf{X}_*(T)$.
\qed

\subsection{Appendix: classical stacks of $G$-bundles}\label{ssec:BunG}

\subsubsection{}

Let $k$ be a commutative ring, $C$ a smooth proper curve over $k$, and $G$ a smooth affine group scheme over $C$.
We write $\frBun_G(C)$ for the classical Artin stack of $G$-bundles over $C$; recall that it is locally of finite type over $k$ with affine diagonal.

\subsubsection{}\label{sssec:BunG diamond}

For $G$ defined over $k$, we denote by $\frBun_G^{\diamond}(C)$ the substack whose $S$-points are $G$-bundles $F$ on $C \times S$ such that for every geometric point $s$ of $S$ and every finite-dimensional $V \in \Rep(G_s)^\heartsuit$, the vector bundle $F_s \times^{G_s} V$ is slope-semistable of degree zero.
In this subsection, we will prove the following result:

\begin{theorem}\label{thm:frBunGdiamond}
	Let $k$ be an algebraically closed field, $E$ an elliptic curve, and $G$ a smooth affine linearly reductive group scheme over $k$.
	Then we have:

	{\rm (i)} $\frBun_G^{\diamond}(E)$ is of finite type over $k$.

 	{\rm (ii)} $\frBun_G^{\diamond}(E)$ admits a good moduli space, which is separated and of finite type over $k$.
\end{theorem}

\subsubsection{}

We will also show:

\begin{proposition}\label{prop:frBunGdiamond vs ss tors top type}
  Let $k$ be an algebraically closed field, $E$ an elliptic curve, and $G$ a connected reductive smooth affine algebraic group over $k$.
  Then $\frBun_G^{\diamond}(E)$ is the substack of $\frBun_G(E)$ consisting of $G$-bundles which are semistable and of torsion topological type.
\end{proposition}

We recall that, for $G$ a smooth affine group scheme over $C$ with connected reductive fibres, a $G$-bundle $F$ on $C$ is \emph{semistable} if for every parabolic subgroup scheme $P$ of the group scheme $\uAut_C^G(F)$ of $G$-equivariant automorphisms, $\operatorname{Lie}(P)$ is of degree $\le 0$ (following Behrend, see e.g. \cite[\S 3.2]{ORT}).
In characteristic zero, this is equivalent to slope-semistability of $\operatorname{Ad}_G(F) := F \times^{G} \operatorname{Lie}(G)$ (see \cite[Lem.~3.8]{ORT}).
We write $\frBun_G^{\ss}(C)$ for the substack of semistable $G$-bundles.

\begin{lemma}\label{lem:adjoint semistability}
  Let $H$ be a smooth affine group scheme over $C$ with connected reductive fibres, and $F$ an $H$-bundle.
  If $\operatorname{Ad}_H(F)$ is slope-semistable of degree zero, then $F$ is semistable.
\end{lemma}
\begin{proof}
  For every parabolic subgroup scheme $P\subseteq\uAut_C^H(F)$, its Lie algebra is a vector subbundle of $\operatorname{Ad}_H(F)$, hence has degree $\le 0$.
\end{proof}

For $G$ a connected reductive algebraic group over $k$ (assumed algebraically closed as above), the connected components of the stack $\frBun_G(C)$ are in bijection with $\pi_1(G)$, the quotient of the group of cocharacters $X_*(T)$ of a maximal torus $T \subset G$ by the lattice of coroots of $G$; see \cite[\S\S 1.4, 2.1]{Fratila}.
The \emph{topological type} $\delta(F) \in \pi_1(G)$ of a $G$-bundle $F$ is the index of the component containing $F$.

\subsubsection{}\label{sssec:frBunGdiamond open}

We first note that $\frBun_G^{\diamond}(E)$ is an \emph{open} substack of $\frBun_G(E)$:

\begin{proposition}\label{prop:frBunGdiamond open}
  Let $k$ be an ordinary commutative ring, $E$ an elliptic curve over $S=\spc(k)$, and $G$ an affine fppf linearly reductive group scheme over $S$.
  Then $\frBun_G^\diamond(E)$ is an open substack of $\frBun_G(E)$.
\end{proposition}

\begin{lemma}\label{lem:diamond one rep}
  Let $k$ be an ordinary commutative ring, $E$ an elliptic curve over $k$, and $G$ a linearly reductive group scheme over $k$ admitting a faithful finite projective representation $V_0$.
  For a $G$-bundle $P$ on $E$ the following are equivalent:

  (i) $P \times^G V_0$ is fibrewise semistable of degree zero;

  (ii) $P \times^G V$ is fibrewise semistable of degree zero for every finite projective $V \in \Rep(G)^\heartsuit$;

  (iii) $P \in \frBun_G^\diamond(E)$.
\end{lemma}
\begin{proof}
  Each condition is fibrewise and stable under base change, so we may assume that $k$ is an algebraically closed field; (iii) $\Rightarrow$ (ii) $\Rightarrow$ (i) is then trivial.
  For (i) $\Rightarrow$ (iii), note that the conditions ``semistable'' and ``degree zero'' are both closed under finite direct sums and direct summands, as well as under tensor products and duals (as follows from the tensor equivalence of Theorem~\ref{thm:HA vs QCoh checkA}).
  Since $V_0$ is faithful, every finite-dimensional $V \in \Rep(G)^\heartsuit$ is a direct summand of a polynomial expression $W$ in $V_0$ and $V_0^\vee$ (see the proof of Lemma~\ref{lem:BG pluperfect lin red}).
  Then (i) implies that $P \times^G W$ is semistable of degree zero (since $V \mapsto P \times^G V$ is a tensor functor), and hence so is its direct summand $P \times^G V$.
\end{proof}

\begin{proof}[Proof of Proposition~\ref{prop:frBunGdiamond open}]
  Openness is Nisnevich-local on $S$, so by \cite[Cor.~6.2(2)]{AHRLocalStructure} we may assume that $G$ is a closed subgroup of $\GL(V_0)$ for a finite projective representation $V_0$.
  By Lemma~\ref{lem:diamond one rep}, $\frBun_G^\diamond(E)$ is the preimage of $\frBun_{\GL(V_0)}^{\ss,0}(E)$ under the ``extension of structure group'' map $\frBun_G(E) \to \frBun_{\GL(V_0)}(E)$.
  Since $\frBun_{\GL(V_0)}^{\ss,0}(E) \subseteq \frBun_{\GL(V_0)}(E)$ is open, so is its preimage.
\end{proof}

\subsubsection{}\label{sssec:BunG decomposition}

We begin working towards the proof of Theorem~\ref{thm:frBunGdiamond}; unless otherwise stated, all notation is as in that statement.

Write $G^\circ$ for the identity component, and consider the finite group $\Gamma := G/G^\circ$.
Since $G$ is linearly reductive, $\Gamma$ is linearly reductive.
Extension of structure group defines a morphism
\begin{equation}\label{eq:BunG to BunGamma}
  \frBun_G(E) \to \frBun_\Gamma(E),
  \qquad F \mapsto F/G^\circ.
\end{equation}
For a $\Gamma$-torsor $P$ on $E$, consider the fibre
\[
  \UU_P := \frBun_G^\diamond(E)
  \times_{\frBun_\Gamma(E)} \{P\},
\]
whose points are pairs $(F,\alpha)$ with $F \in \frBun_G^\diamond(E)$ and $\alpha : F/G^\circ \simeq P$ an isomorphism of $\Gamma$-bundles.
It is equipped with an action of the group $H_P$ of $\Gamma$-equivariant automorphisms of $P$ over $E$, by $h\cdot(F,\alpha)=(F,h\circ\alpha)$, and we have the canonical decomposition
\[
  \frBun_G^\diamond(E)
  \simeq \coprod_{[P]} \UU_P/H_P,
\]
where $[P] \in H^1_{\text{\'et}}(E, \Gamma)$ ranges over isomorphism classes of $\Gamma$-bundles on $E$ (see \cite[Prop.~3.1]{ORT}).
Note that $H_P$ is finite and linearly reductive, since it can be regarded as a subgroup of $\Gamma$.

\subsubsection{}\label{sssec:U_P0 = BunGcirc}

For the trivial $\Gamma$-torsor $P_0 := E \times \Gamma$, we have
\[
  H_{P_0} = \Gamma, \qquad
  \UU_{P_0} \simeq \frBun_{G^\circ}^\diamond(E).
\]
Indeed, for $F \in \frBun_G(E)$ a trivialization of $F/G^\circ$ is equivalent to a reduction of structure group to a $G^\circ$-bundle $Q \in \frBun_{G^\circ}(E)$, and $F \in \frBun_G^\diamond(E)$ if and only if $Q \in \frBun_{G^\circ}^\diamond(E)$: one direction follows from the identity
\[
  F \times^G V \simeq Q \times^{G^\circ} \operatorname{Res}_{G^\circ}^G V
\]
for all finite-dimensional $V \in \Rep(G)^\heartsuit$, while conversely every finite-dimensional $G^\circ$-representation $W$ is a direct summand of $\operatorname{Res}_{G^\circ}^G\operatorname{Ind}_{G^\circ}^G W$ (which is finite-dimensional since $G/G^\circ$ is finite).

\subsubsection{Proof of Proposition~\ref{prop:frBunGdiamond vs ss tors top type}}
\label{sssec:proof of frBunGdiamond vs ss tors top type}

We assume $G$ connected reductive as in the statement.
Then the character--cocharacter pairing induces a perfect pairing
\[
  X^*(G) \times \bigl(\pi_1(G)/\pi_1(G)_{\mathrm{tors}}\bigr)
  \to \mathbf Z,
\]
see e.g. \cite[\S2.1.2]{Schieder} with $P=G$.
For every character $\chi$ of $G$, with corresponding one-dimensional representation $L_\chi$, we have for every $F \in \frBun_G(E)$,
\[
  \deg(F \times^G L_\chi) = \langle \chi, \delta(F) \rangle
\]
by \cite[Prop.~3.2.5(a)]{Schieder}.
Thus $\delta(F)$ is torsion if and only if every associated character line bundle has degree zero.

Now for a $G$-bundle $F$, the above shows that the condition $F \in \frBun_G^\diamond(E)$ implies $F$ has torsion topological type $\delta(F)$.
It also implies that $\operatorname{Ad}_G(F)$ is slope-semistable of degree zero, so $F$ is semistable by Lemma~\ref{lem:adjoint semistability}.

Conversely, let $F$ be a semistable $G$-bundle with torsion $\delta(F)$.
Since the degree of a vector bundle is defined via its determinant, the formation of which commutes with associated bundles, the first paragraph above gives
\[
  \deg(F \times^G V)
  = \deg(F \times^G \det(V))
  = 0
\]
for all finite-dimensional $V \in \Rep(G)^\heartsuit$.
Since $Z(G)$ is diagonalizable, we have a weight decomposition
\[
  V = \bigoplus_{\chi \in X^*(Z(G))} V_\chi,
\]
whose summands are $G$-invariant because $Z(G)$ is central.
By \cite[Cor.~2.10]{Fratila}, each $F \times^G V_\chi$ is slope-semistable.
Each summand has degree zero by the identity above, so their direct sum $F \times^G V$ is slope-semistable of degree zero.
Thus $F \in \frBun_G^\diamond(E)$.
\qed

\subsubsection{}

Fix a $\Gamma$-torsor $P$ on $E$, and choose a pair $(F_0,\alpha_0) \in \UU_P$ (assuming $\UU_P$ nonempty).
We denote by $G_P$ the group scheme over $E$,
\begin{equation}\label{eq:GP automorphisms}
  G_P := \ker\bigl(\uAut_E^G(F_0) \to \uAut_E^\Gamma(P)\bigr).
\end{equation}
Since it is \'etale-locally on $E$ isomorphic to $G^\circ$, it is smooth affine with connected reductive fibres.
Moreover, the fibre of \eqref{eq:BunG to BunGamma} identifies with the classical stack of $G_P$-bundles (see \cite[Lem.~3.3]{ORT}):
\begin{equation}\label{eq:BunGP fibre}
  \frBun_G(E) \times_{\frBun_\Gamma(E)} \{P\}
  \simeq \frBun_{G_P}(E),
\end{equation}
where a pair $(F,\alpha)$ is sent to the $G_P$-torsor of $G$-equivariant isomorphisms $F_0 \to F$ that intertwine $\alpha_0$ and $\alpha$ (with $G_P$ acting by precomposition).
By Proposition~\ref{prop:frBunGdiamond open}, this identifies $\UU_P$ with an open substack of $\frBun_{G_P}(E)$.

\subsubsection{}\label{sssec:G^circ reduction}

Throughout the following, we fix a connected finite \'etale cover $f:E'\twoheadrightarrow E$ and a trivialization $f^*P\simeq E'\times\Gamma$ (e.g. take $E'$ to be a connected component of $P$); note that $E'$ is of genus $1$ by Riemann--Hurwitz.
For $(F,\alpha)$ in the fibre \eqref{eq:BunGP fibre}, the inverse image of the identity section under the map
\[
  f^*F \to f^*(F/G^\circ)
  \simeq f^*P \simeq E'\times\Gamma
\]
defines a $G^\circ$-reduction $F'$ of $f^*F$ over $E'$.
That is, $F'$ is a $G^\circ$-torsor on $E'$ with
\begin{equation}\label{eq:component torsor reduction}
  F'\times^{G^\circ}G\simeq f^*F,
  \qquad
  f^*(F\times^G V)\simeq F'\times^{G^\circ}\operatorname{Res}_{G^\circ}^G V
\end{equation}
for every finite-dimensional $V\in\Rep(G)^\heartsuit$.

The group scheme of automorphisms of $F' \in \frBun_{G^\circ}(E')$ is given by
\[
  \uAut_{E'}^{G^\circ}(F')
  \simeq \ker\bigl(\uAut_{E'}^G(f^*F)\to\uAut_{E'}^\Gamma(f^*P)\bigr).
\]
Indeed, a $G$-equivariant automorphism of $f^*F$ inducing the identity on $f^*P$ preserves the inverse image $F'$ of the identity section, and hence restricts to a $G^\circ$-equivariant automorphism of $F'$.
Conversely, every such automorphism of $F'$ extends uniquely to $f^*F\simeq F'\times^{G^\circ}G$, inducing the identity on $f^*P$.

\subsubsection{}\label{sssec:r map}

Applying (\ref{sssec:G^circ reduction}) to the fixed pair $(F_0,\alpha_0)$ yields a $G^\circ$-reduction $F_0'$ of $f^*F_0$ over $E'$, with $\uAut_{E'}^{G^\circ}(F_0') \simeq f^*G_P$.
In particular, there is a canonical isomorphism of stacks
\[
  \frBun_{f^*G_P}(E') \xrightarrow{\sim}
  \frBun_{G^\circ}(E'),
  \qquad Q \mapsto Q \times^{f^*G_P} F'_0.
\]
We consider the composite
\begin{equation}\label{eq:BunGP pullback twist}
  r_{F_0'}:
  \frBun_{G_P}(E) \to
  \frBun_{f^*G_P}(E') \xrightarrow{\sim}
  \frBun_{G^\circ}(E'),
\end{equation}
sending $Q\mapsto f^*Q \mapsto f^*Q\times^{f^*G_P}F_0'$.
Given a pair $(F,\alpha)$ in \eqref{eq:BunGP fibre}, $r_{F_0'}$ sends the corresponding $Q \in \frBun_{G_P}(E)$ to the reduction $F'$.
By Behrend, \eqref{eq:BunGP pullback twist} is affine of finite presentation (see \cite[Lem.~3.12]{ORT}).

\subsubsection{}

Let $(F,\alpha)$ be a pair in \eqref{eq:BunGP fibre}, $Q \in \frBun_{G_P}(E)$ the corresponding $G_P$-bundle, and $F'=r_{F_0'}(Q)$ the reduction.
We have
\begin{equation}\label{eq:component torsor adjoint}
  \operatorname{Ad}_G(F) \simeq \operatorname{Ad}_{G_P}(Q),
  \qquad
  \operatorname{Ad}_{G^\circ}(F')\simeq f^*\operatorname{Ad}_{G_P}(Q).
\end{equation}
The left-hand isomorphism follows by identifying $\uAut_E^{G_P}(Q)$ with the kernel of $\uAut_E^G(F)\to\uAut_E^\Gamma(P)$, and taking Lie algebras (since $\uAut_E^\Gamma(P)$ is \'etale over $E$).
The right-hand isomorphism follows from \eqref{eq:component torsor reduction} applied to $V=\operatorname{Lie}(G)$, since $\operatorname{Lie}(G^\circ)=\operatorname{Lie}(G)$.

In particular, if $(F,\alpha)\in\UU_P$, then $\operatorname{Ad}_{G_P}(Q)$ is slope-semistable of degree zero, and hence by Lemma~\ref{lem:adjoint semistability} $Q$ is semistable.
We thus have the inclusion
\[
  \UU_P\subseteq\frBun_{G_P}^{\ss}(E).
\]

\subsubsection{}\label{sssec:U_P = r inv diamond}

For the finite \'etale cover $f:E'\twoheadrightarrow E$, a vector bundle $W$ on $E$ is slope-semistable of degree zero if and only if $f^*W$ is so (see \cite[Lem.~2.1(ii--iii)]{Weissmann}).

For a pair $(F,\alpha)$ in \eqref{eq:BunGP fibre}, the argument of (\ref{sssec:U_P0 = BunGcirc}) shows in view of the identity \eqref{eq:component torsor reduction} that
\[
  F\in\frBun_G^\diamond(E)
  \quad\Longleftrightarrow\quad
  F'\in\frBun_{G^\circ}^\diamond(E').
\]
Thus $\UU_P$ is identified under \eqref{eq:BunGP fibre} with the substack of $Q \in \frBun_{G_P}(E)$ for which $r_{F'_0}(Q) \in \frBun_{G^\circ}^\diamond(E')$:
\begin{equation}\label{eq:BunGP diamond inverse image}
  \UU_P = r_{F_0'}^{-1}\bigl(\frBun_{G^\circ}^\diamond(E')\bigr).
\end{equation}

\subsubsection{}

Let $Q \in \frBun_{G_P}(E)$ be arbitrary.
We claim that $Q$ is semistable if and only if $F'=r_{F_0'}(Q)$ is semistable:
\begin{itemize}
	\item
  In characteristic zero, this follows from \eqref{eq:component torsor adjoint}, (\ref{sssec:U_P = r inv diamond}), and the criterion of \cite[Lem.~3.8]{ORT}.

  \item
  In positive characteristic, $G^\circ$ is a torus by Nagata's theorem (cf. \cite[Thm.~1.2]{HallRydh}).
  Then $\operatorname{Ad}_{G^\circ}(F')$ is trivial, so \eqref{eq:component torsor adjoint} shows that $\operatorname{Ad}_{G_P}(Q)$ is slope-semistable of degree zero.
  Thus both $Q$ and $F'$ are \emph{always} semistable by Lemma~\ref{lem:adjoint semistability}.
\end{itemize}

Combining this observation with \eqref{eq:BunGP diamond inverse image} and Proposition~\ref{prop:frBunGdiamond vs ss tors top type}, we find that for $Q\in\frBun_{G_P}^{\ss}(E)$,
\[
  Q\in\UU_P
  \quad\Longleftrightarrow\quad
  \delta\bigl(r_{F_0'}(Q)\bigr)\in\pi_1(G^\circ)_{\mathrm{tors}}.
\]
Since topological type is locally constant in families, $\UU_P$ is an open and closed substack of $\frBun_{G_P}^{\ss}(E)$.

\subsubsection{}\label{sssec:Z^ss gms char 0}

Let $H$ be a smooth affine group scheme over $E$ with connected and linearly reductive geometric fibres.
We next claim that for every connected component $Z$ of $\frBun_H(E)$, the semistable locus $Z^{\ss} := Z \cap \frBun_H^{\ss}(E)$ is of finite type and admits a separated good moduli space.
The same then holds for every open and closed substack $U \subseteq Z^{\ss}$, by \cite[Lem.~4.14]{Alper}.

In characteristic zero, the finite type claim follows from the proof of \cite[Lem.~8.2]{AHLH}, and the proof of \cite[Lem.~3.14]{ORT} shows that $Z^{\ss}$ admits a projective good moduli space.

Suppose now that $k$ is of characteristic $p>0$.
Then $H$ is a torus by Nagata's theorem (cf. \cite[Thm.~1.2]{HallRydh}), and by \cite[Exp.~X, Thm.~5.16]{SGA3} there exists a connected finite \'etale cover $f: E' \twoheadrightarrow E$ along with an isomorphism $f^*H\simeq T\times E'$ over $E'$, where $T = \GGm^r$.
By Lemma~\ref{lem:adjoint semistability}, every $H$-bundle is semistable, since its adjoint bundle is $\operatorname{Lie}(H)$ and becomes trivial on $E'$.
Thus $Z^{\ss}=Z$.

The pullback map $f^*:\frBun_H(E)\to\frBun_T(E') \simeq \PIC(E')^{r}$ is affine of finite presentation by Behrend (see \cite[Lem.~3.12]{ORT}).
Since $Z$ is connected, its image lies in a component $$Z' := \PIC_{d_1}(E') \times \cdots \times \PIC_{d_r}(E') \subset \PIC(E')^{r}$$ for some $(d_1,\ldots,d_r) \in \bZ^{r}$, where $\PIC_d(E')$ is the stack of degree $d$ line bundles.
This is a $T$-gerbe over $M'$, a product of Picard schemes (which can be identified with $(E')^{r}$); in particular, $Z'$ is of finite type and admits $M'$ as a good moduli space (since $T$ is linearly reductive, see \cite[Prop.~4.7(ii)]{Alper}).
As $Z$ is closed in $(f^*)^{-1}(Z')$, it is of finite type and affine over $Z'$; thus by \cite[Lem.~4.14]{Alper} it admits a good moduli space which is affine over $M'$, hence separated over $k$.

\subsubsection{Proof of Theorem~\ref{thm:frBunGdiamond}(i)}\label{sssec:proof frBunGdiamond i}

We prove that $\frBun_G^\diamond(E)$ is of finite type.
In view of the decomposition of (\ref{sssec:BunG decomposition}), and since $H^1_{\text{\'et}}(E,\Gamma)$ is finite (see e.g. \cite[Prop.~3.1]{ORT}), it will suffice to show that each $\UU_P$ is of finite type.
Fix $P$ with $\UU_P$ nonempty and choose a connected finite \'etale cover $f : E' \twoheadrightarrow E$ as in (\ref{sssec:G^circ reduction}).
By the description of $\UU_P$ given in \eqref{eq:BunGP diamond inverse image}, together with the fact that the map $r_{F'_0}$ is of finite presentation (\ref{sssec:r map}), we reduce to checking that $\frBun_{G^\circ}^\diamond(E')$ is of finite type.

For this, note that the abelian group $\pi_1(G^\circ)$ is finitely generated, so its torsion subgroup is finite.
Hence by Proposition~\ref{prop:frBunGdiamond vs ss tors top type} and (\ref{sssec:Z^ss gms char 0}), $\frBun_{G^\circ}^\diamond(E')$ is a finite union of semistable loci of connected components, each of finite type.
The claim follows.
\qed

\subsubsection{Proof of Theorem~\ref{thm:frBunGdiamond}(ii)}\label{sssec:proof frBunGdiamond ii}

We show that $\frBun_G^\diamond(E)$ admits a good moduli space which is separated of finite type.
As in (\ref{sssec:proof frBunGdiamond i}) it will suffice to show that the same holds for each (nonempty) piece $\UU_P/H_P$.
Since $H_P$ is finite and linearly reductive, Lemma~\ref{lem:gms descent} reduces us to showing that $\UU_P$ admits a separated good moduli space.

Since $\UU_P$ is quasi-compact by (i), and the connected components of the locally noetherian stack $\frBun_{G_P}(E)$ are open, $\UU_P$ meets only finitely many such components $Z_1,\ldots,Z_n$.
We have a decomposition
\[
  \UU_P=\coprod_{j=1}^n\bigl(\UU_P\cap Z_j^{\ss}\bigr),
\]
where each intersection is open and closed in $Z_j^{\ss}$ since $\UU_P$ is open and closed in $\frBun_{G_P}^{\ss}(E)$.
By (\ref{sssec:Z^ss gms char 0}), each intersection admits a separated good moduli space $M_{P,j}$.
Taking their finite disjoint union gives a good moduli space $\UU_P\to M_P := \coprod_j M_{P,j}$ with $M_P$ separated.
\qed

\section{Properties of loop spaces}\label{sec:loop properties}

\subsection{Artinness}

\begin{corollary}\label{cor:loop Artin}
  Let $\YY$ be an Artin stack.
  Then $\Loop_A(\YY)$ is Artin in each of the following cases:
  \begin{enumerate}
    \item $A = \GGa$ and $k$ is a $\bQ$-algebra;
    \item $A = \GGm$;
    \item $A$ is an abelian variety over $k$, and $\YY$ is pluperfect fppf-locally on $k$.
  \end{enumerate}
\end{corollary}
\begin{proof}
  In cases (i) and (ii) this follows from Theorem~\ref{thm:LGa LGm}.
  Case (iii) is Theorem~\ref{open in map pluperfect}.
\end{proof}

\begin{corollary}\label{cor:loop Artin VarG}
  Let $A = \GGa$ (with $k$ a $\bQ$-algebra), $A=\GGm$, or $A$ an abelian variety over $k$.
  Then for every affine fppf group scheme $G$ over $k$ and every $X \in \Var^G$, the $G$-equivariant loop space $\Loop_A^G(X)$ is Artin.
\end{corollary}
\begin{proof}
	By definition, $X/G$ is perfect and pluperfect fppf-locally on $\spc(k)$, so we may apply Corollary~\ref{cor:loop Artin}.
\end{proof}

\subsection{Representability and properness}

Let $f : \XX \to \YY$ be a morphism of Artin stacks.
In this subsection we establish some global properties of the map $\Loop_A(f) : \Loop_A(\XX) \to \Loop_A(\YY)$ induced by $f$.

\subsubsection{}

We first observe that the classical truncation of $\Loop_A(\XX)$ only depends on the classical truncation of $\XX$:

\begin{proposition}\label{prop:Loop cl}
Let $\XX$ be an Artin stack, and consider the canonical morphism $\Loop_A(\XX_\cl) \to \Loop_A(\XX)$.
\begin{enumerate}
	\item For $A=\GGa$ (with $k$ a $\bQ$-algebra) or $\GGm$, the above is an isomorphism on classical truncations.
	\item For $A$ an abelian variety, the same holds for $\XX=X/G$ with $X \in \Var^G$ and $G$ an affine fppf group scheme over $k$.
\end{enumerate}
\end{proposition}
\begin{proof}
In the first case this follows from the descriptions given in Theorem~\ref{thm:LGa LGm}.
Suppose $A$ is an abelian variety.
Since the claim is fppf-local on $\spc(k)$, we may assume that $\XX = X/G$ is pluperfect (since $X \in \Var^G$).
Note first that the analogous claim for the morphism
\[
  \MAP(\check{A}, \XX_\cl) \to \MAP(\check{A}, \XX)
\]
holds: for any classical affine $T$ over $k$, every morphism $\check{A}_T \to \XX$ factors uniquely through $i : \XX_\cl \hookrightarrow \XX$ ($\check{A}_T$ is classical because $\check{A}$ is flat over $k$).

Since $\XX$ is pluperfect, $\QCOH(\XX)$ admits a finite set $\{F_j\}$ of perfect tensor generators.
Since $i : \XX_\cl \hookrightarrow \XX$ is affine, $\{i^*F_j\}$ is still a set of perfect tensor generators.
Then the proof of Theorem~\ref{open in map pluperfect} shows that a morphism $\sigma \in \MAP(\check{A}, \XX_\cl)(T)$ is a loop if and only if $\sigma^* (i^* F_j) \in \INDPERF^{\mathrm{hmg}}(\check{A}_T)$ for all $j$.
Applying this observation again to $i\circ \sigma : \check{A}_T \to \XX$ shows that this is the same as the condition that its image in $\MAP(\check{A}, \XX)$ is a loop.
\end{proof}

\subsubsection{}

We next relate the fibres of the maps $\Loop_A(f)$ and $f$.

\begin{lemma}\label{lem:loop fibre}
  Let $f : \XX \to \YY$ be a representable morphism of Artin stacks.
  Given an affine scheme $S$ and a loop $\gamma : S \to \Loop_A(\YY)$ with $y = \ev_0(\gamma) : S \to \YY$ (\ref{sssec:ev_0}), consider the canonical map
  \[ \ev_0 : \Loop_A(\XX) \times_{\Loop_A(\YY)} \{\gamma\} \to \XX \times_{\YY} \{y\} =: \XX_y. \]
  \begin{enumerate}
  \item $A = \GGa$ (with $k$ a $\mathbf{Q}$-algebra): the map above is a closed immersion.
  \item $A = \GGm$: If $f$ is separated (resp. quasi-separated), then the above map is a closed immersion (resp. quasi-compact and locally of finite type, and a monomorphism on classical truncations).
  \item $A = E$ (an elliptic curve over $k$): If $\XX = X/G$ and $\YY = Y/G$ for $X, Y \in \Var^G$, and $f$ is quasi-projective on classical truncations, then the above map is a closed immersion.
  \end{enumerate}
\end{lemma}

\subsubsection{Proof for $A=\GGa$}

By Theorem~\ref{thm:LGa LGm}, $\Loop_{\GGa}(\ZZ) \simeq \mathbf{V}_{\ZZ}(\LL_{\ZZ}[1])$ is the stack of cosections of $\LL_{\ZZ}[1]$ over $\ZZ$ (see \cite[App.~B]{weilres} for notation).
Thus an additive loop $\gamma : S \to \Loop_{\GGa}(\YY)$ amounts to an $S$-point $y$ of $\YY$ together with a cosection $v : y^*\LL_{\YY}[1] \to \OO_S$.
Factoring $\Loop_A(f)$ as $\mathbf{V}_{\XX}(\LL_{\XX}[1]) \to \XX \times_{\YY} \mathbf{V}_{\YY}(\LL_{\YY}[1]) \to \mathbf{V}_{\YY}(\LL_{\YY}[1])$, we have:
\[ \Loop_A(\XX) \times_{\Loop_A(\YY)} \{\gamma\} \simeq \mathbf{V}_{\XX_y}(\LL_{\XX}[1]|_{\XX_y}) \times_{\mathbf{V}_{\XX_y}(y^*\LL_{\YY}[1]|_{\XX_y}),\, v} \XX_y. \]
From the exact triangle $\LL_{\XX/\YY} \to f^*\LL_{\YY}[1] \to \LL_{\XX}[1]$ restricted to $\XX_y$, we have by \cite[Rem.~B.3]{weilres},
\[\mathbf{V}_{\XX_y}(\LL_{\XX}[1]|_{\XX_y}) \simeq \mathbf{V}_{\XX_y}(y^*\LL_{\YY}[1]|_{\XX_y}) \times_{\mathbf{V}_{\XX_y}(\LL_{\XX/\YY}|_{\XX_y}),\, 0} \XX_y.\]
Combining the two isomorphisms above, we identify the fibre as
\[
  \Loop_A(\XX) \times_{\Loop_A(\YY)} \{\gamma\} \simeq \XX_y \times_{w, \mathbf{V}_{\XX_y}(\LL_{\XX/\YY}|_{\XX_y}), 0} \XX_y,
\]
i.e., the zero locus of the cosection
\[\LL_{\XX/\YY}|_{\XX_y} \xrightarrow{\partial} f^*\LL_{\YY}[1] \xrightarrow{v|_{\XX_y}} \OO_{\XX_y}.\]
Since $\LL_{\XX_y/S}$ is connective (as $f$ is representable), the map to $\XX_y$ is a closed immersion.

\subsubsection{Proof for $A = \GGm$}
By Theorem~\ref{thm:LGa LGm}, $\Loop_{\GGm}(\ZZ) \simeq \ZZ \times_{\ZZ \times \ZZ} \ZZ$ naturally in $\ZZ$.
Thus a multiplicative loop $\gamma$ in $\YY$ amounts to a point $y$ of $\YY$ together with a homotopy $\alpha : y \simeq y$; by base change, this determines an automorphism
\[f^*\alpha : \XX_y \simeq \XX_y.\]
By distributivity of fibres over limits, the fibre of $\Loop_A(f)$ at $\gamma$ is
\[
  \Loop_A(\XX) \times_{\Loop_A(\YY)} \{\gamma\}
  \simeq \XX_y \times_{\XX_y \times_S \XX_y} \XX_y,
\]
where the fibred product on the right is formed with respect to the diagonal $\Delta_{\XX_y/S}$ and the graph of $f^*\alpha$.
Hence the claim follows by base change from the analogous properties of $\Delta_{\XX_y/S}$.

\subsubsection{Proof for $A = E$}

\begin{lemma}\label{lem:Loop_E P(V)}
	Let $A = E$ be an elliptic curve over $k$.
	Let $\YY$ be a global quotient stack and $g : \mathcal{P} := \mathbf{P}_{\YY}(V) \to \YY$ for a vector bundle $V$ on $\YY$.
	For a classical affine scheme $S$ and a morphism $\sigma : \check{E}_S \to \mathcal{P}$, denote by $L_\sigma$ the line bundle and $u_\sigma : \sigma^*g^*V \twoheadrightarrow L_\sigma$ the surjection classified by $\sigma$.
	Then the following conditions are equivalent:
	\begin{enumerate}
	  \item\label{item:Loop_E P(V)/loop} $\sigma$ is a loop in $\mathcal{P}$, i.e., lies in $\Loop_E(\mathcal{P}) \subset \MAP(\check{E}, \mathcal{P})$.
		\item\label{item:Loop_E P(V)/hmg} $L_\sigma$ is homogeneous, and $\gamma : \check{E}_S \xrightarrow{\sigma} \mathcal{P} \to \YY$ is a loop in $\YY$.
	\end{enumerate}
\end{lemma}
\begin{proof}
  \itemref{item:Loop_E P(V)/loop} $\implies$ \itemref{item:Loop_E P(V)/hmg} is clear from the definitions.
  For \itemref{item:Loop_E P(V)/hmg} $\implies$ \itemref{item:Loop_E P(V)/loop}:
  Choose a quasi-affine morphism $\YY \to BGL_n$, and let $V_{\YY}$ be the pullback of the standard representation $V_{\mathrm{std}}$.
  The induced morphism $(g, \OO(1)) : \mathcal{P} \to \YY \times B\GGm$ is quasi-affine: its base change along $\YY \to \YY \times B\GGm$ is the complement of the zero section in $\mathbf{V}_{\YY}(V)$.
  The morphism $B(GL_n \times \GGm) \to BGL_{n+1}$ induced by the block inclusion is affine, as $GL_{n+1}/(GL_n \times \GGm)$ is the complement of the incidence divisor in $\mathbf{P}^n \times \mathbf{P}^{n,\vee}$.
  Hence the composite $\pi : \mathcal{P} \to \YY \times B\GGm \to B(GL_n \times \GGm) \to BGL_{n+1}$ is quasi-affine.
  Moreover, $\sigma^*\pi^*V_{\mathrm{std}} \simeq \gamma^*V_{\YY} \oplus L_\sigma$ is homogeneous ($\gamma^*V_{\YY}$ is homogeneous because $\gamma$ is a loop by assumption).
  Thus we may apply Lemma~\ref{lem:std suff for hmg} to conclude that $\sigma$ is a loop.
\end{proof}

We first reduce to the case where $\XX$ and $\YY$ are classical.
The fibre $\Loop_E(\XX) \times_{\Loop_E(\YY)} \{\gamma\}$ is an open substack of the Weil restriction of $\XX \times_{\YY} \check{E}_S$ along $\check{E}_S \to S$ (by Theorem~\ref{open in map pluperfect}), hence an algebraic space by \cite[Thm.~19.1.0.1]{SAG}.
It is therefore enough to show that $\ev_0$ is a closed immersion on classical truncations.
For this we may assume by Proposition~\ref{prop:Loop cl} that $\XX$, $\YY$ and $S$ are classical.

Since $\YY$ has the resolution property, $f$ factors as a quasi-compact immersion $i : \XX \hookrightarrow \mathcal{P} := \mathbf{P}_{\YY}(V)$, for a vector bundle $V$ on $\YY$, followed by the projection $g : \mathcal{P} \to \YY$.
To prove Lemma~\ref{lem:loop fibre} for $A=E$, it will suffice to show:
\begin{enumerate}[label=(\alph*)]
\item\label{item:loop fibre E/proj} The map $\ev_0 : \Loop_E(\mathcal{P}) \times_{\Loop_E(\YY)} \{\gamma\} \to \mathcal{P} \times_{\YY} \{y\}$ is a closed immersion.
\item\label{item:loop fibre E/imm} For every classical affine $T$, a loop $\sigma : T \to \Loop_E(\mathcal{P})$ factors through $\Loop_E(\XX)$ if and only if $\ev_0(\sigma) : T \to \mathcal{P}$ factors through $\XX$.
\end{enumerate}
Indeed, \ref{item:loop fibre E/imm} implies that $\ev_0 : \Loop_E(\XX) \times_{\Loop_E(\YY)} \{\gamma\} \to \XX \times_{\YY} \{y\}$ is the base change of the closed immersion in \ref{item:loop fibre E/proj} (up to classical truncation).

\subsubsection*{Proof of \ref{item:loop fibre E/proj}.}

For a classical affine $T$ over $S$, $T$-points of $\Loop_E(\mathcal{P}) \times_{\Loop_E(\YY)} \{\gamma\}$ are pairs $(L,u)$ where $L$ is a homogeneous line bundle on $\check{E}_T$ and $u : (\gamma^*V)_T \twoheadrightarrow L$ is a surjection (Lemma~\ref{lem:Loop_E P(V)}).
Replacing $E_S \to S$ by its base change $E_T \to T$, we take $T=S$ to simplify the notation.
Recall the exact equivalence of abelian categories $\Phi_{E_S/S} : \Hcat_{E_S/S}^{\heartsuit} \to \IndPerf^{\mathrm{hmg}/S}(\check{E}_S)^\heartsuit$ (Remark~\ref{rem:HA vs QCoh checkA dzbl}).
Since $\gamma$ is a loop, $\gamma^*V$ is a homogeneous vector bundle, as is the kernel of $u$.
Thus by Lemma~\ref{lem:Coh_0^r wedge finite}, $\Phi_{E_S/S}$ identifies such pairs $(L,u)$ with pairs $(F', v)$ where $F' \in \Hcat_{E_S/S}^\heartsuit$ is flat over $S$ with fibres of length one and $v : F := \Phi_{E_S/S}^{-1}(\gamma^*V) \twoheadrightarrow F'$ is surjective.

The map $\ev_0$ sends $(L,u)$ to the $S$-point of $\mathcal{P} \times_{\YY} \{y\}$ classified by the pair $(e^*L, e^*u)$ where $e : S \hookrightarrow \check{E}_S$ is the identity section.
Under $\Phi_{E_S/S}$, $e^*$ corresponds to $p_*$ for $p : E_S \to S$.
In other words, $\ev_0$ sends $(F', v) \mapsto (p_*F', p_*v : p_*(F) \twoheadrightarrow p_*(F'))$.
Since $F$ is an $\OO_D$-module for some $D \subseteq E_S$ finite over $S$, and $p_*$ determines an equivalence $\QCOH(D)^\heartsuit \simeq \MOD_{p_*\OO_{D}}(\QCOH(S))^\heartsuit$, it follows that $\ev_0$ is a monomorphism.

The $p_*\OO_D$-linear quotients of $p_*F = y^*V$ form a closed subscheme $Z \subseteq \mathcal{P} \times_{\YY} \{y\} = \mathbf{P}_S(y^*V)$.
Indeed, let $q : y^*V \twoheadrightarrow \OO(1)$ be the universal quotient with kernel $K$, and $b_1, \ldots, b_r$ generators of the finite $\OO_S$-module $p_*\OO_D$; then $q$ is $p_*\OO_D$-linear exactly on the common zero locus of the sections $q \circ b_i|_K$ of $K^\vee(1)$.
Moreover, every point in $Z$ lies in the image of $\ev_0$: a point of $Z$ is a $p_*\OO_D$-linear surjection $q : y^*V \twoheadrightarrow N$ onto a line bundle, necessarily of the form $q = p_*v$ for a quotient $v : F \twoheadrightarrow F'$ of $\OO_{D}$-modules; since $p_*F' \simeq N$ is a line bundle and $p|_D$ is affine, $F'$ is flat over $S$ with fibres of length one, so $(F',v)$ is a point of the fibre with $\ev_0(F',v) = q$.
Thus $\ev_0$ is an isomorphism onto the closed subscheme $Z$.
\qed

\subsubsection*{Proof of \ref{item:loop fibre E/imm}.}
The ``only if'' direction is clear.
For the converse, factor the quasi-compact immersion $i$ as $\XX \hookrightarrow \bar{\XX} \hookrightarrow \mathcal{P}$.

Suppose first that the $T$-point $\ev_0(\sigma)$ in $\mathcal{P}$ factors through $\bar{\XX}$.
Since $g$ is projective and $\YY$ has the resolution property, there exists a vector bundle $W$ on $\YY$ and a surjection $g^*W \otimes \OO(-m) \twoheadrightarrow \mathcal{I}_{\bar{\XX}}$ for $m \gg 0$, so that $\bar{\XX}$ is the zero locus of the induced map $\phi : g^*W \otimes \OO(-m) \to \OO_{\mathcal{P}}$.
Let $\sigma : \check{E}_T \to \mathcal{P}$ be a loop, so that $\sigma$ factors through $\bar{\XX}$ if and only if $\sigma^*\phi = 0$.
The source $(g\sigma)^*W \otimes L_\sigma^{\otimes -m}$ of $\sigma^*\phi$ is a homogeneous vector bundle (Lemma~\ref{lem:Loop_E P(V)}), so as in the proof of \ref{item:loop fibre E/proj}, $\sigma^*\phi \simeq \Phi_{E_T/T}(\psi)$ for a morphism $\psi$ of $\OO_D$-modules, where $D \subseteq E_T$ is finite over $T$, and $e^*\sigma^*\phi \simeq p_*\psi$.
Since $p_*$ is faithful on $\QCOH(D)^\heartsuit$, $\sigma^*\phi = 0$ if and only if $\ev_0(\sigma)^*\phi = e^*\sigma^*\phi = 0$, i.e., if and only if $\ev_0(\sigma)$ factors through $\bar{\XX}$.
In this case, the factorization $\bar\sigma : \check{E}_T \to \bar{\XX}$ is a loop by Lemma~\ref{lem:std suff for hmg}: for $\pi : \mathcal{P} \to BGL_{n+1}$ the quasi-affine morphism constructed in the proof of Lemma~\ref{lem:Loop_E P(V)}, $\bar\sigma^*(\pi|_{\bar{\XX}})^*V_{\mathrm{std}} \simeq \sigma^*\pi^*V_{\mathrm{std}}$ is homogeneous because $\sigma$ is a loop.

Now suppose $\ev_0(\sigma)$ factors through $\XX$, so that $\sigma$ factors through $\bar{\XX}$ by the above.
Let $Z := \bar{\XX} \setminus \XX$ with its reduced structure, and choose as above a map $\phi_Z : g^*W_Z \otimes \OO(-m_Z) \to \OO_{\mathcal{P}}$ with image $\mathcal{I}_Z$.
Then $\sigma$ factors through $\XX$ if and only if $\sigma^{-1}(Z) = \emptyset$, i.e., if and only if $\sigma^*\phi_Z$ is surjective.
As above, $\sigma^*\phi_Z \simeq \Phi_{E_T/T}(\psi_Z)$ for a morphism $\psi_Z$ of $\OO_D$-modules with $e^*\sigma^*\phi_Z \simeq p_*\psi_Z$.
Since $\Phi_{E_T/T} : \Hcat_{E_T/T}^{\heartsuit} \to \IndPerf^{\mathrm{hmg}/T}(\check{E}_T)^\heartsuit$ is an exact equivalence (Remark~\ref{rem:HA vs QCoh checkA dzbl}), $\sigma^*\phi_Z$ is surjective if and only if $\psi_Z$ is.
Since $p_*$ is exact and conservative on $\QCOH(D)^\heartsuit$, this holds if and only if $p_*\psi_Z = e^*\sigma^*\phi_Z = \ev_0(\sigma)^*\phi_Z$ is surjective, i.e., if and only if $\ev_0(\sigma)$ factors through $\XX$.
The resulting map $\check{E}_T \to \XX$ is then a loop by Lemma~\ref{lem:std suff for hmg}, as $\pi|_{\XX}$ is quasi-affine.
\qed

\begin{corollary}\label{mapsrepr}
Let $A$ be $\GGm$, $\GGa$ with $k$ a $\mathbf{Q}$-algebra, or an elliptic curve $E$ over $k$, and let $G$ be a linearly reductive group scheme over $k$.
For every morphism $f : X \to Y$ in $\Var^G$, we have:
\begin{enumerate}
	\item
	The morphism $\Loop_A^G(f) : \Loop_A^G(X) \to \Loop_A^G(Y)$ is schematic and separated of finite type.

	\item
	The induced morphism $\Loop_A^G(f)/\check{A} : \Loop_A^G(X)/\check{A} \to \Loop_A^G(Y)/\check{A}$ is representable and separated of finite type.

	\item
  If $f$ is proper, then $\Loop_A^G(f)$ and $\Loop_A^G(f)/\check{A}$ are proper.
\end{enumerate}
\end{corollary}
\begin{proof}
It is formal that representability, separatedness, properness, and finite type descend to the quotient $\Loop_A^G(f)/\check{A}$, so it will suffice to prove the claims about $\Loop_A^G(f)$.
The claims are Nisnevich-local on $\spc(k)$, so by Lemma~\ref{lem:VarG projective embedding} we may assume there is a $G_\cl$-equivariant quasi-compact immersion $i : X_\cl \hookrightarrow \mathbf{P}(W)$ for a finite projective representation $W$ of $G_\cl$.
Then $f_\cl/G_\cl : X_\cl/G_\cl \to Y_\cl/G_\cl$ is quasi-projective, since $(i, f_\cl) : X_\cl \hookrightarrow \mathbf{P}(W) \times Y_\cl \simeq \mathbf{P}_{Y_\cl}(W_{Y_\cl})$ is a $G_\cl$-equivariant quasi-compact immersion over $Y_\cl$.

Let $T$ be an affine scheme and $\gamma : T \to \Loop_A^G(Y)$ with $y := \ev_0(\gamma) : T \to Y/G$ (see (\ref{sssec:ev_0}) for notation).
By Lemma~\ref{lem:loop fibre}, the fibre of $\Loop_A^G(f)$ over $\gamma$ is a closed substack of $(X/G) \times_{Y/G} \{y\}$.
The latter is a separated scheme of finite type over $T$, and proper if $f$ is, since its classical truncation is quasi-projective over $T_\cl$.
The claim follows.
\end{proof}

\subsection{Tangent complexes and smoothness}

\begin{theorem}\label{tangent} Let $\YY$ be a stack admitting a tangent complex\footnote{%
  i.e., admitting a cotangent complex which is perfect (whose dual is $\TT_{\YY}$ by definition)
} $\TT_{\YY}$. Then $\Loop_A(\YY)$ admits a tangent complex, given by
$$\TT_{\Loop_A(\YY)} = \Theta^A_{\YY}(\TT_{\YY}),$$
where $\Theta^A_{\YY} : \PERF(\YY) \to \PERF(\Loop_A(\YY))$ is the functor defined in (\ref{sssec:Theta}).
Moreover, if $A$ is of dimension $1$ and $\TT_{\YY}$ has amplitude $[a,b]$, then $\TT_{\Loop_A(\YY)}$ has amplitude $[a-1,b]$.\footnote{%
  using homological grading conventions for Tor-amplitude
} \end{theorem}

\begin{proof} Since $\Loop_A(\YY) \hookrightarrow \MAP(\check{A}, \YY)$ is a monomorphism (\ref{sssec:loops as maps}), $\TT_{\Loop_A(\YY)}$ is the restriction of the tangent complex of $\MAP(\check{A}, \YY)$. For the mapping stack the formula $\TT_{\MAP(\check{A},\YY)} \simeq \mathrm{pr}_* \eval^* \TT_{\YY}$ is standard, see e.g.\ \cite[\S 2.1]{PTVV} or \cite[Prop.~19.1.4.1]{SAG}. Restricting to $\Loop_A(\YY)$ gives $\TT_{\Loop_A(\YY)} = \Theta^A_{\YY}(\TT_{\YY})$. It remains to show that the functor $\Theta^A_{\YY} = \pr_* \eval^*$ has amplitude $[-1, 0]$. This reduces to computing the amplitude of $\OO(\check{A})$ (using the base change and projection formulas).

By Proposition~\ref{prop:fun check A}, $\OO(\check{A}) \simeq \RA^\vee$, where $\RA \simeq \OO_S \oplus \frav[1]$ additively (\ref{sssec:funcheckA add split}). This has amplitude $[0, 1]$, so $\OO(\check{A})$ has amplitude $[-1, 0]$. \end{proof}

\begin{theorem}\label{reltangent} Let $f:\XX\rightarrow\YY$ be a morphism of stacks admitting a perfect relative cotangent complex $\LL_{\XX/\YY}$. Then $\Loop_A(f):\Loop_A(\XX)\rightarrow\Loop_A(\YY)$ has perfect relative tangent complex\footnote{%
  where $\TT_{\XX/\YY}$ is the relative tangent complex, dual to $\LL_{\XX/\YY}$
}
\[ \TT_{\Loop_A(\XX)/\Loop_A(\YY)} \simeq \Theta^A_\XX(\TT_{\XX/\YY}). \]
In particular, if $f$ is étale, then $\Loop_A(f)$ is formally étale; if $f$ is differentially smooth, then $\Loop_A(f)$ is formally differentially quasi-smooth\footnote{%
  i.e., its relative cotangent complex is perfect of amplitude $[0,1]$
}. \end{theorem}
\begin{proof}
As in the proof of Theorem~\ref{tangent}, $\TT_{\Loop_A(\XX)/\Loop_A(\YY)}$ is the restriction of $\TT_{\MAP(\check{A},\XX)/\MAP(\check{A},\YY)}$.
Since $\MAP(\check{A},\XX)$ is the Weil restriction of $\XX\times_{\YY}(\check{A}\times\MAP(\check{A},\YY))$ along $\pr : \check{A}\times\MAP(\check{A},\YY) \to \MAP(\check{A},\YY)$, \cite[Prop.~19.1.4.3]{SAG} gives $\TT_{\MAP(\check{A},\XX)/\MAP(\check{A},\YY)} \simeq \pr_*\eval^*\TT_{\XX/\YY}$.
\end{proof}

\begin{corollary}\label{cor:shifted symplectic}
  Let $k$ be a $\mathbf{Q}$-algebra.
  Let $\YY$ be an Artin stack equipped with an $n$-shifted symplectic structure in the sense of \cite{PTVV}.
  Assume that $\YY$ satisfies $(\textup{Map}\textsubscript{A})$.
  Then $\Loop_A(\YY)$ carries a canonical $(n-1)$-shifted symplectic structure.
\end{corollary}
\begin{proof}
  The AKSZ construction of \cite[Thm.~2.5]{PTVV} equips $\MAP(\check{A}, \YY)$ with an $(n-1)$-shifted symplectic structure.
  Since $\Loop_A(\YY) \hookrightarrow \MAP(\check{A}, \YY)$ is an open immersion by $(\textup{Map}\textsubscript{A})$, the restriction of this symplectic structure to $\Loop_A(\YY)$ is an $(n-1)$-shifted symplectic structure.
\end{proof}

\subsection{Mayer--Vietoris and codescent}

\begin{proposition}\label{prop:fin pres}
  Let $A$ be $\GGa$ (over a $\bQ$-algebra), $\GGm$, or an abelian variety over $k$.
  Let $f:\XX\to\YY$ be a morphism of Artin stacks satisfying $(\textup{Map}\textsubscript{A})$.
  If $f$ is locally of finite presentation, then so is $\Loop_A(f)$.
\end{proposition}
\begin{proof}
  The abelian case is Proposition~\ref{prop:fin pres abvar}.
  The $\GGa$ and $\GGm$ cases follow from the descriptions in Theorem~\ref{thm:LGa LGm}.
  For example, for $A=\GGa$, factor the map as
  \[ \TT[-1]\XX\longrightarrow\XX\times_{\YY}\TT[-1]\YY\longrightarrow\TT[-1]\YY, \]
  where the second map is a base change of $f$.
  The first map is a base change of the zero section of the stack of cosections $\mathbf{V}_{\XX}(\LL_{\XX/\YY})$ (see \cite[Rem.~B.3]{weilres}); by \cite[Prop.~B.1]{weilres}, the latter is locally of finite presentation since $\LL_{\XX/\YY}$ is perfect (see \cite[Prop.~17.4.2.3(1)]{SAG}), hence so is its zero section.
\end{proof}

\begin{lemma}\label{lem:L(f) open/et}
Let $A$ be as in Proposition~\ref{prop:fin pres}, and let $\YY$ be a stack.
\begin{enumerate}
  \item\label{item:L(f) open} If $\mathcal{U}\hookrightarrow \YY$ is a quasi-compact open immersion, then $\Loop_A(\mathcal{U}) \hookrightarrow \Loop_A(\YY)$ is a quasi-compact open immersion.
  \item\label{item:L(f) et} If $f : \YY \to \XX$ is an étale morphism of stacks satisfying $(\textup{Map}\textsubscript{A})$, then $\Loop_A(f) : \Loop_A(\YY) \to \Loop_A(\XX)$ is étale.
\end{enumerate}
\end{lemma}
\begin{proof}
(\ref{item:L(f) open}):
We have the commutative square \[\begin{tikzcd}
\Loop_A(\mathcal{U}) \arrow{r} \arrow[swap]{d}& \MAP(\check{A}, \mathcal{U}) \arrow{d}\\
\Loop_A(\YY)\arrow{r} & \MAP(\check{A}, \YY).
\end{tikzcd}\]
Since $\QCoh(\YY)\rightarrow\QCoh(\mathcal{U})$ is a tensor localization and $\HA\rightarrow\QCoh(\check{A})$ is fully faithful (Theorem~\ref{thm:HA vs QCoh checkA}), Lemma~\ref{Cart} implies that this square is cartesian.
It remains to check that $\MAP(\check{A},\mathcal{U}) \to \MAP(\check{A}, \YY)$ is a quasi-compact open immersion.
For $A=\GGa$ or $A=\GGm$, Theorem~\ref{thm:LGa LGm} shows that it is a base change of $\UU \hookrightarrow \YY$.
For $A$ an abelian variety, it is a quasi-compact open immersion by the proof of \cite[Prop.~19.1.1.1]{SAG}.

(\ref{item:L(f) et}):
Since $f$ is étale, $\MAP(\check{A},f)$ is formally étale (by the lifting criterion), hence so is $\Loop_A(f)$.
By Proposition~\ref{prop:fin pres}, $\Loop_A(f)$ is locally of finite presentation.
Hence $\Loop_A(f)$ is étale.
\end{proof}

\begin{corollary}\label{cor:Kun open}
  Let $A$ be as in Proposition~\ref{prop:fin pres} and $\mathcal{U} \hookrightarrow \YY$ a quasi-compact open immersion.
  Then for every morphism of stacks $\XX \to \YY$, the canonical morphism
  $$ \Loop_A(\mathcal{U} \times_{\YY} \XX) \longrightarrow \Loop_A(\mathcal{U}) \times_{\Loop_A(\YY)} \Loop_A(\XX) $$
  is an isomorphism in $\check{A}\circlearrowright \Stk$.
\end{corollary}
\begin{proof}
  By the proof of Lemma~\ref{lem:L(f) open/et}(\ref{item:L(f) open}), applied to $\mathcal{U} \hookrightarrow \YY$ and to $\mathcal{U} \times_{\YY} \XX \hookrightarrow \XX$, both sides are identified with $\MAP(\check{A}, \mathcal{U}) \times_{\MAP(\check{A}, \YY)} \Loop_A(\XX)$, since $\MAP(\check{A}, -)$ preserves fibred products.
\end{proof}

\begin{proposition}\label{eff epi}\leavevmode
\begin{enumerate}
\item
Let $\YY$ be a stack and $\{\UU_i\}_{i \in I}$ a finite cover by quasi-compact open substacks.
If $\YY$ is perfect fppf-locally on $S=\spc(k)$ and satisfies $(\textup{Map}\textsubscript{A})$, then $\{\Loop_A(\UU_i)\}$ is a cover of $\Loop_A(\YY)$ by quasi-compact open substacks.

\item
Let $\YY$ be a stack and $\{\XX_i\}_{i \in I}$ a finite Nisnevich cover of $\YY$, i.e., $p_i : \XX_i \to \YY$ are étale morphisms which are jointly surjective on $K$-points for every field extension $K/k$.
If $\YY$ is perfect fppf-locally on $\spc(k)$ and satisfies $(\textup{Map}\textsubscript{A})$, then $\{\Loop_A(\XX_i)\}$ is a Nisnevich cover of $\Loop_A(\YY)$.
That is, the morphisms $\Loop_A(p_i) : \Loop_A(\XX_i) \to \Loop_A(\YY)$ are \'etale and jointly surjective on $K$-points for every field extension $K/k$.
\end{enumerate}
\end{proposition}

\begin{proof} By Lemma~\ref{lem:L(f) open/et}, each map $\Loop_A(\UU_i)\rightarrow\Loop_A(\YY)$ is a quasi-compact open immersion. We show that the induced map $\coprod_{i}\Loop_A(\UU_{i})\rightarrow\Loop_A(\YY)$ is surjective on $K$-points for every field extension $K/k$; in particular, it is an effective epimorphism.

Given a $K$-point of $\Loop_A(\YY)$ classifying a tensor functor $\phi : \QCoh(\YY)\rightarrow \stHA(\spc(K)) = \Hcat_{A_K/K}$, we must show that $\phi$ factors through some $\QCoh(\UU_{i})$. The restrictions $\QCoh(\YY) \to \QCoh(\UU_i)$ are tensor localizations (since $\UU_i$ are quasi-compact opens), and they are jointly conservative (by faithfully flat descent), so $\{\QCoh(\UU_i)\}_i$ forms a cover of $\QCoh(\YY)$ in the sense of Definition~\ref{cover}. Since $\Hcat_{A_K/K}$ is an integral domain (Corollary~\ref{HE int}), the claim follows from Lemma~\ref{local}.

For part (ii), we first treat the case of an elementary Nisnevich cover.
Suppose given a Nisnevich square of stacks satisfying $(\textup{Map}\textsubscript{A})$,
  \[\begin{tikzcd}
      \mathcal{V} \ar{r}\ar{d}
      & \XX \ar{d}{p} \\
      \mathcal{U} \ar{r}{j}
      & \YY,
  \end{tikzcd}\]
  i.e., $j$ is a quasi-compact open immersion with reduced complement $\mathcal{Z}$, and $p$ is an étale morphism inducing an isomorphism $p^{-1}(\mathcal{Z}) \simeq \mathcal{Z}$.
  We claim that $\Loop_A(\mathcal{U}) \amalg \Loop_A(\XX) \to \Loop_A(\YY)$ is jointly surjective on $K$-points.
  It will suffice to show that the induced square
  \[\begin{tikzcd}
      \Loop_A(\mathcal{V}) \ar{r}\ar{d}
      & \Loop_A(\XX) \ar{d}{\Loop_A(p)} \\
      \Loop_A(\mathcal{U}) \ar{r}{\Loop_A(j)}
      & \Loop_A(\YY)
  \end{tikzcd}\]
  is Nisnevich.
  By Corollary~\ref{cor:Kun open}, it is cartesian.
  By Lemma~\ref{lem:L(f) open/et}, $\Loop_A(j)$ is a quasi-compact open immersion and $\Loop_A(p)$ is étale.
  By Theorem~\ref{thm:Kun affine}, $\Loop_A(p^{-1}(\mathcal{Z})) \simeq \Loop_A(\mathcal{Z}) \times_{\Loop_A(\YY)} \Loop_A(\XX) \to \Loop_A(\mathcal{Z})$ is invertible.
  Note that $\Loop_A(\UU) \times_{\Loop_A(\YY)} \Loop_A(\ZZ) \simeq \Loop_A(\emptyset)$ is empty by Corollary~\ref{cor:Kun open}.

  To conclude that the square is Nisnevich, it remains to show that $\Loop_A(\ZZ)$ is complementary to $\Loop_A(\UU)$.
  Let us show that any $K$-point $\gamma$ of $\Loop_A(\YY)$ outside $\Loop_A(\mathcal{U})$ lifts to $\Loop_A(\mathcal{Z})$.
  By the cartesian square in the proof of Lemma~\ref{lem:L(f) open/et}(\ref{item:L(f) open}), the underlying map $\gamma : \check{A}_K \to \YY$ does not factor through $\UU$, so $\check{A}_K \times_{\YY} \ZZ$ is nonempty.
  We first show that $\gamma$ factors through $\ZZ$:
  \begin{itemize}
	\item
  For $A = \GGa$ or $\GGm$, $\check{A}_K$ is classical, so $\gamma$ factors through the classical truncation of $\YY$.
  Then, its restriction along the effective epimorphism $\spc(K) \twoheadrightarrow \check{A}_K$ lands in $\ZZ$ (because $\spc(K)$ is reduced), hence so does $\gamma$.

  \item
  For $A$ an abelian variety, $\check{A}_K \times_{\YY} \ZZ \hookrightarrow \check{A}_K$ is a nonempty closed immersion with structure sheaf $\gamma^*\OO_{\ZZ}$.
  The defining ideal of its classical truncation is thus $\mathcal{J} \neq \OO_{\check{A}_K}$.
  In fact, since $\gamma^*\OO_{\ZZ} \in \Hcat_{A_K/K}$ has flat cohomology sheaves (Theorem~\ref{thm:HA vs QCoh checkA}) and $\check{A}_K$ is integral, we have $\mathcal{J} = 0$; hence $\gamma$ factors through $\ZZ$ (as $\check{A}_K$ is classical).
  \end{itemize}
  It remains to see that the resulting factorization $\gamma : \check{A}_K \to \ZZ$ is a loop.
  As $i : \ZZ \to \YY$ is affine, the objects $i^*F$, $F \in \QCoh(\YY)$, generate $\QCoh(\ZZ)$ under colimits.
  Since the $\gamma^*i^*F$ lie in $\Hcat_{A_K/K}$, because $\gamma$ is a loop of $\YY$, and $\Hcat_{A_K/K} \subseteq \QCoh(\check{A}_K)$ is closed under colimits (Theorem~\ref{thm:HA vs QCoh checkA}), it follows that $\gamma^*G \in \Hcat_{A_K/K}$ for every $G \in \QCoh(\ZZ)$.
  This shows the claim, and hence that $\Loop_A(-)$ preserves elementary Nisnevich covers.

  The general case follows, since any Nisnevich cover of a quasi-compact stack with affine diagonal can be refined to a finite sequence of elementary Nisnevich squares (cf.\ \cite[Rem.~5.5]{Rydh11}).
\end{proof}

\begin{theorem}[Zariski codescent]\label{Zarcosh}\leavevmode
  \begin{enumerate}
    \item Let $\YY$ be a stack, and $\{\UU_i\}_{i \in I}$ a finite cover by quasi-compact open substacks.
    If $\YY$ is perfect fppf-locally on $S=\spc(k)$ and satisfies $(\textup{Map}\textsubscript{A})$, then the canonical map
    $$\operatorname*{colim}_{[n] \in \mathbf{\Delta}^{\mathrm{op}}} \coprod_{(i_0,\ldots,i_n) \in I^{n+1}} \Loop_A(\UU_{i_0} \cap \cdots \cap \UU_{i_n}) \xrightarrow{\;\sim\;} \Loop_A(\YY)$$
    is invertible.

    \item If $\YY$ is quasi-compact, perfect fppf-locally on $S=\spc(k)$, and satisfies $(\textup{Map}\textsubscript{A})$, then the assignment $\UU \mapsto \Loop_A(\UU)$ satisfies codescent for the Zariski topology on $\YY$.
  \end{enumerate}
\end{theorem}
\begin{proof}
  We show (i); part (ii) follows because quasi-compactness ensures any Zariski cover has a finite refinement.
  By Proposition~\ref{eff epi} and Lemma~\ref{lem:L(f) open/et}, $\Loop_A(\UU_i) \to \Loop_A(\YY)$ form a jointly surjective family of quasi-compact open immersions, so $\coprod_i \Loop_A(\UU_i) \to \Loop_A(\YY)$ is an effective epimorphism.
  Thus, $\Loop_A(\YY)$ is isomorphic to the geometric realization of its \v{C}ech nerve.
  It remains to observe that by Corollary~\ref{cor:Kun open},
  $$\Loop_A(\UU_{i_0} \cap \cdots \cap \UU_{i_n}) \;\simeq\; \Loop_A(\UU_{i_0}) \times_{\Loop_A(\YY)} \cdots \times_{\Loop_A(\YY)} \Loop_A(\UU_{i_n}),$$
  and that $\Loop_A(-)$ preserves finite disjoint unions: indeed, for a finite family of stacks $\{\mathcal{V}_\alpha\}$, Lemma~\ref{lem:L(f) open/et} makes each $\Loop_A(\mathcal{V}_\alpha) \hookrightarrow \Loop_A(\coprod_\beta \mathcal{V}_\beta)$ a quasi-compact open immersion, Corollary~\ref{cor:Kun open} gives $\Loop_A(\mathcal{V}_\alpha) \times_{\Loop_A(\coprod_\beta \mathcal{V}_\beta)} \Loop_A(\mathcal{V}_{\alpha'}) \simeq \Loop_A(\emptyset) = \emptyset$ for $\alpha \ne \alpha'$, and Proposition~\ref{eff epi} shows that they cover.
\end{proof}

\begin{theorem}[Nisnevich codescent]\label{Niscosh}
  Let $\{p_i : \UU_i \to \YY\}_{i \in I}$ be a finite Nisnevich covering family.
  If each $\UU_i$ is perfect fppf-locally on $S=\spc(k)$ and satisfies $(\textup{Map}\textsubscript{A})$, then the canonical map
  $$\operatorname*{colim}_{[n] \in \mathbf{\Delta}^{\mathrm{op}}} \coprod_{(i_0,\ldots,i_n) \in I^{n+1}} \Loop_A(\UU_{i_0} \times_{\YY} \cdots \times_{\YY} \UU_{i_n}) \xrightarrow{\;\sim\;} \Loop_A(\YY)$$
  is invertible.
\end{theorem}
\begin{proof}
  By Proposition~\ref{eff epi}(ii), the maps $\Loop_A(p_i)$ are jointly surjective on $K$-points, and by Lemma~\ref{lem:L(f) open/et}(\ref{item:L(f) et}) each $\Loop_A(p_i)$ is étale.
  By K\"unneth (Theorem~\ref{Kun}), $\Loop_A$ sends the \v{C}ech nerve of $\{p_i\}$ levelwise to the \v{C}ech nerve of $\{\Loop_A(p_i)\}$ (perfectness ensures K\"unneth applies at each level).
  Hence the left-hand side is the geometric realization of the \v{C}ech nerve of the jointly surjective family $\{\Loop_A(p_i)\}$, which is isomorphic to the target.
\end{proof}

\part{Loop cohomology}\label{part:cohomology}

We now arrive at the main objective of this paper, the construction of cohomology theories associated to the $1$-dimensional commutative group scheme $A$.
Given a $G$-variety $X$, the structure sheaf of the $A$-loop space $\Loop_A^G(X)$, together with its equivariance for the loop-rotation action of $\check{A}$, will give rise to the \emph{$A$-loop cohomology} $\HH_A^G(X)$, valued in the category $\Rot^G(A)$ of $G$-equivariant $A$-\emph{rotation modules}.
The rotation module structure will then be used to extract two further theories, the \emph{periodic $A$-loop cohomology theory} $\Hcoh_A^G(-)$ and, for $G$ diagonalizable, the equivariant refinement $\K_A^G$ of $\K_A$ (\ref{sssec:Bott periodization}), by periodization in two dual ways.
In the elliptic case, $\Hcoh_E^G(-)$ is our definition of equivariant elliptic cohomology (for a fixed elliptic curve), and the $\K_E^G$ assemble into equivariant elliptic cohomology in the sense of Lurie and Gepner--Meier (Subsect.~\ref{ssec:HH on XAG}).

\section{Loop cohomology}\label{sec:cohomology}

The $A$-loop cohomology $\HH_A^G(X)$ will be the pushforward of the structure sheaf along the equivariant $A$-loop space (Definition~\ref{def:HH_A^G}), and periodic $A$-loop cohomology $\Hcoh_A^G(-)$ will be its Euler periodization (Definition~\ref{def:H_A^G}).
In the additive case $\HH_{\GGa}^G(X)$ is the Hodge-filtered derived de Rham complex of $X/G$ (Subsect.~\ref{ssec:dR comparison}). In the multiplicative case $\HH_{\GGm}^G(X)$ and $\Hcoh_{\GGm}^G(X)$ are ``delocalized'' Hochschild homology and periodic cyclic homology, respectively (Subsect.~\ref{ssec:TP/KU}).

Throughout the section, we fix a $1$-shifted Cartier dual pair $(A, \check{A}, \jmath)$, and an affine fppf group scheme $G$, over the connective $E_\infty$-ring $k$.

\subsection{Loop sheaves and rotation modules}\label{ssec:sheaves}

\subsubsection{}

Let $\YY$ be a stack and recall the $A$-loop space $\Loop_A(\YY)$, with its $\check{A}$-action by loop rotation.
The \emph{category of $A$-loop sheaves on $\YY$} is defined as
\[
  \Shv_A(\YY) \;:=\; \INDPERF^{\check{A}}(\Loop_A(\YY)),
\]
the category of ind-perfect complexes on the quotient $\Loop_A(\YY)/\check{A}$ (see Subsect.~\ref{ssec:indperf}).
We denote by
\[
  \Shvper_A(\YY) := \INDPERF^{\check{A}}(\Loop_A(\YY))[\alpha_A^{-1}]
\]
the full subcategory of \emph{periodic $A$-loop sheaves}, i.e., those which are Euler-periodic as in (\ref{sssec:Euler periodization}).

The category $\Shv_A(\YY)$ is closed symmetric monoidal, with unit object $\un_\YY$ (= the structure sheaf of $\Loop_A(\YY)/\check{A}$), tensor product $(-) \otimes_{\YY} (-)$, and internal Hom $\hhoms_{\YY}(-, -)$.
We omit the subscripts when there is no ambiguity.

By construction, $\Shv_A(\YY)$ is compactly generated.
We say that an $A$-loop sheaf is \emph{constructible} if it is compact, i.e., lies in the full subcategory
\[
  \Shv_{A,\mathrm{c}}(\YY) := \PERF^{\check{A}}(\Loop_A(\YY)).
\]

\subsubsection{}

Let $f : \XX \to \YY$ be a morphism of stacks, and denote by
\[\phi := \Loop_A(f)/\check{A} : \Loop_A(\XX)/\check{A} \to \Loop_A(\YY)/\check{A}\]
the induced morphism on loop rotation quotients.
The functors $f^{*} := \phi^*_{\INDPERF}$, $f_{*} := \phi_*^{\INDPERF}$ form an adjoint pair
\[
  f^* : \Shv_A(\YY) \rightleftarrows \Shv_A(\XX) : f_*.
\]

\subsubsection{Equivariant setting}

Let $G$ be an affine fppf group scheme.
For $X \in \Stk^G$, we may specialize the above to $\YY = X/G$.
We write
\[
  \Shv_A^G(X) := \Shv_A(X/G), \qquad
  \Shv_A^{G,\mathrm{per}}(X) := \Shvper_A(X/G),
\]
and $\Shv_{A,\mathrm{c}}^G(X) := \Shv_{A,\mathrm{c}}(X/G)$.

If $f : X \to Y$ is a $G$-equivariant morphism, we write simply
\[
  f^* := (f/G)^* : \Shv_A^G(Y) \to \Shv_A^G(X),
\]
where $f/G : X/G \to Y/G$ is the induced morphism, and similarly for $f_*$.

\subsubsection{Rotation modules}

For $X = *$ we call
\[
  \Rot^G(A) := \Shv_A^G(*) = \INDPERF^{\check{A}}(\Loop_A^G(*))
\]
the category of \emph{$G$-equivariant $A$-rotation modules}.
We denote its full subcategory of Euler-periodic objects by $\Rot^{G,\mathrm{per}}(A) := \Shv_A^{G,\mathrm{per}}(*)$.

For $G=1$ the trivial group, we will omit the $^G$ decoration.
Recall from Corollary~\ref{cor:IndPerf checkA} that we have
\[
  \Rot(A) \simeq \INDPERF(B\check{A}) \simeq \Ind(\MOD_{\RA}^{\mathrm{perf}/k}),
\]
the ind-completion of the category of modules over the rotation algebra $\RA$ (Subsect.~\ref{ssec:RA}) with perfect underlying $k$-module.

For example, for $A = \GGm$ we have $\R_{\GGm} \simeq k[S^1]$ (\ref{sssec:OO(Gm[-1])}), so $\Rot(\GGm)$ is the ind-completion of the category of $k[S^1]$-modules perfect over $k$, i.e. perfect $k$-modules with $S^1$-action (Example~\ref{ex:Tate multiplicative}).

\subsubsection{}

If $f : X \to Y$ is a proper morphism in $\Var^G$, so that $\phi$ is proper representable (Corollary~\ref{mapsrepr}), then we also set
\[
  f^{\dag} := \phi^\dag_{\INDPERF} : \Shv_A^G(Y) \to \Shv_A^G(X),
\]
see Subsect.~\ref{ssec:indperf dualizing complex}.
In this case we also define the \emph{relative dualizing complex} of $f$ as
\[
  \omega_f := f^\dag(\un) \simeq \omega_{\phi}^{\INDPERF} \quad \in \Shv_A^G(X).
\]
By Proposition~\ref{prop:f^! formula indperf}, $f^{\dag}(\cF) \simeq f^*(\cF) \otimes \omega_f$ for $\cF \in \Shv_{A,\mathrm{c}}^G(Y)$.
Construction~\ref{constr:trace indperf} yields for every $\cF \in \Shv_{A,\mathrm{c}}^G(Y)$ the natural morphism
\begin{equation}\label{eq:trace maps}
  \operatorname{tr}^A_f : f_* f^\dag (\cF) \longrightarrow \cF
\end{equation}
in $\Shv_A^G(Y)$.

\subsubsection{Exchange transformations}\label{sssec:exchange}

Consider a cartesian square in $\Var^G$
\[\begin{tikzcd}
  X' \ar{r}{g}\ar[swap]{d}{p} & Y' \ar{d}{q} \\
  X \ar{r}{f} & Y
\end{tikzcd}\]
with $f$ proper.
Applying $\Loop_A^G(-)/\check{A}$, the induced square of stacks is cartesian by Theorem~\ref{Kun}, and Construction~\ref{constr:Ex f^* f^dag} gives the exchange transformation
\[
\operatorname{Ex}^{*,\dag} : p^* f^\dag \to g^\dag q^*
\]
of functors $\Shv_A^G(Y) \to \Shv_A^G(X')$.

\subsubsection{Induction and restriction functors}\label{sssec:res ind}

Given a homomorphism $H \to G$ of fppf group schemes, $H$ acts on $X \in \Stk^G$ by restriction, and there is a canonical map $f : X/H \to X/G$ which induces an adjoint pair
\[
  (\operatorname{res}^G_H, \operatorname{ind}^G_H) := (f^*, f_*), \qquad
  \operatorname{res}^G_H : \Shv_A^G(X) \rightleftarrows \Shv_A^H(X) : \operatorname{ind}^G_H.
\]

\subsubsection{Functoriality in $A$}\label{sssec:Shv func A}

Let $h : (A,\check{A},\jmath) \to (A',\check{A}',\jmath')$ be a morphism of $1$-shifted Cartier dual pairs.
For every stack $\YY$, the correspondence
\[\begin{tikzcd}
  & \Loop_A(\YY)/\check{A}' \ar[swap]{ld}{\beta}\ar{rd}{\alpha} &
  \\
  \Loop_A(\YY)/\check{A} & & \Loop_{A'}(\YY)/\check{A}',
\end{tikzcd}\]
where $\alpha$ is induced by $A \to A'$ and $\beta$ is the base change of $q : B\check{A}' \to B\check{A}$, gives rise to the canonical functor
\[
  \operatorname{res}^{A'}_A := \operatorname{res}_h := \beta_* \alpha^* : \Shv_{A'}(\YY) \to \Shv_A(\YY).
\]

For any morphism $f : \XX \to \YY$, we have the exchange transformation
\begin{equation}
  \operatorname{res}_A^{A'} \bigl(f_*\cF\bigr) \to f_* \bigl( \operatorname{res}_A^{A'} (\cF) \bigr)
\end{equation}
for every $\cF \in \Shv_{A'}(\XX)$.
For example, take $f : X/G \to *$ with $X \in \Var^G$; for $\YY=*$ the correspondence above has $\alpha = \id_{B\check{A}'}$ and $\beta=q$, and this becomes a canonical map
\begin{equation}\label{eq:push to res A A'}
	q_* f_*(\cF) \to f_* \bigl( \operatorname{res}_A^{A'} (\cF) \bigr),
	\qquad \cF \in \Shv_{A'}^G(X),
\end{equation}
in $\Rot(A)$.

\subsubsection{Weights}\label{sssec:Shv wt=chi}

Suppose that $G$ is a diagonalizable group scheme over $k$, with group of characters $X^*(G)$.
Recall from Subsect.~\ref{ssec:LA(BG) diag} the canonical correspondence of $\check{A}$-equivariant stacks
\[
  \begin{tikzcd}
    & \frX_A(G) \times BG \ar[swap]{ld}{p_{A,G}}\ar{rd}{q_{A,G}} &
    \\
    \Loop_A^G(*) & & \frX_A(G),
  \end{tikzcd}
\]
where $p_{A,G}$ is an isomorphism when $k$ is a $\bQ$-algebra, and on classical truncations.
Given a character $\chi \in X^*(G)$, recall the $\mathcal{P}_{-\chi}$-equivariant line bundle $\OO(-\chi)$ of (\ref{sssec:O(chi)}).
We define the functor
\begin{equation}\label{eq:wt=chi}
  (-)^{\mathrm{wt}=\chi} :=
  q_{A,G,*} \bigl( p_{A,G}^*(-) \otimes \OO(-\chi) \bigr)
\end{equation}
on $\Rot^G(A) = \INDPERF^{\check{A}}(\Loop_A^G(*))$.
More precisely, it is the composite
\begin{equation*}
  \begin{aligned}
    \INDPERF^{\check{A}}(\Loop_A^G(*))
    &\xrightarrow{\mathmakebox[2.8em]{p_{A,G}^*}}
    \INDPERF^{\check{A}}(\frX_A(G) \times BG)
    \\
    &\xrightarrow{\mathmakebox[2.8em]{\operatorname{res}^{\check{A}}_{\mathcal{P}_{-\chi}}}}
    \INDPERF^{\mathcal{P}_{-\chi}}(\frX_A(G) \times BG)
    \\
    &\xrightarrow{\mathmakebox[2.8em]{- \otimes \OO(-\chi)}}
    \INDPERF^{\mathcal{P}_{-\chi}}(\frX_A(G) \times BG)
    \\
    &\xrightarrow{\mathmakebox[2.8em]{q_{A,G,*}}}
    \INDPERF^{\mathcal{P}_{-\chi}}(\frX_A(G)),
  \end{aligned}
\end{equation*}
where $\operatorname{res}^{\check{A}}_{\mathcal{P}_{-\chi}}$ restricts the equivariance along the projection $\mathcal{P}_{-\chi} \to \check{A} \times \frX_A(G)$.

By construction, \eqref{eq:wt=chi} is a lift of the usual weight-$\chi$ functor
\[
  \INDPERF(\Loop_A^G(*)) \xrightarrow{p_{A,G}^*}
  \INDPERF(\frX_A(G) \times BG) \xrightarrow{(-)^{\mathrm{wt}=\chi}}
	\INDPERF(\frX_A(G)),
\]
see (\ref{sssec:QCohBG diag}).\footnote{%
  Here, $\frX_A(G)$ and $\frX_A(G) \times BG$ are perfect (since $\frX_A(G)$ is a quasi-compact separated scheme and $G$ is linearly reductive), so $\INDPERF(\YY) = \QCOH(\YY)$ in both cases.
}

\subsubsection{Weight zero} \label{sssec:Shv wt=0}

For $\chi = 0$ the extension $\mathcal{P}_0 \simeq \check{A} \times \GGm \times \frX_A(G)$ is canonically split, since $\ev_A^{\mathrm{str}}$ is strictly multiplicative in the $A$-variable, and the weight-zero functor canonically identifies with
\begin{multline}
	(-)^{\mathrm{wt}=0} \simeq
	q_{A,G,*} p_{A,G}^* :
	\Rot^G(A) \to
	\INDPERF^{\check{A}}(\frX_A(G))
	\\
	\simeq \INDPERF(B\check{A} \times \frX_A(G)).
\end{multline}
This is lax symmetric monoidal (since $q_{A,G,*}$ is), hence preserves commutative algebra objects.
We have:

\begin{proposition}\label{prop:eqper trivialization}
The canonical map
\[
  \OO_{\frX_A(G)}
  \to
  \bigl(\OO_{\Loop_A(BG)}\bigr)^{\mathrm{wt}=0}
\]
is an isomorphism of commutative algebras in $\INDPERF^{\check{A}}(\frX_A(G))$, where the source carries the trivial $\check{A}$-equivariant structure.
Moreover, $(\OO_{\Loop_A(BG)})^{\mathrm{wt}=\chi} \simeq 0$ for all $\chi\ne 0$.
\end{proposition}

\begin{proof}
The map in question is the unit map $\OO \to q_*q^*(\OO)$, where $q : (\frX_A(G) \times BG)/\check{A} \to \frX_A(G)/\check{A} \simeq B\check{A} \times \frX_A(G)$ is the map induced on quotient stacks by $q_{A,G}$.
On perfect complexes, $q^*$ is the totalization of the cosimplicial diagram of functors
\[
  q_\bullet^* : \PERF(\check{A}^\bullet \times \frX_A(G)) \to
  \PERF(\check{A}^\bullet \times \frX_A(G) \times BG),
\]
each of which is identified with the inclusion of weight-zero objects (\ref{sssec:QCohBG diag}).
In particular, $q^*$ is fully faithful on $\PERF$, and hence on $\INDPERF$ by \cite[Prop.~5.3.5.11(1)]{LurieHTT}.

Similarly, for $\chi \ne 0$ consider the map $q_\chi : (\frX_A(G) \times BG)/\mathcal{P}_{-\chi} \to \frX_A(G)/\mathcal{P}_{-\chi}$ induced by $q_{A,G}$; by definition, $(\OO_{\Loop_A(BG)})^{\mathrm{wt}=\chi} \simeq q_{\chi,*}\OO(-\chi)$.
For every $F \in \PERF(\frX_A(G)/\mathcal{P}_{-\chi})$ we have
\[
  \Map_{\INDPERF}(F, q_{\chi,*}\OO(-\chi)) \simeq
  \Map_{\PERF}(q_\chi^*F, \OO(-\chi)) \simeq 0,
\]
since $q_\chi^*F$ has weight zero and $\OO(-\chi)$ has weight $-\chi$.
The same thus holds for every ind-perfect complex $F$, i.e., $q_{\chi,*}\OO(-\chi) \simeq 0$.
\end{proof}

\subsection{$A$-loop cohomology}

Given $X \in \Stk^G$, we will write $a_X : X \to *$ for the structure morphism.

\begin{definition}\label{def:HH_A^G}
The \emph{$G$-equivariant $A$-loop cohomology} of $X \in \Stk^G$ is the $G$-equivariant $A$-rotation module
\[
  \HH^G_A(X) \;:=\; a_{X,*}(\OO) \in \Rot^G(A).
\]
More generally, for $\cF \in \Shv_A^G(X)$, the \emph{$G$-equivariant $A$-loop cohomology of $X$ with coefficients in $\cF$} is
\[
  \HH_A^G(X; \cF) \;:=\; a_{X,*}(\cF).
\]
\end{definition}

\begin{definition}\label{def:BM homology}
The \emph{$G$-equivariant Borel--Moore $A$-loop homology of $X$} is the $G$-equivariant $A$-rotation module
\[
  \HH_A^{G,BM}(X) \;:=\; \HH_A^G(X; \omega_{X}) \;:=\; a_{X,*} a_X^! (\un).
\]
\end{definition}

\begin{definition}\label{def:H_A^G} The \emph{periodic $G$-equivariant $A$-loop cohomology of $X$} is the periodization $$\Hcoh_A^G(X) \;:=\; \HH_A^G(X)[\alpha_A^{-1}] \in \Rot^{G,\mathrm{per}}(A).$$
Similarly, for $\cF \in \Shv_A^G(X)$, $\Hcoh_A^G(X; \cF) := \HH_A^G(X; \cF)[\alpha_A^{-1}]$.
\end{definition}

\subsection{Additive case: Hodge/de Rham}\label{ssec:dR comparison}

We now explain that in characteristic zero, $\HH_{\GGa}^G(X) \in \Rot^G(\GGa)$ is a refinement of the derived de Rham complex of $X/G$, together with its Hodge filtration.

Let $k$ be a regular noetherian $\bQ$-algebra, $G$ an affine fppf group scheme over $k$ such that $BG$ is perfect\footnote{%
  e.g. satisfying the condition of Theorem~\ref{thm:pluperfect quotient}(i)
}, and $X \in \Var^G$.

Under the isomorphism $\Loop_{\GGa}^G(*) \simeq \TT[-1]BG \simeq \mathfrak{g}/G$ of Example~\ref{exam:equivariant loops Ga Gm}, the category $\Rot^G(\GGa)$ of additive $G$-equivariant rotation modules is equivalent to $\INDPERF^{\check{\GG}_a}(\mathfrak{g}/G)$, where $\check{\GG}_a = B\widehat{\GG}_a$ is the $1$-shifted Cartier dual of $\GGa$.

\subsubsection{Underlying quasi-coherent sheaf}

We first forget rotation-equivariance via the functor $\sfe^* : \INDPERF^{\check{\GG}_a}(\mathfrak{g}/G) \to \QCOH(\mathfrak{g}/G)$ (see \ref{sssec:forget rotation equivariance}).\footnote{%
  Since $\mathfrak{g}/G$ is perfect (Theorem~\ref{thm:pluperfect quotient}(i)), $\INDPERF(\mathfrak{g}/G) \simeq \QCOH(\mathfrak{g}/G)$.
}
By Theorem~\ref{thm:LGa LGm} and Lemma~\ref{lem:sfe base change}, $\sfe^* \HH_{\GGa}^G(X)$ is identified with the pushforward of the structure sheaf along
\[
  \pi_* : \QCOH(\TT[-1](X/G)) \to \QCOH(\TT[-1]BG),
\]
where $\pi$ is induced by the projection $f : X/G \to BG$.
Since $\TT[-1]BG \to BG$ is affine, $\QCOH(\TT[-1]BG)$ is the category of $\Sym^*_{\OO_{BG}}(\mathbf{L}_{BG}[1])$-modules in $\QCOH(BG)$, and as such
\[ \sfe^* \HH_{\GGa}^G(X) \simeq f_*\Sym^*_{\OO_{X/G}}(\mathbf{L}_{X/G}[1]) \simeq \bigoplus_{n \ge 0} f_*\Omega^n_{X/G}[n], \]
where $\Omega^n_{X/G} := \wedge^n_{\OO_{X/G}}\mathbf{L}_{X/G}$.
The module structure on $f_*\Sym^*_{\OO_{X/G}}(\mathbf{L}_{X/G}[1])$ is via the algebra map $\Sym^*_{\OO_{BG}}(\mathbf{L}_{BG}[1]) \to f_*\Sym^*_{\OO_{X/G}}(\mathbf{L}_{X/G}[1])$.
Note that $f_*$ commutes with direct sums since $f$ is representable by qcqs algebraic spaces.

Note that $\Omega^*_{X/G} := \bigoplus_{n \ge 0} \Omega^n_{X/G}[n]$ is the underlying graded module of the derived de Rham complex of $X/G$.

\subsubsection{Loop rotation}

After induction along $G \twoheadrightarrow 1$ (\ref{sssec:res ind}), the resulting object $\operatorname{ind}_G^1\HH_{\GGa}^G(X) \simeq \HH_{\GGa}(X/G)$ (Example~\ref{exam:ind_G^1 HH}) has a residual $\check{\GG}_a$-equivariant structure:
\[\HH_{\GGa}(X/G) \in \Rot(\GGa) \simeq \INDPERF(B\check{\GG}_a).\]
Via the Koszul duality equivalence $\INDPERF(B\check{\GG}_a) \simeq \MOD_{k[u]}$ of Corollary~\ref{cor:IndPerf BcheckA Kosz}, we identify this with a $k[u]$-module, where $u$ has degree $-2$.
By \cite[Thm.~1.3]{BNloop}, this $u$-action is the de Rham differential, so
\[ \HH_{\GGa}(X/G) \simeq \bigl(\Omega^*_{X/G}[u],\; u \cdot d_{\mathrm{dR}}\bigr) \quad\text{in } \MOD_{k[u]} \]
is the $u$-twisted derived de Rham complex of $X/G$.
Here we have abused notation by writing simply $\Omega^*_{X/G}$ for its global sections.

\subsubsection{Loop rescaling}

Let $\GGm$ act on $A=\GGa$ by weight $-1$ scaling.
This gives rise to the loop rescaling action of $\GGm$ on $\Loop_{\GGa}^G(X)$, which corresponds to weight $+1$ scaling on $\TT[-1](X/G)$.
This places each $\Omega^n_{X/G}[n] \simeq \Sym^n_{\OO_{X/G}}(\mathbf{L}_{X/G}[1])$ in weight $+n$.
Since $\mathrm{Lie}(\GGa)$ has weight $+1$, the cohomology Euler class $u$ has weight $-1$.
Thus:

\begin{proposition}\label{prop:Ga comparison}
  As an object of $\MOD_{k[u]}^{\GGm}$, $\HH_{\GGa}(X/G)$ is the Hodge-filtered $u$-twisted derived de Rham complex of $X/G$:
  \[
    \HH_{\GGa}(X/G) \simeq (\Omega^*_{X/G}[u],\; u \cdot d_{\mathrm{dR}}),
  \]
  with $\Omega^n_{X/G}$ in weight $+n$ and $u$ in weight $-1$.
  The weight-zero piece expands as $\bigoplus_{n \ge 0} \Omega^n_{X/G}[n] \cdot u^n$, with $u \cdot d_{\mathrm{dR}}\colon \Omega^n[n] \cdot u^n \to \Omega^{n+1}[n+1] \cdot u^{n+1}$ taking each weight-zero generator to the next.
\end{proposition}

\begin{corollary}\label{cor:Ga comparison periodic}
There is a canonical isomorphism
\[
  \Hcoh_{\GGa}(X/G) \simeq (\Omega^*_{X/G}[u, u^{-1}],\; u \cdot d_{\mathrm{dR}}),
\]
identifying the periodic additive $G$-equivariant cohomology of $X$ with the $2$-periodized derived de Rham complex.
\end{corollary}

\begin{proof}
The weight-zero summand of $\Hcoh_{\GGa}(X/G)$ is given by
\[ \bigl(\HH_{\GGa}(X/G) \otimes_{k[u]} k[u, u^{-1}]\bigr)^{\GGm = 0} \simeq \bigoplus_{n \ge 0} \Sym_{\OO_{X/G}}^n(\mathbf{L}_{X/G}[1]) \cdot u^n, \]
since $\Sym^n$ has weight $+n$, $u^j$ has weight $-j$, and weight $0$ forces $n = j \geq 0$.
Multiplication by $u$ shifts weights by $-1$, so we have
$$\Hcoh_{\GGa}(X/G) \simeq \Hcoh_{\GGa}(X/G)^{\GGm = 0} \otimes_k k[u, u^{-1}] \simeq (\Omega^*_{X/G}, d_{\mathrm{dR}}) \otimes_k k[u, u^{-1}]$$
as claimed.
\end{proof}

\subsubsection{Shearing}

Recall that \emph{shearing} defines a symmetric monoidal equivalence
\[
  \MOD_{k[u]}^{\GGm} \;\xrightarrow{\;\sim\;}\; \MOD^{\GGm}_{k[t]},
\]
where $t$ is in degree $0$ and weight $-1$, preserving weights and shifting each weight-$w$ piece by $[-2w]$.\footnote{%
  The even shift is so that the Koszul sign rule on tensor products is unchanged (whence the symmetric monoidality).
}
See e.g. \cite[\S 2.2.1]{BLR}.
In particular, applied to $\Omega^*_{X/G} = \bigoplus_{n \ge 0} \Omega^n_{X/G}[n]$, with $\Omega^n_{X/G}[n]$ in weight $+n$, it gives
\[
  \Omega^{*,\Rightarrow}_{X/G} := \bigoplus_{n \ge 0} \Omega^n_{X/G}[-n].
\]
We thus get:

\begin{corollary}\label{cor:Ga comparison sheared}
As an object of $\MOD_{k[t]}^{\GGm}$, $\HH_{\GGa}(X/G)$ is identified with the Rees construction of the Hodge-filtered derived de Rham complex of $X/G$:
\[
  \HH_{\GGa}(X/G) \simeq (\Omega^{*,\Rightarrow}_{X/G}[t],\; d_{\mathrm{dR}}),
\]
with $\Omega^n_{X/G}$ in weight $+n$ and $t$ in weight $-1$.
The weight-zero piece expands as $\bigoplus_{n \ge 0} \Omega^n_{X/G}[-n] \cdot t^n$, with $d_{\mathrm{dR}}\colon \Omega^n[-n] \cdot t^n \to \Omega^{n+1}[-n-1] \cdot t^{n+1}$ taking each weight-zero generator to the next.
\end{corollary}

In other words, under the equivalence between $\MOD_{k[t]}^{\GGm}$ and the category of filtered objects in $\MOD_k$, $\HH_{\GGa}(X/G)$ with its $\GGm$-equivariant structure encodes the Hodge-filtered de Rham complex of $X/G$.

\subsubsection{}

To summarize, $\HH_{\GGa}^G(X)$ encodes a one-parameter family of cohomology theories:
\[\begin{tikzcd}[column sep=tiny, row sep=small]
  & \HH_{\GGa}^G(X) = \text{Hodge} \arrow[dl, squiggly, "u\,=\,0"'] \arrow[dr, squiggly, "u\,=\,1"] & \\
  \Omega^*_{X/G} = \text{Dolbeault} & & (\Omega^*_{X/G}, d_{\mathrm{dR}}) = \text{de Rham}
\end{tikzcd}\]
where the degree $-2$ parameter $u$ can be converted to the usual degree $0$ parameter $t$ by shearing.

We comment on the relation to Simpson's Hodge stack, which is another geometric incarnation of the Hodge-filtered de Rham complex (see \cite{Simp} or \cite{BhattFGauges}).
The $\GGa$-loop space $\Loop_{\GGa}(X/G)$ appears in Simpson's picture as the \emph{Dolbeault} stack,
\[ \widehat{\Loop}_{\GGa}(X/G) \simeq \widehat{\TT}[-1](X/G) \simeq (X/G)_{\mathrm{Dol}}, \]
where the ``hats'' indicate formal completion (along constant loops / the zero section).
The stack $\widehat{\TT}[-1](X/G)$ itself encodes only the Dolbeault complex $\Omega^*_{X/G} := \bigoplus_{n \ge 0} \Omega^n_{X/G}[n]$.
The two approaches can be summed up as follows:
\begin{itemize}
  \item Simpson constructs a one-parameter family $(X/G)_{\mathrm{Hod}} \to \mathbf{A}^1/\GGm$, with coordinate $t$, with $(X/G)_{\mathrm{Dol}}$ at $t=0$.
  There is an equivalence between $\QCOH(\aff/\GGm)$ and the category of filtered objects in $\MOD_k$, and the Hodge-filtered de Rham complex is then encoded by functions on $(X/G)_{\mathrm{Hod}}$.
  The de Rham complex itself is encoded by functions on the $t=1$ fibre, the de Rham stack $(X/G)_{\mathrm{dR}}$.
  \item We consider the loop rotation action of $\check{\GG}_a = B\GGhat_a$ on $\Loop_{\GGa}^G(X)$, and the loop rescaling action of $\GGm$ on $\Loop_{\GGa}(X/G)$.
  The Hodge-filtered de Rham complex is then encoded by functions $\OO(\Loop_{\GGa}^G(X))$, equipped with the induced $\GGm$-equivariant and $\OO(B\check{\GG}_a) \simeq k[u]$-module structures.
  The de Rham complex itself is encoded by inverting $u$ (the Euler periodization) and taking the weight $0$ part.
\end{itemize}

\subsection{Multiplicative case: periodic cyclic homology}
\label{ssec:TP/KU}

\subsubsection{}

Let $A = \GGm$ over a connective $E_\infty$-ring $k$.
For an affine fppf group scheme $G/k$ and $X \in \Var^G$, the projection $\Loop_A^G(X) \to \Loop_A^G(*)$ is identified by Theorem~\ref{thm:LGa LGm} with the canonical map
\[ \pi : L(\XX) \simeq \XX \times_{\XX \times \XX} \XX \to BG \times_{BG \times BG} BG \simeq G/G, \]
where $L(\XX)$ is the $S^1$-loop space of $\XX := X/G$.
This map is equivariant for the loop rotation action of $\check{\GG}_m \simeq S^1$.

The category $\Rot^G(\GGm)$ of $G$-equivariant multiplicative rotation modules is identified with $\INDPERF^{S^1}(G/G)$.
The $G$-equivariant $\GGm$-loop cohomology of $X$ (Definition~\ref{def:HH_A^G}) is correspondingly identified with the ind-perfect complex
\[ \HH_{\GGm}^{G}(X) \simeq \pi_*(\OO_{L(\XX)}) \quad \text{in} \ \INDPERF^{S^1}(G/G). \]
After induction along $G \twoheadrightarrow 1$ (\ref{sssec:res ind}), it recovers the Hochschild homology of $\XX$,\footnote{%
  Formed relative to the base $k$, i.e. $\operatorname{HH}(\XX/k)$; for $k = \mathbf{S}$ this is $\operatorname{THH}(\XX)$.
}
\[ \operatorname{ind}_G^1\HH_{\GGm}^{G}(X) \simeq \OO(L(\XX)) \simeq \operatorname{HH}(\XX), \]
together with its canonical $S^1$-equivariant structure (using \cite[Thm.~1.2(1), Cor.~4.8]{BFN}).
In particular, $\HH_{\GGm}^{G}(X)$ is a ``delocalization'' of $\operatorname{HH}(\XX)$ in the sense of \cite{BaumBrylinskiMacPherson}.

Periodization with respect to the Euler map $\alpha_{\GGm}$ yields the $G$-equivariant periodic $\GGm$-loop cohomology of $X$,
\[ \Hcoh_{\GGm}^{G}(X) := \HH_{\GGm}^{G}(X)[\alpha_{\GGm}^{-1}] \in \Rot^{G,\mathrm{per}}(\GGm), \]
satisfying
\[ \operatorname{ind}_G^1\Hcoh_{\GGm}^{G}(X) \simeq \operatorname{HH}(\XX)[\alpha_{\GGm}^{-1}] \in \Rot^{\mathrm{per}}(\GGm). \]

Recall from (\ref{sssec:OO(Gm[-1])}) that $\R_{\GGm}$ is identified with the group ring $k[S^1]$.
Consider the localization
\begin{multline*}
  \Ren_{B\check{\GG}_m} : \Rot(\GGm) \simeq \INDPERF(B\check{\GG}_m) \\
  \to \QCOH(B\check{\GG}_m) \simeq \operatorname{Fun}(BS^1, \MOD_{k})
\end{multline*}
of Subsect.~\ref{ssec:indperf}.
By Corollary~\ref{rights}, $\operatorname{ind}_G^1\HH_{\GGm}^{G}(X) \simeq \Ren^R_{B\check{\GG}_m}(\operatorname{HH}(\XX))$, and the periodic $\GGm$-loop cohomology $\Hcoh_{\GGm}^{G}(X)$ recovers some familiar invariants.

For example, the Tate construction of $\operatorname{ind}_G^1\Hcoh_{\GGm}^{G}(X)$, which agrees with that of $\operatorname{ind}_G^1\HH_{\GGm}^{G}(X)$,
\[
  (\operatorname{ind}_G^1\Hcoh_{\GGm}^{G}(X))^{t\check{\GG}_m} \simeq \bigl(\Ren^R_{B\check{\GG}_m}\operatorname{HH}(\XX)\bigr)^{t\check{\GG}_m} \simeq \operatorname{HH}(\XX)^{tS^1}
\]
recovers periodic cyclic homology $\operatorname{HP}(\XX) := \operatorname{HH}(\XX)^{tS^1}$, by Theorem~\ref{thm:Tate vs invariants of Bott}(ii) and Example~\ref{ex:Tate multiplicative}.
For $k = \mathbf{S}$ this is topological periodic cyclic homology $\operatorname{TP}(\XX)$, and for $k$ discrete it agrees with the classical definition via mixed complexes (see \cite[Thm.~2.1]{HoyoisCirc}).

Taking instead Bott-periodized coinvariants (\ref{sssec:Bott periodization}), for $X = \spc(k)$ and $G$ trivial, gives
\[
  \K_{\GGm} :=
    \HH_{\GGm}(*)_{h\check{\GG}_m}[\beta_{\GGm}^{-1}]
    \simeq \lam_{\GGm}[\beta_{\GGm}^{-1}]
    \simeq k \otimes_{\mathbf{S}} \mathrm{KU},
\]
where the last identification comes from Snaith's theorem (see Theorem~\ref{thm:Snaith}).

For $G$ diagonalizable, the equivariant refinement $\K_{\GGm}^G \simeq \OO_{\frX_{\GGm}(G)} \otimes \lam_{\GGm}[\beta_{\GGm}^{-1}]$ of Subsect.~\ref{ssec:HH on XAG} lives on $\frX_{\GGm}(G) = \HOM(X^*(G), \GGm) \simeq G$ (see Subsect.~\ref{ssec:LA(BG) diag}).
Since $G = \spc(k[X^*(G)])$ is affine, $\K_{\GGm}^G$ is determined by its global sections
\[
  \Gamma(G, \K_{\GGm}^G)
  \simeq k[X^*(G)] \otimes_k \lam_{\GGm}[\beta_{\GGm}^{-1}]
  \simeq k \otimes_{\mathbf{S}} \mathrm{KU}[X^*(G)].
\]
For $k = \mathbf{S}$, the group ring $\mathrm{KU}[X^*(G)]$ is the $E_\infty$-ring of $\mathsf{G}$-equivariant complex K-theory of the point, where $\mathsf{G}$ is the compact abelian Lie group with character group $X^*(G)$ (see \cite[\S 4]{GepMei}). Its homotopy ring is $R(\mathsf{G})[\beta_{\GGm}^{\pm 1}]$, where $R(\mathsf{G}) \simeq \bZ[X^*(G)]$ is the representation ring.

\subsection{Loop cohomology on $\frX_A(G)$}
\label{ssec:HH on XAG}

\subsubsection{}

Suppose that $G$ is a diagonalizable group scheme over $k$.
By definition, the $G$-equivariant $A$-loop cohomology of the point is the $E_\infty$-algebra
\[
  \HH_A^G(\ast)
  =
  \OO_{\Loop_A(BG)}
  \in
  \CAlg\bigl(\Rot^G(A)\bigr).
\]
Applying the weight-zero functor of (\ref{sssec:Shv wt=0}), we get the commutative algebra
\[
  \HH_A^G(\ast)^{\mathrm{wt}=0} \simeq
  \OO_{\frX_A(G)}
  \in \CAlg\Bigl( \INDPERF(\frX_A(G) \times B\check{A}) \Bigr)
\]
by Proposition~\ref{prop:eqper trivialization}.

\subsubsection{}

Since $\HH_A^G(\ast)^{\mathrm{wt}=0}$ lies in the essential image of the pullback functor $\mathsf{p}_{\frX}^*$, where $\mathsf{p}_{\frX} : \frX_A(G) \times B\check{A} \to \frX_A(G)$ is the base change of $\sfp : B\check{A} \to *$, we can form its \emph{Bott-periodized coinvariants} in the sense of (\ref{sssec:Bott periodization}).
This is a commutative algebra in $\QCOH(\frX_A(G))$, which we denote
\[
  \K_A^G :=
  \bigl(\HH_A^G(\ast)^{\mathrm{wt}=0}\bigr)_{h\check{A}}[\beta_A^{-1}]
  \simeq \OO_{\frX_A(G)} \otimes \lam_A[\beta_A^{-1}].
\]
This is the Bott-periodized double rotation algebra $\K_A = \lam_A[\beta_A^{-1}] \in \CAlg_k$ (\ref{sssec:Bott periodization}), pulled back along $\frX_A(G) \to \spc(k)$.

\subsubsection{}

Given a closed subgroup $H \subset G$ defined by a surjective homomorphism $X^*(G) = M \twoheadrightarrow N$, there is an induced morphism $i_H : \frX_A(H) \to \frX_A(G)$.
This is a base change of the zero section of $\frX_A(K)$, where $K$ is the diagonalizable group scheme whose character group is the kernel of $M \twoheadrightarrow N$.
Since $\frX_A(K)$ is separated (see \ref{sssec:XAG}), $i_H$ is a closed immersion.

The pushforward of $\K_A^H$ along $i_{H}$ is again a commutative algebra in $\QCOH(\frX_A(G))$,
\begin{equation}\label{eq:KAG}
  \K_A^{G,H} :=
  i_{H,*} (\K_A^H) \simeq
  i_{H,*} (\OO_{\frX_A(H)}) \otimes \lam_A[\beta_A^{-1}],
\end{equation}
where the isomorphism is the projection formula.

\subsubsection{}

We now explain the relationship between this construction and equivariant elliptic cohomology in the sense of Lurie and Gepner--Meier \cite{LurieEllSurv,GepMei}.

Write $\mathsf{G}$ for the compact abelian Lie group with character group $X^*(G)$, and similarly $\mathsf{H} \subset \mathsf{G}$ for a subgroup $H \subset G$ defined by a quotient of $X^*(G)$.

\begin{theorem}[Lurie, Gepner--Meier]\label{thm:KAG functoriality}
  Let $A=\GGa$ (with $k$ a $\bQ$-algebra), $A=\GGm$, or $A=E$ an elliptic curve over $k$.
  Suppose given a nondegenerate periodicity datum for $A$ over an $E_\infty$-ring $R$, classified by a homomorphism $\lam_A[\beta_A^{-1}] \to R$ (see Definition~\ref{def:periodicity datum}, Proposition~\ref{prop:bott per univ prop}).
  The assignment
  \[
    \mathsf{G}/\mathsf{H} \mapsto \K_A^{G,H} \otimes_{\lam_A[\beta_A^{-1}]} R
  \]
  determines a functor $\mathrm{Orb}_{\mathsf{G}}^{\mathrm{op}} \to \QCOH(\frX_{A_R}(G_R))$.
  Moreover, it extends uniquely to a limit-preserving functor
  \begin{equation}\label{eq:GM cohomology}
    \Kcoh_A^G : (\mathrm{Spc}^{\mathsf{G}})^{\mathrm{op}} \to \QCOH(\frX_{A_R}(G_R))
  \end{equation}
  on the category of $\mathsf{G}$-spaces.
\end{theorem}

We regard \eqref{eq:GM cohomology} as a cohomology theory for $\mathsf{G}$-spaces.
For $A=E$, it is precisely $\mathsf{G}$-equivariant elliptic cohomology in the sense of \cite{LurieEllSurv,GepMei}.
The category $\mathrm{Spc}^{\mathsf{G}}$ of \emph{$\mathsf{G}$-spaces} is equivalent to presheaves on $\mathrm{Orb}_{\mathsf{G}}$ (Elmendorf's theorem, cf. \cite[Constr.~6.4]{GepMei}), so \eqref{eq:GM cohomology} is equivalent to its restriction to $\mathrm{Orb}_{\mathsf{G}}$ by the universal property of presheaves.
The \emph{orbit category} $\mathrm{Orb}_{\mathsf{G}}$ has objects labelled $\mathsf{G}/\mathsf{H}$ and, since $\mathsf{G}$ is abelian, the mapping groupoid
\[
  \Map(\mathsf{G}/\mathsf{H}, \mathsf{G}/\mathsf{K})
\]
is the homotopy type of $\mathsf{G}/\mathsf{K}$ when $\mathsf{H} \subset \mathsf{K}$, and is empty otherwise.

\subsubsection{Proof of Theorem~\ref{thm:KAG functoriality}}

For $A=E$ an elliptic curve (and $R$ having noetherian connective cover), this is \cite[Constr.~6.4]{GepMei} (following \cite{LurieEllSurv}).\footnote{%
  Indeed, \cite[Ex.~6.2]{GepMei} gives a natural identification of $\frX_{E_R}(H_R) \simeq \HOM(N, E_R)$ with $\mathrm{Ell}(B\mathsf{H})$, hence identifies the functor we construct here with \cite[Constr.~6.4]{GepMei}.
}
The following is a variation on their argument.

Let $\mathcal{C}$ denote the category of prestacks over $R_0$, i.e. functors from the category $\CAlg_{R_0}$ of connective $E_\infty$-$R_0$-algebras to the category of groupoids.
Then $A_{R_0}$ defines a commutative group object of $\mathcal{C}$, and the given periodicity datum for $A$ is a preorientation of $A_{R_0}$ in the sense of \cite[Def.~7.2.1]{LurieEll2} and \cite[Def.~3.6]{GepMei}.
By \cite[Prop.~2.16, Constr.~3.13, Prop.~3.15]{GepMei} and (\ref{sssec:XAG}), this data determines a functor
\[
  \mathrm{Orb}_{\mathsf{G}} \to \mathcal{C}_{/\frX_{A_{R_0}}(G_{R_0})},
  \quad
  \mathsf{G}/\mathsf{H} \mapsto \bigl(i_H : \frX_{A_{R_0}}(H_{R_0}) \to \frX_{A_{R_0}}(G_{R_0})\bigr).
\]
This takes values in affine morphisms, and by applying the functor $f \mapsto f_*(\OO)$ and then base change along $\Specnc(R) \to \spc(R_0)$, further gives rise to a functor $\mathrm{Orb}_{\mathsf{G}} \to \QCOH(\frX_{A_R}(G_R))^{\mathrm{op}}$.
By \eqref{eq:KAG} and the base change formula for $i_{H,*}$, it is given on objects by the assignment
\[
  \mathsf{G}/\mathsf{H}
  \;\mapsto\;
  i_{H_R,*}\OO_{\frX_{A_R}(H_R)} \simeq \K_A^{G,H} \otimes_{\lam_A[\beta_A^{-1}]} R.
  \qed
\]

\begin{theorem}[Gepner--Meier]\label{thm:KAG genuine}
  Let the notation be as in Theorem~\ref{thm:KAG functoriality} and assume that $R$ has noetherian connective cover.
	If $A=E$ is an elliptic curve and $G=\GGm$ (so that $\mathsf{G}=U(1)$), then the functor \eqref{eq:GM cohomology} extends to a limit-preserving functor
  \[ (\operatorname{Spt}^{\mathsf{G}})^{\mathrm{op}} \to \QCOH(\frX_{E_R}(G_R)) \]
  on the category of genuine $\mathsf{G}$-spectra.
\end{theorem}

\begin{proof}
  The given nondegenerate periodicity datum is an orientation of $E_R$ in the sense of \cite[Def.~7.2.7]{LurieEll2} by Lemma~\ref{lem:Bott vs Lurie}, hence in the sense of \cite[Def.~5.17]{GepMei} by \cite[Prop.~4.3.23]{LurieEll2}.
  Identifying \eqref{eq:GM cohomology} with $\mathrm{Ell}_{\mathsf{G}}$ as in the proof of Theorem~\ref{thm:KAG functoriality}, the claim is \cite[Prop.~9.2]{GepMei} for compact $\mathsf{G}$-spectra; this is ind-extended to all genuine $\mathsf{G}$-spectra as in \cite[Constr.~9.3]{GepMei}.
\end{proof}

\begin{remark}\label{rem:GM general base}
  As with the rest of the material in this section, the constructions and arguments above go through over a general base $S$, $A \to S$ a relative $1$-dimensional commutative group scheme, and so on.
  The nondegenerate periodicity data for $A$ over the nonconnective $E_\infty$-algebra $R$ (i.e. the homomorphism $\lam_{A/S}[\beta_{A/S}^{-1}] \to R$) can be replaced by a nonconnective spectral DM stack\footnote{%
    satisfying the mild finiteness conditions in \cite[\S 1, Conventions]{GepMei}
  } $\mathcal{P}$ over $\Specnc_S(\lam_{A/S}[\beta_{A/S}^{-1}])$.
  For example, we may take the universal case of $S = \Mell$, $A$ the universal elliptic curve over $\Mell$, and $\mathcal{P} = \Mell^{\mathrm{per},\mathrm{nd}}$ the moduli stack of elliptic curves with nondegenerate periodicity data (\ref{sssec:proof of thm:Otop}), i.e. Lurie's moduli stack $\mathcal{M}_{\mathrm{Ell}}^{\mathrm{or}}$ of oriented elliptic curves \cite{LurieEll2}.
\end{remark}

\begin{example}\label{exam:KAG G=1}
  Let $G = 1$ be the trivial group, so that $\frX_A(G) = \spc(k)$ and $\K_A \simeq \lam_A[\beta_A^{-1}]$, and take the universal $R=\lam_A[\beta_A^{-1}]$.
  Then $\Kcoh_A$ \eqref{eq:GM cohomology} is the unique limit-preserving functor $\operatorname{Spc}^{\mathrm{op}} \to \Mod_{\K_A}$ sending the point to $\K_A$.
  For $A = \GGa$ over a $\bQ$-algebra $k$, $\K_{\GGa} \simeq k[\beta^{\pm 1}]$ (\ref{sssec:Bott classical}); for $A = \GGm$ over $\mathbf{S}$, $\K_{\GGm} \simeq \mathrm{KU}$ (Theorem~\ref{thm:Snaith}).

  For $A=\mathcal{E}$ the universal elliptic curve over $S=\Mell$ as in the remark above, $R$ is replaced by $\mathcal{P} = \Mell^{\mathrm{per},\mathrm{nd}}$, so that $\K_{\mathcal{E}/\Mell} \simeq \OO_{\mathcal{P}} = \lam_{\mathcal{E}/\Mell}[\beta_{\mathcal{E}/\Mell}^{-1}]$; taking its global sections over the classical truncation $\Mell^{\heartsuit}$ yields the $E_\infty$-ring spectrum $\mathrm{TMF}$ by Theorem~\ref{thm:Otop}.
\end{example}

\section{Formal structure of loop cohomology}\label{sec:coh formal}

As in \S \ref{sec:cohomology}, we fix a $1$-shifted Cartier dual pair $(A, \check{A}, \jmath)$ and an affine fppf group scheme $G$ over $k$.
We establish some of the formal structure of loop cohomology, including pullbacks, cup products, the induction axiom, and functoriality in $A$.
We also construct change-of-groups isomorphisms comparing $\HH_A^G(-)$ with $\HH_A^H(-)$, for a subgroup $H \subset G$.

\subsection{Pullbacks}\label{ssec:coh pullbacks}

Let $f : X \to Y$ be a morphism in $\Stk^G$.
For every coefficient $\cF \in \Shv_A^G(Y)$, the adjunction unit
\[
  \cF \to f_* f^*(\cF)
\]
induces, after pushforward to $\Loop_A^G(*)$, the \emph{pullback map}
\[ f^* : \HH_A^G(Y; \cF) \to \HH_A^G(X; f^*\cF). \]
We also write $f^* : \Hcoh_A^G(Y; \cF) \to \Hcoh_A^G(X; f^*\cF)$ for the induced map on periodizations.

In particular, $\HH_A^G(-)$ determines a functor
\[\HH_A^G : (\Stk^G)^{\mathrm{op}} \longrightarrow \Rot^G(A),\]
and similarly for $\Hcoh_A^G(-)$.

\subsection{Products}\label{ssec:coh products}

Let $X \in \Stk^G$ and $\cF_1, \cF_2 \in \Shv_A^G(X)$.
Writing $a : X \to *$ for the structure morphism, the lax symmetric monoidal structure on $a_* : \Shv_A^G(X) \to \Shv_A^G(*) = \Rot^G(A)$ encodes maps
\[
  a_*(\cF_1) \otimes a_*(\cF_2) \to a_*(\cF_1\otimes\cF_2),
\]
which define the \emph{cup product} maps
\[
  \cup : \HH_A^G(X;\cF_1)\otimes \HH_A^G(X;\cF_2) \longrightarrow \HH_A^G(X;\cF_1\otimes\cF_2).
\]
These are compatible with pullbacks (Subsect.~\ref{ssec:coh pullbacks}).
They equip $\HH_A^G(X)$ with the structure of a commutative algebra, and each $\HH_A^G(X;\cF)$ with the structure of an $\HH_A^G(X)$-module.

Since periodization is symmetric monoidal, this structure passes to $\Hcoh_A^G(-)$.

\subsection{Induction}\label{ssec:coh induction}

Let $H \subset G$ be a normal fppf subgroup scheme such that $G/H$ is an affine fppf group scheme.
Denote by $\operatorname{ind}_G^{G/H} : \Rot^G(A) \to \Rot^{G/H}(A)$ the functor of induction along $G \twoheadrightarrow G/H$ (\ref{sssec:res ind}).

\begin{proposition}\label{prop:coh induction}
For every $X \in \Stk^G$, there is a canonical isomorphism
\[ \HH_A^{G/H}(X/H) \simeq \operatorname{ind}^{G/H}_G \HH_A^G(X) \]
in $\Rot^{G/H}(A)$.
Similarly, $\Hcoh_A^{G/H}(X/H) \simeq \operatorname{ind}^{G/H}_G \Hcoh_A^G(X)$.
\end{proposition}
\begin{proof}
The structural morphism $a_{X/H} : X/H \to *$ is $(G/H)$-equivariant, and the induced map on quotient stacks factors as
\[
  (X/H)/(G/H) \simeq
  X/G \to
  BG \to
  B(G/H).
\]
The claim now follows by definition of the induction functor, and functoriality of pushforwards in $\Shv_A(-)$.
\end{proof}

\begin{example}\label{exam:ind_G^1 HH}
	For every affine fppf group scheme $G$ and every $X \in \Stk^G$, we have
	\[
	  \operatorname{ind}_G^1 \HH_A^G(X) \simeq \HH_A(X/G)
	\]
	in $\Rot(A)$, and similarly for $\Hcoh_A^G(X)$.
\end{example}

\begin{example}\label{exam:ind_H^G HH}
	Let $H \subset G$ be a closed fppf subgroup scheme.
	For every $Y \in \Stk^H$, identifying $(G \times^H Y)/G \simeq Y/H$ and its structure map with $Y/H \to BH \to BG$, we have
	\[
	  \operatorname{ind}^G_H \HH_A^H(Y) \simeq \HH_A^G(G \times^H Y)
	\]
	in $\Rot^G(A)$, and similarly for $\Hcoh_A^H(Y)$.
\end{example}

\subsection{Change of groups}

Let $k$ be a $\bQ$-algebra, $A$ as in Corollary~\ref{mapsrepr}, $G$ a reductive group scheme over $k$, and $H \subset G$ a closed fppf subgroup scheme.

\begin{theorem}\label{thm:change of groups}
Let $X \in \Var^G$ and let $H\hookrightarrow G$ as above.
\begin{enumerate}
\item If $G/H$ is smooth and $G$-projective, then the canonical map $$\operatorname{res}^G_H\HH_A^G(X)\longrightarrow \HH_A^H(X)$$ is an isomorphism.

\item If $H$ is reductive or satisfies the hypothesis of Theorem~\ref{thm:pluperfect quotient}(i), and $X$ is smooth and proper, then the canonical map $$\operatorname{res}^G_H\HH_A^G(X)\longrightarrow \HH^H_A(X)$$ is an isomorphism.
\end{enumerate}
Moreover, the same statements hold for periodic $A$-loop cohomology $\Hcoh_A^G(-)$. \end{theorem}

\begin{proof}
Consider the cartesian square
$$\begin{tikzcd}
X/H\arrow{r}{r'} \arrow[swap]{d}{a'} & X/G\arrow{d}{a} \\
BH\arrow{r}{r} & BG.
\end{tikzcd}$$
In both cases the canonical map is the exchange transformation
\[
  \operatorname{Ex}^*_* : r^* a_*(\un) \to a'_* r'^*(\un) \simeq a'_*(\un),
\]
which is by definition the corresponding exchange transformation in $\INDPERF$ for the induced square of loop rotation quotients.
That the latter is cartesian follows by Theorem~\ref{Kun}, since $X/G$, $BG$ and $BH$ are perfect fppf-locally on $\spc(k)$: the first since $X \in \Var^G$, the second by Theorem~\ref{thm:pluperfect quotient}(ii), and the third by Lemma~\ref{lem:qproj implies variety} applied to $G/H$ in case (i), resp. by Theorem~\ref{thm:pluperfect quotient} applied to $H$ in case (ii).

(i) We can write $\Loop_A(r) = \Loop_A^G(G/H \to *)$, where $G/H \in \Var^G$ by Lemma~\ref{lem:qproj implies variety}.
Since $G/H$ is smooth and proper, $\Loop_A(r)/\check{A}$ is proper, representable and quasi-smooth (Corollary~\ref{mapsrepr} and Theorem~\ref{reltangent}), and $\Loop_A(a)/\check{A}$ is qcqs and representable (Corollary~\ref{mapsrepr}).
The desired base change isomorphism is thus Theorem~\ref{thm:bcindperf}(i).

(ii) Since $X$ is smooth proper, $\Loop_A(a)/\check{A}$ is proper, representable and quasi-smooth (Corollary~\ref{mapsrepr} and Theorem~\ref{reltangent}).
The desired base change isomorphism is thus Theorem~\ref{thm:bcindperf}(ii).
\end{proof}

\subsection{Functoriality in $A$}\label{ssec:HH func A}

Let $h : (A,\check{A},\jmath) \to (A',\check{A}',\jmath')$ be a morphism of $1$-shifted Cartier dual pairs, and recall the functor $\operatorname{res}_A^{A'} : \Rot^G(A') \to \Rot^G(A)$ of (\ref{sssec:Shv func A}).

\begin{proposition}\label{prop:HH func A}
  Assume that the map $q : B\check{A}' \to B\check{A}$ induces an isomorphism $\OO_{B\check{A}} \simeq q_*(\OO_{B\check{A}'})$ in $\QCOH(B\check{A})$ after any base change.
  Then for every $X \in \Stk^G$, there is a canonical map in $\Rot^G(A)$,
  \[
    \operatorname{res}_A^{A'} \HH_{A'}^G(X) \to \HH_{A}^G(X),
  \]
  natural with respect to pullbacks in $X$.
\end{proposition}

\begin{proof}
Consider the diagram of commutative squares
\[\begin{tikzcd}
  \Loop_A^G(X)/\check{A} \ar[swap]{d}{\pi}
  & \Loop_A^G(X)/\check{A}' \ar{d}{\tau}\ar[swap]{l}{\beta_X}\ar{r}{\alpha_X}
  & \Loop_{A'}^G(X)/\check{A}' \ar{d}{\pi'}
  \\
  \Loop_A^G(*)/\check{A}
  & \Loop_A^G(*)/\check{A}' \ar[swap]{l}{\beta}\ar{r}{\alpha}
  & \Loop_{A'}^G(*)/\check{A}',
\end{tikzcd}\]
where we write $\pi$ and $\tau$ for the quotients of $\pi_X^{A,G}$ by $\check{A}$ and $\check{A}'$, respectively, and $\pi'$ for the quotient of $\pi_X^{A',G}$ by $\check{A}'$ (see Definition~\ref{def:Loop_A^G} and (\ref{sssec:Shv func A}) for the notation).
Since $\HH_{A'}^G(X) := \pi'^{\INDPERF}_*(\OO_{\Loop_{A'}^G(X)/\check{A}'})$, the counit of $\pi'$ gives a canonical map
\[
  \tau^*_{\INDPERF} \alpha^*_{\INDPERF} \HH_{A'}^G(X) \simeq
  (\alpha_X)^*_{\INDPERF} \pi'^{*}_{\INDPERF} \HH_{A'}^G(X) \xrightarrow{\mathrm{counit}}
  \OO_{\Loop_A^G(X)/\check{A}'},
\]
which by adjunction corresponds to a map
\[
  \alpha^*_{\INDPERF} \HH_{A'}^G(X) \longrightarrow
  \tau_*^{\INDPERF} \OO_{\Loop_A^G(X)/\check{A}'}.
\]
By the assumption, the unit $\OO \to (\beta_X)_*\OO$ is invertible in $\QCOH$, and hence by Corollary~\ref{rights} in $\INDPERF$.
Applying $\beta_*^{\INDPERF}$ to the preceding map therefore gives
\[\begin{aligned}
  \operatorname{res}_{A}^{A'} :
  \beta_*^{\INDPERF}\alpha^*_{\INDPERF} \bigl( \HH_{A'}^G(X) \bigr)
  &\longrightarrow \beta_*^{\INDPERF}\tau_*^{\INDPERF} (\OO) \\
  &\simeq \pi_*^{\INDPERF}(\beta_X)_*^{\INDPERF} (\OO)
  \simeq \HH_A^G(X).
\end{aligned}\]
This is natural in $X$ by construction.
\end{proof}

\begin{remark}
With $q : B\check{A}' \to B\check{A}$ as above, and $\operatorname{ind}_G^1 : \Rot^G(-) \to \Rot(-)$ induction along $G \twoheadrightarrow 1$, Proposition~\ref{prop:HH func A} yields the map
\[
  q_* \operatorname{ind}_G^1 \HH_{A'}^G(X)
  \xrightarrow{\eqref{eq:push to res A A'}}
  \operatorname{ind}_G^1 \operatorname{res}_A^{A'} \HH_{A'}^G(X)
  \xrightarrow{\ref{prop:HH func A}} \operatorname{ind}_G^1 \HH_{A}^G(X)
\]
in $\Rot(A)$, which by Example~\ref{exam:ind_G^1 HH} is equivalently a canonical map
\[
  q_* \HH_{A'}(X/G) \to \HH_A(X/G)
\]
in $\Rot(A)$.
\end{remark}

\begin{example}[Adams operations]\label{exam:Adams}
Let $n\geq 1$ be a positive integer invertible in the base $k$.
Multiplication by $n$ defines a morphism of $1$-shifted Cartier dual pairs $([n],[n]) : (A, \check{A}) \to (A, \check{A})$.
For $A=\GGa$, multiplication by $n$ on $\check{A}$ is an isomorphism.
For $A=\GGm$ and $A=E$, the map $B[n]:B\check{A}\to B\check{A}$ is a gerbe banded by $\mathbf{Z}/n\mathbf{Z}$ and $\check{E}[n]$, respectively.
These groups are finite étale of invertible order, so the unit $\OO_{B\check{A}} \to B[n]_*\OO_{B\check{A}}$ is invertible in $\QCOH$ after every base change.
We refer to the map given by Proposition~\ref{prop:HH func A} in this case as the $n$th \emph{Adams operation}.
\end{example}

\section{Properties of loop cohomology}\label{sec:coh properties}

\subsection{Mayer--Vietoris and descent}\label{ssec:coh descent}

\begin{theorem}\label{thm:coh MV}
Let $X \in \Var^G$ and let $\{U_i\}_{i\in I}$ be a finite cover of $X$ by $G$-invariant quasi-compact opens.
Then, writing $U_J := \bigcap_{j\in J} U_j$, we have:
$$\HH_A^G(X)\xrightarrow{\sim} \lim_{\emptyset\neq J\subseteq I}\HH_A^G(U_J).$$
More generally, $\HH_A^G(-)$ satisfies descent for finite Nisnevich covers $\{U_i \to X\}_{i\in I}$ in $\Var^G$.
Moreover, the same statements hold for $\Hcoh_A^G(-)$.
\end{theorem}

\begin{proof}
By Proposition~\ref{eff epi}(i), the $\Loop_A^G(U_i)$ form an open cover of $\Loop_A^G(X)$, with intersections $\bigcap_{j\in J}\Loop_A^G(U_j) \simeq \Loop_A^G(U_J)$ for every $\emptyset\neq J\subseteq I$, by Corollary~\ref{cor:Kun open}.
The same holds on loop rotation quotients.
By Zariski descent for $\QCOH$, we obtain
$$ (\pi_X^G)_*^{\QCOH}\OO \simeq \lim_{\emptyset\neq J\subseteq I}(\pi_{U_J}^G)_*^{\QCOH}\OO $$
in $\QCOH^{\check{A}}(\Loop_A^G(*))$.
Applying the functor $\Ren^R$ of Subsect.~\ref{ssec:indperf}, which preserves limits, acts as the identity on perfect complexes, and commutes with pushforwards (see Corollary~\ref{rights}), this gives
$$ \HH_A^G(X)\simeq \lim_{\emptyset\neq J\subseteq I}\HH_A^G(U_J). $$
For Nisnevich descent one argues similarly, applying Theorem~\ref{Niscosh} to the Nisnevich cover $\{U_i/G\to X/G\}$.
\end{proof}

\subsection{K\"unneth formulas}\label{ssec:coh Kunneth}

\begin{theorem}\label{thm:H Kun}
Suppose $k$ is a $\bQ$-algebra and $G$ is reductive.
\begin{enumerate}
\item \label{item:IntKun} \emph{Interior K\"unneth formula}.
Suppose $X \in \Var^G$ is smooth and proper.
Then for every $Y \in \Var^G$, we have:
$$\HH_A^G(X\times Y)\simeq \HH_A^G(X) \otimes \HH_A^G(Y).$$

\item \label{item:ExtKun} \emph{Exterior K\"unneth formula}.
Let $H$ be a reductive group scheme over $k$, and let $X \in \Var^G$ and $Y \in \Var^H$ be smooth and proper.
Then the canonical map
\[
  \HH_A^G(X)\boxtimes\HH_A^H(Y) \longrightarrow \HH_A^{G\times H}(X\times Y)
\]
is an isomorphism.
\end{enumerate}
Moreover, the same statements hold for $\Hcoh_A^G(-)$.
\end{theorem}

\begin{proof}(i) By Theorem~\ref{Kun}, we have the cartesian square
$$\begin{tikzcd}
\Loop_A^G(X\times Y)/\check{A}
  \arrow{r}{\widetilde{\pi}_X}
  \arrow{d}{\widetilde{\pi}_Y}
  & \Loop^G_A(Y)/\check{A}\arrow{d}{\pi_Y} \\
\Loop^G_A(X)/\check{A}\arrow{r}{\pi_X}
  & \Loop^G_A(*)/\check{A}.
\end{tikzcd}$$
Since $X$ is smooth and proper, $\pi_X$ is proper, representable and quasi-smooth (Corollary~\ref{mapsrepr} and Theorem~\ref{reltangent}).
Using the projection formula (Lemma~\ref{lem:indperf push easy}) and base change (Theorem~\ref{AppMain}(ii)), we have:
\begin{align*}
\HH_A^G(X)\otimes\HH_A^G(Y)
&\simeq \pi_{Y,*}^{\INDPERF}(\pi_Y)^*_{\INDPERF}\pi_{X,*}^{\INDPERF}\OO \\
&\simeq \pi_{Y,*}^{\INDPERF}\widetilde{\pi}_{X,*}^{\INDPERF}(\widetilde{\pi}_Y)^*_{\INDPERF}\OO \\
&\simeq \pi_{X\times Y,*}^{\INDPERF}\OO
\simeq \HH_A^G(X\times Y).
\end{align*}

(ii) By Theorem~\ref{Kun}, we have
$$\Loop_A^{G\times H}(X\times Y) \simeq \Loop_A^G(X)\times\Loop_A^H(Y),$$
compatibly with loop rotation.
Under this identification it will suffice to show the identity
\[
  (\pi_X^G \times \pi_Y^H)_*(\OO) \simeq \pi_{X,*}^G(\OO) \boxtimes \pi_{Y,*}^H(\OO)
\]
in $\INDPERF$.
The structure maps $\pi_X^G$, $\pi_Y^H$ and $\pi_{X\times Y}^{G \times H}$ are proper, representable and quasi-smooth (Corollary~\ref{mapsrepr} and Theorem~\ref{reltangent}).
Hence by Lemma~\ref{lem:pushindperf}, the claim follows from the analogous identity in $\QCOH$.
\end{proof}

\subsection{Homotopy invariance}\label{ssec:homotopy invariance}

In characteristic zero, we show that $\Hcoh_A^G(-)$ is homotopy invariant.
Note that in positive characteristic, homotopy invariance fails even for $A=\GGm$ and $G$ trivial (see Remark~\ref{rem:HI char p}).

\begin{theorem}\label{thm:homotopy invariance}
Suppose $k$ is a $\mathbf{Q}$-algebra.
If $X$ is a qcqs scheme with an action of an affine algebraic group $G$, and $p : V \to X$ is the total space of a $G$-equivariant vector bundle, then the map $$p^* : \Hcoh_A^G(X)\rightarrow \Hcoh_A^G(V)$$ is invertible.
\end{theorem}

\begin{remark}
In fact, we will see that the statement holds more generally when $p : Y \to X$ is a morphism in $\Var^G$ with a section $s$, and there exists a morphism $h : Y \times \aff \to Y$ over $X$ with $h \circ \iota_1 \simeq \id$ and $h \circ \iota_0 \simeq s \circ p$, where $\iota_t$ denotes the inclusion at $t \in \aff$.
\end{remark}

\begin{remark}\label{rem:HI char p}
Theorem~\ref{thm:homotopy invariance} does not hold outside characteristic zero.
For example, take $A=\GGm$ and $G$ trivial over $k=\mathbf{F}_p$.
Subsect.~\ref{ssec:TP/KU} identifies $\HH_{\GGm}(*) \to \HH_{\GGm}(\aff)$ with the inclusion of the first summand $k$ in
\[
  \HH_{\GGm}(\aff)
  \simeq \operatorname{THH}(k[t]/k)
  \simeq k\otimes_{\mathbf{S}}\operatorname{THH}(\mathbf{S}[t])
  \simeq k\oplus\bigoplus_{n\ge1}
  k\otimes\Sigma^\infty_+(S^1/C_n),
\]
where the last isomorphism is \cite[Prop.~11.3]{BMS2}.
We claim that its cofibre does not vanish after Euler periodization (which commutes with cofibres and direct summands).
By Theorem~\ref{thm:Tate vs invariants of Bott}, it suffices to show that the Tate construction $(M_p)^{t\check{\GG}_m}$ of the weight $p$ component $M_p := k\otimes\Sigma^\infty_+(S^1/C_p)$ is nonzero.
Here $(M_p)^{t\check{\GG}_m}$ agrees with the $S^1$-Tate construction by Example~\ref{ex:Tate multiplicative} (since $M_p$ has perfect underlying $k$-module $k\oplus k[1]$), which by \cite[Thm.~4.10]{GreenleesMay} and \cite[Prop.~3.8]{GreenleesMayTate} is $(M_p)^{tS^1} \simeq k^{tC_p}[1]$.
This is nonzero: $\pi_0(k^{tC_p})$ is the Tate cohomology group $\widehat{H}^0(C_p; k) \simeq k/pk \simeq \mathbf{F}_p$.
\end{remark}

\subsubsection{Homotopy invariance in $\Shv_A$}

\begin{lemma}\label{lem:a1 weight splitting}
Suppose $k$ is a $\mathbf{Q}$-algebra.
Then the structure map $\pi : \Loop_A(\aff)/\check{A} \to B\check{A}$ is affine, and the unit $\OO_{B\check{A}} \to \pi^{\QCOH}_*(\OO)$ induces a canonical decomposition
\[
  \pi^{\QCOH}_*(\OO) \simeq \OO_{B\check{A}} \oplus P
\]
in $\QCOH(B\check{A})$.
Moreover, the map $\id_{P} \otimes \alpha_A : P \to P \otimes \OO_{B\check{A}}\Thom{1}$ in $\QCOH(B\check{A})$ is null-homotopic, where $\alpha_A$ is the Euler map (\ref{sssec:Euler map}).
\end{lemma}
\begin{proof}
By Corollary~\ref{cor:hkr0}, $\Loop_A(\aff) \simeq \mathbf{V}_{\aff}(\LL_{\aff} \otimes \frav[1])$ is an affine scheme, and the $\check{A}$-action factors through $\check{A} \to B\mathbf{V}(\fra)$.
Hence $\pi$ fits in the cartesian square
\[\begin{tikzcd}
\Loop_A(\aff)/\check{A} \arrow{r} \arrow[swap]{d}{\pi} & \Loop_A(\aff)/B\mathbf{V}(\fra) \arrow{d} \\
B\check{A} \arrow{r} & B^2\mathbf{V}(\fra).
\end{tikzcd}\]
Since every morphism from an affine $T$ to $B^2\mathbf{V}(\fra)$ is null-homotopic, the base change of $\pi$ along any $T \to B\check{A}$ is $T \times \Loop_A(\aff) \to T$, so $\pi$ is affine.
Since $k$ is a $\mathbf{Q}$-algebra, $\RA \simeq k \oplus \frav[1]$ is a trivial square-zero extension (\ref{sssec:funcheckA char 0}).
Under the equivalence $\QCOH(B\check{A}) \simeq \MOD_{\RA}$ of Theorem~\ref{thm:QCoh checkA}, $\pi^{\QCOH}_*(\OO)$ is then the $\RA$-module $M$ with underlying complex $\OO(\Loop_A(\aff)) \simeq k[x] \oplus k[x]dx \otimes \frav[1]$, on which $\frav[1]$ acts by the de Rham differential (see \cite[Prop.~1.2]{BNloop}).
The scaling action of $\GGm$ on $\aff$ commutes with loop rotation and yields a weight decomposition $M = \bigoplus_{n \ge 0} M_n$ in $\MOD_{\RA}$.
\begin{itemize}
\item The weight zero piece is spanned by $1 \in k[x]$, on which $\frav[1]$ acts by zero, i.e., $M_0 \simeq \eps_*(k)$.
Thus the weight zero piece of $\pi^{\QCOH}_*(\OO)$ is the corresponding object $\OO_{B\check{A}}$ (\ref{sssec:eps_lam}).

\item For $n > 0$, $\frav[1]$ maps $x^n$ to $n\,x^{n-1}dx$ with $n$ invertible, so $M_n$ is the free $\RA$-module on $x^n$.
Thus the weight $n$ piece of $\pi^{\QCOH}_*(\OO)$ is the corresponding object $\sfe_\sharp(k \cdot x^n)$ (where $\sfe_\sharp$ is as in the proof of Theorem~\ref{thm:QCoh checkA}).
\end{itemize}
Now consider the unit map $\OO_{B\check{A}} \to \pi^{\QCOH}_*(\OO)$.
Since it is $\GGm$-equivariant with source of weight zero, it factors through $\OO_{B\check{A}}$.
The resulting endomorphism $\OO_{B\check{A}} \to \OO_{B\check{A}}$ is invertible, since on underlying $k$-modules it is the identity of $k = k \cdot 1$.
Thus the stated decomposition holds with $P := \sfe_\sharp(M_+)$, where $M_+ := \bigoplus_{n>0} k \cdot x^n$.
By the projection formula, $\id_P \otimes \alpha_A \simeq \sfe_\sharp(\id_{M_+} \otimes \sfe^*\alpha_A)$, and $\sfe^*\alpha_A$ is null-homotopic by Remark~\ref{rem:sfe^*beta}.
\end{proof}

\begin{proposition}\label{prop:a1 invariance shv}
Suppose $k$ is a $\mathbf{Q}$-algebra.
Let $\YY$ be a stack and $p : \YY \times \aff \to \YY$ the projection.
For every $F \in \Shv_A(\YY)$, the unit $F \to p_* p^* F$ is invertible in $\Shvper_A(\YY)$.
\end{proposition}
\begin{proof}
By the projection formula (Lemma~\ref{lem:indperf push easy}), and since Euler periodization commutes with $-\otimes F$, we reduce to the case $F=\OO$.
By Theorem~\ref{thm:Kun affine} we have the $\check{A}$-equivariant isomorphism $\Loop_A(\YY \times \aff) \simeq \Loop_A(\YY) \times \Loop_A(\aff)$, so $q := \Loop_A(p)/\check{A} : \Loop_A(\YY \times \aff)/\check{A} \to \Loop_A(\YY)/\check{A}$ is a base change of the map $\pi : \Loop_A(\aff)/\check{A} \to B\check{A}$ of Lemma~\ref{lem:a1 weight splitting}.
The latter, together with the base change formula for $\pi_*$ in $\QCOH$ (see \cite[Prop.~6.3.4.1]{SAG}), identifies the unit $\OO \to q^{\QCOH}_*(\OO)$ with the summand inclusion
\[
   \OO_{\Loop_A(\YY)/\check{A}} \longrightarrow \OO_{\Loop_A(\YY)/\check{A}} \oplus P,
\]
where $P$ is the pullback along $\Loop_A(\YY)/\check{A} \to B\check{A}$ of the second summand in Lemma~\ref{lem:a1 weight splitting}.
Applying $\Ren^R$ (Subsect.~\ref{ssec:indperf}), which preserves finite sums and the object $\OO$, Corollary~\ref{rights} identifies the unit $\OO \to p_*(\OO)$ in $\Shv_A(\YY)$ with the summand inclusion $\OO \to \OO \oplus \Ren^R(P)$.
It thus remains to show that $\Ren^R(P)$ vanishes after Euler periodization.
This follows from the fact that the Euler map $\alpha_A : \Ren^R(P) \to \Ren^R(P)\Thom{1}$ is null-homotopic, since it is identified via Lemma~\ref{perftens} with $\Ren^R$ of the null-homotopic map $\id \otimes \alpha_A : P \to P \otimes \OO\Thom{1}$.
\end{proof}

\begin{theorem}\label{thm:homotopy invariance shv}
Suppose $k$ is a $\mathbf{Q}$-algebra.
Let $a : \YY \to \XX$ be a morphism of stacks with a section $s$.
If there exists a morphism
\[
  h : \YY \times \aff \to \YY,
  \qquad \text{with } h \circ \iota_1 \simeq \id,
  \quad h \circ \iota_0 \simeq s \circ a
\]
over $\XX$, then for every $F \in \Shv_A(\XX)$, the unit $F \to a_*a^*(F)$ is invertible in $\Shvper_A(\XX)$.
\end{theorem}

\begin{proof}
As in the proof of Proposition~\ref{prop:a1 invariance shv}, the projection formula reduces us to the case $F = \OO$.
Proposition~\ref{prop:a1 invariance shv}, applied to $\YY$ and the projection $p : \YY \times \aff \to \YY$, shows that the unit $\OO \to p_*(\OO)$ is invertible after Euler periodization.
Applying $a_*$, which commutes with Euler periodization by the projection formula (Lemma~\ref{lem:indperf push easy}), the induced map $a_*(\OO) \to (a \circ p)_*(\OO)$ is also invertible.

For a stack $\ZZ$ with a morphism $f : \ZZ \to \XX$, write
\[
  \Hcoh(\ZZ) := \Hcoh_A(\ZZ_{/\XX}) := f_*(\OO)[\alpha_A^{-1}]
  \qquad \text{in } \Shvper_A(\XX).
\]
The claim is that $a^* : \Hcoh(\XX) \to \Hcoh(\YY)$ is invertible.
The first paragraph shows that $p^* : \Hcoh(\YY) \to \Hcoh(\YY \times \aff)$ is invertible.
By functoriality, $\iota_0^* \simeq (p^*)^{-1} \simeq \iota_1^*$, hence
\[ \id \simeq \iota_1^* \circ h^* \simeq \iota_0^* \circ h^* \simeq a^*\circ s^*. \]
Since also $s^* \circ a^* \simeq (a \circ s)^* \simeq \id$ by functoriality, the claim follows.
\end{proof}

\subsubsection{Proof of Theorem~\ref{thm:homotopy invariance}}

Apply Theorem~\ref{thm:homotopy invariance shv} with $a$ the projection $p/G : V/G \to X/G$, $s$ the zero section, $h$ the scaling action, and $F=\un$ the unit object.
\qed

\section{Twists, pushforwards, and Euler classes}
\label{sec:pushforwards}

In this section we discuss pushforwards in $A$-loop cohomology.
We fix a linearly reductive group scheme $G$ over $k$ and restrict our attention to $\HH_A^G(-)$ as a functor on $\Var^G$.

For a proper morphism of smooth $G$-varieties, the trace of the relative dualizing complex gives a Gysin push-forward $f_*$ (Subsect.~\ref{ssec:pushforwards}). Applied to the zero section of a vector bundle $V$, it produces the Euler class $e_A(V) : \vartheta^!_A(V) \to \un$, multiplicative under direct sums (Subsect.~\ref{ssec:euler}).

\subsection{Twists}\label{ssec:twists}

Recall the discussion of dualizing complexes in $\Shv_A^G(-)$ in Subsect.~\ref{ssec:sheaves}.

\begin{definition}\label{def:theta^! twist}
Let $X \in \Var^G$.
For every $G$-equivariant vector bundle $V \to X$, we write
\[
  \vartheta^!_A(V) := \omega_s \quad \in \Shv_A^G(X),
\]
for the dualizing complex of the zero section $s : X \to V$ (see Subsect.~\ref{ssec:sheaves}).
\end{definition}

\subsubsection{}

Recall from (\ref{sssec:Theta}) the $\check{A}$-equivariant perfect complexes
\[ \Theta_{X/G}^A(F) \in \PERF^{\check{A}}(\Loop_A^G(X)), \qquad F \in \PERF(X/G). \]
When $k$ is a $\bQ$-algebra\footnote{%
  or, if we were to work in the animated rather than spectral setting
}, we can form its graded determinant,
\[
  \vartheta_A(F) := \det(\Theta_{X/G}^A(F)),
\]
which is an $\check{A}$-equivariant shifted line bundle on $\Loop_A^G(X)$.

\begin{proposition}\label{prop:theta^! vs theta}
	Let $k$ be a $\bQ$-algebra.
	Assume that $A=\GGm$ or $A=E$ is an elliptic curve over $k$.
	\begin{enumerate}
	  \item If $X,Y \in \Var^G$ are smooth, then for every $G$-equivariant proper morphism $f : X \to Y$, there is a canonical isomorphism $\omega_f \simeq \vartheta_A(\TT_{f})^{-1}$.\footnote{%
			 Here $\TT_f$ is regarded with its $G$-equivariant structure, i.e., it stands for the relative tangent complex of $f/G : X/G \to Y/G$.
		}
	  \item If $X \in \Var^G$ is smooth, then for every $G$-equivariant vector bundle $V \to X$, there is a canonical isomorphism $\vartheta^!_A(V) \simeq \vartheta_A(V)$.\footnote{%
			 We tacitly identify $F \in \PERF(X/G)_{[0,0]}$ with the vector bundle $V = \mathbf{V}(F^\vee) = \spc(\Sym^*(F^\vee))$.
		}
	\end{enumerate}
\end{proposition}
\begin{proof}
	The second claim follows from the first.
	When $X$ and $Y$ are smooth, $\phi := \Loop_A^G(f)/\check{A}$ is a morphism between Artin stacks that are quasi-smooth over the common base $\Loop_A^G(*)/\check{A}$ (by Theorem~\ref{reltangent}).
	Therefore, the main result of \cite{cohdual} computes the relative dualizing complex $\omega_\phi^{\QCOH}$ as $\det(\LL_\phi)$.
	In particular, it is perfect, hence also $\omega_\phi \simeq \det(\LL_\phi)$ (see Subsect.~\ref{ssec:indperf dualizing complex}).
	The claim now follows from the description of (the dual of) $\LL_\phi$ given in Theorem~\ref{reltangent}.
\end{proof}

\subsubsection{Base change and composition}

\begin{construction}\label{constr:vartheta^! funct}
Let $f : X' \to X$ be a morphism in $\Var^G$ and $V \to X$ a $G$-equivariant vector bundle.
Write $s$ and $s'$ for the zero sections of $V \to X$ and $f^*V \to X'$, and $f' : f^*V \to V$ for the base change of $f$.
The canonical morphism $f^* \vartheta_A^!(V) \to \vartheta_A^!(f^*V)$ in $\Shv_A^G(X')$ is defined as the exchange transformation
\[
  \operatorname{Ex}^{*,\dag} : f^* s^\dag (\un) \to s'^\dag f'^* (\un)
\]
of Construction~\ref{constr:Ex f^* f^dag}.
\end{construction}

\begin{construction}\label{constr:vartheta^! sum}
Let $V$ and $W$ be $G$-equivariant vector bundles over $X \in \Var^G$.
Denote by $s$ the zero section of $p : V \to X$.
We define $\vartheta^!_A(V) \otimes \vartheta^!_A(W) \to \vartheta^!_A(V \oplus W)$ as the composite
\[
  \vartheta^!_A(V) \otimes s^*p^*\vartheta^!_A(W)
  \xrightarrow{\mathrm{can}} \vartheta^!_A(V) \otimes s^*\vartheta^!_A(p^*W)
  \to \vartheta^!_A(V \oplus W),
\]
where the first map is Construction~\ref{constr:vartheta^! funct}, and the second is Construction~\ref{constr:omega composite indperf}.
\end{construction}

\subsection{Push-forwards}
\label{ssec:pushforwards}

Let $f : X \to Y$ be a proper morphism in $\Var^G$.

\begin{construction}\label{push}
Let $\cF \in \Shv_{A,\mathrm{c}}^G(Y)$.
Pushing forward the trace map $\operatorname{tr}^A_f : f_*(\omega_f \otimes f^*\cF) \simeq f_* f^\dag(\cF) \to \cF$ \eqref{eq:trace maps} along $Y \to *$ gives the \emph{push-forward} along $f$ with coefficients in $\cF$:
\[
  f_* : \HH_A^G(X; \omega_f \otimes f^*\cF) \longrightarrow \HH_A^G(Y; \cF).
\]
We also write $f_* : \Hcoh_A^G(X; \omega_f \otimes f^*\cF) \longrightarrow \Hcoh_A^G(Y; \cF)$ for the induced map on periodizations.
\end{construction}

\begin{example}
	Let $k$ be a $\bQ$-algebra, and $A = \GGm$ or $E$.
	Then by Proposition~\ref{prop:theta^! vs theta}, for a proper morphism $f : X \to Y$ in $\Var^G$ with $X, Y$ smooth, and the coefficient $\cF=\OO$, the push-forward becomes
	\[
	  f_* : \HH_A^G(X; \vartheta_A(\TT_f)^{-1}) \longrightarrow \HH_A^G(Y).
	\]
\end{example}

\begin{proposition}\label{prop:push properties}
Let $f : X \to Y$ be a proper morphism in $\Var^G$.
\begin{enumerate}
\item \emph{Functoriality.}
If $g : Y \to Z$ is a proper morphism in $\Var^G$, then for every $\cF \in \PERF^{\check{A}}(\Loop_A^G(Z))$ the composite
\[
  \HH_A^G(X; (gf)^\dag\cF) \xrightarrow{f_*} \HH_A^G(Y; g^\dag\cF) \xrightarrow{g_*} \HH_A^G(Z; \cF)
\]
is homotopic to $(gf)_*$.

\item \emph{Projection formula.}
For all $\cF, \cG \in \PERF^{\check{A}}(\Loop_A^G(Y))$, the square
\[\begin{tikzcd}
  \HH_A^G(X; \omega_f \otimes f^*\cF) \otimes \HH_A^G(Y; \cG) \ar{r}{f_* \otimes \id}\ar{d}{\cup \circ (\id \otimes f^*)}
  & \HH_A^G(Y; \cF) \otimes \HH_A^G(Y; \cG) \ar{d}{\cup}
  \\
  \HH_A^G(X; \omega_f \otimes f^*(\cF \otimes \cG)) \ar{r}{f_*}
  & \HH_A^G(Y; \cF \otimes \cG)
\end{tikzcd}\]
commutes.
\end{enumerate}
Moreover, the same holds on $\Hcoh_A^G(-)$.
\end{proposition}
\begin{proof}
Both claims follow from the analogous property of the counit $\phi_*^{\QCOH}\phi^!_{\QCOH}(\cF) \to \cF$, by applying $\Ren^R$ (and Lemma~\ref{perftens} and Corollary~\ref{rights}) and the (lax symmetric monoidal) pushforward to $\Loop_A^G(*)$.
\end{proof}

\begin{remark}\label{rem:push not perfect}
Suppose $\cF = \Ren^R(F)$ where $F \in \QCOH^{\check{A}}(\Loop_A^G(Y))$ is not perfect.
Then we can extend the above construction by
\[
  f_* : \HH_A^G(X; \Ren^R\phi^!_{\QCOH}F) \to \HH_A^G(Y; \cF).
\]
This was used implicitly in the formulation of Proposition~\ref{prop:push properties}(i).
\end{remark}

\subsection{Euler classes}\label{ssec:euler}

We now construct the Euler class $e_{A}(V)$ for $V$ a $G$-equivariant vector bundle on $X$.

\begin{construction}
Let $p : V \to X$ be a $G$-equivariant vector bundle over $X$ with zero section $s : X \to V$.
Using the trace map $\operatorname{tr}_s^A : s_*s^\dag(\un) \to \un$ \eqref{eq:trace maps} and the unit map $\un \to s_*\un$, we get the morphism
\[
  \vartheta^!_A(V) = s^\dag(\un) \simeq p_* s_* s^\dag(\un)
  \xrightarrow{\operatorname{tr}^A_s} p_*\un
  \xrightarrow{p_* \mathrm{unit}} p_* s_*\un \simeq \un
\]
in $\Shv_A^G(X)$.
We call the composite $e_A(V) : \vartheta^!_A(V) \to \un$ the \emph{Euler class} of $V$, and also regard it as a map
\[
  e_A(V) : \HH_A^G(X; \vartheta^!_A(V)) \to \HH_A^G(X)
\]
in $\Rot^G(A)$ by pushing forward.
By periodization, it also descends to a map on $\Hcoh_A^G(-)$ which we denote the same way.
\end{construction}

\begin{remark}
The map $e_A(V) : \HH_A^G(X; \vartheta^!_A(V)) \to \HH_A^G(X)$ is the composite
\[
  \HH_A^G(X; \vartheta^!_A(V)) \xrightarrow{s_*}
  \HH_A^G(V) \xrightarrow{s^*}
  \HH_A^G(X).
\]
\end{remark}

\begin{remark}
When $\vartheta^!_A(V)$ is invertible, $e_A(V) \otimes \vartheta^!_A(V)^{-1}$ is a map $\OO \to \vartheta^!_A(V)^{-1}$, i.e.\ a class $e_A(V) \in \HH_A^G(X; \vartheta^!_A(V)^{-1})$.
\end{remark}

\begin{remark}
	If $k$ is a $\bQ$-algebra and $X$ is smooth, $\vartheta^!_A(V) \simeq \vartheta_A(V)$ (Proposition~\ref{prop:theta^! vs theta}).
	Thus in that case, $e_A(V)$ can be regarded as a class in (i.e., section of) $\HH_A^G(X; \vartheta_A(V)^{-1})$.
\end{remark}

\begin{proposition}\label{prop:euler properties}
Let $V$ be a $G$-equivariant vector bundle over $X \in \Var^G$.
\begin{enumerate}
\item \emph{Naturality.}
For every morphism $f : X' \to X$ in $\Var^G$, there is a commutative diagram
\[\begin{tikzcd}[column sep=small]
  \HH_A^G(X; \vartheta^!_A(V)) \ar{r}{f^*}\ar[swap]{d}{e_A(V)}
  & \HH_A^G(X'; f^*\vartheta^!_A(V)) \ar{r}{\mathrm{can}}
  & \HH_A^G(X'; \vartheta^!_A(f^*V)) \ar{d}{e_A(f^*V)}
  \\
  \HH_A^G(X) \ar{rr}{f^*}
  &
  & \HH_A^G(X')
\end{tikzcd}\]
where the map labelled $\mathrm{can}$ comes from Construction~\ref{constr:vartheta^! funct}.

\item \emph{Sum formula.}
For every $G$-equivariant vector bundle $W$ over $X$, there is a commutative diagram
\[\begin{tikzcd}[column sep=large]
  \HH_A^G(X; \vartheta^!_A(V)) \otimes \HH_A^G(X; \vartheta^!_A(W)) \ar{r}{e_A(V) \otimes e_A(W)}\ar[swap]{d}{\cup}
  & \HH_A^G(X) \otimes \HH_A^G(X) \ar{dd}{\cup}
  \\
  \HH_A^G(X; \vartheta^!_A(V) \otimes \vartheta^!_A(W)) \ar[swap]{d}{\mathrm{can}}
  &
  \\
  \HH_A^G(X; \vartheta^!_A(V \oplus W)) \ar{r}{e_A(V \oplus W)}
  & \HH_A^G(X)
\end{tikzcd}\]
where the map labelled $\mathrm{can}$ comes from Construction~\ref{constr:vartheta^! sum}.
\end{enumerate}
\end{proposition}

\begin{proof}
Standard diagram chases that we leave to the reader.
\end{proof}

\appendix

\section{$E_\infty$-group stacks}
\label{app:strictification}

\subsection{Definitions}

Let $S$ be a stack over $k$.
An \emph{$E_\infty$-group stack} over $S$ is a grouplike $E_\infty$-monoid object in stacks over $S$, in the sense of \cite[Def.~1.3.8]{LurieEll1}.

We denote by $\GrpEoo_S$ the category of $E_\infty$-group stacks over $S$, and by $\CGrp_S$ the category of commutative group stacks over $S$, taken in the strict sense of Subsect.~\ref{ssec:intro/conv}.
There is a forgetful functor $\CGrp_S \to \GrpEoo_S$, which identifies $\CGrp_S$ with the category of $\bZ$-module objects in $\GrpEoo_S$ (see \cite[Constr.~1.3.9, Ex.~1.2.14]{LurieEll1}).
We write $\HOM_S^{E_\infty}(-,-)$ for the internal hom in $\GrpEoo_S$, and $\HOM_S(-,-)$ for the internal hom in $\CGrp_S$.

For example, $\GGm$ is a commutative group scheme, while $\GL_1$ is only an $E_\infty$-group scheme (cf. Example~\ref{ex:strict units comparison}).

\subsection{Strictification}

The forgetful functor $\CGrp_S\to\GrpEoo_S$ admits left and right adjoints.
The left adjoint is the free $\bZ$-module functor, and the right adjoint is the \emph{strictification} functor
\[
  (-)^{\mathrm{str}}: \GrpEoo_S\to \CGrp_S,
\]
given by
\[
  X^{\mathrm{str}}(T)
  =
  \Hom_{\GrpEoo}
  \bigl(\bZ, X(T)\bigr)
\]
for $X\in\GrpEoo_S$ and affine $T$ over $S$.

Unravelling definitions, for $X\in\GrpEoo_S$ and $Y\in\CGrp_S$ we have a natural isomorphism of group-theoretic mapping stacks
\begin{equation}\label{eq:str adjunction}
  \HOM_S\bigl(Y,X^{\mathrm{str}}\bigr)
  \simeq
  \HOM_S^{E_\infty}(Y,X),
\end{equation}
and a natural isomorphism
\[
  (X^{\mathrm{str}})_{\mathrm{cl}}
  \simeq
  (X_{\mathrm{cl}})^{\mathrm{str}}.
\]

\begin{example}\label{ex:strict units comparison}
  We denote by $\GGm$ the commutative group scheme $\spc(k[\mathbf{Z}])$ over $k$, and by $\GL_1$ the $E_\infty$-group scheme representing the functor which sends a connective $E_\infty$-$k$-algebra $R$ to its $E_\infty$-group of units $R^\times$.
  There is a canonical homomorphism of $E_\infty$-groups $\GGm \to \GL_{1}$, which is an isomorphism if and only if $k$ is a $\bQ$-algebra, and which is always an isomorphism on classical truncations.
  By \cite[Constr.~1.6.10 and Rems.~1.6.11--1.6.12]{LurieEll2}, we have
  \[
    \GGm \simeq (\GL_{1})^{\mathrm{str}}.
  \]
\end{example}

\subsection{Strictification in low degrees and characteristic zero}

\begin{lemma}\label{lem:strictification stable defect}
Let $X$ be an $E_\infty$-group such that $\pi_iX$ is a $\bQ$-vector space for every $i\geq2$.
Then the canonical map $X^{\mathrm{str}} \to X$ is an isomorphism on $\pi_i$ for $i\geq 1$ and on $\pi_0$ induces an isomorphism
\[
  \pi_0(X^{\mathrm{str}}) \simeq \ker\bigl(\eta:\pi_0X\to\pi_1X\bigr),
\]
where $\eta\in\pi_1(\mathbf{S})$ is the stable Hopf element.
\end{lemma}

\begin{proof}
When $X$ is $1$-truncated, this follows from the proof of \cite[Prop.~3.8]{CarmeliStrictPicard}.
In general, consider the fibre $X_{\geq2}$ of $X\to X_{\leq1}$ and form the commutative diagram of fibre sequences
\[
  \begin{tikzcd}
    (X_{\geq2})^{\mathrm{str}} \ar[r] \ar[d] & X^{\mathrm{str}} \ar[r] \ar[d] & (X_{\leq1})^{\mathrm{str}} \ar[d] \\
    X_{\geq2} \ar[r] & X \ar[r] & X_{\leq1}.
  \end{tikzcd}
\]
Since $X_{\geq 2}$ is rational by assumption and $\mathbf S\to H\bZ$ is invertible after rationalization, the left-hand vertical map is invertible.
In degrees $i\geq2$, $\pi_i$ of both right-hand terms vanish, so the middle vertical morphism is identified on $\pi_i$ with the left-hand morphism.
In degrees $i \in \{0,1\}$, $\pi_i$ of both left-hand terms vanish, so the middle vertical morphism is identified on $\pi_i$ with the right-hand morphism, and the claim follows from the $1$-truncated case.
\end{proof}

\begin{corollary}\label{cor:strictification classical strict source}
Let $Y\in\CGrp_S$ be a commutative group stack over $S$.
If $Y(T)$ is $1$-truncated for every classical affine $T$ over $S$, then $Y\to Y^{\mathrm{str}}$ induces an isomorphism on classical truncations.
\end{corollary}

\begin{proof}
By \cite[Ex.~1.2.14]{LurieEll1}, $Y$ is a $\bZ$-module object in $\GrpEoo_S$; in particular, for every $T$ over $S$, the connective spectrum associated to $Y(T)$ is an $H\bZ$-module.
It follows that the action of $\eta\in\pi_1\mathbf S$ on its homotopy groups factors through $\pi_1H\bZ=0$.
Thus for every classical affine $T$ over $S$, Lemma~\ref{lem:strictification stable defect} applied to the underlying $E_\infty$-group of $Y(T)$ shows that the counit map $Y^{\mathrm{str}}(T) \to Y(T)$ is invertible.
The triangle identity identifies its inverse with the underlying $E_\infty$-map of the adjunction unit $Y(T)\to Y^{\mathrm{str}}(T)$, so the latter is also invertible.
Allowing $T$ to vary over classical affines over $S$, it follows that $Y \to Y^{\mathrm{str}}$ is invertible on classical truncations.
\end{proof}

The same argument also shows:

\begin{corollary}\label{cor:strictification strict source}
Let $Y \in \CGrp_S$ be a commutative group stack over $S$.
Suppose that for every affine $T$ over $S$, $\pi_i Y(T)$ is a $\bQ$-vector space for every $i\geq2$.
Then the canonical map $Y\to Y^{\mathrm{str}}$ is invertible.
\end{corollary}

\bigskip
{\footnotesize \parskip 0pt \parindent 0pt
\textsc{Emile Bouaziz}\\[0.2em]
Shanghai Institute for Mathematics and Interdisciplinary Sciences\\
Fudan University, Block A, International Innovation Plaza\\
No.~657 Songhu Road, Yangpu District, Shanghai, China

\medskip
\textsc{Adeel A. Khan}\\[0.2em]
Institute of Mathematics, Academia Sinica\\
Astronomy-Mathematics Building, No.~1, Sec.~4, Roosevelt Road\\
Taipei 106319, Taiwan\\[0.4em]
National Center for Theoretical Sciences, National Taiwan University\\
Cosmology Building, No.~1, Sec.~4, Roosevelt Road\\
Taipei 106319, Taiwan\par}

\end{document}